\documentclass{amsart}[10pt]
\usepackage{amssymb}\usepackage{euscript}
\usepackage[all]{xy}

\newtheorem{theorem}{Theorem}[section]
\newtheorem{lemma}[theorem]{Lemma}
\newtheorem{proposition}[theorem]{Proposition}
\newtheorem{corollary}[theorem]{Corollary} 
\newtheorem{definition}[theorem]{Definition}
\newtheorem{lemmadef}[theorem]{Lemma-Definition} 
\newtheorem{propdef}[theorem]{Proposition-Definition} 
\theoremstyle{definition}

\newtheorem{conjecture}[theorem]{Conjecture}  
\newtheorem{remark}[theorem]{Remark}

\newcommand{\HH}{{\mathfrak{H}}}\newcommand{\on}{\operatorname}   
\newcommand{\eps}{\varepsilon}

\newcommand{\wh}{\widehat} 
\newcommand{\ul}{\underline}

\newcommand{\cC}{{\mathcal C}}
\newcommand{\cF}{{\mathcal F}}
\newcommand{\cD}{{\mathcal D}}
\newcommand{\cG}{{\mathcal G}}

\newcommand{\CC}{{\mathbb{C}}}\newcommand{\QQ}{{\mathbb{Q}}}
\newcommand{\RR}{{\mathbb{R}}}
\newcommand{\ZZ}{{\mathbb{Z}}}\newcommand{\NN}{{\mathbb N}}

\newcommand{\kk}{{\mathbf k}}

\renewcommand{\b}{{\mathfrak{b}}}
\newcommand{\f}{{\mathfrak{f}}}\newcommand{\h}{{\mathfrak{h}}}

\renewcommand{\t}{{\mathfrak{t}}}\renewcommand{\u}{{\mathfrak{u}}}
\newcommand{\vv}{{\mathfrak{v}}}

\newcommand{\G}{{\mathfrak{G}}}
\newcommand{\SL}{{\mathfrak{sl}}}\newcommand{\GL}{{\mathfrak{gl}}}

\renewcommand{\r}{{\mathfrak{r}}}
\newcommand{\g}{{\mathfrak{g}}}\newcommand{\grt}{{\mathfrak{grt}}}
\newcommand{\gt}{{\mathfrak{gt}}}

\begin{document}

\title[Elliptic associators]{Elliptic associators}

\author{Benjamin Enriquez}
\address{Institut de Recherche Math\'ematique Avanc\'ee, UMR 7501, 
Universit\'e de Strasbourg et CNRS, 
7 rue Ren\'e Descartes, 67000 Strasbourg, France}
\email{b.enriquez@math.unistra.fr}

\begin{abstract}
We construct a genus one analogue of the theory of associators and the 
Grothen\-dieck-Teichm\"uller group. The analogue of the Galois action on the profinite 
braid groups is an action of the arithmetic fundamental group of a moduli space 
$M_{1,\vec1}^{\QQ}$ of elliptic curves on the profinite braid groups in genus one. 
This action factors through an explicit profinite group $\widehat{\on{GT}}_{ell}$, which 
admits an interpretation in terms of decorations of braided monoidal categories. 
We relate this group to a prounipotent group scheme $\on{GT}_{ell}(-)$. We construct 
a torsor over this group, the scheme of elliptic associators. An explicit family of elliptic 
associators is constructed, based on our earlier work with Calaque and Etingof 
on the universal KZB connexion. The existence of elliptic associators enables 
one to show that the Lie algebra of $\on{GT}_{ell}(-)$ is isomorphic to a graded Lie 
algebra, on which we obtain several results: semidirect product structure; 
explicit generators. This existence also allows one to compute the Zariski closure of 
the mapping class group $B_{3}$ in genus one in the automorphism groups 
of the prounipotent completions of braid groups in genus one. The analytic 
study of the family of elliptic associators produces relations between MZVs 
and iterated integrals of Eisenstein series.
\end{abstract}
\maketitle

\section*{Introduction}

The theory of associators and of the Grothendieck-Teichm\"uller group has been 
developed by Drinfeld (\cite{Dr:Gal}) in relation with certain problems of quantum groups. 
This theory was based on several previous pieces of work: on the one hand, the 
approach proposed by Grothendieck to the study of $G_\QQ$, the absolute Galois 
group of $\QQ$, via its action on the Teichm\"uller tower in genus zero, and in 
particular on the profinite completions of the braid groups (\cite{G:Esq}); on the 
other hand, rational homotopy theory, in particular the computation by 
Kohno of the prounipotent completions of the pure braid groups, based 
on the study of a particular connexion on the configuration spaces of 
the plane, which may be identified with a universal version of the 
Knizhnik-Zamolodchikov (KZ) connection.  

The main actors of associator theory are: a profinite group $\widehat{\on{GT}}$, of 
categorical origin, containing $G_{\QQ}$; pro-$l$, proalgebraic variants of 
this group, and the associated Lie algebra ${\mathfrak{gt}}$; a principal homogeneous 
space, the space of associators, which enables one to prove that ${\mathfrak{gt}}$ is 
isomorphic to a graded Lie algebra ${\mathfrak{grt}}$; a particular associator, the 
KZ associator, whose study allows one both to derive a system of 
relations between multizeta values (MZVs) and a collection of generators 
for ${\mathfrak{grt}}$. The theory of associators is therefore related to the theory of 
MZVs and motives (\cite{An}); it allows one to exhibit conditions satisfied 
by elements of motivic Lie algebras. 

The purpose of the present work is to construct the analogous theory in genus one. The 
analogue of the Galois action on the braid groups in genus zero is the action of the arithmetic 
fundamental group of the moduli space of elliptic curves $M_{1,1}^{\QQ}$ on the profinite 
completions of braid groups in genus one, studied in \cite{G:LM,O}. The analogue of the 
rational homotopy part is the computation of the prounipotent completion of braid 
groups in genus one, first obtained by Bezrukavnikov using minimal model theory, and later 
rederived in \cite{CEE} using an analogue of the KZ connection, the universal KZB 
connection (this connection was independently obtained in \cite{LR}). A new feature of the 
KZB connection is its horizontal part (related to variation of the elliptic modulus), which 
corresponds to an extension of the holonomy Lie algebra ${\mathfrak t}_{1,n}$ by a 
Lie algebra of derivations $\langle\delta_{2n},n \geq -1\rangle$. 

Our construction of the genus one analogue of Grothendieck-Teichm\"uller theory consists 
of several steps. We first construct a genus one analogue of the theory of braided monoidal 
categories (BMCs). This enables us to define the genus one analogue 
$\widehat{\on{GT}}_{ell}$ of $\widehat{\on{GT}}$, which is a profinite group 
containing $\pi_1(M_{1,\vec1}^{\QQ})$. We construct the pro-$l$ and proalgebraic variants 
of this group; the associated Lie algebra is denoted ${\mathfrak{gt}}_{ell}$. We 
construct a torsor 
under this proalgebraic group: the scheme of elliptic associators. We present two 
constructions of elliptic associators: (a) we define an explicit map from the set of associators 
to its elliptic analogue; (b) the KZB connection gives rise to a map $\tau\to e(\tau)$
from the Poincar\'e 
half-plane to the set of elliptic associators. We study the properties of this map: differential 
system, modular behavior, behavior at infinity; this shows in particular that the constructions 
(a) and (b) are related to each other by suitable specializations and limiting procedures. The 
existence of elliptic associators then enable us to construct an isomorphism between 
${\mathfrak{gt}}_{ell}$ and an explicit Lie algebra ${\mathfrak{grt}}_{ell}$. 
We prove several results on ${\mathfrak{grt}}_{ell}$: (a) ${\mathfrak{grt}}_{ell}$ 
is a semidirect product of ${\mathfrak{grt}}$ by a Lie algebra ${\mathfrak r}_{ell}$, 
which is therefore acted upon by ${\mathfrak{grt}}$; (b) we construct explicit generators 
of ${\mathfrak r}_{ell}$. 

Beside these results, which may be viewed as internal to the theory, our work leads to the 
following results: 

(a) the outer action of the arithmetic fundamental group of $M_{1,\vec1}^{\QQ}$ on the 
$\QQ_{l}$-points 
of the prounipotent completions of the braid group in genus one with various numbers of 
strands factors through the action of the group of $\QQ_{l}$-points of one and the same 
proalgebraic group, which is $\on{GT}_{ell}(-)$; 

(b) the mapping class group in genus one, which is isomorphic to the group $B_3$ of braids 
with three strands, naturally acts on the pure braid groups in genus one. We compute the 
Zariski closure of $B_3$ in the automorphism group of their prounipotent completions, in 
terms of the Lie algebra $\langle\delta_{2n},n \geq -1\rangle$;

(c) the study of the above-mentioned map from the Poincar\'e half-plane to the space of 
elliptic associators leads to relations between MZVs and iterated integrals of Eisenstein 
series.  

This paper is organized as follows. In Section 1, we define the genus one counterpart 
of the notion of braided monoidal category. This enables us to define the group 
$\widehat{\on{GT}}_{ell}$ in Section 2, as well as its pro-$l$ and prounipotent 
variants. In Section 3, we introduce the space of elliptic associators, prove its 
nonemptiness and study its torsor structure. This leads us to the definition of the group
scheme $\on{GRT}_{ell}(-)$ in Section 4; we prove the announced results on its Lie 
algebra ${\mathfrak{grt}}_{ell}$: isomorphism with ${\mathfrak{gt}}$, 
generators, semidirect product structure. In Section 5, we carry out the computation of 
Zariski closure of $B_{3}$ explained above. In Section 6, we introduce the map 
$\tau\mapsto e(\tau)$ and study its properties. We define the iterated integrals of 
Eisenstein series in Section 7, and prove there their relations with MZVs. In Section 8, 
we recall the relations between $G_{\QQ}$, $\widehat{\on{GT}}$ and the 
Teichm\"uller groupoid in genus zero, and generalize these results to genus one. 
Section 9 raises a question on the structure of the kernel ${\mathfrak r}_{ell}$
of a natural morphism ${\mathfrak{grt}}_{ell}\to {\mathfrak{grt}}$, and its relation 
with a transcendence conjecture on the KZ associator (which is related to the 
Grothendieck period conjecture); namely, it is shown that an affirmative 
answer to both questions imply the same (also conjectural) statement on the behavior of 
certain isomorphisms arising from associators (see Propositions \ref{prop:77} and 
\ref{prop:78}). 

Let us now mention some works and projects related to the present work. Hain and 
Matsumoto construct a theory of ``mixed elliptic motives'' (\cite{HM:ell}). This gives rise to a 
proalgebraic $\QQ$-group scheme $G_{MEM}(-)$, equipped with a morphism 
$G_{MEM}(-)\to G_{MTM}(-)$. One may expect a commutative diagram from this 
morphism to $\on{GT}_{ell}(-)\to \on{GT}(-)$. The Lie algebra $\langle\delta_{2n},
n\geq -1\rangle$ is a Lie subalgebra of the graded version of the kernels of both 
morphisms and was studied in Pollack's Ph.D. thesis (\cite{Po}). On the other hand, 
Brown and Levin develop a parallel theory of elliptic motives (\cite{BL}); the elliptic 
multiple zeta values arising from this theory could be related to the family 
$\tau\mapsto e(\tau)$ of elliptic associators studied here.  

The author expresses his thanks to P. Etingof, H. Furusho, R. Hain, 
and L. Schneps for useful discussions related to this work, as well as 
to D. Calaque and P. Etingof for collaboration in \cite{CEE}.

\section{Elliptic structures over braided monoidal categories}
\label{sec:1}

In \cite{CEE}, we introduced a notion of elliptic structure over a braided monoidal
category (BMC) $\cC$. It consists in a category ${\mathcal E}$, a functor 
${\mathcal E}\to \cC$, 
and additional data. In this section, we introduce a variant of this notion, which 
consists in a category $\tilde \cC$, a functor $\cC\to\tilde\cC$, and additional data. 
The two definitions can be related by adjunction, as will be explained in forthcoming 
joint work with P.~Etingof. As is the notion from \cite{CEE}, the variant 
presented here  is related with elliptic braid groups in the same way as BMCs 
are related to usual braid groups.

\subsection{Definition}

Let $(\cC,\otimes,\beta_{\ldots},a_{\ldots},{\bf 1})$ 
be a braided monoidal category (see e.g.~\cite{Ka}). Here
$\otimes:\cC\times\cC$ is the tensor product,  
$\beta_{X,Y} : X\otimes Y\to Y \otimes X$ 
and $a_{X,Y,Z}: (X\otimes Y)\otimes Z\to 
X\otimes (Y\otimes Z)$ are the braiding and associativity 
isomorphisms and ${\bf 1}$ is the unit object. They satisfy in particular 
the pentagon and hexagon identities
$$
a_{X,Y,Z\otimes T}a_{X\otimes Y,Z,T} = (\on{id}_X\otimes a_{Y,Z,T})
a_{X,Y\otimes Z,T}(a_{X,Y,Z}\otimes\on{id}_T), 
$$
$$
(\on{id}_Y \otimes \beta_{X,Z}^{\pm}) a_{Y,X,Z}
(\beta_{X,Y}^{\pm} \otimes \on{id}_Z) = a_{Y,Z,X} 
\beta_{X,Y\otimes Z}^{\pm} a_{X,Y,Z}, 
$$
where
$\beta_{X,Y}^+=\beta_{X,Y}$, $\beta_{X,Y}^-=\beta_{Y,X}^{-1}$. 

\begin{definition}
An elliptic structure over the braided monoidal category $\cC$ is a set 
$(\tilde\cC,F,A_{\ldots}^+,A^-_{\ldots})$, where $\cC$ is a 
category, $F:\cC\to\tilde\cC$ is a 
functor\footnote{For $\cC$ a category, $\on{Ob}\cC$ is its class of 
objects; for $X,Y\in\on{Ob}\cC$, $\on{Iso}_\cC(X,Y)\subset 
\cC(X,Y)$ are the sets of isomorphisms and morphisms
$X\to Y$; $\on{Aut}_\cC(X)=\on{Iso}_\cC(X,X)$.}, 
and $A^\pm_{\ldots}$ are natural\footnote{Natural means 
that if $\varphi\in \cC_0(X,X')$, $\psi\in\cC_0(Y,Y')$, 
then $A^\pm_{X',Y'} F(\varphi\otimes\psi) = F(\varphi\otimes\psi) 
A^\pm_{X,Y}$} assignments $(\on{Ob}\cC)^2\ni (X,Y)\mapsto A^\pm_{X,Y}
\in \on{Aut}_{\tilde\cC}(F(X\otimes Y))$, such that: 
\begin{equation} \label{def:ell:str:1}
\alpha^\pm_{Z,X,Y}\alpha^\pm_{Y,Z,X}\alpha^\pm_{X,Y,Z}
=\on{id}_{(X\otimes Y)\otimes Z}, 
\end{equation}
where $\alpha^\pm_{X,Y,Z} = F(\beta^\pm_{X,Y\otimes Z}) A^\pm_{X,Y\otimes Z}
F(a_{X,Y,Z})$, 
\begin{eqnarray} \label{def:ell:str:3}
F(\beta_{Y,X}\beta_{X,Y}\otimes\on{id}_Z) = 
&& \nonumber
\Big( F(a^{-1}_{X,Y,Z}) A^-_{X,Y\otimes Z} F(a_{X,Y,Z}), 
\\ && 
F\big( (\beta_{X,Y}^{-1}\otimes\on{id}_Z)a^{-1}_{Y,X,Z} \big) 
(A^+_{Y,X\otimes Z})^{-1} F \big( 
a_{Y,X,Z}(\beta_{Y,X}^{-1}\otimes\on{id}_Z)\big)\Big),  
\end{eqnarray}
(identities\footnote{In (\ref{def:ell:str:3}) and later, 
we set $(g,h):= ghg^{-1}h^{-1}$.} 
in $\on{Aut}_{\tilde\cC}(F((X\otimes Y)\otimes Z))$), for any 
$X,Y,Z\in \on{Ob}\cC$, and
\begin{equation} \label{def:ell:str:4}
A^\pm_{{\bf 1},X} = \on{id}_{F({\bf 1}\otimes X)} \quad
{\text for\ any\ }X\in\on{Ob}\cC. 
\end{equation}
\end{definition}

Dropping associativity constraints and the functor $F$ (which can be put in automatically), 
the two first conditions mean that the cycles 
$$
\xymatrix{
\scriptstyle{X\otimes Y\otimes Z} \ar[r]^{A^{\pm}_{X,Y\otimes Z}}& 
\scriptstyle{X\otimes Y\otimes Z} \ar[r]^{\beta^{\pm}_{X,Y\otimes Z}}& 
\scriptstyle{Y\otimes Z\otimes X} \ar[d]^{A^{\pm}_{Y,Z\otimes X}} \\
\scriptstyle{Z\otimes X\otimes Y}\ar[u]^{\beta^{\pm}_{Z,X\otimes Y}} & 
\scriptstyle{Z\otimes X\otimes Y} \ar[l]^{A^{\pm}_{Z,X\otimes Y}} & 
\scriptstyle{Y\otimes Z\otimes X}\ar[l]^{\beta^{\pm}_{Y,Z\otimes X}}
}
\quad \on{and} \quad 
\xymatrix{
\scriptstyle{Y\otimes X\otimes Z} \ar[r]^{A_{Y,X\otimes Z}}& 
\scriptstyle{Y\otimes X\otimes Z}\ar[r]^{\beta_{YX}\otimes\on{id}_{Z}}& 
\scriptstyle{X\otimes Y\otimes Z} \ar[d]^{B_{X,Y\otimes Z}^{-1}}\\
\scriptstyle{X\otimes Y\otimes Z}\ar[u]^{\beta_{YX}^{-1}\otimes\on{id}_{Z}} & & 
\scriptstyle{X\otimes Y\otimes Z}\ar[d]^{\beta_{YX}^{-1}\otimes\on{id}_{Z}} \\
\scriptstyle{X\otimes Y\otimes Z}\ar[u]^{B_{X,Y\otimes Z}}  & 
\scriptstyle{Y\otimes X\otimes Z}\ar[l]^{\beta_{XY}^{-1}\otimes \on{id}_{Z}} & 
\scriptstyle{Y\otimes X\otimes Z}\ar[l]^{A^{-1}_{Y,X\otimes Z}}}$$
are identity morphisms, where $A_{\cdots} = A^{+}_{\cdots}$, $B_{\cdots} = 
A^{-}_{\cdots}$. 

A morphism $(\cC,\tilde\cC,F,A^\pm_{\ldots})\to (\cC',\tilde\cC',F',
A^{\prime\pm}_{\ldots})$ is then the data of a tensor functor
$\cC\stackrel{\varphi}{\to} \cC'$ and a functor 
$\tilde\cC\stackrel{\tilde\varphi}{\to}\tilde\cC'$, 
such that $\begin{matrix} \cC&\to &\cC' \\ \downarrow & & \downarrow 
\\ \tilde\cC& \to & \tilde\cC'\end{matrix}$ commutes, and 
$\tilde\varphi(A_{X,Y}^\pm) = A^{\prime\pm}_{\varphi(X),\varphi(Y)}$. 

\begin{remark} By setting $Z = {\bf 1}$, the axioms 
(\ref{def:ell:str:1})-(\ref{def:ell:str:4}) imply 
\begin{equation} \label{for:quot}
F(\beta^\pm_{Y,X})A^\pm_{Y,X} 
F(\beta^\pm_{X,Y}) A^\pm_{X,Y} = \on{id}_{F(X\otimes Y)}, \quad 
F(\beta_{Y,X}\beta_{X,Y}) = (A^-_{X,Y},A^+_{X,Y}),
\end{equation} 
which in their turn imply 
\begin{equation} \label{3bis}
A^{\pm}_{X,{\bf 1}} = \on{id}_{F(X\otimes{\bf 1})}.
\end{equation} 
Taking these identities, (\ref{def:ell:str:4}) and the hexagon identities 
into account, axiom (\ref{def:ell:str:1}) can be replaced by 
\begin{equation} \label{rel:variant}
A^\pm_{X\otimes Y,Z} = F((\beta^\pm_{Y,X}\otimes \on{id}_Z)a^{-1}_{Y,X,Z}) 
A^\pm_{Y,X\otimes Z} F(a_{Y,X,Z}(\beta^\pm_{X,Y}\otimes \on{id}_Z)
a^{-1}_{X,Y,Z}) A^\pm_{X,Y\otimes Z} F(a_{X,Y,Z}).  
\end{equation}
\end{remark}

\subsection{Relation with elliptic braid groups}

For $n\geq 1$, the reduced pure elliptic braid group on $n$ strands $P_{1,n}$
is the fundamental group of the reduced configuration space 
$\overline{\on{Cf}}_{n}(T) := \on{Cf}_{n}(T)/T$, where
$\on{Cf}_{n}(T) = T^{n} - ($diagonals) is the configuration space of 
$n$ points on the topological torus $T := \RR^{2}/\ZZ^{2}$, on which $T$
acts diagonally. The reduced elliptic braid group $B_{1,n}$ is the fundamental group
of the quotient $\overline{\on{Cf}}_{[n]}(T):= \overline{\on{Cf}}_{n}(T)/S_{n}$. 
We then have an exact sequence 
$$
1\to P_{1,n}\to B_{1,n}\to S_{n}\to 1. 
$$
These definitions are extended by $P_{1,0} = B_{1,0}=\{1\}$. 

The group $B_{1,n}$ ($n\geq 1$) can be presented by generators 
$\sigma_i$ ($i=1,\ldots,n-1$), $X_1^\pm$, and relations 
$$
(\sigma_1^{\pm1}X_1^\pm)^2 = 
(X_1^\pm\sigma_1^{\pm1})^2, \quad (X_1^\pm,\sigma_i)=1 \text{\ for\ } 
i=2,\ldots,n-1,\quad (X_1^-,(X^+_2)^{-1})=\sigma_1^2, 
\quad X_1^\pm\cdots X_n^\pm=1, 
$$
\begin{equation} \label{Artin}
(\sigma_i,\sigma_j)=1 \text{\ for\ }|i-j|>1, \quad 
\sigma_i\sigma_{i+1}\sigma_i=\sigma_{i+1}\sigma_i\sigma_{i+1}
\text{\ for\ }i=1,\ldots ,n-2,  
\end{equation}
where $X_{i+1}^\pm = \sigma_i^{\pm1}X_i^\pm\sigma_i^{\pm1}$
for $i=1,\ldots ,n-1$ (\cite{Bi}). In particular, $P_{1,1}=B_{1,1}=\{1\}$, 
and $P_{1,2}$ is the free group with two generators $X_{1}^{\pm}$. 

The braid group $B_n$ on $n$ strands  ($n\geq 1$) is presented by generators 
$\sigma_i$, $i=1,\ldots ,n-1$ and the Artin relations (\ref{Artin}). Its definition 
is extended to $n=0$ by $B_{0} = \{1\}$. There is a unique morphism
$B_n\to B_{1,n}$ such that $\sigma_i\mapsto \sigma_i$. 
If $\cC$ is a braided monoidal category and $X\in \on{Ob}\cC$, then 
there is a unique group morphism $\varphi:B_n\to\on{Aut}_\cC(X^{\otimes n})$
($X^{\otimes n}$ is defined by right parenthesization, so 
$X^{\otimes n}=X\otimes X^{\otimes n-1}$), such that 
$$
\sigma_i\mapsto a_i ((\on{id}_{X^{\otimes i-1}}\otimes\beta_{X,X})
\otimes\on{id}_{X^{\otimes n-i-1}}) a_i^{-1},  
$$
where $a_i : (X^{\otimes i-1}\otimes X^{\otimes 2})\otimes 
X^{\otimes n-i-1} \to X^{\otimes n}$ is the associativity constraint. 

\begin{proposition} \label{prop:gp:morphism}
If $(\tilde\cC,F,A^\pm_{\ldots})$ is an elliptic structure over $\cC$
and $X\in\on{Ob}\cC$, then there is a unique group morphism 
$B_{1,n}\to \on{Aut}_{\tilde\cC}(F(X^{\otimes n}))$, 
such that 
$$
X_1^\pm\mapsto A^\pm_{X,X^{\otimes n-1}}, \quad 
\sigma_i\mapsto F(\varphi(\sigma_i)). 
$$
\end{proposition}

{\em Proof.} Let us check that $(\sigma_1X_1^+)^2=(X_1^+\sigma_1)^2$, i.e., 
$(X_1^+,X_2^+)=1$ is preserved (for simplicity, we omit the associativity 
constraints). By naturality, $(\beta_{X,Y}\otimes\on{id}_Z)A^+_{X\otimes Y,Z}
=A^+_{Y\otimes X,Z}(\beta_{X,Y}\otimes\on{id}_Z)$. Plugging in this equality 
the relation (\ref{rel:variant}) and its analogue with $X,Y$ exchanged, we obtain 
$$
(A^+_{X,Y\otimes Z},F(\beta_{Y,X}\otimes\on{id}_Z)
A^+_{Y,X\otimes Z}F(\beta_{X,Y}\otimes\on{id}_Z))=1;
$$ if we set $Y:= X$, 
$Z:= X^{\otimes n-2}$, this says that $(X^+_1,X^+_2)=1$ is preserved. 
Similarly, one proves that (\ref{def:ell:str:1}) with sign $-$ implies that 
$(\sigma_1^{-1}X_1^-)^2 = (X_1^-\sigma_1^{-1})^2$ is preserved. 
(\ref{def:ell:str:3}) immediately implies that 
$(X_1^-,(X_2^+)^{-1})=\sigma_1^2$ is preserved. The naturality assumption
implies that $(X_1^{\pm},\sigma_i)=1$ ($i>1$) is preserved. 
One shows by induction that the image of $X_{k}^{\pm}$ is
$F(\beta_{X^{\otimes k-1},X}^{\pm}\otimes\on{id}_{X^{\otimes n-k}})
A^{\pm}_{X,X^{\otimes n-1}}F(\beta^{\pm}_{X,X^{\otimes k-1}}\otimes 
\on{id}_{X^{\otimes n-k}})$, therefore the image of $X_{1}^{\pm}
\cdots X_{k}^{\pm}$ is $A^{\pm}_{X^{\otimes k},X^{\otimes n-k}}$. 
It follows that the image of $X_1^\pm\cdots X^\pm_n$
is $A^\pm_{X^{\otimes n},{\bf 1}}$, which is $\on{id}_{F(X^{\otimes n})}$ 
by (\ref{3bis}). 
Finally, as $B_n\to\on{Aut}_{\cC}(F(X^{\otimes n}))$ is a 
group morphism, the Artin relations are preserved. 
\hfill\qed\medskip 

\subsection{Universal elliptic structures}

Let ${\bf PaB}$ be the braided monoidal category of parenthesized braids
(see \cite{JS,Ba}). Its set of objects is  ${\bf Par}:= \sqcup_{n\geq 0}{\bf Par}_n$, 
where ${\bf Par}_n=\{$parenthe\-si\-za\-tions of the word $\bullet\ldots\bullet$
of length $n\}$, so ${\bf Par}_0 = \{{\bf 1}\}$, ${\bf Par}_1=
\{\bullet\}$, ${\bf Par}_2=\{\bullet\bullet\}$, ${\bf Par}_3
=\{(\bullet\bullet)\bullet,\bullet(\bullet\bullet)\}$, etc.  
For $O,O'\in{\bf Par}$, we set $|O|:=$ the integer such that 
$O\in{\bf Par}_{|O|}$, and 
$\cC_0(O,O'):= 
  \left\{ \begin{array}{ll}
         B_{|O|} & \on{\ if\ } |O|=|O'|\\
        \emptyset & \on{otherwise}\end{array} \right.$. 
The composition is the 
product in $B_{|O|}$. The tensor product is defined 
at the level of objects, as the juxtaposition, and at the 
level of morphisms, by the group morphism $B_n\times B_{n'}
\to B_{n+n'}$, $(\sigma_i,1)\mapsto\sigma_i$, $(1,\sigma_j)
\mapsto\sigma_{n+j}$. We set $a_{O,O',O''}:=1\in B_{|O|+|O'|+|O''|}
= {\bf PaB}((O\otimes O')\otimes O'',O\otimes(O'\otimes O''))$
and $\beta_{O,O'}:=\sigma_{|O|,|O'|}\in B_{|O|+|O'|} = 
{\bf PaB}(O\otimes O',O'\otimes O)$, where $\sigma_{n,n'}:= 
(\sigma_n\cdots \sigma_1)(\sigma_{n+1}\cdots \sigma_2)\cdots 
(\sigma_{n+n'-1}\cdots 
\sigma_{n'})\in B_{n+n'}$. 

Let now ${\bf PaB}_{ell}$ be the category with 
the same objects, ${\bf PaB}_{ell}(O,O'):= 
  \left\{ \begin{array}{ll}
         B_{1,|O|} & \on{\ if\ } |O|=|O'|\\
        \emptyset & \on{otherwise}\end{array} \right.$
and whose product is the composition in $B_{1,|O|}$. Let $F:{\bf PaB}\to
{\bf PaB}_{ell}$ be the functor induced by the identity at the 
level of objects, and by $B_n\to B_{1,n}$, 
$\sigma_i\mapsto\sigma_i$ at the level of morphisms. 
For $O,O'\in{\bf Par}$, set $A^\pm_{O,O'}:= X^\pm_1\cdots X^\pm_{|O|}
\in B_{1,|O|+|O'|}$. Then $({\bf PaB}_{ell},F,A^\pm_{\ldots})$ is an 
elliptic structure over ${\bf PaB}$. Indeed, relations (\ref{def:ell:str:1}) 
and (\ref{def:ell:str:3})
for objects $O,O',O''$ are consequences of the identities $(\sigma_{2}^{\pm1}
\sigma_{1}^{\pm1}X_{1}^{\pm})^{3}=1$ and $(X_{1}^{-},(\sigma_{1}X_{1}^{+}
\sigma_{1})^{-1}) = \sigma_{1}^{2}$ in $P_{1,3}$ under the morphism 
$P_{1,3}\to P_{1,|O|+|O'|+|O''|}$ induced by the replacement of the 
first (resp., second, third) strand by $|O|$ (resp., $|O'|,|O''|$) consecutive strands. 

The pair $({\bf PaB},\bullet)$ has the following universal property: 
for any pair $(\cC,M)$, where $\cC$ is a braided monoidal 
category and $M\in\on{Ob}\cC$, there exists a unique tensor 
functor $\varphi_0: {\bf PaB}\to\cC$, such that $F(\bullet)=M$. 
Proposition \ref{prop:gp:morphism} immediately implies that this property 
extends as follows. 

\begin{proposition}
If $\tilde\cC$ is an elliptic structure over $\cC$, then there
exists a unique morphism $({\bf PaB},{\bf PaB}_{ell})\to(\cC,
\tilde\cC)$, extending $\varphi$. 
\end{proposition}

\section{The elliptic Grothendieck-Teichm\"uller group} \label{sec:2}

In this section, we introduce the group $\on{GT}_{ell}$ of 
universal automorphisms of elliptic structures over BMCs, which we call the 
elliptic Grothendieck-Teichm\"uller group. We compute the ``naive'' version of this 
group, and then introduce its variants (profinite, pro-$l$, proalgebraic) by playing on 
the classes of considered BMCs. We study the relations between these groups and 
the corresponding variants of $\on{GT}$; we construct in particular, in the various 
frameworks, a section of the natural morphism $\on{GT}_{ell}\to \on{GT}$. 
This shows that $\on{GT}_{ell}$ and its variants have semidirect product structures.

\subsection{Reminders about ${\ul{\on{GT}}}$ and its variants}

According to \cite{Dr:Gal}, $\ul{\on{GT}}$ is the set of pairs
$(\lambda,f)\in(1+2\ZZ)\times F_2$, $F_{2}$ being the free group with 
generators $X$ and $Y$, such that 
\begin{equation} \label{hex:f}
f(X_3,X_1)X_3^m f(X_2,X_3)X_2^m f(X_1,X_2)
X_1^m = 1,\quad m = {{\lambda-1}\over 2},  \quad X_1X_2X_3=1, 
\end{equation}
$$
f(Y,X)=f(X,Y)^{-1}, \quad 
\partial_3(f)\partial_1(f)=\partial_0(f)\partial_2(f)\partial_4(f),   
$$
where\footnote{$P_n = \on{Ker}(B_n\to S_n,\sigma_i\mapsto (i,i+1))$
is the pure braid group on $n$ strands.} 
$\partial_i:F_2 \subset P_3\to P_4$ are simplicial 
morphisms. It is equipped with a semigroup structure
with $(\lambda,f)(\lambda',f')=(\lambda'',f'')$, with 
$$\lambda'':=\lambda\lambda', \quad f''(X,Y):=
f(f'(X,Y)X^{\lambda'}f'(X,Y)^{-1},Y^{\lambda'})f'(X,Y).
$$

One defines similarly semigroups $\wh{\ul{\on{GT}}}$, $\ul{\on{GT}}_l$, 
$\ul{\on{GT}}(\kk)$ by replacing in the above definition 
$(\ZZ,F_2)$ by their profinite, pro-$l$, $\kk$-prounipotent 
versions (where $\kk$ is a $\QQ$-ring). We then have morphisms
$\ul{\on{GT}}\hookrightarrow \wh{\ul{\on{GT}}}\to
\ul{\on{GT}}_l\hookrightarrow \ul{\on{GT}}(\QQ_l)$ and 
$\ul{\on{GT}}\to \ul{\on{GT}}(\kk)$ for any $\kk$. 

$\ul{\on{GT}}$ acts on $\{$braided monoidal categories (BMCs)$\}$ 
by $(\lambda,f)(\cC_0,\beta_{\ldots},a_{\ldots}):= (\cC_0,
\beta'_{\ldots},a'_{\ldots})$, where
$$
\beta'_{X,Y} := \beta_{X,Y}(\beta_{Y,X}\beta_{X,Y})^m, \quad  
a'_{X,Y,Z}:= a_{X,Y,Z}f(\beta_{YX}\beta_{XY}\otimes\on{id}_Z,
a_{X,Y,Z}^{-1}(\on{id}_X\otimes\beta_{ZY}\beta_{YZ})a_{X,Y,Z}).
$$ 
Similarly, $\wh{\ul{\on{GT}}}$ (resp., $\ul{\on{GT}}_l$, 
$\ul{\on{GT}}(\kk)$) act on $\{$BMCs $\cC_0$ such that 
$\on{Aut}_{\cC_0}(X)$ is finite for any $X\in\on{Ob}\cC_0\}$
(resp., such that the image of $P_n\to\on{Aut}_{\cC_0}(X_1\otimes
\cdots
\otimes X_n)$ is an $l$-group, is contained in a unipotent 
group).

\subsection{The semigroup ${\ul{\on{GT}}}_{ell}$ and its variants}

Let us define $\ul{\on{GT}}_{ell}$ as the set of all $(\lambda,f,g_\pm)$, where
$(\lambda,f)\in\ul{\on{GT}}$, $g_\pm\in F_2$ are such that 
\begin{equation} \label{def:GTell:1}
(\sigma_2^{\pm 1}\sigma_1^{\pm 1}(\sigma_1\sigma_2^2\sigma_1)^{\pm m}
g_\pm(X_1^+,X_1^-) f(\sigma_1^2,\sigma_2^2))^3=1, 
\end{equation}
\begin{equation} \label{def:GTell:3}
u^2 = (g_-,u^{-1}g_+^{-1}u^{-1}) 
\end{equation}
(identities in $B_{1,3}$) where $u=f(\sigma_1^2,
\sigma_2^2)\sigma_1^\lambda f(\sigma_1^2,\sigma_2^2)^{-1}$, 
$g_\pm=g_\pm(X_1^+,X_1^-)$. 

If $\cC$ is a BMC and $(\tilde\cC,F,A^{\pm}_{\cdots})$ is an elliptic structure 
over $\cC$, then $(\tilde\cC,F,A^{\prime\pm}_{\cdots})$ is an elliptic
structure over $\cC'$, where 
\begin{equation} \label{transfo:cat}
\cC':= (\lambda,f)*\cC, \quad
A^{\prime\pm}_{X,Y} = g_{\pm}(A^{+}_{X,Y},A^{-}_{X,Y})
(\in\on{Aut}_{\cC}(X\otimes Y)). 
\end{equation}
The following statement is then the analogue of equations (\ref{for:quot}). 
\begin{lemma}
The conditions (\ref{def:GTell:1}), (\ref{def:GTell:3}) imply the identities 
\begin{equation} \label{GTell:B12}
(\sigma_1^{\pm\lambda}g_\pm(X_1^+,X_1^-))^2=1, \quad 
\sigma_1^{2\lambda}=(g_-(X_1^+,X_1^-),g_+(X_1^+,X_1^-))
\end{equation}
in $B_{1,2}$. 
\end{lemma}

{\em Proof.} Let $\sigma_{\pm}:= \sigma_{2}^{\pm1}\sigma_{1}^{\pm1}
(\sigma_{1}\sigma_{2}^{2}\sigma_{1})^{\pm m}$, $g_{\pm}:=g_{\pm}(X_{1}^{+},
X_{1}^{-})$, $f:= f(\sigma_{1}^{2},\sigma_{2}^{2})$, then the first equation of 
(\ref{def:GTell:1}) is rewritten as
$\on{Ad}(\sigma_{\pm})^{-1}(g_{\pm}f)\cdot g_{\pm}f \cdot
\on{Ad}(\sigma_{\pm})(g_{\pm}f) = \sigma_{\pm}^{-3}$, an identity in 
$P_{1,3}$. There is a unique morphism $P_{1,3}\to P_{1,2}$, corresponding to 
the erasing of the third point, i.e. to the map $\on{Cf}_{3}(T)\to \on{Cf}_{2}(T)$, 
$(x_{1},x_{2},x_{3})\mapsto (x_{1},x_{2})$. It is given by $X_{1}^{\pm}
\mapsto X_{1}^{\pm}$, $X_{2}^{\pm}\mapsto 1$, $X_{3}^{\pm}\mapsto 
(X_{1}^{\pm})^{-1}$, $\sigma_{1}^{2}\mapsto \sigma_{1}^{2}$, 
$\sigma_{2}^{2}\mapsto 1$, $(\sigma_{1}\sigma_{2})^{3}\mapsto \sigma_{1}^{2}$. 
The image of the above identity in $P_{1,3}$ by this morphism is the identity 
$g_{\pm}(X_{1}^{+},X_{1}^{-})\cdot \on{Ad}(\sigma_{1}^{\pm\lambda})
(g_{+}(X_{1}^{+},X_{1}^{-})) = \sigma_{1}^{\mp 2\lambda}$ in $P_{1,2}$, 
which is equivalent to the first equations of (\ref{GTell:B12}). The same morphism
similarly takes (\ref{def:GTell:3}) to the last equation of (\ref{GTell:B12}). 
\hfill \qed\medskip 

For $(\lambda,f,g_{\pm}),(\lambda',f',g'_{\pm})\in \ul{\on{GT}}_{ell}$, 
we set 
$$(\lambda,f,g_\pm)
(\lambda',f',g'_\pm):=(\lambda'',f'',g''_\pm), \quad \text{where} \quad 
g''_\pm(X,Y) = g_\pm(g'_+(X,Y),g'_-(X,Y)).
$$ 

\begin{proposition} \label{prop:GTell:sg}
This defines a semigroup structure on $\ul{\on{GT}}_{ell}$. 
We have a semigroup inclusion $\ul{\on{GT}}_{ell}\subset 
\ul{\on{GT}}\times\on{End}(F_2)^{op}$, $(\lambda,f,g_\pm)
\mapsto ((\lambda,f),\theta_{g_\pm})$, where 
$\theta_{g_\pm} = (X\mapsto g_+(X,Y), Y\mapsto g_-(X,Y))$. 
\end{proposition}

{\em Proof.} We first prove: 

\begin{lemma} \label{lemma:gp:pties}
If $(\lambda,f,g_\pm)\in\ul{\on{GT}}_{ell}$, then there is a unique 
endomorphism of $B_{1,3}$, such that 
$$
\sigma_1\mapsto \tilde\sigma_1 := f(\sigma_1^2,\sigma_2^2)\sigma_1^\lambda
f(\sigma_1^2,\sigma_2^2)^{-1}, \quad \sigma_2\mapsto 
\tilde\sigma_2 := \sigma_2^\lambda, \quad X_1^\pm\mapsto
g_\pm(X_1^+,X_1^-). 
$$
For any $\lambda'\in 2\ZZ+1$, we then have 
\begin{equation} \label{*}
f(\sigma_1^2,\sigma_2^2)\sigma_2^{\pm 1}\sigma_1^{\pm 1}
(\sigma_1\sigma_2^2\sigma_1)^{\pm{{\lambda\lambda'-1}\over 2}}
= \tilde\sigma_2^{\pm 1}\tilde\sigma_1^{\pm 1}(\tilde\sigma_1\tilde\sigma_2^2
\tilde\sigma_1)^{\pm{{\lambda'-1}\over 2}}.  
\end{equation}
\end{lemma}

{\em Proof.} Recall that we have an elliptic structure 
$({\bf PaB},{\bf PaB}_{ell},F,A^\pm_{\ldots})$. Applying 
$(\lambda,f,g_\pm)$ to it, we get an elliptic structure
$(\ul{{\bf PaB}},\ul{{\bf PaB}}^{ell},F,\ul A^\pm_{\ldots})$.
An endomorphism of $B_{1,3}$ is given by the composition 
$$
B_{1,3}\to \on{Aut}_{\ul{{\bf PaB}}^{ell}}
(\bullet(\bullet\bullet))\simeq B_{1,3}, 
$$ 
where the first morphism arises from the elliptic structure of 
$\ul{{\bf PaB}}^{ell}$, and the second morphism arises from 
the isofunctor $\ul{{\bf PaB}}^{ell} \simeq {\bf PaB}^{ell}$. One checks 
that this endomorphism of $B_{1,3}$ is given by the above formulas. 

We now prove (\ref{*}).  The hexagon identity implies
$$(\sigma_2^2)^m f(\sigma_1^2,\sigma_2^2) (\sigma_1^2)^m
f((\sigma_1^2\sigma_2^2)^{-1},\sigma_1^2)(\sigma_1^2\sigma_2^2)^{-m}
f(\sigma_2^2,(\sigma_1^2\sigma_2^2)^{-1})=1.
$$ 
Now since $(\sigma_1^2\sigma_2^2)^{-1}\equiv \sigma_1\sigma_2^2
\sigma_1^{-1} \equiv \sigma_2^{-1}\sigma_1^2\sigma_2$ mod $Z(B_3)$, 
since $f(a,b)=f(a',b')$ for any group $G$ and any $a,a',b,b'\in G$
with $a\equiv a'$, $b\equiv b'$ mod $Z(G)$ (as $f\in F_{2}' = (F_{2},F_{2})$), 
and by the duality identity, 
this is rewritten 
$$
(\sigma_2^2)^m f(\sigma_1^2,\sigma_2^2)
\sigma_1^{2m+1} f(\sigma_1^2,\sigma_2^2)^{-1} \sigma_1^{-1}
(\sigma_1^2\sigma_2^2)^{-m}\sigma_2^{-1} 
f(\sigma_1^2,\sigma_2^2)^{-1} \sigma_2=1,
$$ 
which yields (\ref{*}) with $(\pm,\lambda')=(+,1)$. 

(\ref{*}) with $\pm = +$ then follows from 
\begin{equation} \label{transfo:sigmas}
\tilde\sigma_1\tilde\sigma_2^2\tilde\sigma_1 = 
(\sigma_1\sigma_2^2\sigma_1)^\lambda, 
\end{equation}
which is proved as follows. The hexagon identity (\ref{hex:f})
implies that if $X_1X_2X_3$ commutes with all the $X_i$, then 
$$
f(X_3,X_1) X_3^m f(X_2,X_3) X_2^m f(X_1,X_2) = (X_2X_3)^m. 
$$
Applying this to $X_1 = \sigma_2^2$, $X_2 = \sigma_1
\sigma_2^2\sigma_1^{-1}$, $X_3 = \sigma_1^2$, and using 
$\sigma_{1}^{2} = \on{Ad}(\sigma_{2}^{-1}\sigma_{1}^{-1})(\sigma_{2}^{2})$,  
this implies 
\begin{align*}
& \tilde\sigma_1\tilde\sigma_2 
f(\sigma_{1}^{2},\sigma_2^2)
= f(\sigma_1^2,\sigma_2^2)(\sigma_1^2)^m f(\sigma_1\sigma_2^2\sigma_1^{-1},\sigma_2^2)
(\sigma_1\sigma_2^2\sigma_1^{-1})^m f(\sigma_2^2,\sigma_1\sigma_2^2\sigma_1^{-1})\sigma_1\sigma_2
\\ & = (\sigma_1\sigma_2^2\sigma_1)^m\sigma_1\sigma_2. 
\end{align*}
Using the same identity with $X_{1} = \sigma_{2}^{2}$, $X_{2} = \sigma_{1}^{2}$, 
$X_{3} = \sigma_{1}^{-1}\sigma_{2}^{2}\sigma_{1}$, one proves similarly that 
$$
f(\sigma_2^2,\sigma_{1}^{2})\tilde\sigma_2\tilde\sigma_1
= \sigma_2\sigma_1(\sigma_1\sigma_2^2\sigma_1)^m. 
$$
The product of these identities yields (\ref{transfo:sigmas}). 

Each side of (\ref{*}) with
$\pm=-$ identifies with the same side of (\ref{*}) with $\pm=+$
and $\lambda'$ replaced by $-\lambda'$. This implies (\ref{*}) 
with $\pm=-$. 
\hfill \qed\medskip 

{\em End of proof of Proposition \ref{prop:GTell:sg}.}
It suffices to prove that $(\lambda'',f'',g''_\pm)\in \ul{\on{GT}}_{ell}$, 
i.e., that it satisfies conditions (\ref{def:GTell:1}) and (\ref{def:GTell:3}). 

Condition  (\ref{def:GTell:1}) is expressed as follows
$$
\big(\sigma_{2}^{\pm1}\sigma_{1}^{\pm1}
(\sigma_{1}\sigma_{2}^{2}\sigma_{1})^{\pm {{\lambda\lambda'-1}\over2}}
g_{\pm}(g'_{+}(X_{1}^{+},X_{1}^{-}),g'_{-}(X_{1}^{+},X_{1}^{-}))
f(\on{Ad}(f'(\sigma_{1}^{2},\sigma_{2}^{2}))
(\sigma_{1}^{2\lambda'}),\sigma_{2}^{2\lambda'})f'(\sigma_{1}^{2},
\sigma_{2}^{2})\big)^{3}=1, $$
i.e., according to (\ref{*}), as follows
$$
\big(g_{\pm}(g'_{+}(X_{1}^{+},X_{1}^{-}),
g'_{-}(X_{1}^{+},X_{1}^{-})) f(\tilde\sigma_{1}^{\prime 2},
\tilde\sigma_{1}^{\prime 2}) \tilde\sigma_{2}^{\prime\pm1}
\tilde\sigma_{1}^{\prime\pm1}(\tilde\sigma_{1}^{\prime}
\tilde\sigma_{2}^{\prime 2}\tilde\sigma_{1}^{\prime})^{\pm m}
\big)^{3}=1, 
$$
where $\tilde\sigma'_{1}, \tilde\sigma'_{2}$ are the analogue of 
$\tilde\sigma_{1},\tilde\sigma_{2}$ from Lemma \ref{lemma:gp:pties}
with $(\lambda',f')$ instead of $(\lambda,f)$. The latter identity is the image of identity 
(\ref{def:GTell:1}) satisfied by $(\lambda,f,g_{\pm})$ by the endomorphism
of $B_{1,3}$ attached to $(\lambda',f',g'_{\pm})$ by Lemma \ref{lemma:gp:pties}. 

Condition (\ref{def:GTell:3}) is the image of identity (\ref{def:GTell:3}) satisfied by 
$(\lambda,f,g_{\pm})$ under the endomorphism of $B_{1,3}$ attached to  
$(\lambda',f',g'_{\pm})$ by Lemma \ref{lemma:gp:pties}. 
\hfill \qed\medskip 

The operation $(\lambda,f,g_\pm)(\cC,\tilde\cC,F,A^\pm_{\ldots}):= 
(\cC',\tilde\cC,F,A^{\prime\pm}_{\ldots})$, where $\cC',
A^{\prime\pm}_{\ldots}$ are as in (\ref{transfo:cat}) defines an action of 
$\ul{\on{GT}}_{ell}$ on $\{(\cC,\tilde\cC,F,A^\pm_{\ldots}) | \cC$ is a BMC, 
$(\tilde\cC,F,A^\pm_{\ldots})$ is an elliptic structure over it$\}$. 

As before, we define semigroups $\underline{\wh{\on{GT}}}_{ell}$, 
$\underline{\on{GT}}^{ell}_l$, $\underline{\on{GT}}(\kk)$ 
by replacing in the definition of $\underline{\on{GT}}_{ell}$, 
$(\underline{\on{GT}},F_2)$ by\footnote{For $G$ a group (other than 
$\on{GT}$, $\on{GT}_{ell}$ or $R_{ell}$), $\widehat G$ is its profinite 
completion. If $G$ is a free or pure (elliptic) braid group, $G_l$, $G(\kk)$
are its pro-$l$, $\kk$-prounipotent completions. Here $G(-)$ is the
prounipotent $\QQ$-group scheme associated to $G$; it is
characterized by $\on{Hom}_{groups}(G,U(\QQ)) \simeq 
\on{Hom}_{gp\ schemes}(G(-),U)$ for any unipotent group scheme $U$. 
If $G = B_n$ or $B_{1,n}$, then $G_l:= P_l *_P G$, 
$G(\kk):= P(\kk)*_P G$, where $P = \on{Ker}(G\to S_n)$ and $*_{P}$
denotes the amalgamated product over the group $P$.} 
$(\underline{\widehat{\on{GT}}},
\widehat F_2)$, $(\underline{\on{GT}}_l,(F_2)_l)$, 
$(\underline{\on{GT}}(\kk), F_2(\kk))$. They act on 
the sets of pairs $(\cC,\tilde\cC)$, such that $\cC$ satisfies
the same conditions as above, together with: 
$\on{Aut}_{\tilde\cC}(F(X))$ is finite for any 
$X\in \on{Ob}\cC$ (resp., the image of 
$P_{1,n}\to \on{Aut}_{\tilde\cC}(F(X_1\otimes\cdots\otimes X_n))$
is an $l$-group, is contained in a unipotent group). We have morphisms
$\underline{\on{GT}}_{ell} \hookrightarrow
\underline{\widehat{\on{GT}}}_{ell} \to 
\underline{\on{GT}}^{ell}_l \hookrightarrow 
\underline{\on{GT}}_{ell}(\QQ_l)$ and $\ul{\on{GT}}_{ell}
\to \ul{\on{GT}}_{ell}(\kk)$
compatible with the similar `non-elliptic' morphisms. 

\subsection{Computation of ${\ul{\on{GT}}}_{ell}$}

Recall that the braid group $B_{3}$ is presented by generators $\Psi_{\pm}$
and relations $\Psi_{+}\Psi_{-}\Psi_{+} = \Psi_{-}\Psi_{+}\Psi_{-}$
($\Psi_{\pm}$ are the $\sigma_{1},\sigma_{2}$ of the standard presentation
and are used in order to avoid confusion with previous notation). Its center
$Z(B_{3})$ is isomorphic to $\ZZ$ and generated by $(\Psi_{+}\Psi_{-})^{3}$. 
There is a central exact sequence 
$$
1\to 2Z(B_{3})\to B_{3}\to \on{SL}_{2}(\ZZ)\to 1, 
$$
given by $\Psi_{+}\mapsto \bigl( \begin{smallmatrix} 
1 & 1 \\ 0 & 1 \end{smallmatrix}\bigr)$, 
$\Psi_{-}\mapsto \bigl( \begin{smallmatrix} 
1 & 0 \\ -1 & 1 \end{smallmatrix}\bigr)$. 

\begin{proposition}
Let $\tilde B_{3}$ be the group generated by $\Psi_{\pm},\varepsilon$
and relations 
$$
\Psi_{+}\Psi_{-}\Psi_{+} = \Psi_{-}\Psi_{+}\Psi_{-}, \quad 
\varepsilon\Psi_{+}\varepsilon\Psi_{-}=1, \quad \varepsilon^{2}=1. 
$$
There is an exact sequence $1\to B_{3}\to\tilde B_{3}\to\ZZ/2\to 1$, 
where $\tilde B_{3}\to\ZZ/2$ is given by $\Psi_{\pm}\mapsto 1$, $\varepsilon
\to -1$. There is also a (non-central) exact sequence $1\to 2Z(B_{3})\to\tilde B_{3}
\to \on{GL}_{2}(\ZZ)\to 1$, where $\tilde B_{3}\to\on{GL}_{2}(\ZZ)$ extends
$B_{3}\to\on{SL}_{2}(\ZZ)$ by $\varepsilon\mapsto 
\mapsto \bigl( \begin{smallmatrix} 0 & 1 \\ 1 & 0 
\end{smallmatrix}\bigr)$. All these morphisms fit in the diagram  
$$
\xymatrix{
 &  & 1\ar[d] & 1\ar[d]  & \\
 1\ar[r] & 2Z(B_{3})\ar[r]\ar@{=}[d] & B_{3}\ar[r]\ar[d] 
 & \on{SL}_{2}(\ZZ)\ar[r]\ar[d] & 1 \\
 1\ar[r] & 2Z(B_{3})\ar[r] & \tilde B_{3}\ar[r]\ar[d] & 
 \on{GL}_{2}(\ZZ)\ar[r]\ar[d] & 1 \\
 &  & \ZZ/2\ar@{=}[r]\ar[d]  & \ZZ/2\ar[d] &  \\
 & & 1 & 1 & 
}$$
\end{proposition}
The proof is straightforward.

\begin{proposition} \label{prop:sg:morph}
1) There is a unique semigroup morphism 
$\tilde B_{3}\to \underline{\on{GT}}_{ell}$, 
such that: 
$$
\Psi_{+}\mapsto (\lambda,f,g_+,g_-)=(1,1,g_+(X,Y)=X,g_-(X,Y)=YX), 
$$
$$
\Psi_{-}\mapsto (\lambda,f,g_{+},g_{-}) = (1,1,g_{+}(X,Y) = XY^{-1}, 
g_{-}(X,Y) = Y), $$
$$
\varepsilon\mapsto(\lambda,f,g_+,g_-)=(-1,1,g_+(X,Y)=Y,g_-(X,Y)=X),
$$
It fits in a commutative diagram 
$$
\begin{matrix}
\tilde B_{3} 
& \to&\underline{\on{GT}}_{ell} \\
\downarrow  & & \downarrow\\
\ZZ/2 &\to & \underline{\on{GT}} 
\end{matrix}
$$

2) The horizontal maps in this diagram are isomorphisms. 
\end{proposition}

{\em Proof.} Set $X_{i}:= X_{i}^{+}$, $Y_{i}:= X_{i}^{-}$. 
Using the commutation of $\sigma_{2}$ with $X_{1}$ and the braid relation 
between $\sigma_{1}$ and $\sigma_{2}$, one obtains 
$(\sigma_{2}\sigma_{1}X_{1})^{3} = X_{3}X_{2}X_{1}=1$ (relation in 
$B_{1,3}$). In the same way, $(\sigma_{2}^{-1}\sigma_{1}^{-1}Y_{1}X_{1})^{3}$
expresses as an element of $P_{1,3}$ as 
$$
Y_{3}X_{3}\sigma_{2}^{-1}\sigma_{1}^{-2}\sigma_{2}^{-1}
Y_{2}X_{2}\sigma_{1}^{-2}Y_{1}X_{1}. 
$$
Since $(Y_{1},X_{2}^{-1}) = \sigma_{1}^{2}$, $X_{2}\sigma_{1}^{-2}Y_{1}$
can be replaced by $Y_{1}X_{2}$; in the resulting expression, $Y_{2}Y_{1}X_{2}X_{1}$
can then be replaced by $Y_{3}^{-1}X_{3}^{-1}$. The above expression is therefore equal
to 
$$
Y_{3}X_{3}\sigma_{2}^{-1}\sigma_{1}^{-2}\sigma_{2}^{-1}Y_{3}^{-1}X_{3}^{-1}.
$$
One has $(Y_{3}^{-1},X_{3}^{-1}) = (Y_{2}Y_{1},X_{3}^{-1}) = 
Y_{1}(Y_{2},X_{3}^{-1})Y_{1}^{-1}(Y_{1},X_{3}^{-1})$; one computes 
$(Y_{2},X_{3}^{-1}) = \sigma_{2}^{2}$, $(Y_{1},X_{3}^{-1}) = 
\sigma_{2}^{-1}\sigma_{1}^{2}\sigma_{2}$, which implies 
that $(Y_{3}^{-1},X_{3}^{-1}) = \sigma_{2}\sigma_{1}^{2}\sigma_{2}$ and 
therefore that 
$$(\sigma_{2}^{-1}\sigma_{1}^{-1}Y_{1}X_{1})^{3}=1
\quad \text{(equality in }B_{1,3}).
$$ 
Finally, $(Y_{1}X_{1},\sigma_{1}^{-1}X_{1}^{-1}\sigma_{1}^{-1}) = 
(Y_{1}X_{1},X_{2}^{-1}) = (Y_{1},X_{2}^{-1}) = \sigma_{1}^{2}$ (equality
in $B_{1,3}$), where the second equality uses the commutation of $X_{1}$ and 
$X_{2}$. All this implies that $(1,1,X,YX)\in\underline{\on{GT}}_{ell}$. 

If $(\lambda,f,g_{+},g_{-}) = (-1,1,Y,X)$, then $m=-1$, therefore
$$(\sigma_{2}^{\pm1}\sigma_{1}^{\pm1}(\sigma_{1}\sigma_{2}^{2}
\sigma_{1})^{\pm m}g_{\pm}(X_{1}^{+},X_{1}^{-})f(\sigma_{1}^{2},
\sigma_{2}^{2}))^{3} = (\sigma_{2}^{\mp1}\sigma_{1}^{\mp1}X_{1}^{\mp})^{3}
=1 \quad \text{(relation in }B_{1,3}). 
$$
The relation $u^{2} = (g_{-},u^{-1}g_{+}^{-1}u^{-1})$ follows from $\sigma_{1}^{-2}
=(X_{1},Y_{2}^{-1})$ (relation in $B_{1,3}$). All this implies that 
$(-1,1,Y,X)\in\underline{\on{GT}}_{ell}$. 

One checks that $(1,1,XY^{-1},Y) = (-1,1,Y,X)(1,1,X,YX)^{-1}(-1,1,Y,X)$, 
therefore 
$$
(1,1,XY^{-1},Y)\in \underline{\on{GT}}_{ell}.
$$ 
Finally, one checks that
the relations between $\Psi_+,\Psi_{-}$ and $\varepsilon$ are also satisfied 
by their images in $\underline{\on{GT}}_{ell}$.  All this proves 1). 

Let us prove 2). 
The bijectivity of $\ZZ/2 \to \underline{\on{GT}}$
is proved in \cite{Dr:Gal}, Proposition 4.1. Set $\underline{R}_{ell}:= 
\on{Ker}(\underline{\on{GT}}_{ell}\to \underline{\on{GT}})$, 
then the commutativity of the above diagram implies that its 
upper map restricts to a morphism $B_{3}\to
\underline{R}_{ell}$, and we need to prove that it is bijective. 
According to the second identity in (\ref{GTell:B12}), 
$\underline{R}_{ell}\subset \{(g_+,g_-)\in(F_2)^2| 
(g_-(X,Y),g_+(X,Y))=(Y,X)\}$. We now recall some results 
due to Nielsen. 

\begin{theorem} \label{thm:nielsen} (\cite{N})
1) The morphism $\on{Out}(F_2)\to\on{GL}_2(\ZZ)$ 
induced by abelianization is an isomorphism. 

2) $\on{Im}(\on{Aut}(F_2)\to (F_2)^2) = \{(g_+,g_-)\in(F_2)^2|
\exists k\in F_2, \exists \epsilon\in\{\pm 1\}, 
(g_-(X,Y),g_+(X,Y)) = k(Y,X)^\epsilon k^{-1}\}$, where the 
map $\on{Aut}(F_2)\to (F_2)^2$ is $\theta\mapsto(\theta(X),\theta(Y))$. 
\end{theorem} 

The bijectivity of $B_{3} \to \underline{R}_{ell}$, 
together with the equality $\underline{R}_{ell} = \{(g_+,g_-)\in (F_2)^2|
(g_-(X,Y),$ $g_+(X,Y))=(Y,X)\}$ are then proven in the following
corollary to Theorem \ref{thm:nielsen}: 

\begin{corollary} \label{cor:nielsen}
We have bijections 
$$
B_{3}\to \on{Aut}_{(X,Y)}(F_2)
\to \{(g_+,g_-)\in(F_2)^2|(g_-(X,Y),g_+(X,Y)) = (Y,X)\},  
$$
where $\on{Aut}_{(X,Y)}(F_2) = \{\theta\in\on{Aut}(F_2)|\theta((X,Y))
=(X,Y)\}$, the first map is as in Proposition 
\ref{prop:sg:morph} and the second map 
is $\theta\mapsto(\theta(X),\theta(Y))$. 
\end{corollary}

{\em Proof of Corollary \ref{cor:nielsen}.} The bijectivity 
of the second map follows from the injectivity of $\on{Aut}(F_2)
\to (F_2)^2$, $\theta\mapsto (\theta(X),\theta(Y))$ and from 
Theorem \ref{thm:nielsen}, 2). Let us now prove the bijectivity of 
the map $B_{3}\to\on{Aut}_{(X,Y)}(F_2)$. 
The kernel of $B_{3}\to\on{Aut}_{(X,Y)}(F_2)$
is contained in $\on{Ker}(B_{3}\to\on{Aut}_{(X,Y)}(F_2)
\to\on{Out}(F_2)\to\on{GL}_2(\ZZ)) = \langle(\Psi_{+}\Psi_{-})^{6}\rangle$. 
On the other hand, $B_{3}\to\on{Aut}_{(X,Y)}(F_2)$ takes 
$(\Psi_{+}\Psi_{-})^{6}$ to\footnote{Here and later, $\on{Ad}(g)$ 
is the inner automorphism $x\mapsto gxg^{-1}$.} 
$\on{Ad}((X,Y)^{-1})$, so the restriction of 
$B_{3}\to\on{Aut}_{(X,Y)}(F_2)$ to $\langle(\Psi_{+}\Psi_{-})^{6}\rangle$
is injective. It follows that $B_{3}\to\on{Aut}_{(X,Y)}(F_2)$ is injective. 

Let us now show that $B_{3}\to\on{Aut}_{(X,Y)}(F_2)$
is surjective. We have a commutative diagram 
$$
\xymatrix{
\on{Aut}(F_{2})\ar[r] & \on{Out}(F_{2})\ar[r]^{\sim} & \on{GL}_{2}(\ZZ)\\
\on{Aut}_{(X,Y)}(F_{2})\ar[rr]\ar@{^{(}->}[u] & & \on{SL}_{2}(\ZZ) 
\ar@{^{(}->}[u]}
$$
where the isomorphism follows from Theorem \ref{thm:nielsen}, 1), and the bottom 
map is given by abelianization. 
It follows that $\on{Ker}(\on{Aut}_{(X,Y)}(F_2)\to 
\on{SL}_2(\ZZ)) = \on{Ker}(\on{Aut}_{(X,Y)}(F_2)\to\on{Out}(F_2))
= \on{Aut}_{(X,Y)}(F_2)\cap \on{Inn}(F_2) = \{\theta\in 
\on{Aut}(F_2)|\exists k\in F_2, \theta = \on{Ad}(k)$ and 
$k$ commutes with $(X,Y)\}$. The subgroup of $F_{2}$ generated by 
$k$ and $(X,Y)$ is abelian and, according to \cite{N:danish:1921}, free, 
and therefore isomorphic to $\ZZ$. If $(X,Y)$ is a power of an element $h$ of $F_{2}$, 
then the sum of the degrees of $h$ in $X$ and in $Y$ is zero, and comparing coefficients in 
$[\on{log}X,\on{log}Y]$ in $\on{log}(X,Y)$ and $\on{log}h$ in 
the Lie algebra of the prounipotent completion of $F_{2}$, one sees that $h$ is 
$(X,Y)$ or its inverse, therefore $(X,Y)$ is not the power of an element of $F_{2}$
other that itself or its inverse. All this implies that $k$ should be a power of $(X,Y)$, 
therefore 
$\on{Ker}(\on{Aut}_{(X,Y)}(F_2)\to \on{SL}_2(\ZZ)) 
= \langle \on{Ad}(X,Y)\rangle = \langle (\Psi_{+}\Psi_{-})^{6}
\rangle$. On the other hand, as the composition $B_{3}\to 
\on{Aut}_{(X,Y)}(F_{2})\to \on{SL}_{2}(\ZZ)$ is surjective, so is the morphism 
$\on{Aut}_{(X,Y)}(F_{2})\to \on{SL}_{2}(\ZZ)$. All this implies that there is 
an exact sequence 
$$
1\to \langle(\Psi_{+}\Psi_{-})^{6}\rangle \to \on{Aut}_{(X,Y)}(F_{2})\to 
\on{SL}_{2}(\ZZ)\to 1. 
$$
Let us denote this exact sequence as $1\to K\to G\to H\to 1$, and let 
$G':= \on{Im}(B_{3}\to G) \subset G$.
To prove that $G' = G$, it suffices to prove that $\on{Im}(G'\subset G\to H) = H$
and that $G'\supset K$. The first statement follows from the surjectivity to 
$B_{3}\to\on{Aut}_{(X,Y)}(F_{2})\to\on{SL}_{2}(\ZZ)$, 
while the second statement follows from the fact that $(\Psi_{+}\Psi_{-})^{6}
\in \on{Im}(B_{3}\to \on{Aut}_{(X,Y)}(F_{2}))$. \hfill \qed \qed\medskip 

\begin{remark}
As $\on{GL}_2(\ZZ)$ is the nonoriented mapping class group of 
the topological torus, we have a morphism $\on{GL}_{2}(\ZZ)\to 
\on{Out}(B_{1,n})$, obtained by applying mapping class group elements to 
elliptic braids; its target is an outer automorphism group because 
the mapping class group does not preserve a base point of the elliptic configuration 
space. This morphism lifts to a morphism 
\begin{equation} \label{act:B3:ellbraids}
\tilde B_{3}\to \on{Aut}(B_{1,n}),
\end{equation} 
given by $\Psi_{+}\mapsto (X_{1}\mapsto X_{1}, Y_{1}\mapsto Y_{1}X_{1}, 
\sigma_{i}\mapsto\sigma_{i})$, $\Psi_{-}\mapsto (X_{1}\mapsto 
X_{1}Y_{1}^{-1},\sigma_{i}\mapsto\sigma_{i})$, $\varepsilon\mapsto 
(X_{1}\leftrightarrow Y_{1}, \sigma_{i}\mapsto \sigma_{i}^{-1})$. It is 
such that $(\Psi_{+}\Psi_{-})^6\mapsto ($conjugation by the image of $z\in P_n
\to B_{1,n})$, where $z$ is a generator of $Z(P_n) \simeq \ZZ$. 
The assignment $\{$elliptic structures over BMCs$\}\to\{$representations of
$B_{1,n}\}$ is then $\tilde B_{3}$-equivariant. 
\end{remark}

\begin{remark} 
The morphisms $\tilde B_{3}\to \on{Aut}(B_{1,n})$
and $\ul{\on{GT}}_{ell}\to \on{End}(B_{1,3})$ from 
Lemma \ref{lemma:gp:pties} admit a common 
generalization to a morphism $\ul{\on{GT}}_{ell}\to \on{End}(B_{1,n})$, 
taking $(\lambda,f,g_{+},g_{-})$ to the endomorphism 
$X_{1}\mapsto g_{+}(X_{1},Y_{1})$, $Y_{1}\mapsto 
g_{-}(X_{1},Y_{1})$, $\sigma_{i}\mapsto 
\on{Ad}(f(\sigma_{i}^{2},\sigma_{i+1}\cdots\sigma_{n-1}^{2}\cdots
\sigma_{i+1}))(\sigma_{i}^{\lambda})$; this corresponds to the identification 
of $B_{1,n}$ with $\on{Aut}_{{\bf PaB}_{ell}}(\bullet(\bullet\cdots
(\bullet\bullet)))$. This morphism extends to the various setups (profinite, etc.). 
\end{remark}

\subsection{The semigroup scheme ${\ul{\on{GT}}}_{ell}(-)$}

For $\kk$ a $\QQ$-ring, we set\footnote{The kernel of a morphism of semigroups
with unit is the preimage of the unit of the target semigroup; it is again a semigroup 
with unit.} $R_{ell}(\kk):= \on{Ker}(
\ul{\on{GT}}_{ell}(\kk)\to \ul{\on{GT}}(\kk))$. 
The assignments $\kk\mapsto \ul{\on{GT}}_{(ell)}(\kk), 
R_{ell}(\kk)$ are functors $\{\QQ$-rings$\}\to\{$semigroups$\}$, 
i.e., semigroup schemes over $\QQ$. 

\begin{proposition} \label{prop15}
We have a commutative diagram of morphisms of
semigroup schemes
$$
\begin{matrix}
R_{ell}(-) & \to & \ul{\on{GT}}_{ell}(-) &\to& \ul{\on{GT}}(-) \\
\downarrow & & \downarrow & & \downarrow \\
\on{SL}_2(-) &\to&  \on{M}_2(-) & \stackrel{\on{det}}\to & {\mathbb A}^{1}(-)
\end{matrix}
$$
where $\ul{\on{GT}}(\kk)\to\kk$ is $(\lambda,f)\mapsto \lambda$
and $\ul{\on{GT}}_{ell}(\kk)\to \on{M}_2(\kk)$ is $(\lambda,f,g_\pm)
\mapsto \bigl( \begin{smallmatrix} \alpha_+ & \beta_+ \\ \alpha_- & 
\beta_- \end{smallmatrix}\bigr)$, where $\on{log}g_\pm(X,Y)=
\alpha_\pm \on{log} X+\beta_\pm\on{log}Y$ mod\footnote{Recall that 
$F_2(\kk)=\on{exp}(\hat\f_2^\kk)$, where $\hat\f_2^\kk$ is the topologically
free $\kk$-Lie algebra in two generators $\on{log}X$ and $\on{log}Y$.} 
$[\hat\f_2^\kk,\hat\f_2^\kk]$. 
\end{proposition}

{\em Proof.} It suffices to show that the right square is commutative, 
which follows by abelianization from the second part of (\ref{GTell:B12}). 
\hfill \qed\medskip 

Recall that\footnote{If $S$ is a semigroup with unit, $S^\times$
is the group of its invertible elements.} $\on{GT}(\kk) = 
\ul{\on{GT}}(\kk)^\times$. We set 

\begin{definition}
$\on{GT}_{ell}(\kk):= \ul{\on{GT}}_{ell}(\kk)^\times$. 
\end{definition}

\begin{proposition}
1) $\on{GT}_{ell}(\kk) = \ul{\on{GT}}_{ell}(\kk)\times_{\on{M}_2(\kk)}
\on{GL}_2(\kk)$ (Cartesian product in the category of proalgebraic varieties). 

2) $R_{ell}(\kk)$ is a group. 
\end{proposition}

{\em Proof.} Let $(\lambda,f,g_\pm)\in\ul{\on{GT}}_{ell}(\kk)$
be invertible as an element of $\ul{\on{GT}}(\kk)\times
\on{End}(F_2(\kk))^{op}$, with inverse $(\lambda',f',g'_{\pm})$. 
Then the endomorphism of Lemma \ref{lemma:gp:pties} attached to 
$(\lambda,f,g_{\pm})$ is an automorphism of $B_{1,3}(\kk)$. 
The identities $(\sigma_{2}^{\pm1}\sigma_{1}^{\pm1}X_{1}^{\pm})^{3}=1$,
$\sigma_{1}^{2}=(X_{1}^{-},(X_{2}^{+})^{-1})$ in $B_{1,3}(\kk)$ are the 
images by this automorphism of the identities expressing that $(\lambda',f',g'_{\pm})$
belongs to $\underline{\on{GT}}_{ell}(\kk)$. It follows that $(\lambda',f',g'_{\pm})
\in\underline{\on{GT}}_{ell}(\kk)$. The element $(\lambda,f,g_{\pm})$ 
is invertible iff the image of $(\lambda,f,g_\pm)\in\ul{\on{GT}}_{ell}(\kk)\to 
\on{M}_2(\kk)$ lies in $\on{GL}_2(\kk)$. All this proves 1). 2) is then immediate. 
\hfill \qed\medskip 

Recall that for any $\QQ$-ring $\kk$, $\on{GT}_1(\kk)= 
\on{Ker}(\ul{\on{GT}}(\kk)\to\kk)$. We also set 
$$
\on{GT}_{I_2}^{ell}(\kk):= \on{Ker}(\ul{\on{GT}}_{ell}(\kk)
\to \on{M}_2(\kk)), \quad 
R_{I_2}^{ell}(\kk) := \on{Ker}(R_{ell}(\kk)\to\on{SL}_2(\kk)). 
$$
Then $\kk\mapsto \on{GT}_{I_2}^{(ell)}(\kk)$, $R_{I_2}^{ell}(\kk)$
are $\QQ$-group schemes. It is known that $\on{GT}_1(-)$ is prounipotent.

\begin{proposition}
The group schemes $\on{GT}_{I_2}^{ell}(-)$ and $R_{I_2}^{ell}(-)$
are prounipotent. 
\end{proposition}

{\em Proof.} $\on{GT}_{I_2}^{ell}(\kk)\subset \on{GT}_1(\kk)\times 
\on{Aut}_{I_2}(F_2(\kk))^{op}$, where $\on{Aut}_{I_2}(F_2(\kk)) = \on{Ker}
(\on{Aut}(F_2(\kk))\to \on{GL}_2(\kk))$; $\kk\mapsto \on{Aut}_{I_2}(F_2(\kk))$
is prounipotent, so $\kk\mapsto \on{GT}_{I_{2}}^{ell}(\kk)$ is prounipotent as the 
subgroup of a prounipotent group scheme. The same argument implies that 
$R_{I_2}^{ell}(-)$ is prounipotent. \hfill \qed\medskip 

\begin{proposition} We have exact sequences 
$1\to R_{I_2}^{ell}(\kk)\to R_{ell}(\kk)\to \on{SL}_2(\kk)\to 1$
and $1\to \on{GT}_{I_2}^{ell}(\kk)\to \on{GT}_{ell}(\kk)\to 
\on{GL}_2(\kk)\to 1$. 
\end{proposition}

{\em Proof.} We need to prove that $R_{ell}(\kk)\to\on{SL}_2(\kk)$
is surjective. Set $G(\kk):= \on{Im}(R_{ell}(\kk)\to\on{SL}_2(\kk))$, 
then $\kk\mapsto G(\kk)$ is a group subscheme of $\on{SL}_2$. 
We have two morphisms ${\mathbb G}_a\to R_{ell}(-)$, extending  
$\ZZ\to B_{3}$, $1\mapsto \Psi_\pm$
in the sense that 
$$
\begin{matrix}
B_{3} 
&\to & R_{ell}(\kk)\\
\uparrow & & \uparrow \\
\ZZ &\to & {\mathbb G}_{a}(\kk) 
\end{matrix}
$$
commutes; then ${\mathbb G}_a\to R_{ell}\to \on{SL}_2$
are the morphisms $t\mapsto \bigl( \begin{smallmatrix} 1 & t \\ 0
 & 1 \end{smallmatrix}\bigr), \bigl( \begin{smallmatrix} 1 & 0 \\ -t
 & 1 \end{smallmatrix}\bigr)$. So the Lie algebra of $G(-)$
contains both $\bigl( \begin{smallmatrix} 0 & 1 \\ 0
 & 0 \end{smallmatrix}\bigr)$ and $\bigl( \begin{smallmatrix} 0 & 0 \\ 1
 & 0 \end{smallmatrix}\bigr)$, hence is equal to $\SL_2$, so $G=\on{SL}_2$. 

Let us now prove that $\on{GT}_{ell}(\kk)\to\on{GL}_2(\kk)$ 
is surjective. Set $\tilde G(\kk):= \on{Im}(\on{GT}_{ell}(\kk)\to 
\on{GL}_2(\kk))$, then $\on{SL}_2 \subset \tilde G(-)\subset
\on{GL}_2$. We will construct in Section \ref{sect:section} a
semigroup scheme morphism $\ul{\on{GT}}(-)\stackrel{\sigma}{\to}
\ul{\on{GT}}_{ell}(-)$, such that 
$$
\begin{matrix}
\ul{\on{GT}}(-) & \to & {\mathbb A}^{1}(-)\\
\scriptstyle{\sigma}\downarrow & & \downarrow \\
\ul{\on{GT}}_{ell}(-) & \to & \on{M}_2(-)
\end{matrix}
$$
commutes, where ${\mathbb A}^{1}\to \on{M}_2$ is $t\mapsto 
\bigl( \begin{smallmatrix} t & 0 \\ 0
& 1 \end{smallmatrix}\bigr)$. Then $\tilde G(-)$
contains the image of ${\mathbb G}_m\to \on{GT}(-)
\stackrel{\sigma}{\to} \on{GT}_{ell}(-)\to\on{GL}_2$, where
${\mathbb G}_m\to \on{GT}(-)$ is a section of $\on{GT}(-)\to 
{\mathbb G}_m$ (see \cite{Dr:Gal}), which is the image of 
${\mathbb G}_m\to\on{GL}_2$, $t\mapsto \bigl( \begin{smallmatrix} t & 0 \\ 0
& 1 \end{smallmatrix}\bigr)$. Then $\on{Lie}(\tilde G) = \GL_2$, 
so $\tilde G=\on{GL}_2$, as wanted. \hfill \qed\medskip 

\subsection{The Zariski closure $\langle B_3\rangle \subset R_{ell}(-)$}

Recall that we have a group morphism $B_3=R_{ell}\to R_{ell}(\QQ)$. 
The Zariski closure $\langle B_3 \rangle \subset 
R_{ell}(-)$ is then the subgroup scheme defined as 
$$ 
\langle B_3 \rangle  :=\bigcap_{
\begin{smallmatrix} 
G\subset R_{ell}(-)\ 
\on{subgroup}\ \on{scheme}| \\
G(\QQ) \supset \on{Im}(B_3
\to R_{ell}(\QQ))
\end{smallmatrix} } G 
$$
Let us compute the Lie algebra\footnote{Recall that the Lie 
algebra of a $\QQ$-group scheme $G$ is $\on{Ker}(G(\QQ[\eps]/(\eps^2))\to
G(\QQ))$.} inclusion $\on{Lie}\langle B_3\rangle \subset \on{Lie}R_{ell}(-)$. 
First, $\on{Lie}R_{ell}(-)$ is a Lie subalgebra of
$$
\on{Lie}\on{Aut}(F_2(-))^{op} \simeq 
\on{Lie}\on{Aut}(\hat{\mathfrak f}_{2})^{op} \simeq 
(\on{Der}\hat\f_{2})^{op} \simeq \hat\f_{2}^{2}, 
$$
where: 

$\bullet$ $\hat\f_2 := \hat\f_2^{\QQ}$ is the Lie algebra freely 
generated by $\xi:= \on{log}X$ and $\eta := \on{log}Y$; 

$\bullet$ the first map is based on the isomorphism 
$F_{2}(\kk)\simeq \on{exp}(\hat{\mathfrak f}_{2}^{\kk})$; 

$\bullet$ the Lie algebra structure on $\hat\f_{2}^{2}$ is given by 
$$
[(\alpha,\beta),(\alpha',\beta')]:= 
(D_{\alpha',\beta'}(\alpha),D_{\alpha',\beta'}(\beta))
-(D_{\alpha,\beta}(\alpha'),D_{\alpha,\beta}(\beta')) ,
$$ 
where $D_{\alpha,\beta}\in \on{Der}(\hat\f_2)$ is given by 
$\xi\mapsto\alpha$, $\eta\mapsto\beta$; 

$\bullet$ the last isomorphism $(\on{Der}\hat\f_{2})^{op} \simeq 
\hat\f_{2}^{2}$ has inverse $(\alpha,\beta)\mapsto D_{\alpha,\beta}$. 

\begin{lemma} \label{descr:R:ell}
$\on{Lie}R_{ell}(-)\subset \on{Lie}\on{Aut}(F_{2}(-))^{op}$ identifies with the 
set of $(\alpha,\beta)\in\hat\f_{2}^{2}$ such that 
$$
\tilde\alpha(X_{1},Y_{1})+\tilde\alpha(X_{2}\sigma_{1}^{-2},Y_{2})
+\tilde\alpha(X_{3}\sigma_{1}^{-1}\sigma_{2}^{-2}\sigma_{1}^{-1},Y_{3})
=0, $$
$$
\tilde\beta(X_{1},Y_{1})+\tilde\beta(X_{2},Y_{2}\sigma_{1}^{2})
+\tilde\beta(X_{3},Y_{3}\sigma_{1}\sigma_{2}^{2}\sigma_{1})=0, 
$$
$$
(\on{Ad}X_{2}^{-1}-1)\tilde\beta(X_{1},Y_{1}) 
+(1-\on{Ad}Y_{1}^{-1})\tilde\alpha(X_{2}\sigma_{1}^{-2},
\sigma_{1}^{2}Y_{2}) = 0$$
(relations in $\on{Lie}P_{1,3}(-)$). Here $\tilde\alpha(X_{1},Y_{1}),\ldots$
are the images of the elements 
$$
\tilde\alpha(e^{\xi},e^{\eta}) := {{1- e^{-\on{ad}\xi}}\over{\on{ad}\xi}}
(\alpha(\xi,\eta)), \quad 
\tilde\beta(e^{\xi},e^{\eta}) := {{1- e^{-\on{ad}\eta}}\over{\on{ad}\eta}}
(\beta(\xi,\eta)) 
$$
of $\hat\f_{2}$ by the morphism $\hat\f_{2}\to \on{Lie}P_{1,3}(-)$, 
$\xi\mapsto \on{log}X_{1}$, $\eta\mapsto \on{log}Y_{1}$, etc., and
$X_{i}:= X_{i}^{+}$, $Y_{i}:= X_{i}^{-}$ (elements of $P_{1,}$). 

The above relations imply the relations
$$
\tilde\alpha(X_{1},Y_{1}) + \tilde\alpha(X_{1}^{-1}\sigma_{1}^{-2},
Y_{1}^{-1}) = 0, \quad 
\tilde\beta(X_{1},Y_{1}) + \tilde\beta(X_{1}^{-1},Y_{1}^{-1}
\sigma_{1}^{2}) = 0, $$
$$
(\on{Ad}X_{1}-1)\tilde\beta(X_{1},Y_{1}) + (1-\on{Ad}Y_{1}^{-1})
\tilde\alpha(X_{1}^{-1}\sigma_{1}^{-2},\sigma_{1}^{2}Y_{1}^{-1}) = 0$$
in $\on{Lie}P_{1,2}(-)$. 
\end{lemma}

{\em Proof.} $(\alpha,\beta)\in\hat\f_{2}^{2}\simeq (\on{Der}\hat\f_{2})^{op}$
induces the infinitesimal automorphism of $F_{2}(\QQ)$ given by 
$X\mapsto g_{+}(X,Y) = X(1+\epsilon\tilde\alpha(X,Y))$, 
$Y\mapsto g_{-}(X,Y) = Y(1+\epsilon\tilde\beta(X,Y))$, where $\epsilon^{2}=0$.
The condition that $(1,1,g_{+},g_{-})$ belongs to $R_{ell}(\QQ[\epsilon]
/(\epsilon^{2}))$ linearizes as follows
$$
\big(\on{id}+ \on{Ad}(\sigma_{2}\sigma_{1}X_{1}) 
+ \on{Ad}(\sigma_{2}\sigma_{1}X_{1})^{2} \big)
(\tilde\alpha(X_{1},Y_{1}) ) = 0, 
$$
$$
\big(\on{id}+ \on{Ad}(\sigma_{2}^{-1}\sigma_{1}^{-1}Y_{1}) 
+ \on{Ad}(\sigma_{2}^{-1}\sigma_{1}^{-1}Y_{1})^{2} \big)
(\tilde\beta(X_{1},Y_{1}) ) = 0, 
$$
$$
\big( Y_{1}(1+\epsilon\tilde\beta(X_{1},Y_{1})),
(1-\epsilon\tilde\alpha(\sigma_{1}^{-2}X_{2},Y_{2}\sigma_{1}^{2}))
X_{2}^{-1}\big) = (Y_{1},X_{2}^{-1}), 
$$
which are equivalent to the announced identities using the relations in $P_{1,3}$ : 
$(X_{i},X_{j}) = (Y_{i},Y_{j})
=1$, 
$$
(Y_{1},X_{1}) = \sigma_{1}\sigma_{2}^{2}\sigma_{1}, \quad
(Y_{1},X_{2}^{-1}) = \sigma_{1}^{2} = 
(Y_{2}^{-1},X_{1}) , \quad 
(Y_{1},X_{3}^{-1}) = \sigma_{2}^{-1}\sigma_{1}^{2}\sigma_{2}, $$
$$
(X_{1},Y_{3}^{-1}) = \sigma_{2}\sigma_{1}^{-2}\sigma_{2}^{-1}, \quad 
(Y_{2},X_{3}^{-1}) = \sigma_{2}^{2} = (Y_{3}^{-1},X_{2}), 
\quad (Y_{3}^{-1},X_{3}^{-1}) = 
\sigma_{2}\sigma_{1}^{2}\sigma_{2}.  $$
\hfill \qed\medskip 


We now compute $\on{Lie}\langle B_3\rangle\subset 
\on{Lie}R_{ell}(-)$. 

\begin{lemma}
Let $u_+:= (0,{{\on{ad}\xi_-}\over{1-e^{-\on{ad}\xi_-}}}(\xi_+))$, 
$u_-:= ({{\on{ad}\xi_+}\over{1-e^{-\on{ad}\xi_+}}}(\xi_-),0)$ in 
$\on{Lie}\on{Aut}(F_{2}(-))^{op}\simeq \hat\f_{2}^{2}$, then 
$u_\pm \in \on{Lie}\langle B_3\rangle$.  
\end{lemma} 

{\em Proof.} We have morphisms ${\mathbb G}_a\to \langle 
B_3\rangle \subset \on{Aut}(F_2(-))^{op}$, 
extending $\ZZ\to B_3$, 
$1\mapsto\Psi_\pm^{\pm 1}$. The corresponding morphisms
$(\kk,+)\to\on{Aut}(F_2(\kk))^{op}$ are 
$t\mapsto (X\mapsto X,Y\mapsto YX^t)$ and 
$t\mapsto (X\mapsto XY^t,Y\mapsto Y)$. The equality 
$XY^{t} = e^{\xi}e^{t\eta} = \on{exp}(\xi + t{{\on{ad}\xi} 
\over {1-e^{-\on{ad}\xi}}}(\eta))$, valid for $t^{2}=0$, and the similar 
equality for $YX^{t}$, imply that the associated Lie algebra morphisms 
are $\QQ\to \on{Lie}\on{Aut}(F_2(-))^{op}$, $1\mapsto u_\pm$, 
which proves that $u_\pm\in\on{Lie}\langle B_3\rangle$. 
\hfill\qed\medskip 

\begin{proposition} \label{prop:Lie<SL_2>}
$\on{Lie}\langle B_3\rangle\subset 
\on{Lie}\on{Aut}(F_2(-))^{op} \simeq \hat\f_2^2$
is the smallest closed Lie subalgebra containing $u_+$ and $u_-$. 
In particular, the image of $\on{Lie}\langle B_{3}\rangle$
by the morphism $\on{Der}(\hat\f_{2})^{op}\to{\mathfrak{gl}}_{2}$
induced by the abelianization map $\f_{2}\to\QQ^{2}$ is ${\mathfrak{sl}}_{2}$. 
\end{proposition}

We first prove: 

\begin{lemma} \label{lemma:Levi}
Let $G$ be a proalgebraic group 
over $\QQ$ fitting in $1\to U\to G\to G_0\to 1$, where
$G_0$ is semisimple and $U$ is prounipotent. Let $0\to 
\u\to \g\to \g_0\to 0$ be the corresponding exact sequence of
Lie algebras. Then $H\mapsto \on{Lie}H$ sets up a bijection 
$\{$proalgebraic subgroups $H\subset G$, such that $\on{Im}(H\subset G
\to G_0) = G_0\} \stackrel{\sim}{\to}
\{$closed Lie subalgebras $\h\subset \g$, 
such that $\on{Im}(\h\subset\g\to \g_0)=\g_0\}$. 
\end{lemma}

{\em Proof.} If $H$ is in the first set, then we have an exact sequence 
$1\to H\cap U\to H\to G_0\to 1$, where $H\cap U$ is necessarily 
prounipotent, hence connected, which implies that $H$ is connected. 
According to \cite{TY}, Prop. 24.3.5, ii), if $\tilde G$ is an 
algebraic group, then the map $\{$connected algebraic subgroups of 
$\tilde G\}\to\{$Lie subalgebras of $\on{Lie}\tilde G\}$ defined by 
taking Lie algebras, is injective. Applying this to the algebraic quotients of 
$G$, one derives the injectivity of the map $H\mapsto \on{Lie}H$. 

Let us prove its surjectivity. Let $\h$ belong to the second set. 

First note that according to the Levi-Mostow decomposition (\cite{Mo}, 
\cite{BS} prop. 5.1), there exists a section $\tilde\sigma:G_0\to G$
of $G\to G_0$. We denote by $\sigma:\g_0\to\g$ its infinitesimal. 
Any section of $\g\to\g_0$ is then conjugate to $\sigma$
by an element of $U(\QQ)$. 

Then we have an exact sequence $0\to \h\cap\u\to \h\to\g_0\to 0$; 
applying the Levi decomposition theorem for Lie algebras, we obtain a section 
$\tau:\g_0\to\h$ of $\h\to\g_0$. Now the composite map 
$\g_0\stackrel{\tau}{\to}\h\hookrightarrow\g$ is a section of $\g\to\g_0$, 
hence of the form $\on{Ad}(x)\circ\sigma$, where $x\in U(\QQ)$. 
If $\vv:= \h\cap\u$, we then have $[\on{Ad}(x)(\sigma(\g_0)),\vv]
\subset\vv$. 

Let then $V\subset U$ be the subgroup with Lie algebra $\vv$; if
we set $H:= V \cdot \on{Ad}(x)(\tilde\sigma(G_0)) 
= \on{Ad}(x)(\tilde\sigma(G_0)) \cdot V$, then $H$ is in the first set, 
and has Lie algebra $\h$. \hfill \qed\medskip 

{\em Proof of Proposition \ref{prop:Lie<SL_2>}.}
Let $\on{Lie}(u_+,u_-)$ be the smallest closed Lie 
subalgebra of $$\on{Lie}\on{Aut}(F_2(-))^{op}$$ containing
$u_+$ and $u_-$. Then $\on{Lie}\langle B_3 \rangle\supset 
\on{Lie}(u_+,u_-)$. Apply now Lemma \ref{lemma:Levi}
with $G=R_{ell}(-)$, $G_0 = \on{SL}_2$. The map 
$\g\to\g_0 = {\mathfrak{sl}}_{2}$ is such that $u_+\mapsto 
\bigl( \begin{smallmatrix} 0 & 0 \\ 1 & 0
\end{smallmatrix}\bigr)$ and $u_-\to 
\bigl( \begin{smallmatrix} 0 & 1 \\ 0 & 0
\end{smallmatrix}\bigr)$, 
so if $\h:= \on{Lie}(u_+,u_-)$, then $\on{Im}(\h\subset
\g\to\g_0)=\g_0$. Let then $H\subset R_{ell}(-)$ be the proalgebraic 
subgroup corresponding to $\h$ by Lemma \ref{lemma:Levi}; 
then $\langle B_3\rangle\supset H$. 
On the other hand, we have group morphisms ${\mathbb G}_a\to H$ corresponding to 
$\QQ\to\h$, $1\mapsto u_\pm$, whose versions over $\QQ$ 
are $(\QQ,+)\to H(\QQ)\subset \on{Aut}(F_2(\QQ))^{op}$, 
$t\mapsto \Psi_\pm^t$. Setting $t=1$, we obtain 
$H(\QQ)\ni\Psi_\pm$, and as 
$\Psi_+,\Psi_-$ generate $B_3$, 
$H(\QQ)\supset B_{3}$. So 
$\langle B_3\rangle=H$. Taking Lie 
algebras, we obtain Proposition \ref{prop:Lie<SL_2>}. 
\hfill \qed\medskip

\begin{remark} Let $d_\pm:= [[u_+,u_-],u_\pm]\pm 2 u_\pm$. 
Then for any $(\alpha_\pm,\beta)\in\NN^3$, 
$$
x_{\alpha_\pm,\beta}^\pm:= \on{ad}(u_+)^{\alpha_+}
\on{ad}(u_-)^{\alpha_-}\on{ad}([u_+,u_-])^\beta(d_\pm)
\in \on{Ker}(\on{Lie}\langle B_{3}\rangle\to\SL_2).
$$ 
Then $\on{Ker}(\on{Lie}\langle B_3\rangle\to\SL_2)$
is topologically generated by these elements, more precisely, 
it is equal to $\{\sum_{n\geq 1}P_n((x^\pm_{\alpha_\pm,\beta})) |
(P_n)_n\in\prod_{n\geq 1}\f_n\}$, where $\f_n$ is the 
part of degree $n$ of the free Lie algebra with generators indexed 
by $\NN^3\times\{\pm\}$ (each generator having degree 1). Then 
$\on{Lie}\langle B_3\rangle = 
\on{Ker}(\on{Lie}\langle B_3\rangle\to\SL_2)
\oplus \on{Span}_\QQ(u_+,u_-,[u_+,u_-])$. 
\end{remark}

\subsection{A morphism ${\ul{\on{GT}}}\to{\ul{\on{GT}}}_{ell}$ and its variants}
\label{sect:section}

We now construct a section of the semigroup morphism 
${\ul{\on{GT}}}_{ell}\to{\ul{\on{GT}}}$ and of its variants.  

\begin{proposition} \label{prop:lift}
There exists a unique semigroup morphism 
${\ul{\on{GT}}}\to{\ul{\on{GT}}}_{ell}$, defined by 
$(\lambda,f)\mapsto (\lambda,f,g_\pm)$, where 
$$
g_+(X,Y)= f(X,(Y,X))X^\lambda f(X,(Y,X))^{-1}, 
$$
$$
g_-(X,Y)=(Y,X)^{{{\lambda-1}\over 2}}
f(YX^{-1}Y^{-1},(Y,X))Y f(X,(Y,X))^{-1}. 
$$
The same formulas define semigroup morphisms 
$\widehat{{\ul{\on{GT}}}}\to\widehat{{\ul{\on{GT}}}}_{ell}$, 
$\ul{\on{GT}}_l\to \ul{\on{GT}}^{ell}_l$, and a semigroup 
scheme morphism ${\ul{\on{GT}}}(-)\to {\ul{\on{GT}}}_{ell}(-)$,
compatible with the natural maps between the various versions of 
$\ul{\on{GT}}_{(ell)}$. 

There are commutative diagrams 
$$
\xymatrix{
{\ul{\on{GT}}} \ar[r]\ar[d]_{\sim} & {\ul{\on{GT}}}_{ell} \ar[d]_{\sim} \\
\ZZ/2 \ar[r] & \tilde B_{3} 
} \quad \text{and }\quad
\xymatrix{
{\ul{\on{GT}}}(-) \ar[r]\ar[d] & {\ul{\on{GT}}}_{ell}(-) \ar[d] \\
{\mathbb A}^{1} \ar[r] & \on{M}_2 } 
$$
where the bottom morphisms are 
$\bar 1\mapsto \varepsilon\Psi_{+}\Psi_{-}\Psi_{+}$
and $\lambda\mapsto \bigl( \begin{smallmatrix} \lambda & 0 \\ 0 & 1
\end{smallmatrix}\bigr)$. 
\end{proposition}

{\em Proof.} 1) As the center $Z(B_{n+1})$ of $B_{n+1}$ is contained in the pure braid 
group $P_{n+1}:=\on{Ker}(B_{n+1}\to S_{n+1})$, the morphism $B_{n+1}\to 
S_{n+1}$ descends to a morphism $B_{n+1}/Z(B_{n+1})\to S_{n+1}$. 
Identify $S_{n+1}\simeq\on{Perm}(\{0,\ldots,n\})$ and let $S_n\subset
S_{n+1}$ be $\{\sigma|\sigma(0)=0\}$. The Cartesian product 
$$
(B_{n+1}/Z(B_{n+1}))\times_{S_{n+1}} S_{n}
$$
then identifies with the quotient $(B_{n+1}\times_{S_{n+1}}S_{n})/Z(B_{n+1})$
relative to the sequence of inclusions $Z(B_{n+1})\subset B_{n+1}\times_{S_{n+1}}
S_{n} \subset B_{n+1}$. The middle subgroup identifies with a type B braid group
and is generated by $\sigma_{0}^{2},\sigma_{1},\ldots,\sigma_{n-1}$, where
the generators of $B_{n+1}$ are labeled $\sigma_{0},\ldots,\sigma_{n-1}$. 
Using the presentation of the type B group, one proves that there is a unique morphism 
$B_{n+1}\times_{S_{n+1}}S_{n}\to B_{n+1}$, such that $\sigma_{0}^{2}\mapsto 
X_{1}^{+}$, $\sigma_{i}\mapsto \sigma_{i}$ ($i>1$). Moreover,  
this morphism takes a generator of $Z(B_{n+1}) \simeq \ZZ$ to $X_{1}^{+}\cdots
X_{n}^{+} = 1\in B_{1,n}$. It follows that it factors though a morphism 
$(B_{n+1}\times_{S_{n+1}}S_{n})/Z(B_{n+1})\to B_{1,n}$, i.e., 
\begin{equation} \label{morph:B:B}
(B_{n+1}/Z(B_{n+1}))\times_{S_{n+1}}S_n\to B_{1,n}. 
\end{equation}
This morphism admits the following interpretation. 
If $X$ is a topological additive group, let $C_{[n]}(X):= \on{Inj}([n],X)/S_{n}$, 
where $\on{Inj}$ means the space of injections, $[n]:= \{1,\ldots,n\}$,
and $\overline C_{[n]}(X):= C_{[n]}(X)/X$, where $X$ acts by addition of 
a constant function. We then have the identifications 
$$
\pi_{1}(C_{[n]}(\CC^{\times}))\simeq B_{n+1}\times_{S_{n+1}} S_{n}, \quad 
\pi_{1}(\overline C_{n}(\CC^{\times})) \simeq 
(B_{n+1}\times_{S_{n+1}}S_{n})/Z(B_{n+1}), 
$$
$$
\pi_{1}(\overline C_{n}(\CC^{\times}/q^{\ZZ})) \simeq B_{1,n}, $$
where $q$ is a real number with $0<q<1$. The canonical projection $\CC^{\times}
\to \CC^{\times}/q^{\ZZ}$ then induces a group morphism $\pi_{1}(\overline C_{n}
(\CC^{\times}))\to \pi_{1}(\overline C_{n}(\CC^{\times}/q^{\ZZ}))$, which 
turns out to coincide with (\ref{morph:B:B}). 

Any $(\lambda,f)\in\ul{\on{GT}}$ induces an endomorphism $F_{\lambda,f}$ of 
${\bf PaB}$, such that for any object $O$ and $z\in {\bf PaB}(O)$, 
$$
F_{\lambda,f}(z) = z^{\lambda}
$$
if $z$ corresponds to an element of $Z(B_{|O|})$. 

The element $\sigma_{2}\sigma_{1}\sigma_{0}^{2}
\in B_{4}$ corresponds to
\begin{equation} \label{elt:PaB}
(\on{id}_{\bullet}\otimes\beta_{\bullet,\bullet\bullet})
a_{\bullet,\bullet,\bullet\bullet}
(\beta^{2}_{\bullet,\bullet}\otimes\on{id}_{\bullet\bullet})
a_{\bullet,\bullet,\bullet\bullet}^{-1}
(\on{id}_{\bullet}\otimes a_{\bullet,\bullet,\bullet})
\in {\bf PaB}(\bullet((\bullet\bullet)\bullet)).
\end{equation}
The image of (\ref{elt:PaB}) by this endomorphism 
is the product of the images of its factors, namely 
$$
F_{\lambda,f}(\on{id}_{\bullet}\otimes\beta_{\bullet,\bullet\bullet}) = 
\on{id}_{\bullet}\otimes \beta_{\bullet,\bullet\bullet}
(\beta_{\bullet\bullet,\bullet}\beta_{\bullet,\bullet\bullet})^{m}
\in{\bf PaB}(\bullet(\bullet(\bullet\bullet)),\bullet((\bullet\bullet)
\bullet))\leftrightarrow 
\sigma_{2}\sigma_{1}(\sigma_{1}\sigma_{2}^{2}\sigma_{1})^{m}\in B_{4}, 
$$
\begin{align*}
F_{\lambda,f}(a_{\bullet,\bullet,\bullet\bullet}) = 
a_{\bullet,\bullet,\bullet\bullet}f(\beta_{\bullet,\bullet}^{2}
\otimes\on{id}_{\bullet\bullet},a^{-1}(\on{id}_{\bullet}\otimes
\beta_{\bullet\bullet,\bullet}\beta_{\bullet,\bullet\bullet})a)
\in  & {\bf PaB}( (\bullet\bullet)(\bullet\bullet),
\bullet(\bullet(\bullet\bullet)) ) \\
 & \leftrightarrow 
f(\sigma_{0}^{2},\sigma_{1}\sigma_{2}^{2}\sigma_{1})\in B_{4}, 
\end{align*}
$$
F_{\lambda,f}(\beta_{\bullet,\bullet}^{2}\otimes \on{id}_{\bullet\bullet})
 = \beta_{\bullet,\bullet}^{2\lambda}\otimes\on{id}_{\bullet\bullet}
 \in {\bf PaB}((\bullet\bullet)(\bullet\bullet))\leftrightarrow 
 \sigma_{0}^{2\lambda}\in B_{4}, 
$$
$$
F_{\lambda,f}(\on{id}_{\bullet}\otimes a_{\bullet,\bullet,\bullet})
 = \on{id}_{\bullet}\otimes a_{\bullet,\bullet,\bullet}
 f(\beta_{\bullet\bullet}^{2}\otimes\on{id}_{\bullet},
 a^{-1}(\on{id}_{\bullet}\otimes\beta_{\bullet\bullet}^{2})a)
\in {\bf PaB}(\bullet((\bullet\bullet)\bullet),
\bullet(\bullet(\bullet\bullet)))
\leftrightarrow f(\sigma_{1}^{2},\sigma_{2}^{2})\in B_{4}. 
$$
Therefore 
$$
F_{\lambda,f}((\ref{elt:PaB}))\in {\bf PaB}(\bullet((\bullet\bullet)\bullet))
\leftrightarrow 
\sigma_2\sigma_1(\sigma_1\sigma_2^2\sigma_1)^m
f(\sigma_0^2,\sigma_1\sigma_2^2\sigma_1)\sigma_0^{2\lambda}
f^{-1}(\sigma_0^2,\sigma_1\sigma_2^2\sigma_1) 
f(\sigma_1^2,\sigma_2^2) \in B_{4}. 
$$
Now $(\sigma_{2}\sigma_{1}\sigma_{0}^{2})^{3}$ generates $Z(B_{4})$, 
therefore
$$
(\ref{elt:PaB})^{3} \in{\bf PaB}(\bullet((\bullet\bullet)
\bullet))\leftrightarrow (\sigma_{2}\sigma_{1}\sigma_{0}^{2})^{3}
\in Z(B_{4}). 
$$
It follows that $F_{\lambda,f}((\ref{elt:PaB})^{3}) = (\ref{elt:PaB})^{3\lambda}$. 
The image of this equality in $B_{4}$ is 
$$
\Big(\sigma_{2}\sigma_{1}(\sigma_{1}\sigma_{2}^{2}\sigma_{1})^{m}
f(\sigma_{0}^{2},\sigma_{1}\sigma_{2}^{2}\sigma_{1})
\sigma_{0}^{2\lambda}f^{-1}(\sigma_{0}^{2},\sigma_{1}\sigma_{2}^{2}
\sigma_{1}) f(\sigma_{1}^{2},\sigma_{2}^{2})\Big)^{3} 
= (\sigma_{2}\sigma_{1}\sigma_{0}^{2})^{3\lambda}. 
$$
As $Z(B_{4})$ is in the kernel of $B_{4}\times_{S_{4}}S_{4}\to B_{1,3}$, 
the image of the left hand side of this equality under this morphism is $1\in B_{1,3}$. 
It follows that 
\begin{equation} \label{a'}
\Big( \sigma_2\sigma_1 (\sigma_1\sigma_2^2\sigma_1)^m
f(X_1^+,(X_1^-,X_1^+)) (X_1^+)^\lambda f^{-1}(X_1^+,
(X_1^-,X_1^+)) f(\sigma_1^2,\sigma_2^2)\Big)^3=1
\end{equation}
in $B_{1,3}$.
This means that identity (\ref{def:GTell:1}) is satisfied with $\pm=+$. 

2) We show that $g_-$ satisfies (\ref{def:GTell:1})
with $\pm=-$, i.e., 
$$
\Big(\sigma_2^{-1}\sigma_1^{-1}(\sigma_1\sigma_2^2\sigma_1)^{-m}
g_-(X_1,Y_{1})f(\sigma_1^2,\sigma_2^2)\Big)^3=1
$$
in $B_{1,3}$ (we set $X_{i}:=X_{i}^{+}$, $Y_{i}:=X_{i}^{-}$). 
Substituting the given expression for $g_{-}(X_{1},Y_{1})$, 
using $(Y_{1},X_{1}) = \sigma_{1}\sigma_{2}^{2}\sigma_{1}$, 
the identities $(\sigma_{2}^{-1}\sigma_{1}^{-1}) Y_{1}X_{1}^{-1}Y_{1}^{-1} 
= X_{3}^{-1}(\sigma_{2}^{-1}\sigma_{1}^{-1})$, 
$(\sigma_{2}^{-1}\sigma_{1}^{-1})\sigma_{1}\sigma_{2}^{2}\sigma_{1}
= \sigma_{2}\sigma_{1}^{2}\sigma_{2}(\sigma_{2}^{-1}\sigma_{1}^{-1})$, 
and 
after a suitable conjugation, this equality is equivalent to 
\begin{equation} \label{labex}
\Big(Y_{1}f^{-1}(X_{1},\sigma_{1}\sigma_{2}^{2}\sigma_{1})
f(\sigma_{1}^{2},\sigma_{2}^{2})f(X_{3}^{-1},\sigma_{2}\sigma_{1}^{2}
\sigma_{2})\sigma_{2}^{-1}\sigma_{1}^{-1}\Big)^{3}=1. 
\end{equation}

As $f\in F_{2}'$, $f(a\alpha,b) = f(a,b)$ if $\alpha$ commutes with both $a$ and $b$. 
In particular, $\sigma_{1}^{2}$ commutes (in $B_{4}$) with both $\sigma_{0}
\sigma_{1}^{2}\sigma_{0}$ and $\sigma_{2}\sigma_{1}^{2}\sigma_{2}$. 
It follows that $f(\sigma_{0}\sigma_{1}^{2}\sigma_{0},\sigma_{2}\sigma_{1}^{2}
\sigma_{2}) = f((\sigma_{1}^{2}\sigma_{0})^{2},\sigma_{2}\sigma_{1}^{2}
\sigma_{2})$. Since $(\sigma_{1}^{2}\sigma_{0})^{2}
 = (\sigma_{1}\sigma_{0}^{2})^{2}$, $f(\sigma_{0}\sigma_{1}^{2}\sigma_{0},
 \sigma_{2}\sigma_{1}^{2}\sigma_{2}) = f((\sigma_{1}\sigma_{0}^{2})^{2},
 \sigma_{2}\sigma_{1}^{2}\sigma_{2})$. Substituting this identity in the 
 pentagon identity 
 $$
 f(\sigma_{1}^{2},\sigma_{2}^{2})f(\sigma_{0}\sigma_{1}^{2}\sigma_{0},
 \sigma_{2}\sigma_{1}^{2}\sigma_{2})f(\sigma_{0}^{2},\sigma_{1}^{2})
 = f(\sigma_{0}^{2},\sigma_{1}\sigma_{2}^{2}\sigma_{1})
 f(\sigma_{1}\sigma_{0}^{2}\sigma_{1},\sigma_{2}^{2})
 $$
in $P_{4}:= \on{Ker}(B_{4}\to S_{4})$, taking the image of the resulting 
identity by the morphism $P_{4}\subset B_{4}\times_{S_{4}}S_{3}\to B_{1,3}$, 
and using the identity $X_{2}X_{1} = X_{3}^{-1}$ in $B_{1,3}$, one obtains
$$
f(\sigma_{1}^{2},\sigma_{2}^{2})f(X_{3}^{-1},\sigma_{2}\sigma_{1}^{2}
\sigma_{2})f(X_{1},\sigma_{1}^{2}) = f(X_{1},\sigma_{1}\sigma_{2}^{2}
\sigma_{1}) f(X_{2},\sigma_{2}^{2}) 
$$
(identity in $B_{1,3}$). Using this identity, (\ref{labex}) is equivalent to 
$$
(Y_{1}A\sigma_{2}^{-1}\sigma_{1}^{-1})^{3}=1,
$$
where $A := f(X_{2},\sigma_{2}^{2})f^{-1}(X_{1},\sigma_{1}^{2})$. 
Using $Y_{3} = \sigma_{2}^{-1}\sigma_{1}^{-1}Y_{1}\sigma_{1}^{-1}
\sigma_{2}^{-1}$, $Y_{2} = \sigma_{2}\sigma_{1}\sigma_{2}^{-1}\sigma_{1}^{-1}
Y_{1}\sigma_{2}^{-1}\sigma_{1}^{-1}$, the latter identity is equivalent to 
\begin{equation} \label{idex}
Y_{1} A Y_{3} (\sigma_{2}\sigma_{1}A\sigma_{1}^{-1}\sigma_{2}^{-1})
Y_{2} (\sigma_{1}\sigma_{2}A\sigma_{2}^{-1}\sigma_{1}^{-1})=1.
\end{equation}
As $Y_{1}X_{2} = (X_{2}\sigma_{1}^{-2})Y_{1}$, $Y_{1}\sigma_{2}^{2}
=\sigma_{2}^{2}Y_{1}$, $X_{1}Y_{3}=Y_{3}(\sigma_{2}\sigma_{1}^{2}
\sigma_{2}^{-1}X_{1})$, $\sigma_{1}^{2}Y_{3}=Y_{3}\sigma_{1}^{2}$, 
and $Y_{1}Y_{3}=Y_{2}^{-1}$, 
$$
Y_{1}AY_{3} = f(X_{2}\sigma_{1}^{-2},\sigma_{2}^{2})Y_{2}^{-1}
f^{-1}(\sigma_{2}\sigma_{1}^{2}\sigma_{2}^{-1}X_{1},\sigma_{1}^{2}).
$$
As $\on{Ad}(\sigma_{2}\sigma_{1})(X_{2}) 
= \sigma_{2}\sigma_{1}^{2}\sigma_{2}^{-1}
X_{1}$, $\on{Ad}(\sigma_{2}\sigma_{1})(\sigma_{2}^{2})
= \sigma_{1}^{2}$, $\on{Ad}(\sigma_{2}\sigma_{1})(X_{1})
= X_{3}\sigma_{2}^{-1}\sigma_{1}^{-2}\sigma_{2}^{-1}$, 
$\on{Ad}(\sigma_{2}\sigma_{1})(\sigma_{1}^{2}) = 
\sigma_{2}\sigma_{1}^{2}\sigma_{2}^{-1}$, 
$$
\sigma_{2}\sigma_{1}A\sigma_{1}^{-1}\sigma_{2}^{-1} = 
f(\sigma_{2}\sigma_{1}^{2}\sigma_{2}^{-1}X_{1},\sigma_{1}^{2})
f^{-1}(X_{3}\sigma_{2}^{-1}\sigma_{1}^{-2}\sigma_{2}^{-1},
\sigma_{2}\sigma_{1}^{2}\sigma_{2}^{-1}). 
$$
As $\on{Ad}(\sigma_{1}\sigma_{2})(X_{2}) = X_{3}\sigma_{2}^{-1}
\sigma_{1}^{-2}\sigma_{2}$, $\on{Ad}(\sigma_{1}\sigma_{2})(\sigma_{2}^{2}) 
= \sigma_{2}\sigma_{1}^{2}\sigma_{2}^{-1}$, $\on{Ad}(\sigma_{1}\sigma_{2})
(X_{1}) = X_{2}\sigma_{1}^{-2}$, $\on{Ad}(\sigma_{1}\sigma_{2})
(\sigma_{1}^{2}) = \sigma_{2}^{2}$, 
$$
\sigma_{1}\sigma_{2} A \sigma_{2}^{-1}\sigma_{1}^{-1} = 
f(X_{3}\sigma_{2}^{-1}\sigma_{1}^{-2}\sigma_{2},\sigma_{2}
\sigma_{1}^{2}\sigma_{2}^{-1}) f^{-1}(X_{2}\sigma_{1}^{-2},\sigma_{2}^{2}).
$$
Taking these equalities into account and after simplification and 
conjugation, (\ref{idex}) is equivalent to 
$$
Y_{2}^{-1}f^{-1}(X_{3}\sigma_{2}^{-1}\sigma_{1}^{-2}\sigma_{2}^{-1},
\sigma_{2}\sigma_{1}^{2}\sigma_{2}^{-1}) Y_{2} f(X_{3}\sigma_{2}^{-1}
\sigma_{1}^{-2}\sigma_{2},\sigma_{2}\sigma_{1}^{2}\sigma_{2}^{-1})=1,  
$$
which follows from $Y_{2}^{-1} \cdot X_{3}\sigma_{2}^{-1}\sigma_{1}^{-2}
\sigma_{2}^{-1}\cdot Y_{2} = X_{3}\sigma_{2}^{-1}\sigma_{1}^{-2}\sigma_{2}$
and from the fact that $Y_{2}$ commutes with $\sigma_{2}\sigma_{1}^{2}
\sigma_{2}^{-1}$. 

3) Since $\sigma_{1}^{2}$ commutes with both $\sigma_{2}\sigma_{1}^{2}
\sigma_{2}$ and $\sigma_{0}\sigma_{1}^{2}\sigma_{0}$ and since $f\in F_{2}'$, 
one has 
$$
f(\sigma_{2}\sigma_{1}^{2}\sigma_{2},\sigma_{0}\sigma_{1}^{2}\sigma_{0}) 
= f(\sigma_{2}\sigma_{1}^{2}\sigma_{2},(\sigma_{0}\sigma_{1}^{2})^{2})
$$
(equality in $B_{4}$). Since $(\sigma_{0}\sigma_{1}^{2})^{2} \equiv (\sigma_{2}\sigma_{1}\sigma_{0}^{2}\sigma_{1}
\sigma_{2})^{-1}$ mod $Z(P_{4})$ and $f\in F_{2}'$, one has
$$
f(\sigma_{2}\sigma_{1}^{2}\sigma_{2},(\sigma_{0}\sigma_{1}^{2})^{2}) =  
f(\sigma_{2}\sigma_{1}^{2}\sigma_{2},
(\sigma_{2}\sigma_{1}\sigma_{0}^{2}\sigma_{1}
\sigma_{2})^{-1})
$$ 
(equality in $B_{4}$). Plugging these equalities in the pentagon equation 
$$
f(\sigma_1\sigma_0^2\sigma_1,\sigma_2^2)f^{-1}(\sigma_0^2,\sigma_1^2)
f(\sigma_2\sigma_1^2\sigma_2,\sigma_0\sigma_1^2\sigma_0)
f^{-1}(\sigma_1^2,\sigma_2^2) f(\sigma_0^2,\sigma_1\sigma_2^2\sigma_1)
=1 
$$
(in $B_{4}$) and multiplying by 
$f^{-1}(\sigma_0^2,\sigma_1\sigma_2^2\sigma_1)$ from the right, one obtains
$$
f(\sigma_1\sigma_0^2\sigma_1,\sigma_2^2)f^{-1}(\sigma_0^2,\sigma_1^2)
f(\sigma_{2}\sigma_{1}^{2}\sigma_{2},
(\sigma_{2}\sigma_{1}\sigma_{0}^{2}\sigma_{1}
\sigma_{2})^{-1})
f^{-1}(\sigma_1^2,\sigma_2^2) = f^{-1}(\sigma_0^2,\sigma_1\sigma_2^2\sigma_1)  
$$
(in $B_{4}$). As $\sigma_{2}$ commutes with both $\sigma_{0}^{2}$ and $\sigma_{1}
\sigma_{2}^{2}\sigma_{1}$, the right side of this equality, and therefore also its
left side, commutes with $\sigma_{2}^{\lambda}$. It follows that the equality also 
holds with the left side replaced by its conjugation of  $\sigma_{2}^{\lambda}$; 
multiplying the resulting equality by $f(\sigma_{0}^{2},\sigma_{1}\sigma_{2}^{2}
\sigma_{1})$ from the right, and using the identities $\sigma_{2}\sigma_{1}^{2}
\sigma_{2} = \sigma_{2}^{-1}\sigma_{1}^{-1}(\sigma_{1}\sigma_{2}^{2}
\sigma_{1})\sigma_{1}\sigma_{2}$, 
$(\sigma_{2}\sigma_{1}\sigma_{0}^{2}\sigma_{1}\sigma_{2})^{-1}
= \sigma_{2}^{-1}\sigma_{1}^{-1}(\sigma_{1}\sigma_{2}^{2}\sigma_{1}
\sigma_{0}^{2})^{-1}\sigma_{1}\sigma_{2}$, one obtains 
\begin{align} \label{interm:1}
\sigma_{2}^{\lambda}f(\sigma_1\sigma_0^2\sigma_1,\sigma_2^2)
f^{-1}(\sigma_0^2,\sigma_1^2) \sigma_{2}^{-1}\sigma_{1}^{-1}
f(\sigma_{1}\sigma_{2}^{2}\sigma_{1},
(\sigma_{1}\sigma_{2}^{2}\sigma_{1}\sigma_{0}^{2})^{-1})
\sigma_{1}\sigma_{2}
f^{-1}(\sigma_1^2,\sigma_2^2) \sigma_{2}^{-\lambda}
& \nonumber f(\sigma_0^2,\sigma_1\sigma_2^2\sigma_1) \\
 & =1.   
\end{align}
On the other hand, $\sigma_1^{-1}\sigma_2^2\sigma_1
=\sigma_2\sigma_1^2\sigma_2^{-1}\equiv (\sigma_2^2\sigma_1^2)^{-1}$
mod $Z(P_3)$ (equalities in $P_{3}$); together with $f\in F_{2}'$, this implies 
$f((\sigma_{2}^{2}\sigma_{1}^{2})^{-1},\sigma_{2}^{2}) = f(\sigma_{2}
\sigma_{1}^{2}\sigma_{2}^{-1},\sigma_{2}^{2})$ 
and $f(\sigma_{1}^{2},(\sigma_{2}^{2}\sigma_{1}^{2})^{-1}) 
= f(\sigma_{1}^{2},\sigma_{1}^{-1}\sigma_{2}^{2}\sigma_{1})$ 
(equalities in $P_{3}$).  Plugging these equalities in the hexagon equation
$$
1 = (\sigma_{2}^{2})^{m}f((\sigma_{2}^{2}\sigma_{1}^{2})^{-1},\sigma_{2}^{2})
(\sigma_{2}^{2}\sigma_{1}^{2})^{-m}f(\sigma_{1}^{2},(\sigma_{2}^{2}
\sigma_{1}^{2})^{-1})(\sigma_{1}^{2})^{m}f(\sigma_{2}^{2},\sigma_{1}^{2}) 
$$
(in $P_{3}$), 
using the equalities $f(\sigma_{2}\sigma_{1}^{2}\sigma_{2}^{-1},\sigma_{2}^{2})
= \sigma_{2}f(\sigma_{1}^{2},\sigma_{2}^{2})\sigma_{2}^{-1}$, 
$f(\sigma_{1}^{2},\sigma_{1}^{-1}\sigma_{2}^{2}\sigma_{1}) = 
\sigma_{1}^{-1}f(\sigma_{1}^{2},\sigma_{2}^{2})\sigma_{1}$, 
$(\sigma_{2}^{2})\sigma_{2} = \sigma_{2}^{\lambda}$,  
$\sigma_{1}(\sigma_{1}^{2})^{m} = \sigma_{1}^{\lambda}$, 
multiplying by $\sigma_{1}\sigma_{2}f^{-1}(\sigma_{1}^{2},\sigma_{2}^{2})
\sigma_{2}^{-\lambda}$ from the left and using $\sigma_{1}(\sigma_{2}^{2}
\sigma_{1}^{2})^{-m}\sigma_{1}^{-1} = (\sigma_{1}\sigma_{2}^{2}
\sigma_{1})^{{1-\lambda}\over 2}$, one obtains 
$$
\sigma_1\sigma_2 f^{-1}(\sigma_1^2,\sigma_2^2)\sigma_2^{-\lambda}
= (\sigma_1\sigma_2^2\sigma_1)^{{{1-\lambda}\over 2}} 
f(\sigma_1^2,\sigma_2^2) \sigma_1^\lambda f^{-1}(\sigma_1^2,\sigma_2^2). 
$$
Plugging this equality in (\ref{interm:1}), one obtains 
\begin{eqnarray*}
&& \sigma_2^\lambda f(\sigma_1\sigma_0^2\sigma_1,\sigma_2^2)
f^{-1}(\sigma_0^2,\sigma_1^2)\sigma_2^{-1}\sigma_1^{-1}
f(\sigma_1\sigma_2^2\sigma_1,\sigma_0^{-2}\sigma_1^{-1}\sigma_2^{-2}
\sigma_1^{-1}) 
\\ &&
(\sigma_1\sigma_2^2\sigma_1)^{{{1-\lambda}\over 2}}
f(\sigma_1^2,\sigma_2^2) \sigma_1^\lambda
f^{-1}(\sigma_1^2,\sigma_2^2) f(\sigma_0^2,\sigma_1\sigma_2^2\sigma_1) = 1 
\end{eqnarray*}
(in $B_{4}\times_{S_{4}}S_{3}$). Taking the image of this equality under 
$B_{4}\times_{S_{4}}S_{3}\to B_{1,3}$ and multiplying the resulting 
equality by $\sigma_{2}f^{-1}(X_{2},\sigma_{2}^{2})\sigma_{2}^{-\lambda}$ 
from the left, one obtains 
\begin{eqnarray*}
\sigma_2 f^{-1}(X_2,\sigma_2^2)\sigma_2^{-\lambda}
& = & \sigma_1^{-1} f^{-1}(X_2\sigma_1^{-2},\sigma_2^2)
f((Y_1,X_1),X_1^{-1}(Y_1,X_1)^{-1})
(Y_1,X_1)^{{{1-\lambda}\over 2}} \\ && 
f(\sigma_1^2,\sigma_2^2)
\sigma_1^\lambda f^{-1}(\sigma_1^2,\sigma_2^2) f(X_1,(Y_1,X_1))
\end{eqnarray*}
(in $B_{1,3}$). 
As both $\sigma_{2}$ and $X_{2}$ commute with $X_{1}$, the left side of 
this equality commutes with $X_{1}^{\lambda}$, therefore so does its right side. 
Expressing the equality of $X_{1}^{\lambda}$ with its conjugate by the right side, 
and conjugating the resulting equality, one obtains 
\begin{align} \label{iben:bis}
& \on{Ad}\Big(f((Y_{1},X_{1}),X_{1}^{-1}(Y_{1},X_{1})^{-1})
(Y_{1},X_{1})^{{1-\lambda}\over 2}f(\sigma_{1}^{2},\sigma_{2}^{2})
\sigma_{1}^{\lambda}
f^{-1}(\sigma_{1}^{2},\sigma_{2}^{2}) f(X_{1},(Y_{1},X_{1}))\Big)
((X_{1})^{\lambda}) \nonumber \\
 & = \on{Ad}\Big(f(X_{2}\sigma_{1}^{-2},\sigma_{2}^{2})\sigma_{1}\Big)
((X_{1})^{\lambda}) 
\end{align}
(in $B_{1,3}$). 

The equality 
\begin{align*}
& a_{\bullet\bullet,\bullet,\bullet}
(\beta_{\bullet,\bullet\bullet}\beta_{\bullet\bullet,\bullet}
\otimes\on{id}_{\bullet})
a_{\bullet\bullet,\bullet,\bullet}^{-1}
\\
 & = 
 a_{\bullet,\bullet,\bullet\bullet}^{-1}
(\on{id}_{\bullet}\otimes a_{\bullet,\bullet,\bullet})
(\on{id}_{\bullet}\otimes (\beta_{\bullet,\bullet}\otimes\on{id}_{\bullet}))
(\on{id}_{\bullet}\otimes 
a_{\bullet,\bullet,\bullet}^{-1})
a_{\bullet,\bullet,\bullet\bullet}
\cdot (\beta_{\bullet,\bullet}^{2}\otimes\on{id}_{\bullet\bullet}) \cdot
\\
 & 
a^{-1}_{\bullet,\bullet,\bullet\bullet}
(\on{id}_{\bullet}\otimes 
a_{\bullet,\bullet,\bullet})
(\on{id}_{\bullet}\otimes (\beta_{\bullet,\bullet}\otimes\on{id}_{\bullet}))
(\on{id}_{\bullet}\otimes a_{\bullet,\bullet,\bullet}^{-1})
a_{\bullet,\bullet,\bullet\bullet}
\end{align*}
in ${\bf PaB}((\bullet\bullet)(\bullet\bullet))$ follows from the fact that 
both sides correspond to the element $\sigma_{1}\sigma_{0}^{2}
\sigma_{1}\in B_{4}$. Applying the automorphism $F_{\lambda,f}$ to this 
equality, one obtains an equality in
${\bf PaB}((\bullet\bullet)(\bullet\bullet))$, which translates into 
the equality 
\begin{align*}
& f(\sigma_{1}\sigma_{0}^{2}\sigma_{1},\sigma_{2}^{2})
(\sigma_{1}\sigma_{0}^{2}\sigma_{1})^{\lambda}
f(\sigma_{1}\sigma_{0}^{2}\sigma_{1},\sigma_{2}^{2})
\\
 & = f(\sigma_{0}^{2},\sigma_{1}\sigma_{2}^{2}\sigma_{1})^{-1}
f(\sigma_{1}^{2},\sigma_{2}^{2})
\sigma_{1}^{\lambda}
f(\sigma_{1}^{2},\sigma_{2}^{2})^{-1}
f(\sigma_{0}^{2},\sigma_{1}\sigma_{2}^{2}\sigma_{1})
\cdot \sigma_{0}^{2\lambda}\cdot \\
 & 
\cdot f(\sigma_{0}^{2},\sigma_{1}\sigma_{2}^{2}\sigma_{1})^{-1}
f(\sigma_{1}^{2},\sigma_{2}^{2})
\sigma_{1}^{\lambda}
f(\sigma_{1}^{2},\sigma_{2}^{2})^{-1}
f(\sigma_{0}^{2},\sigma_{1}\sigma_{2}^{2}\sigma_{1})
\end{align*}
in $B_{4}\times_{S_{4}}S_{3}\subset B_{4}$. 

As $(Y_{1},X_{1}) = \sigma_{1}\sigma_{2}^{2}\sigma_{1}$ (relation  in $B_{1,3}$), 
the image of this equality in $B_{1,3}$ is 
$$
f(X_{2},\sigma_{2}^{2})X_{2}^{\lambda} f(X_{2},\sigma_{2}^{2})^{-1}
= f(X_{1},(Y_{1},X_{1}))^{-1} u g_{+} u f(X_{1},(Y_{1},X_{1})),  
$$
where 
$$
u:= f(\sigma_{1}^{2},\sigma_{2}^{2})\sigma_{1}^{\lambda}
f(\sigma_{1}^{2},\sigma_{2}^{2})^{-1}, \quad
g_{+}:= f(X_{1},(Y_{1},X_{1}))X_{1}^{2\lambda}f(X_{1},(Y_{1},X_{1}))^{-1}
$$
(elements of $B_{1,3}$). Conjugating by $Y_{1}$ and using 
$Y_{1}X_{2}Y_{1}^{-1} = X_{2}\sigma_{1}^{-2}$, 
$Y_{1}\sigma_{1}^{2}Y_{1}^{-1} = \sigma_{1}^{2}$, one obtains 
$$
f(X_{2}\sigma_{1}^{-2},\sigma_{2}^{2})(X_{2}\sigma_{1}^{-2})^{\lambda}
f(X_{2}\sigma_{1}^{-2},\sigma_{2}^{2})^{-1} = 
Y_{1} f(X_{1},(Y_{1},X_{1}))^{-1} ug_{+}u f(X_{1},(Y_{1},X_{1})) 
Y_{1}^{-1}.  
$$ 
As $X_{2}\sigma_{1}^{-2} = \sigma_{1}X_{1}\sigma_{1}^{-1}$, 
the left side of this identity identifies with the right side of (\ref{iben:bis}). 
Combining these identities, one gets 
$$
\on{Ad}\Big( f((Y_{1},X_{1}),Y_{1}X_{1}^{-1}Y_{1}^{-1})
(Y_{1},X_{1})^{{1-\lambda}\over 2}\Big)(ug_{+}u^{-1}) = 
\on{Ad}\Big( Y_{1}f(X_{1},(Y_{1},X_{1}))^{-1}\Big)(ug_{+}u), 
$$
which gives after conjugation 
$$
ug_{+}u^{-1} = g_{-}ug_{+}ug_{-}^{-1}, 
$$
where $g_{-}:= g_{-}(X_{1},Y_{1})$, which is equivalent to 
$u^{2} = (ug_{+}^{-1}u^{-1},g_{-}^{-1})$, so the pair $(g_{+},g_{-})$ 
defined in the statement of the Proposition satisfies (\ref{def:GTell:3}). 

4) The fact that $\ul{\on{GT}}\to\ul{\on{GT}}_{ell}$, $(\lambda,f)
\mapsto (\lambda,f,g_\pm)$ is a morphism of semigroups follows
from the identity $(g_-(X,Y),g_+(X,Y))=(Y,X)^\lambda$. It is 
straightforward to check the commutativity of the first diagram; 
the second diagram follows from $((\lambda,f)\in\ul{\on{GT}}(\kk))
\Rightarrow (\on{log}f\in [\hat\f_2^\kk,\hat\f_2^\kk])$. 

5) The arguments used in the case of $\ul{\on{GT}}_{(ell)}$
extend {\it mutatis} to their profinite, pro-$l$ and prounipotent 
versions. \hfill \qed\medskip 

\begin{remark}
There are compatible group morphisms ${\on{GT}}
\to\on{Aut}(R_{ell})$, ${\on{GT}}_l\to\on{Aut}(R_{l}^{ell})$, 
${\on{GT}}(\kk)\to\on{Aut}(R_{ell}(\kk))$ (where 
$R_l^{ell} = \on{Ker}(\ul{\on{GT}}^{ell}_l \to \ul{\on{GT}}_l)$), 
defined by $(\lambda,f)\mapsto\theta_{\lambda,f}:=$ conjugation 
by the image of $(\lambda,f)\mapsto(\lambda,f,g_\pm)$ from 
Proposition \ref{prop:lift}. One computes  
$\theta_{\lambda,f}(\Psi_+) = \Psi_+^{1/\lambda}$
and $\theta_{\lambda,f}((\Psi_+\Psi_-)^3)
=((\Psi_+\Psi_-)^{3/\lambda}$, where 
$(\Psi_+\Psi_-)^3$ is a generator of 
$Z(B_3)=\ZZ$ and $(\Psi_{+}\Psi_{-})^{3(1+2m)} = 
(\Psi_{+}\Psi_{-})^{3}(\Psi_{+}\Psi_{-})^{6m} = 
(\Psi_{+}\Psi_{-})^{3}\on{Ad}(Y,X)^{m}$. \hfill \qed\medskip 
\end{remark}

The semigroup scheme morphism from Proposition \ref{prop:lift} 
restricts to a group scheme morphism, which yields an action of 
$\on{GT}(-)$ on $R_{ell}(-)$. The group scheme $\on{GT}_{ell}(-)$
has then a semidirect product structure, fitting in the 
diagram 
$$
\begin{matrix}
\on{GT}_{ell}(-) & \simeq & R_{ell}(-)\rtimes \on{GT}(-)\\
\downarrow & & \downarrow \\
\on{GL}_2 &\simeq  & \on{SL}_2\rtimes {\mathbb G}_m 
\end{matrix}
$$
where the bottom map is induced by ${\mathbb G}_m\to\on{GL}_2$, 
$\lambda\mapsto \bigl( \begin{smallmatrix} \lambda & 0 \\ 0 & 1
\end{smallmatrix}\bigr)$.

\section{Elliptic associators} \label{sec:3}

In this section, we introduce the notion of elliptic associator. This notion yields 
particular elliptic structures over BMCs. It gives rise to a scheme of elliptic 
associators, which appears to be a torsor under the 
action of the group scheme $\on{GT}_{ell}(-)$. We construct a morphism of torsors
from the scheme of associators to its elliptic analogue, which enables us to establish 
the existence of rational elliptic associators. 

\subsection{Lie algebras $\t_n$ and $\t_{1,n}$}

Let $\kk$ be a $\QQ$-ring. 
If $S$ is a finite set, we define $\t_S^\kk$ as the $\kk$-Lie 
algebra with generators $t_{ij}$, $i\neq j\in S$ and relations 
$t_{ji}=t_{ij}$, $[t_{ij},t_{ik}+t_{jk}]=0$ for $i,j,k$
distinct, $[t_{ij},t_{kl}]=0$ for $i,j,k,l$ distinct. 
We define $\hat\t_S^\kk$ as its degree completion, where
$\on{deg}(t_{ij})=1$. 

For $S'\supset D_\phi\stackrel{\phi}{\to}S$ a partially defined 
map, there is a unique Lie algebra morphism $\t_S^\kk\to 
\t_{S'}^\kk$, $x\mapsto x^\phi$, defined by 
$(t_{ij})^\phi:= \sum_{i'\in\phi^{-1}(i),j'\in\phi^{-1}(j)}
t_{i'j'}$. Then $S\mapsto \t_S^\kk$ is a contravariant 
functor $($finite sets, partially defined maps$)\to\{$Lie algebras$\}$. 

We also define $\t_{1,S}^\kk$ as the $\kk$-Lie algebra 
with generators $x_i^\pm$, $i\in S$ and relations 
$\sum_{i\in S}x_i^\pm=0$, $[x_i^\pm,x_j^\pm]=0$ for 
$i\neq j$, $[x_i^+,x_j^-]=[x_j^+,x_i^-]$
for $i\neq j$, $[x_k^\pm,[x_i^+,x_j^-]]=0$ for $i,j,k$
distinct. We then have a Lie algebra morphism 
$\t_S^\kk\to\t_{1,S}^\kk$, $t_{ij}\mapsto 
[x_i^+,x_j^-]$, which we denote by $x\mapsto \{x\}$. 
We will also write $t_{ij} = [x_i^+,x_j^-]$. 
We define $\hat\t_{1,S}^\kk$ as the degree completion of 
$\t_{1,S}^\kk$, where $\on{deg}(x_i^\pm)=1$. 

For $S'\stackrel{\phi}{\to}S$ a map, there is a unique 
Lie algebra morphism $\t_{1,S}^\kk\to \t_{1,S'}^\kk$, 
$x\mapsto x^\phi$, such that $(x_i^\pm)^\phi:= 
\sum_{i'\in \phi^{-1}(i)} x^\pm_{i'}$. Then $S\mapsto 
\t_{1,S}^{\kk}$ is a contravariant functor $($finite sets, maps$)
\to\{$Lie algebras$\}$. By restriction, $S\mapsto \t_S^{\kk}$ may be
viewed as a contravariant functor of the same type, and the morphism
$\t_S^{\kk}\to\t_{1,S}^{\kk}$ is then functorial; i.e., we have 
$\{x\}^\phi = \{x^\phi\}$ for $x\in\t_S$ and any map   
$S'\stackrel{\phi}{\to} S$. 

We set $\t_n^\kk:= \t_{[n]}^\kk$, $\t_{1,n}:= 
\t_{1,[n]}^\kk$, where $[n] = \{1,\ldots,n\}$, and we write 
$x^\phi$ as $x^{I_1,...,I_n}$, where $I_i=\phi^{-1}(i)$ 
for $x\in\t_n^\kk$ or $x\in\t_{1,n}^\kk$. 

\subsection{Elliptic associators} \label{sec:ell:ass}

Recall that the set $\ul{M}(\kk)$ of associators defined 
over $\kk$ is the set of $(\mu,\Phi)\in\kk\times \on{exp}(\hat\f_2^\kk)$, 
such that $\Phi^{3,2,1}=\Phi^{-1}$, 
\begin{equation} \label{hex:Phi}
e^{\mu t_{23}/2}\Phi^{1,2,3}e^{\mu t_{12}/2}\Phi^{3,1,2}
e^{\mu t_{31}/2}\Phi^{2,3,1}=e^{\mu(t_{12}+t_{13}+t_{23})/2}, 
\end{equation}
\begin{equation} \label{pent}
\Phi^{2,3,4}\Phi^{1,23,4}\Phi^{1,2,3}=
\Phi^{1,2,34}\Phi^{12,3,4}, 
\end{equation}
where $\Phi$ is viewed as an element of $\on{exp}(\hat\t_{3}^{\kk})$ via 
the inclusion $\hat\f_2^\kk\subset \hat\t_3^\kk$, $A,B\mapsto t_{12},t_{23}$. 

\begin{definition}
The set $\ul{Ell}(\kk)$ of elliptic associators defined over $\kk$
is the set of quadruples $(\mu,\Phi,A_+,A_{-})$, 
where $(\mu,\Phi)\in \ul{M}(\kk)$ and $A_\pm\in 
\on{exp}(\hat{\t}^\kk_{1,2})$, such that: 
\begin{equation} \label{def:ell:ass:1}
\alpha_\pm^{3,1,2}\alpha_\pm^{2,3,1}\alpha_\pm^{1,2,3} = 1, 
\text{\ where\ }\alpha_\pm = \{e^{\pm\mu(t_{12}+t_{13})/2}\}
A_\pm^{1,23}\{\Phi^{1,2,3}\}, 
\end{equation}
\begin{equation} \label{def:ell:ass:2}
\{e^{\mu t_{12}}\} = \big(\{\Phi\}^{-1}A_-^{1,23}\{\Phi\},
\{e^{-\mu t_{12}/2}(\Phi^{2,1,3})^{-1}\}(A_+^{2,13})^{-1}
\{\Phi^{2,1,3}e^{-\mu t_{12}/2}\} \big).  
\end{equation}
\end{definition}

\begin{remark}
We then have $\{e^{\pm\mu t_{12}/2}\}A_\pm^{2,1}
\{e^{\pm\mu t_{12}/2}\}A_\pm^{1,2}=1$ and 
$\{e^{\mu t_{12}}\} = (A_-,A_+)$; here as in (\ref{def:ell:ass:2}),
the notation $(g,h)$ stands for the group commutator $ghg^{-1}h^{-1}$. 
\end{remark}

Then $\kk\mapsto \ul{M}(\kk),\ul{Ell}(\kk)$ are 
functors $\{\QQ$-rings$\}\to\{$sets$\}$, i.e., $\QQ$-schemes. 
We have an obvious scheme morphism $\ul{Ell}\to\ul{M}$, 
$(\mu,\Phi,A_+,A_{-})\mapsto(\mu,\Phi)$. 

Define also a scheme morphism $\ul{Ell}\to \on{M}_2$
by $(\mu,\Phi,A_+,A_{-})\mapsto 
\bigl( \begin{smallmatrix} u_+ & v_+ \\ u_- & v_-
\end{smallmatrix}\bigr)$, where $u_{\pm},v_{\pm}$ are the 
coefficients arising from $\on{log}A_\pm 
\equiv u_\pm x_1^+ + v_\pm x_1^-$ mod $[\hat\t_{1,2},
\hat\t_{1,2}]$. Then relation (\ref{def:ell:ass:2}) implies that the ldiagram 
$$
\begin{matrix}
\ul{Ell} & \to & \ul{M} \\
\downarrow  & & \downarrow \\
\on{M}_2 & \stackrel{\on{det}}{\to}& {\mathbb A}
\end{matrix}
$$
commutes. 

\subsection{Categorical interpretations}

\begin{definition} (see \cite{Dr:Gal})
An infinitesimally braided monoidal category (IBMC) over $\kk$
is a set $(\cC,\otimes,c_{\ldots},a_{\ldots},U_{\ldots},t_{\ldots})$, such 
that: 

1) $(\cC,\otimes,c_{\ldots},a_{\ldots})$ is a symmetric  
monoidal category (i.e., $c_{Y,X}c_{X,Y}=\on{id}_{X\otimes Y}$); 

2) $\on{Ob}\cC\ni X\mapsto U_X\lhd \on{Aut}_{\cC}(X)$ 
is such that $U_X$ is a $\kk$-prounipotent group, 
and $i U_X i^{-1}=U_Y$ for any $i\in \on{Iso}_{\cC}(X,Y)$;  

3) $(\on{Ob}\cC)^2\ni (X,Y)\mapsto t_{X,Y}\in
\on{Lie}U_{X\otimes Y}$ is a natural assignment; 

4) $t_{Y,X} = c_{X,Y}t_{X,Y} c_{X,Y}^{-1}$ and 
$$
t_{X\otimes Y,Z} = a_{X,Y,Z}(\on{id}_X\otimes
t_{Y,Z})a_{X,Y,Z}^{-1}
+ (c_{Y,X}\otimes\on{id}_Z)a_{Y,X,Z}
(\on{id}_Y\otimes t_{X,Z})
((c_{Y,X}\otimes\on{id}_Z)a_{Y,X,Z})^{-1}. 
$$
\end{definition}

A functor $f:\cC\to\cC'$ between IBMCs is then a tensor functor,  
such that $f(U_X)\subset U'_{f(X)}$ and $f(t_{X,Y})=t'_{f(X),f(Y)}$. 
An example of IMBC is constructed as follows: $\cC = {\bf PaCD}$
is the category with the same objects as ${\bf PaB}$, 
${\bf PaCD}(O,O') := 
  \left\{ \begin{array}{ll}
         \on{exp}(\hat\t_{|O|})\rtimes S_{|O|} & \on{\ if\ } |O|=|O'|\\
        \emptyset & \on{otherwise}\end{array} \right.$, 
$c_{O,O'} = s_{|O|,|O'|}\in S_{|O|+|O'|}
\subset \on{Aut}_{{\bf PaCD}}(O\otimes O')$ is the permutation 
$i\mapsto i+|O'|$ for $i\in[1,|O|]$, $i\mapsto i-|O|$ for $i\in[|O|+1,|O|+|O'|]$, 
$a_{O,O',O''}:= 1$, $U_O = \on{exp}(\hat\t_{|O|}^\kk)\lhd \on{Aut}_{\cC}(O)$,
$t_{O,O'}:= \sum_{i=1}^{|O|}\sum_{i'=|O|+1}^{|O|+|O'|}t_{ii'}$. The pair
$({\bf PaCD},\bullet)$ is initial among pairs (an IBMC, a distinguished
object). 

We then set: 

\begin{definition}
An elliptic structure over the IBMC $\cC$ is a set 
$(\tilde\cC,F,\tilde U_{\ldots},x^\pm_{\ldots})$, where 
$\tilde\cC$ is a category, $F:\cC\to\tilde\cC$ is a functor, 
$\on{Ob}\tilde\cC\ni \tilde X\mapsto\tilde U_{\tilde X}
\lhd\on{Aut}_{\tilde\cC}(\tilde X)$
is the assignment of a $\kk$-prounipotent group, where 
$\tilde i\tilde U_{\tilde X}\tilde i^{-1}=\tilde U_{\tilde Y}$ for 
$\tilde i\in \on{Iso}_{\tilde\cC}(\tilde X,\tilde Y)$ and 
$F(U_X)\subset \tilde U_{F(X)}$, 
and $(\on{Ob}\cC)^2\ni(X,Y)\mapsto x^\pm_{X,Y}
\in \on{Lie}\tilde U_{F(X\otimes Y)}$ is a natural assignment,
such that 
$$
x^\pm_{Y,X} = F(c_{X,Y})x^\pm_{X,Y} F(c_{X,Y}^{-1}), 
\quad x^\pm_{X,{\bf 1}}=0,
$$ 
\begin{align*}
& x^\pm_{X\otimes Y,Z} + F (c_{X,Y\otimes Z}a_{X,Y,Z})^{-1}
x^{\pm}_{Y\otimes Z,X} F(c_{X,Y\otimes Z}a_{X,Y,Z})  
\\
 & + F(a^{-1}_{Z,X,Y}c_{X\otimes Y,Z})^{-1}
x^{\pm}_{Z\otimes X,Y} F(a^{-1}_{Z,X,Y}c_{X\otimes Y,Z}) = 0, 
\end{align*}
$$
F(t_{X,Y}\otimes\on{id}_Z) = [F(a_{X,Y,Z})^{-1}x^+_{X,Y\otimes Z}
F(a_{X,Y,Z}),
F((c_{X,Y}\otimes \on{id}_Z)^{-1}a_{Y,X,Z}) x^-_{Y,X\otimes Z}
F(a_{Y,X,Z}^{-1}(c_{X,Y}\otimes \on{id}_Z))]. 
$$
\end{definition}
Functors between pairs (an IBMC, an elliptic structure over it)
are defined in an obvious way. An elliptic structure over 
${\bf PaCD}$ is defined as follows: $\tilde\cC:= {\bf PaCD}_{ell}$
is the category with the same objects as ${\bf PaB}$, 
${\bf PaCD}_{ell}(O,O') := 
  \left\{ \begin{array}{ll}
         \on{exp}(\hat\t^\kk_{1,|O|})\rtimes S_{|O|} 
& \on{\ if\ } |O|=|O'|\\
        \emptyset & \on{otherwise,}\end{array} \right.$
$\tilde U_O = \on{exp}(\hat\t_{1,|O|}) \lhd 
\on{Aut}_{{\bf PaCD}_{ell}}(O)$, $F$ is induced by the morphism
$\t_n\to\t_{1,n}$, $x\mapsto \{x\}$ and the identity between 
symmetric groups, $x^\pm_{O,O'} = \sum_{i=1}^{|O|} x^\pm_i
\in \on{Lie}\tilde U_{O\otimes O'}$.  The triple $({\bf PaCD},
{\bf PaCD}_{ell},\bullet)$ is universal for triples (an IBMC,
an elliptic structure over it, a distinguished object). 

Let us say that a $\kk$-BMC is a braided monoidal category (BMC)
$\cC$, such that the image of each morphism $P_n\to
\on{Aut}_{\cC}(X_1\otimes\cdots\otimes X_n)$ is contained in 
a $\kk$-prounipotent group. Then each $(\mu,\Phi)\in\ul{M}(\kk)$
gives rise to a map $\{$IBMCs$\}\to\{\kk$-BMCs$\}$, 
$\cC\mapsto (\mu,\Phi)*\cC$, where $(\mu,\Phi)*\cC
:= (\cC,\otimes,\beta_{X,Y}:= c_{X,Y}e^{\mu t_{X,Y}/2},
\tilde a_{X,Y,Z}:= \Phi(a_{X,Y,Z}(t_{X,Y}\otimes\on{id}_Z)a_{X,Y,Z}^{-1},
\on{id}_X\otimes t_{Y,Z})a_{X,Y,Z})$. 

In the same say, a $\kk$-elliptic structure over a $\kk$-BMC
is an elliptic structure, such that the image of each 
morphism $P_{1,n}\to
\on{Aut}_{\tilde\cC}(F(X_1\otimes\cdots\otimes X_n))$ is contained in 
a $\kk$-prounipotent group. Then each $(\mu,\Phi,A_+,A_{-})
\in\ul{Ell}(\kk)$ gives rise to a map $\{$(an IBMC, an elliptic
structure over it)$\}\to\{($a $\kk$-BMC, an elliptic structure 
over it)$\}$, $(\cC,\tilde\cC)\mapsto (\mu,\Phi,A_+,A_{-})*(\cC,\tilde\cC)
= (\cC',\tilde\cC')$, where $\cC'=(\mu,\Phi)*\cC$
and $\tilde\cC '= (\tilde\cC,F,\tilde A^{+}_{X,Y},\tilde A^{-}_{X,Y})$, 
where $\tilde A^\pm_{X,Y}:= A_\pm(x^+_{X,Y},x^-_{X,Y})$.  

\subsection{Action of ${\ul{\on{GT}}}_{ell}(-)$ on ${\ul{Ell}}$}

Recall first that there is an action of $\ul{\on{GT}}(\kk)$
on $\ul{M}(\kk)$, defined by 
$$
(\lambda,f)*(\mu,\Phi):= (\lambda\mu,\Phi(A,B)f(e^{\mu A},\Phi(A,B)^{-1}
e^{\mu B}\Phi(A,B)))=(\mu',\Phi'). 
$$
For $(\lambda,f,g_+,g_{-})\in\ul{\on{GT}}_{ell}(\kk)$ and $(\mu,\Phi,A_+,A_{-})
\in\ul{Ell}(\kk)$, we set 
$$
(\lambda,f,g_+,g_{-})*(\mu,\Phi,A_+,A_{-}):= (\mu',\Phi',A'_+,A'_{-})
$$
where $A'_\pm := g_\pm(A_+,A_-)$. 

\begin{proposition}
This defines an action of $\ul{\on{GT}}_{ell}(\kk)$ on 
$\ul{Ell}(\kk)$. 
\end{proposition}

{\em Proof.} For $g_{ell}\in \ul{\on{GT}}_{ell}(\kk)$, 
and $(\cC,\tilde\cC)\in \{($a $\kk$-BMC, an elliptic structure over it$)\}$, 
we have $g_{ell}*((\mu,\Phi,A_+,A_{-})*(\cC,\tilde\cC)) = 
(g_{ell}*(\mu,\Phi,A_\pm))*(\cC,\tilde\cC)$. When 
$(\cC,\tilde\cC)=({\bf PaCD},{\bf PaCD}_{ell})$, 
$(\mu,\Phi,A_+,A_{-})$ can be recovered uniquely from 
$(\mu,\Phi,A_+,A_{-})*(\cC,\tilde\cC)$, as
$e^{\mu t_{12}}=\beta_{\bullet,\bullet}^2$, 
$\Phi=\tilde a_{\bullet,\bullet,\bullet}$,
and $A_\pm = A^\pm_{\bullet,\bullet}$, which implies that 
the above formula defines an action.  
\hfill \qed\medskip 

\begin{remark}
The actions of $\ul{\on{GT}}(\kk)$ on $\{\kk$-BMCs$\}$
and on $\ul{M}(\kk)$ are compatible, 
in the sense that for $g\in\ul{\on{GT}}(\kk)$, 
$g*((\mu,\Phi)*\cC) = (g*(\mu,\Phi))*\cC$. 
In the same way, the actions of $\ul{\on{GT}}_{ell}(\kk)$ on 
$\{($a $\kk$-BMC, an elliptic structure over it$)\}$
and on $\ul{Ell}(\kk)$ are compatible. 
\end{remark}

\begin{remark} The morphism $\ul{Ell}\to\on{M}_2$
from Section \ref{sec:ell:ass} is compatible with the semigroup 
scheme morphism $\ul{\on{GT}}_{ell}(-)\to\on{M}_2$
from Proposition \ref{prop15}, with the action of $\ul{\on{GT}}_{ell}$
on $\ul{Ell}$, and with the left multiplication action of 
$\on{M}_2$ on itself. 
\end{remark}

\subsection{A morphism ${\ul{M}}\to {\ul{Ell}}$}

The scheme morphism $\ul{Ell}\to\ul{M}$, $(\mu,\Phi,A_{+},A_-)
\to (\mu,\Phi)$ is clearly compatible with the semigroup 
scheme morphism $\ul{\on{GT}}_{ell}(-)\to\ul{\on{GT}}(-)$. 
We now construct a section of this morphism. 

\begin{proposition} \label{prop:sigma}
There is a unique scheme morphism $\sigma:\ul{M}\to\ul{Ell}$, 
$(\mu,\Phi)\to(\mu,\Phi,A_+,A_{-})$, where 
$$
A_+ := \Phi({{\on{ad}x_1}\over{e^{\on{ad}x_1}-1}}(y_2),t_{12})
\cdot e^{\mu{{\on{ad}x_1}\over{e^{\on{ad}x_1}-1}}(y_2)}
\cdot \Phi({{\on{ad}x_1}\over{e^{\on{ad}x_1}-1}}(y_2),t_{12})^{-1}, 
$$
$$
A_- := e^{\mu t_{12}/2}
\Phi({{\on{ad}x_2}\over{e^{\on{ad}x_2}-1}}(y_1),t_{21})e^{x_1}
\Phi({{\on{ad}x_1}\over{e^{\on{ad}x_1}-1}}(y_2),t_{12})^{-1} 
$$
(we set $x_i:=x_i^+,y_i:=x_i^-$). It is compatible with 
the semigroup scheme morphism $\ul{\on{GT}}(-)\to
\ul{\on{GT}}_{ell}(-)$ from Proposition \ref{prop:lift}. 
\end{proposition}

One checks that $\sigma$ fits in a diagram 
$$
\begin{matrix}
\ul{M} & \stackrel{\sigma}{\to}& \ul{Ell}\\
\downarrow  & & \downarrow \\
{\mathbb A} & \to & \on{M}_2 
 \end{matrix}
$$
where the bottom map is $\mu\mapsto \bigl( \begin{smallmatrix} 0 & -\mu 
\\ 1 & 0 \end{smallmatrix}\bigr)$. This diagram is compatible with the last 
diagram of Proposition \ref{prop:lift}. 

{\em Proof.} By \cite{CEE}, Prop. 5.3, $(\mu,\Phi,A_+,A_{-})$
satisfies 
$$
A_\pm^{12,3} = \{e^{\pm\mu t_{12}/2}(\Phi^{-1})^{2,1,3}\}
A_\pm^{2,13}\{\Phi^{2,1,3}e^{\pm\mu t_{12}/2}\Phi^{-1}\}
A_\pm^{1,23}\{\Phi\}, 
$$
and therefore (\ref{def:ell:ass:1}). 

The last identity of {\it loc.~cit.}~can be rewritten as follows
(using the commutation of $\{t_{12}\}$ with $A_+^{12,3}$)
$$
A_-^{2,13}\{\Phi^{2,1,3}\} A_+^{12,3} \{(\Phi^{2,1,3})^{-1}\}
(A_-^{2,13})^{-1} = \{(\Phi^{3,1,2})^{-1}e^{\mu t_{12}/2}\Phi^{3,2,1}
e^{\mu t_{23}} \Phi^{1,2,3}e^{-\mu t_{12}/2}\} A_+^{12,3}\{\Phi^{3,1,2}\}.
$$
Now the hexagon and duality identities imply 
\begin{equation} \label{CEE'}
(\Phi^{3,2,1})^{-1} e^{\mu t_{12}/2} \Phi^{3,2,1} e^{\mu t_{23}}
\Phi^{1,2,3} e^{-\mu t_{12}/2} = e^{-\mu t_{13}/2} \Phi^{2,3,1}
e^{\mu t_{23}/2} (\Phi^{3,2,1})^{-1} e^{\mu t_{3,12}/2}, 
\end{equation}
$$
\Phi^{3,1,2} = e^{\mu t_{3,21}/2} \Phi^{3,2,1} e^{-\mu t_{23}/2}
(\Phi^{2,3,1})^{-1} e^{-\mu t_{13}/2}, 
$$
and 
$$
\Phi^{2,1,3} = e^{\mp\mu t_{13}/2} \Phi^{2,3,1} e^{\mp\mu t_{23}/2}
(\Phi^{3,2,1})^{-1} e^{\pm\mu t_{3,12}/2}, 
$$
so (\ref{CEE'}) is rewritten (using the commutation of $\{t_{13}\}$
with $A_-^{2,13}$)
\begin{eqnarray} \label{CEE''}
\nonumber
\{e^{-\mu t_{13}/2}\} A_-^{2,13} \{\Phi^{2,3,1}e^{-\mu t_{23}/2}
(\Phi^{3,2,1})^{-1} e^{\mu t_{3,12}/2}\} A_+^{12,3}  \{
e^{\mu t_{3,12}/2}\Phi^{3,2,1} e^{-\mu t_{23}/2}(\Phi^{2,3,1})^{-1}\} &&
\\  (A_-^{2,13})^{-1} \{e^{-\mu t_{13}/2}\} && 
\\ \nonumber 
= \{e^{-\mu t_{13}/2}\Phi^{2,3,1}e^{\mu t_{23}/2}(\Phi^{3,2,1})^{-1}
e^{\mu t_{3,21}/2}\} A_+^{12,3} \{e^{\mu t_{3,21}/2}\Phi^{3,2,1}
e^{-\mu t_{23}/2}(\Phi^{2,3,1})^{-1}e^{-\mu t_{13}/2}\}. &&  
\end{eqnarray} 
As $A_+^{2,1} e^{\mu t_{12}/2} A_+ e^{\mu t_{12}/2}=1$, 
we have $e^{\mu t_{3,12}/2} A_+^{12,3} e^{\mu t_{3,12}/2}
=(A_+^{3,12})^{-1}$; using this identity and performing the 
transformation of indices $(1,2,3)\to(3,1,2)$, (\ref{CEE''}) 
yields (\ref{def:ell:ass:2}). So $(\mu,\Phi,A_+,A_{-})\in\ul{Ell}(\kk)$. 
The compatibility of $\sigma:\ul{M}(\kk)\to\ul{Ell}(\kk)$ with the 
semigroup morphism $\ul{\on{GT}}(\kk)\to\ul{\on{GT}}_{ell}(\kk)$
follows from $(A_-,A_+) = e^{\mu t_{12}}$. 
\hfill \qed\medskip 

\subsection{A subscheme $Ell\subset{\underline{Ell}}$  
and its torsor structure under $\on{GT}_{ell}(-)$}

Set $M(\kk) := \{(\mu,\Phi) | \mu\in\kk^\times\}\subset \ul{M}(\kk)$
and $Ell(\kk):= \{(\mu,\Phi,A_\pm) | \mu\in\kk^\times\}\subset 
\ul{Ell}(\kk)$. The actions of $\ul{\on{GT}}_{(ell)}$ restrict
to actions of $\on{GT}(\kk)$ on $M(\kk)$ and $\on{GT}_{ell}(\kk)$
on $Ell(\kk)$. Recall that $M(\kk)$ is a principal homogeneous 
space under the 
action of $\on{GT}(\kk)$ (\cite{Dr:Gal}). Similarly: 

\begin{proposition}
$Ell(\kk)$ is empty or a principal homogeneous space 
under the action of $\on{GT}_{ell}(\kk)$. 
\end{proposition}

{\em Proof.} Assume $Ell(\kk)\neq\emptyset$ and let us show 
that the action is free. If 
$(\lambda,f,g_+,g_{-})\in\on{Stab}(\mu,\Phi,A_+,A_{-})$, then 
by the freeness of the action of $\on{GT}(\kk)$ on 
$M(\kk)$, $(\lambda,f)=1$. Then $A_\pm = g_\pm(A_+,A_-)$. 
Relation (\ref{def:ell:ass:2}) implies that if $a_\pm,b_\pm\in\kk$
are such that $\on{log}A_\pm\equiv a_\pm x_1^++b_\pm x_1^-$
mod degree $\geq 2$ (where $x_1^\pm$ have degree 1), then 
$a_+b_--a_-b_+=\mu$, which implies that $(\on{log}A_+,\on{log}A_-)$
generate $\hat{\t}_{1,2}^\kk$, and therefore that $g_\pm=1$. 

We now prove that the action is transitive. As the action of 
$\on{GT}(\kk)$ on $M(\kk)$ is transitive, and as 
$\on{GT}_{ell}(\kk)\to\on{GT}(\kk)$ is surjective (as the morphism 
defined in Proposition \ref{prop:lift} restricts to a section of it),
it suffices to prove that for any $(\mu,\Phi)\in Ell(\kk)$, the 
action of $R_{ell}(\kk)$ on $\{(A_+,A_-)|(\mu,\Phi,A_+,A_{-})\in
Ell(\kk)\}$ is transitive. If $(A_+,A_-)$ and $(A'_+,A'_-)$
belong to this set, then there is a unique $(g_+,g_-)\in F_2(\kk)^2
\simeq P_{1,2}(\kk)^2$ such that $A'_\pm=g_\pm(A_+,A_-)$.  
Then
$$
\alpha_\pm^{1,2,3} \alpha_\pm^{3,1,2}
\alpha_\pm^{2,3,1}=1, \text{\ where\ }
\alpha_\pm = g_\pm(A_+^{1,23},A_-^{1,23})\{\Phi^{1,2,3}
e^{\pm\mu t_{12,3}/2}\}. 
$$
The canonical morphism $B_{1,3}\to \on{Aut}_{(\mu,\Phi,A_+,A_{-})
*{\bf PaCD}}(\bullet(\bullet\bullet))=\on{exp}(\hat\t_{1,3}^\kk)
\rtimes S_3$ extends to an isomorphism $B_{1,3}(\kk)
\simeq \on{exp}(\hat\t_{1,3}^\kk)\rtimes S_3$, given by 
$X_1^\pm\mapsto A_\pm^{1,23}$, $\sigma_1\mapsto \{\Phi e^{\mu t_{12}/2}\}
(12)\{\Phi\}^{-1}$, $\sigma_2\mapsto \{e^{\mu t_{23}/2}\}(23)$. 
It is such that $\sigma_2^{\pm 1}\sigma_1^{\pm 1}\mapsto 
\{\Phi e^{\pm(\mu/2)t_{3,12}}\}(23)(12)$. The preimage of the above 
identity by this isomorphism then yields $(g_\pm(X_1^+,X_1^-)
\sigma_2^{\pm 1}\sigma_1^{\pm 1})^3=1$. Similarly, the preimage of 
the identity 
$$
\{e^{\mu t_{12}}\} = \big(\{\Phi^{-1}\}
g_-(A_+^{1,23},A_-^{1,23})\{\Phi\},\{e^{-(\mu/2)t_{12}
(\Phi^{2,1,3})^{-1}}\} g_+^{-1}(A_+^{2,13},A_-^{2,13})\{\Phi^{2,1,3}
e^{-(\mu/2)t_{12}}\}\big)
$$ 
yields $\sigma_1^2 =(\sigma_1 g_+^{-1}(X_1^+,X_1^-)
\sigma_1,g_-(X_1^+,X_1^-))$.  \hfill \qed\medskip 

Recall the following definition: 

\begin{definition} \label{def:torsors}
A $\QQ$-torsor is the data of: $\QQ$-group schemes $G,H$, 
a $\QQ$-scheme $X$, commuting left and right actions of $G,H$ on $X$, such that:
for any $\kk$ with $X(\kk)\neq\emptyset$, the action of $G(\kk)$ and $H(\kk)$
on $X(\kk)$ is free and transitive. 
\end{definition}

Morphisms of torsors are then defined in the
obvious way. 

The above $\QQ$-scheme morphisms between $\ul{Ell}$ and $\ul{M}$
restrict to a torsor morphism $Ell\to M$ and a section of it 
$M\stackrel{\sigma}{\to} Ell$, fitting in commutative diagrams
$$
\begin{matrix}
Ell & \to &  M \\
\downarrow & & \downarrow \\
\on{GL}_2 & \stackrel{\on{det}}{\to}& {\mathbb G}_m
\end{matrix}
\text{\ and \ }
\begin{matrix}
M & \to &  Ell \\
\downarrow & & \downarrow \\
{\mathbb G}_m & \stackrel{\mu\mapsto 
\bigl( \begin{smallmatrix} 0 & -1 \\ \mu & 0
\end{smallmatrix}\bigr)}{\to}& \on{GL}_2
\end{matrix}
$$ 

\section{The group $\on{GRT}_{ell}(\kk)$
and isomorphisms of Lie algebras} \label{sec:4}

In this section, we study the group scheme $\on{GRT}_{ell}(-)$ of 
$\on{GT}_{ell}(-)$-automorphisms of the scheme of elliptic associators. We show that 
its Lie algebra ${\mathfrak{grt}}_{ell}$ is graded and equipped with a graded 
morphism ${\mathfrak{grt}}_{ell}\to {\mathfrak{grt}}$. We construct a section 
of this morphism, which brings to light the semidirect product structure of 
${\mathfrak{grt}}_{ell}$. We show that the Lie subalgebra $\mathfrak{sl}_{2}
\subset \on{Der}(\t_{1,2})$ and the derivations $\delta_{2k},k\geq 0$ of 
${\mathfrak t}_{1,2}$ from \cite{CEE} give rise to 
a family of generators of the kernel ${\mathfrak r}_{ell}:= 
\on{Ker}({\mathfrak{grt}}_{ell}\to {\mathfrak{grt}})$. 
The existence of rational elliptic associators enables us to construct an isomorphism 
between the group schemes $\on{GT}_{ell}(-)$ and $\on{GRT}_{ell}(-)$, compatible 
with their semidirect product structures and with their actions on the elliptic 
braid groups and their graded versions. 

\subsection{Reminders about $\on{GRT}(\kk)$}

Let $\kk$ be a $\QQ$-ring. 
Recall (\cite{Dr:Gal}) that $\on{GRT}_1(\kk)$ is defined as the set of all 
$g\in\on{exp}(\hat\f_2^\kk)\subset \on{exp}(\hat\t_3^\kk)$, 
such that: 
$$
g^{3,2,1}=g^{-1}, \quad g^{3,1,2}g^{2,3,1}g^{1,2,3}=1 
\quad \text{(relations in }\on{exp}(\hat\t_{3}^\kk)),
$$
$$
t_{12}+\on{Ad}(g^{1,2,3})^{-1}(t_{23})
+ \on{Ad}(g^{2,1,3})^{-1}(t_{13}) = t_{12}+t_{13}+t_{23}
\quad \text{(relation in }\hat\t_{3}^{\kk}), 
$$
$$
g^{2,3,4}g^{1,23,4}g^{1,2,3}=g^{1,2,34}g^{12,3,4}
\quad \text{(relation in }\on{exp}(\hat\t_{4}^{\kk})). 
$$
This is a group with law $(g_1*g_2)(A,B):= 
g_1(\on{Ad}(g_2(A,B))(A),B)g_2(A,B)$. Note that 
$g\in \on{GRT}_1(\kk)$ gives rise to $\theta_g\in \on{Aut}(\hat\t_3^{\kk})$, 
defined by 
$$
\theta_g : t_{12}\mapsto t_{12}, \quad 
t_{23}\mapsto \on{Ad}(g^{1,2,3})^{-1}(t_{23}), 
\quad 
t_{13}\mapsto \on{Ad}(g^{2,1,3})^{-1}(t_{13}).  
$$
Then $g_1*g_2 = g_1\theta_{g_2}(g_1)$, and 
$\theta_{g_1*g_2} = \theta_{g_2}\theta_{g_1}$, 
so $g\mapsto\theta_g$ is a group antimorphism. 

The group $\kk^\times$ acts on $\on{GRT}_1(\kk)$ by 
$(c\cdot g)(A,B):= g(c^{-1}A,c^{-1}B)$, and 
one sets $\on{GRT}(\kk):= \on{GRT}_1(\kk)\rtimes\kk^\times$. 
$\on{GRT}_1(-)$ is a prounipotent group scheme. 

\subsection{The group $\on{GRT}_{ell}(\kk)$}

Define $\on{GRT}_1^{ell}(\kk)$ as the set of all $(g,u_+,u_{-})$, 
such that $g\in\on{GRT}_1(\kk)$, $u_\pm\in\hat{\t}_{1,2}^\kk$, 
and 
\begin{equation} \label{def:grt:ell:1}
\on{Ad}(g^{1,2,3})^{-1}(u_\pm^{1,23})+
\on{Ad}(g^{2,1,3})^{-1}(u_\pm^{2,13})+u_\pm^{3,12}=0, 
\end{equation}
\begin{equation} \label{def:grt:ell:2}
[\on{Ad}(g^{1,2,3})^{-1}(u_\pm^{1,23}),u_\pm^{3,12}]=0, 
\end{equation}
\begin{equation} \label{def:grt:ell:3}
[\on{Ad}(g^{1,2,3})^{-1}(u_+^{1,23}),
\on{Ad}(g^{2,1,3})^{-1}(u_-^{2,13})]=t_{12}
\end{equation}
(relations in $\hat\t_{1,3}^{\kk}$). 
Set $(g_1,u^1_+,u^{1}_{-})*(g_2,u^2_+,u^{2}_{-}):= (g,u_+,u_{-})$, 
where 
\begin{equation} \label{pdt:upm}
u_\pm(x_1,y_1):= u_\pm^1(u_+^2(x_1,y_1),u_-^2(x_1,y_1))
\end{equation}
(where $\t_{1,2}^\kk$ is viewed as the free Lie algebra 
generated by $x_1,y_1$). 

We first prove: 

\begin{lemma} \label{lemma:charact}
$(g,u_+,u_{-})\in\on{GRT}_1^{ell}(\kk)$ iff there exists an 
automorphism of $\hat\t_{1,3}^\kk$ (henceforth denoted
$\theta_{g,u_\pm}$), such that 
$$
x_1^\pm\mapsto \on{Ad}(g^{1,2,3})^{-1}(u_\pm^{1,23}), \quad
x_2^\pm\mapsto \on{Ad}(g^{2,1,3})^{-1}(u_\pm^{2,13}), \quad
x_3^\pm\mapsto u_\pm^{3,12},
$$
$$
t_{12}\mapsto t_{12}, \quad 
t_{23}\mapsto \on{Ad}(g^{1,2,3})^{-1}(t_{23}), \quad
t_{13}\mapsto \on{Ad}(g^{2,1,3})^{-1}(t_{13}). 
$$
\end{lemma}

{\em Proof.} The condition that the relations 
$x_1^\pm+x_2^\pm+x_3^\pm=0$ (resp., 
$[x_1^\pm,x_3^\pm]=0$, $[x_1^+,x_2^-]=t_{12}$) 
are preserved is equivalent to condition
(\ref{def:grt:ell:1}) (resp., 
(\ref{def:grt:ell:2}), (\ref{def:grt:ell:3})), and 
the relation $[t_{12},x_3^\pm]=0$ is automatically preserved. 
Then the relation $g^{3,1,2}g^{2,3,1}g^{1,2,3}=1$ implies that 
$\theta_{g,u_\pm}(x^{2,3,1})=\on{Ad}(g^{1,2,3})^{-1}
(\theta_{g,u_\pm}(x)^{2,3,1})$ for $x\in\{x_i^{\pm},t_{ij}\}$. 
So the other relations $[t_{ij},x_k^\pm]=0$ are also 
preserved. 
\hfill \qed\medskip 

\begin{proposition}
$\on{GRT}_1^{ell}(\kk)$, equipped with the above product, is a group.  
\end{proposition}

{\em Proof.} The product is that of the group 
$\on{GRT}_1(\kk)\times\on{Aut}(\hat{\t}_{1,2}^\kk)^{op}$, 
so it remains to prove that $\on{GRT}_1^{ell}(\kk)$ is stable 
under the operations of product and inverse. If  
$(g_i,u^i_\pm)\in\on{GRT}_1^{ell}(\kk)$ ($i=1,2$), 
then the action of $\theta_{g_2,u_2^\pm}
\theta_{g_1,u_1^\pm}$ on the generators of $\hat\t_{1,3}^\kk$
is given by the formulas of Lemma \ref{lemma:charact}, with
$g= g_1*g_2$ and $u_\pm$ as in (\ref{pdt:upm}). So $(g,u_\pm)
\in\on{GRT}_1^{ell}(\kk)$, as claimed. Similarly, if $(g,u_\pm)
\in\on{GRT}_1^{ell}(\kk)$, then the action of $\theta_{g,u_\pm}^{-1}$
on the generators of $\hat\t_{1,3}^\kk$
is as in Lemma \ref{lemma:charact}, with $(g,u_\pm)$ replaced by 
$($inverse of $g$ in $\on{GRT}_1(\kk)$, inverse of $(u_+,u_-)$ in 
$\on{Aut}(\hat{\t}_{1,2}^\kk))$, so $(g,u_\pm)$ is invertible.  
\hfill \qed\medskip 

In particular, we have 
\begin{equation} \label{pdt:theta}
\theta_{(g_2,u_2^\pm)}\theta_{(g_1,u_1^\pm)} = 
\theta_{(g_1,u_1^\pm)*(g_2,u_2^\pm)}. 
\end{equation}

The assignments $\kk\mapsto \on{GRT}_1(\kk),  \on{GRT}_1^{ell}(\kk)$ 
are then $\QQ$-group schemes. 

For $(g,u_\pm)\in\on{GRT}_1^{ell}(\kk)$, define $a_\pm,b_\pm\in\kk$
by $u_\pm(x_1,y_1) = a_\pm x_1+b_\pm y_1$ mod $[\hat{\t}_{1,2}^\kk,
\hat{\t}_{1,2}^\kk]$. 

\begin{lemma} \label{lemma:sec:SL2}
1) There is a unique group scheme morphism $\on{GRT}_1^{ell}(-)
\to \on{SL}_2$, $(g,u_\pm)\mapsto 
\bigl( \begin{smallmatrix} a_+ & b_+ \\ a_- & b_-
\end{smallmatrix}\bigr)$. 

2) This morphism has a section $\on{SL}_2\to \on{GRT}_1^{ell}(-)$, 
given by $
\bigl( \begin{smallmatrix} a_+ & b_+ \\ a_- & b_-
\end{smallmatrix}\bigr)\mapsto (1,u_\pm(x_1,y_1) = 
a_\pm x_1+b_\pm y_1)$. 
\end{lemma}

{\em Proof.} 1) $a_+b_--a_-b_+=1$ follows 
from (\ref{def:grt:ell:3}); the morphism property is clear. 
2) is straightforward. 
\hfill \qed\medskip 

We now set $\on{GRT}_{I_2}^{ell}(\kk):= \on{Ker}\big(
\on{GRT}_1^{ell}(\kk)\to \on{SL}_2(\kk)\big)$. This defines
a group scheme $\on{GRT}_{I_2}^{ell}(-)$. 

\begin{lemma} 
$\on{GRT}_{I_2}^{ell}(-)$ is a prounipotent group scheme; 
we have $\on{GRT}_1^{ell}(-) = \on{GRT}_{I_2}^{ell}(-) \rtimes\on{SL}_2$. 
\end{lemma}

{\em Proof.} $\on{GRT}_{I_2}^{ell}(\kk)$ is a subgroup of 
$\on{GRT}_1(\kk)\times \on{Ker}\big( \on{Aut}(\hat{\t}_{1,2}^\kk)
\to \on{GL}_2(\kk)\big)$; the assignment 
$\kk\mapsto ($the latter group$)$ is a prounipotent group 
scheme, hence so is $\on{GRT}_{I_2}^{ell}(-)$. 
The second statement follows from Lemma \ref{lemma:sec:SL2}. 
\hfill\qed\medskip 

The group $\kk^\times$ acts on $\on{GRT}_1^{ell}(\kk)$ by 
$c\cdot (g,u_\pm):= (c\cdot g,c\cdot u_\pm)$, where
$c\cdot g$ is as above, $(c\cdot u_+)(x_1^+,x_1^-):= 
u_+(x_1^+,c^{-1}x_1^-)$, $(c\cdot u_-)(x_1^+,x_1^-):= 
cu_-(x_1^+,c^{-1}x_1^-)$. We then set 
$\on{GRT}_{ell}(\kk):= \on{GRT}_1^{ell}(\kk)\rtimes \kk^\times$. 
Then $\kk\mapsto \on{GRT}_{ell}(\kk)$ is a $\QQ$-group scheme, 
and $\on{GRT}_{ell}(-)=\on{GRT}^{ell}_1(-)\rtimes{\mathbb G}_m$. 

There is a unique group scheme morphism 
$\on{GRT}^{ell}_{1}(-)\to \on{GRT}_{1}(-)$, given by 
$(g,u_\pm)\mapsto g$; it extends to a group scheme morphism 
\begin{equation} \label{proj:GRT}
\on{GRT}_{ell}(-)\to\on{GRT}(-),
\end{equation}
whose restriction to  ${\mathbb G}_m$ is the identity. 

To elucidate the structure of $\on{GRT}_{ell}(-)$, we use the following
statement on iterated semidirect products: 

\begin{lemma} \label{lemma:iterated}
Let $G_i$ be groups ($i=1,2,3$). The following data 
are equivalent: 

(a) actions\footnote{The action of $g_j\in G_j$ on $g_i\in G_i$
is denoted $g_j * g_i\in G_i$.} 
$G_j\to\on{Aut}(G_i)$ for $i<j$, such that 
$g_3*(g_2*g_1) = (g_3*g_2)*(g_3*g_1)$; 

(b) actions $G_j\to\on{Aut}(G_i)$ for $(i,j)=(1,2)$ 
and $(2,3)$,  and an action $G_{23}\to\on{Aut}(G_{12})$
(where $G_{ij}:= G_i\rtimes G_j$), compatible with the 
actions of $G_j$ on $G_i$ for $(i,j)=(1,2)$ or $(2,3)$, 
and with the adjoint action of $G_2$ on itself. 

These equivalent data yield actions $G_3\to \on{Aut}(G_{12})$
and $G_{23}\to\on{Aut}(G_1)$, and we then have a canonical 
isomorphism $(G_1\rtimes G_2)\rtimes G_3\simeq G_1\rtimes 
(G_2\rtimes G_3)$. 
\end{lemma}

{\em Proof.} Straightforward. \hfill \qed\medskip 

We then have an action of $\on{GL}_2$ on $\on{GRT}_1^{ell}(-)$, 
given by $\gamma\cdot (g,u_+,u_-):= (\on{det}\gamma\cdot g,\tilde u_+,
\tilde u_-)$, where 
$\bigl( \begin{smallmatrix} \tilde u_+(x_1,y_1) \\ \tilde u_-(x_1,y_1)
\end{smallmatrix}\bigr) := \gamma^{-1} \bigl( \begin{smallmatrix} 
u_+(\tilde x_1,\tilde y_1) \\ u_-(\tilde x_1,\tilde y_1)
\end{smallmatrix}\bigr)$ and 
$\bigl( \begin{smallmatrix} \tilde x_1 \\ \tilde y_1
\end{smallmatrix}\bigr) := \gamma \bigl( \begin{smallmatrix} 
x_1 \\ y_1 \end{smallmatrix}\bigr)$. It satisfies the 
conditions of Lemma \ref{lemma:iterated}, (b), where: 
$G_1 = \on{GRT}_{I_2}^{ell}(-)$, $G_2 = \on{SL}_2$, 
$G_3 = {\mathbb G}_m$, the isomorphism $G_2\rtimes G_3\simeq
\on{GL}_2$ being given by ${\mathbb G}_m\to\on{GL}_2$, 
$c\mapsto \bigl( \begin{smallmatrix} 
1 & 0 \\ 0 & c \end{smallmatrix}\bigr)$. We have therefore
an isomorphism 
$$
\on{GRT}_{ell}(-)\simeq \on{GRT}_{I_2}^{ell}(-)\rtimes\on{GL}_2,  
$$
where we recall that $\on{GRT}_{I_2}^{ell}(-)$ is prounipotent. 

The morphism (\ref{proj:GRT}) then fits in a commutative diagram 
$$
\begin{matrix}
\on{GRT}_{ell}(-) & \to & \on{GRT}(-)\\
\downarrow & & \downarrow \\
\on{GL}_2 & \stackrel{\on{det}}{\to} & {\mathbb G}_m 
\end{matrix}
$$
as the morphism $G_2\rtimes G_3\to G_3$ coincides with $\on{det}$. 

\subsection{A morphism $\on{GRT}(\kk)\to\on{GRT}_{ell}(\kk)$}

We now construct a section of the morphism (\ref{proj:GRT}).  We first set 
\begin{equation} \label{def:ti}
t_{0i}:= - {{\on{ad}x_i}\over{e^{\on{ad}x_i}-1}}(y_i)\in
\hat\t^\kk_{1,n} 
\text{\ for\ }i\in\{1,...,n\}. 
\end{equation}
For $g = g(A,B)\in\on{exp}(\hat\f_2^\kk)$, 
we set $g^{0,1,2}:= g(t_{01},t_{12})\in\on{exp}(\hat\t_{1,2}^\kk)$
$g^{0,2,1}:= g(t_{02},t_{21})\in\on{exp}(\hat\t_{1,2}^\kk)$. 

\begin{lemmadef} \label{def:alpha:g}
For $g\in\on{exp}(\hat\f_2^\kk)$, there exists 
$\alpha_g\in\on{Aut}(\hat\t_{1,2}^\kk)$, uniquely defined by 
$\alpha_g(x_1) = \on{log}(g^{0,2,1}e^{x_1}(g^{0,1,2})^{-1})$, 
$\alpha_g(t_{01}) = g^{0,1,2} t_{01} (g^{0,1,2})^{-1}$. 
We set 
$$(u_+^g,u_-^g):= (\alpha_g(x_1),\alpha_g(y_1))\in (\hat\t_{1,2}^{\kk})^{2}.$$  
\end{lemmadef}

{\em Proof.} This follows from the fact that $\t_{1,2}^\kk$
is freely generated by $x_1$ and $t_{01}$. \hfill \qed\medskip 

\begin{proposition} \label{prop26}
There exists a unique group morphism 
$\on{GRT}_1(\kk)\to \on{GRT}_1^{ell}(\kk)$, given by 
$g\mapsto (g,u_+^g,u_-^g)$. It is 
compatible with the action of $\kk^\times$, hence extends to 
a group morphism $\on{GRT}(\kk)\to \on{GRT}_{ell}(\kk)$, which is a 
section of (\ref{proj:GRT}) and fits in a commutative diagram 
$$
\begin{matrix}
\on{GRT}(-) & \to & \on{GRT}_{ell}(-)\\
\downarrow  & & \downarrow \\
{\mathbb G}_m & \to & \on{GL}_2, 
\end{matrix}
$$
where the bottom morphism is $c\mapsto
\bigl( \begin{smallmatrix} 
1 & 0 \\ 0 & c \end{smallmatrix}\bigr)$. 
\end{proposition}

{\em Proof.}
We first prove: 

\begin{lemma} $\hat\t_{1,n}^\kk$ admits the following presentation: 
generators $x_i,t_{\alpha\beta}$
($i\in\{1,...,n\}$, $\alpha\neq\beta\in\{0,...,n\}$); one sets 
$X_{i}:= e^{x_{i}}$; relations ($i,j,...$ run over $\{1,...,n\}$
while $\alpha,\beta,...$ run over $\{0,\ldots,n\}$): 
\begin{equation} \label{pres:tell:0}
t_{\beta\alpha}=t_{\alpha\beta} \text{\ for\ }
\alpha\neq\beta,\quad 
[t_{\alpha\beta},t_{\gamma\delta}]=[t_{\alpha\beta},
t_{\alpha\gamma}+t_{\beta\gamma}]=0\text{\ for\ }
\alpha,...,\delta\text{\ all\ different,}
\end{equation}
\begin{equation} \label{pres:tell:1}
\on{log}(X_i,X_j) = \on{log}(\prod_i X_i)=0, 
\end{equation}
\begin{equation} \label{pres:tell:2}
X_i(t_{0j}+t_{ij})X_i^{-1}=t_{0j} \text{\ if\ }
i\neq j, \quad
X_i t_{0i} X_i^{-1} = \sum_{\alpha\neq i}t_{\alpha i}, 
\end{equation}
\begin{equation} \label{pres:tell:3}
X_i t_{jk}X_i^{-1}=t_{jk} \text{\ for\ }i,j,k\text{\ distinct}, 
\quad (X_jX_k)t_{jk}(X_jX_k)^{-1}=t_{jk} \text{\ for\ }i\neq j, 
\end{equation}
\begin{equation} \label{pres:tell:4}
\sum_{0\leq\alpha<\beta\leq n}t_{\alpha\beta}=0. 
\end{equation}
\end{lemma}

{\em Proof.} One first checks that if one defines
$t_{0i}$ as in (\ref{def:ti}), then the above relations
are satisfied; conversely, if one sets $y_i:=
-{{e^{\on{ad}x_i}-1}\over{\on{ad}x_i}}(t_{0i})$,
then the above relations lead to the defining relations of 
$\hat\t_{1,n}^\kk$. \hfill \qed\medskip 

\begin{lemma} \label{lemma28}
Let $(g,u_\pm)\in\on{GRT}_1^{ell}(\kk)$
and $\alpha\in\on{Aut}(\hat\t_{1,2}^\kk)$
be defined by $\alpha(x_1^\pm)=u_\pm$. Then $\theta_{g,u_\pm}
\in\on{Aut}(\hat\t_{1,3}^\kk)$ (see Lemma \ref{lemma:charact})
may be defined by 
$$
\theta_{g,u_\pm} : X_1 \mapsto \on{Ad}(g^{1,2,3})^{-1}(\alpha(X_1)^{1,23}), 
\quad 
X_2 \mapsto \on{Ad}(g^{2,1,3})^{-1}(\alpha(X_1)^{2,13}), 
\quad 
X_3 \mapsto \alpha(X_1)^{3,12}, 
$$
$$
t_{01} \mapsto 
\on{Ad}(g^{1,2,3})^{-1}(\alpha(t_{01})^{1,23}), 
\quad 
t_{02} \mapsto \on{Ad}(g^{2,1,3})^{-1}(\alpha(t_{01})^{2,13}), 
\quad 
t_{03} \mapsto \alpha(t_{01})^{3,12}, 
$$
$$
 t_{12} \mapsto t_{12}, \quad 
t_{23} \mapsto \on{Ad}(g^{1,2,3})^{-1}(t_{23}), 
\quad 
t_{13} \mapsto \on{Ad}(g^{2,1,3})^{-1}(t_{13}).  
$$
\end{lemma}

{\em Proof.} Immediate. \hfill \qed\medskip

\begin{lemma} \label{lemma29}
Let $g\in \on{GRT}_1(\kk)$. There is a unique $\tilde\theta_g\in
\on{Aut}(\hat\t_{1,3}^\kk)$, such that
$$
\tilde\theta_g : X_1\mapsto (g^{1,2,3})^{-1}g^{0,23,1} X_1
(g^{0,1,23})^{-1}g^{1,2,3}, \quad 
X_2\mapsto  (g^{2,1,3})^{-1}g^{0,13,2} X_2
(g^{0,2,13})^{-1}g^{2,1,3}, 
$$
$$
X_3\mapsto g^{0,12,3}X_3(g^{0,3,12})^{-1}, 
$$
\begin{equation} \label{t:ij:1}
t_{01}\mapsto \on{Ad}((g^{1,2,3})^{-1}g^{0,1,23})(t_{01}), \quad
t_{02}\mapsto \on{Ad}((g^{2,1,3})^{-1}g^{0,2,13})(t_{02}), \quad
t_{03}\mapsto \on{Ad}(g^{0,3,12})(t_{03}),  
\end{equation}
\begin{equation} \label{t:ij:2}
t_{12}\mapsto t_{12}, \quad t_{23}\mapsto \on{Ad}(g^{1,2,3})^{-1}(t_{23}), 
\quad t_{13}\mapsto \on{Ad}(g^{2,1,3})^{-1}(t_{13}). 
\end{equation}
\end{lemma}

{\em Proof.} Let us first prove that relations (\ref{pres:tell:0}) and 
(\ref{pres:tell:4}) (for $n=3$) are preserved. In Subsection 
\ref{sec:cat}, we will construct an elliptic IBMC $g*{\bf PaCD}$ with distinguished
object $\bullet$, which gives rise to a functor ${\bf PaCD}\to g*{\bf PaCD}$. 
One derives from there an automorphism $\on{exp}(\hat\t_n^\kk)\rtimes
S_n = \on{Aut}_{{\bf PaCD}}(O)\to \on{Aut}_{g*{\bf PaCD}}(O) =
 \on{exp}(\hat\t_n^\kk)\rtimes
S_n$ for any $O\in{\bf PaCD}(O)$, $|O|=n$. When 
$O=\bullet((\bullet\bullet)\bullet)$, the resulting automorphism of 
$\hat\t_4^\kk$ is given by (\ref{t:ij:1}), (\ref{t:ij:2}). So relations 
(\ref{pres:tell:0}) are preserved. The automorphism necessarily preserves
$Z(\hat{\t}_4^\kk) = \kk\cdot \sum_{\alpha<\beta} t_{\alpha\beta}$, 
so relation (\ref{pres:tell:4}) is also preserved. 

Note for later use that 
\begin{equation} \label{simpl}
\tilde\theta_g(x^{2,3,1}) = \on{Ad}(g^{1,2,3})^{-1}
(\tilde\theta_g(x)^{2,3,1})\text{\ for\ }x\in\{x_i,t_{\alpha\beta}\}.  
\end{equation}

We have 
\begin{eqnarray*}
\tilde\theta_g(X_2) \tilde\theta_g(X_3) & 
= & (g^{2,1,3})^{-1} g^{0,13,2} X_2 (g^{0,2,13})^{-1}
g^{2,1,3} g^{0,21,3} X_3 (g^{0,3,12})^{-1}
\\ &  
= &(g^{2,1,3})^{-1} g^{0,13,2} X_2 g^{02,1,3}
(g^{0,2,1})^{-1} X_3 (g^{0,3,12})^{-1}
\\ & = &(g^{2,1,3})^{-1} g^{0,13,2} g^{0,1,3} X_2 X_3 
(g^{03,2,1})^{-1} (g^{0,3,12})^{-1}
\\ & = &(g^{2,1,3})^{-1} (g^{1,3,2})^{-1}
g^{0,1,32} g^{01,3,2} X_2 X_3 
(g^{0,3,2})^{-1} (g^{0,32,1})^{-1}(g^{3,2,1})^{-1}
\\ & = &
(g^{2,1,3})^{-1} (g^{1,3,2})^{-1}
g^{0,1,32} X_2 X_3 
(g^{0,32,1})^{-1}(g^{3,2,1})^{-1}, 
\end{eqnarray*}
while 
\begin{eqnarray*}
\tilde\theta_g(X_3) \tilde\theta_g(X_2) & 
= &g^{0,21,3} X_3 (g^{0,3,12})^{-1}
(g^{2,1,3})^{-1} g^{0,13,2} X_2 (g^{0,2,13})^{-1}
g^{2,1,3} 
\\ &  
= &g^{0,21,3} X_3 g^{03,1,2} (g^{0,3,1})^{-1}
X_2 (g^{0,2,13})^{-1}g^{2,1,3}
\\ & 
= &g^{0,12,3} g^{0,1,2} X_3 X_2 (g^{02,3,1})^{-1} (g^{0,2,13})^{-1}
g^{2,1,3}
\\ & 
= &(g^{1,2,3})^{-1}g^{0,1,23}g^{01,2,3}X_3X_2 (g^{0,2,3})^{-1}(g^{0,23,1})^{-1}
(g^{2,3,1})^{-1}g^{2,1,3}
\\ & 
= &(g^{1,2,3})^{-1}g^{0,1,23}X_3X_2 (g^{0,23,1})^{-1}
(g^{2,3,1})^{-1}g^{2,1,3}, 
\end{eqnarray*}
which implies $(\tilde\theta_g(X_2),\tilde\theta_g(X_3))=1$. Then 
(\ref{simpl}) implies that 
$(\tilde\theta_g(X_i),\tilde\theta_g(X_j))=1$ for any $i,j$. 

The above computation of $\tilde\theta_g(X_2) \tilde\theta_g(X_3)$
implies that 
\begin{eqnarray*}
&&\tilde\theta_g(X_1)\tilde\theta_g(X_2)
\tilde\theta_g(X_3)
\\ && 
= (g^{1,2,3})^{-1} g^{0,23,1} X_1 (g^{0,1,23})^{-1} g^{1,2,3}
(g^{2,1,3})^{-1}(g^{1,3,2})^{-1}X_2X_3(g^{0,32,1})^{-1}
(g^{3,2,1})^{-1}=1
\end{eqnarray*}
as $X_1X_2X_3=1$. So $X_1X_2X_3=1$ is preserved. 

$\tilde\theta_g(X_3)$ clearly commutes with $\tilde\theta_g(t_{12})$, 
which implies that $X_jt_{jk}X_i^{-1}=t_{jk}$ is preserved in 
view of (\ref{simpl}), as well as $X_jX_kt_{jk}(X_jX_k)^{-1}=t_{jk}$
(as the $X_i$ commute and $X_1X_2X_3=1$). 

Now 
\begin{eqnarray*}
&& \tilde\theta_g(t_{02}+t_{12}) = \on{Ad}((g^{2,1,3})^{-1}
g^{0,2,31})(t_{02})+t_{12}
= \on{Ad}(g^{0,21,3})(\on{Ad}(g^{0,2,1})(t_{02})+t_{12})
\\ && 
= \on{Ad}(g^{0,21,3})(t_{01}+t_{02}+t_{12} - \on{Ad}(g^{0,1,2})(t_{01}))
= t_{12} + \on{Ad}(g^{0,21,3})(t_{01}+t_{02}) 
- \on{Ad}(g^{0,12,3}g^{0,1,2})(t_{01}). 
\end{eqnarray*}
Then 
\begin{eqnarray*}
&& \tilde\theta_g(X_1) \tilde\theta_g(t_{02}+t_{12}) \\
 & &  = 
(g^{1,2,3})^{-1}g^{0,23,1}X_1 \Big( (g^{0,1,23})^{-1} g^{1,2,3} t_{12}
+ (g^{0,1,23})^{-1} g^{1,2,3} g^{0,12,3}(t_{01}+t_{02})(g^{0,12,3})^{-1}
\\ && - (g^{0,1,23})^{-1} g^{1,2,3} g^{0,12,3} g^{0,1,2}t_{01}
(g^{0,1,2})^{-1}(g^{0,12,3})^{-1}\Big) 
\\ && 
= (g^{1,2,3})^{-1}g^{0,23,1}X_1 \Big(g^{01,2,3}(g^{0,1,2})^{-1} t_{12}
(g^{0,12,3})^{-1} + g^{01,2,3}(g^{0,1,2})^{-1} (t_{01}+t_{02})
(g^{0,12,3})^{-1} \\ && - g^{01,2,3} t_{01} (g^{0,1,2})^{-1} (g^{0,12,3})^{-1}
\Big) 
\\ & &
= (g^{1,2,3})^{-1}g^{0,23,1}X_1 g^{01,2,3}(t_{02}+t_{12})
(g^{0,1,2})^{-1}(g^{0,12,3})^{-1}
\\ &  &
= (g^{1,2,3})^{-1}g^{0,23,1} g^{0,2,3} t_{02} X_1 
(g^{0,1,2})^{-1}(g^{0,12,3})^{-1}, 
\end{eqnarray*}
while 
\begin{eqnarray*}
&& \tilde\theta_g(t_{02})\tilde\theta_g(X_1)
= (g^{2,1,3})^{-1} g^{0,2,31} t_{02} (g^{0,2,31})^{-1} g^{2,1,3}
(g^{1,2,3})^{-1} g^{0,23,1} X_1 (g^{0,1,23})^{-1} g^{1,2,3}
\\ && 
= (g^{2,1,3})^{-1} g^{0,2,31} t_{02} g^{02,3,1} (g^{0,2,3})^{-1}
X_1 (g^{0,1,23})^{-1}g^{1,2,3}
\\ && 
= (g^{2,1,3})^{-1} g^{0,2,31} g^{02,3,1} t_{02} X_1
(g^{01,2,3})^{-1}(g^{0,1,23})^{-1}g^{1,2,3}
\\ && 
= g^{3,1,2}g^{2,3,1}g^{0,23,1}g^{0,2,3}t_{02}X_1
(g^{0,1,2})^{-1}(g^{0,12,3})^{-1}, 
\end{eqnarray*}
so the relation $X_1(t_{02}+t_{12})X_1^{-1}=t_{02}$
is preserved. (\ref{simpl}) then implies that the relations
$X_i(t_{0j}+t_{ij})X_i^{-1}=t_{0j}$ are preserved. Together with 
the other relations, these relations imply the relations $X_i t_{0i} X_i^{-1}
= \sum_{\alpha\neq i}t_{\alpha i}$, which are therefore also preserved. 
\hfill \qed \medskip

{\em End of proof of Proposition \ref{prop26}.} If $g\in
\on{GRT}_1(\kk)$, then one checks that the automorphisms 
$\tilde\theta_g$ from Lemma \ref{lemma29} and $\alpha_g
\in\on{Aut}(\hat\t_{1,2}^\kk)$ from Lemma-Definition 
\ref{def:alpha:g} are related in the same way as $\theta_{g,u_\pm}$
and $\alpha$ are in Lemma \ref{lemma28}. It follows that if 
$u_\pm^g$ are as in Lemma-Definition \ref{def:alpha:g}, 
then $(g,u_+^g,u_{-}^{g})\in\on{GRT}_1^{ell}(\kk)$. This defines a 
map $\on{GRT}_1(\kk)\to \on{GRT}_1^{ell}(\kk)$. 

Let us show that $\on{GRT}_1(\kk)\to \on{GRT}_1^{ell}(\kk)$
is a group morphism. In view of (\ref{pdt:theta}), it 
suffices to prove that $\tilde\theta_{g_2}\tilde\theta_{g_1}
=\tilde\theta_{g_1*g_2}$, which can be checked directly, e.g., 
\begin{eqnarray*}
&& \tilde\theta_{g_2}(\tilde\theta_{g_1}(X_1)) = 
\tilde\theta_{g_2}(g_1^{0,2,1}X_1(g^{0,1,2})^{-1})
= \tilde\theta_{g_2}(g_1(t_{02},t_{21}) X_1 g_1^{-1}(t_{01},t_{12}))
\\ && = g_1(\on{Ad}(g_2^{0,2,1})(t_{02}),t_{21})g_2^{0,2,1} 
X_1 (g_2^{0,2,1})^{-1}g_1^{-1}(\on{Ad}(g_2^{0,2,1})(t_{01}),t_{12})
\\ && = (g_1*g_2)^{0,2,1}X_1 ((g_1*g_2)^{0,1,2})^{-1} = 
\tilde\theta_{g_1*g_2}(X_1), 
\end{eqnarray*}
etc. 

Let us prove that $\on{GRT}_1(\kk)\to \on{GRT}_1^{ell}(\kk)$
is compatible with the actions of $\kk^\times$. If $c\cdot 
(g,u_\pm) = (\tilde g,\tilde u_\pm)$, then $\theta_{\tilde g,\tilde u_\pm}$
and $\theta_{g,u_\pm}$ are related by $\theta_{\tilde g,\tilde u_\pm}
= \gamma_c \theta_{g,u_\pm} \gamma_c^{-1}$, where $\gamma_c\in 
\on{Aut}(\hat\t_{1,3}^\kk)$ is given by $\gamma_c(x_i^+)=x_i^+$, 
$\gamma_c(x_i^-)=c^{-1}x_i^-$. It then suffices to prove that 
$\tilde\theta_{\tilde g} = \gamma_c \tilde\theta_g \gamma_c^{-1}$, 
where we recall that $\tilde g(A,B) = g(c^{-1}A,c^{-1}B)$, which 
follows from $\gamma_c(x_i)=x_i$, $\gamma_c(t_{\alpha\beta})=
c^{-1}t_{\alpha\beta}$ for $0\leq \alpha\neq\beta\leq 3$. 

The final commutative diagram follows from 
$$\xymatrix{
\on{GRT}_{1}^{ell}(\kk)\rtimes\kk^{\times} \ar[d]\ar[dr] \\
\on{SL}_{2}(\kk)\rtimes\kk^{\times} \ar[r]^{\sim} & \on{GL}_{2}(\kk)}$$

We set 
\begin{equation} \label{def:Rell:gr}
R_{ell}^{gr}(\kk):= \on{Ker}\big(\on{GRT}_{ell}(\kk)\to\on{GRT}(\kk)
\big).  
\end{equation}
Explicitly, 
\begin{eqnarray} \label{explicit:Rell:gr}
R_{ell}^{gr}(\kk) & = & \{(u_+,u_-)\in
(\hat{\t}_{1,2}^{\kk})^2 | u_\pm^{1,23}+u_\pm^{2,31}+u_\pm^{3,12}=0,
[u_\pm^{1,23},u_\pm^{2,13}]=0, [u_+^{1,23},u_-^{2,13}]=t_{12}\}
\nonumber \\ && \subset \on{Aut}(\hat{\t}_{1,2}^{\kk})^{op}.
\end{eqnarray} 
Then $\kk\mapsto R_{ell}^{gr}(\kk)$ is $\QQ$-group scheme, and 
we have a commutative diagram 
$$
\begin{matrix}
1 &\to & R_{ell}^{gr}(-) & \to & \on{GRT}_{ell}(-) & \to &
\on{GRT}(-) & \to & 1 \\
  &    &  \downarrow     &     & \downarrow        &     &
\downarrow  &     &   \\
1 &\to & \on{SL}_2       & \to & \on{GL}_2  & \stackrel{\on{det}}{\to} &
{\mathbb G}_m & \to & 1 
\end{matrix}
$$
The lift of $\on{GRT}_{ell}(-)\to\on{GL}_2$ restricts to a 
morphism $\on{SL}_2\to R_{ell}^{gr}(-)$, and the structure 
of $R^{gr}_{ell}(-)$ is therefore 
$$
R^{gr}_{ell}(-) = \on{Ker}\big(R_{ell}^{gr}(-)\to\on{SL}_2\big)
\rtimes\on{SL}_2,  
$$
in which the kernel is prounipotent. 

The morphism from Proposition \ref{prop26} enables us to 
define an action of $\on{GRT}(-)$ on $R_{ell}^{gr}(-)$. 
$\on{GRT}_{ell}(-)$ has then the structure of a semidirect 
product, fitting in 
$$
\begin{matrix}
\on{GRT}_{ell}(-) & \simeq & R_{ell}^{gr}(-)\rtimes \on{GRT}(-)\\
\downarrow & & \downarrow \\
\on{GL}_2 & \simeq & \on{SL}_2 \rtimes{\mathbb G}_m 
\end{matrix}
$$
where the bottom morphism is induced by ${\mathbb G}_m\to\on{GL}_2$, 
$c\mapsto \bigl( \begin{smallmatrix} 1 & 0 \\ 0 & c
\end{smallmatrix}\bigr)$.

\subsection{Categorical interpretations} \label{sec:cat}

A left action of $\on{GRT}(\kk)$ on $\{$IMBCs$\}$ is defined
as follows: $g\in\on{GRT}_1(\kk)$ acts on 
$(\cC,c_{\ldots},a_{\ldots},t_{\ldots})$ by only modifying  
$a_{XYZ}$ into $a'_{X,Y,Z}:= a_{XYZ}g(t_{XY}\otimes\on{id}_Z,
a_{XYZ}^{-1}(\on{id}_X\otimes t_{YZ})a_{XYZ})$
and $c\in\kk^\times$ acts by only modifying $t_{XY}$ into 
$ct_{XY}$. 

Similarly, one can show that 
a left action of $\on{GRT}_{ell}(\kk)$ on $\{$(an IBMC,
an elliptic structure over it)$\}$ is defined as follows: 
$(g,u_+,u_{-})\in\on{GRT}_1^{ell}(\kk)$ acts on $(\cC,\tilde\cC)$
as $(g,u_{+},u_-)*(\cC,\tilde\cC):= (g*\cC,\tilde\cC')$, where 
for $\tilde\cC=(\tilde\cC,F,x^\pm_{\ldots})$, we set 
$\tilde\cC'=(\tilde\cC,F,\ul x_{\ldots}^\pm)$, where
$\ul x_{X,Y}^\pm = u^\pm(x^+_{X,Y},x^-_{X,Y})$, and 
$c\in\kk^\times$ acts on $(\cC,\tilde\cC)$ as
$c*(\cC,\tilde\cC):= (c*\cC,\tilde\cC')$, where 
$\tilde\cC'=(\tilde\cC,F,x^+_{X,Y},cx^-_{X,Y})$. 

\subsection{Action of $\on{GRT}_{ell}(\kk)$ on $\ul{Ell}(\kk)$}

Recall that $\on{GRT}(\kk)$ acts on $\ul{M}(\kk)$ from the 
right as follows: for $g\in \on{GRT}_1(\kk)$ and $(\mu,\Phi)
\in \ul{M}(\kk)$, $(\mu,\Phi)*g:= (\mu,\tilde\Phi)$, where 
$$
\tilde\Phi(t_{12},t_{23}) = \Phi(\on{Ad}(g^{1,2,3})(t_{12}),t_{23})
g^{1,2,3}, 
$$
and for $c\in\kk^\times$,$(\mu,\Phi)*c:= (c\mu,c*\Phi)$, 
where $(c*\Phi)(A,B) = \Phi(cA,cB)$. This action is compatible 
with the maps $\{$IBMCs$\}\to\{$BMCs$\}$ induced by elements of 
$\ul{M}(\kk)$: $\Phi*(g*\cC_0) = (\Phi*g)*\cC_0$
for any $\Phi\in \ul{M}(\kk)$, $g\in\on{GRT}(\kk)$
and IBMC $\cC_0$.

For $(g,u_\pm)\in\on{GRT}_1^{ell}(\kk)$ and $(\mu,\Phi,A_\pm)
\in\ul{Ell}(\kk)$, we set $(\mu,\Phi,A_\pm)*(g,u_\pm):= 
(\mu,\tilde\Phi,\tilde A_\pm)$, where 
$$
\tilde A_\pm(x_1,y_{1}):= 
A_\pm(u_+(x_1,y_1),u_-(x_1,y_1)) 
$$
(in other terms, $\tilde A_\pm=\theta(A_\pm)$, where
$\theta\in\on{Aut}(\hat{\t}_{1,2}^\kk)$ is 
$x_1^\pm\mapsto u_\pm(x_1^+,x_1^-)$) and for $c\in\kk^\times$, 
we set $(\mu,\Phi,A_\pm)*c:= (\mu,c*\Phi,c\sharp A_\pm)$, where
$(c\sharp A_\pm)(x_1^+,x_1^-):= A_\pm(x_1^+,cx_1^-)$. 

\begin{proposition}
This defines a right action of $\on{GRT}_{ell}(\kk)$
on $\ul{Ell}(\kk)$, commuting with the 
left action of $\ul{\on{GT}}_{ell}(\kk)$ and compatible
with the right action of $\on{GL}_2(\kk)$ on $\on{M}_2(\kk)$.
\end{proposition}

{\em Proof.} Let us show that $(\mu,\tilde\Phi,\tilde A_\pm)
\in\ul{Ell}(\kk)$. If $\theta\in\on{Aut}(\hat\t_{1,2}^\kk)$
is defined by $\theta(x_1^\pm) = u_\pm$, and $\tilde\theta:= 
\theta_{g,u_\pm}$, then one checks that
$$
\tilde\theta(x^{1,23}) = \on{Ad}(g^{1,2,3})^{-1}(\theta(x)^{1,23}), \quad 
\tilde\theta(x^{2,31}) = \on{Ad}(g^{2,1,3})^{-1}(\theta(x)^{2,31}), \quad 
\tilde\theta(x^{3,12}) = \theta(x)^{3,12} 
$$ 
for any $x\in\hat\t_{1,2}^\kk$. Applying $\tilde\theta$ to 
(\ref{def:ell:ass:1}), one gets  
\begin{eqnarray*}
 \theta(\{e^{\pm\mu t_{12}/2}\}A_\pm)^{3,12}
\tilde\theta(\Phi^{3,1,2})(g^{2,1,3})^{-1}
\theta(\{e^{\pm\mu t_{12}/2}\}A_\pm)^{2,31}
g^{2,1,3}\tilde\theta(\Phi^{2,3,1})(g^{1,2,3})^{-1}
\theta(\{e^{\pm\mu t_{12}/2}\}A_\pm)^{1,23}
&& \\ g^{1,2,3}\tilde\theta(\Phi^{1,2,3})=1. && 
\end{eqnarray*}
Using the identities 
$\tilde\theta(\Phi^{3,1,2})(g^{2,1,3})^{-1} = \tilde\Phi^{3,1,2}$, 
$g^{2,1,3}\tilde\theta(\Phi^{2,3,1})(g^{1,2,3})^{-1} 
= \tilde\Phi^{2,3,1}$, $g^{1,2,3}\tilde\theta(\Phi^{1,2,3}) 
= \tilde\Phi^{1,2,3}$, and $\theta(\{e^{\pm\mu t_{12}/2}\}A_\pm)
=\{e^{\pm\mu t_{12}/2}\}\tilde A_\pm$, one obtains that 
$(\mu,\tilde\Phi,\tilde A_\pm)$
satisfies (\ref{def:ell:ass:1}). 

Applying now $\tilde\theta$ to (\ref{def:ell:ass:2}), one gets 
$$
e^{\mu t_{12}} = \big(\tilde\theta(\Phi)^{-1}g^{-1}\theta(A_-)^{1,23}
g\tilde\theta(\Phi),e^{-\mu t_{12}/2}\tilde\theta(\Phi^{2,1,3})^{-1}
(g^{2,1,3})^{-1}(\theta(A)^{2,13})^{-1} g^{2,1,3}
\tilde\theta(\Phi^{2,1,3})e^{-\mu t_{12}/2}\big).
$$
Using again $g\tilde\theta(\Phi) = \tilde\Phi$ and $g^{2,1,3}
\tilde\theta(\Phi^{2,1,3})=\tilde\Phi^{2,1,3}$, together with 
$\theta(A_\pm) = \tilde A_\pm$, one obtains that $(\mu,\tilde\Phi,
\tilde A_\pm)$ satisfies (\ref{def:ell:ass:2}). 

Similarly, applying the automorphism $x_i^+\mapsto x_i^+$, 
$x_i^-\mapsto cx_i^-$ to identities (\ref{def:ell:ass:1}), 
(\ref{def:ell:ass:2}), one obtains that $(\mu,\Phi,A_\pm)* c$ 
satisfies the same identities, hence belongs to $\ul{Ell}(\kk)$. 
It is then immediate to check that this defines a right action of 
$\on{GRT}_{ell}(\kk)$, commuting with the left action of 
$\ul{\on{GT}}_{ell}(\kk)$. \hfill\qed\medskip 

\begin{proposition}
The action of $\on{GRT}_{ell}(\kk)$ on $\ul{Ell}(\kk)$ restricts to an action on 
$Ell(\kk)\subset \ul{Ell}(\kk)$, which is free and transitive. 
\end{proposition}

{\em Proof.} Given that the action of $\on{GRT}(\kk)$ on 
$M(\kk)$ is free and transitive, it suffices to prove that the 
action of $R_{ell}^{gr}(\kk)$ on $Ell_{(\mu,\Phi)}(\kk):= 
Ell(\kk)\times_{M(\kk)}\{(\mu,\Phi)\}$ is free and transitive for any 
$(\mu,\Phi)\in Ell(\kk)$. 

Recall that $R_{ell}^{gr}(\kk)$ is explicitly described by 
(\ref{explicit:Rell:gr}); its inclusion into 
$\on{Aut}(\hat{\t}_{1,2}^\kk)$ is given by $(u_+,u_-)
\mapsto \theta_{u_+,u_-} = (x_1^\pm\mapsto u_\pm)$. 
On the other hand, $Ell_{(\mu,\Phi)}(\kk) = \{(A_+,A_-)$
satisfying (\ref{def:ell:ass:1}), (\ref{def:ell:ass:2})$\}$. 
Then 
\begin{equation} \label{act:AB}
(A_+,A_-)*(u_+,u_-) = (\theta_{u_\pm}(A_+),
\theta_{u_\pm}(A_-)).
\end{equation} 
Relation (\ref{def:ell:ass:2})
implies that $(A_-,A_+)=e^{\mu t_{12}}$, which together with 
$\mu\in\kk^\times$ implies that $\hat\t_{1,2}^\kk$ is 
generated by $\on{log}A_+,\on{log}A_{-}$. Together with (\ref{act:AB}), 
this implies that the action of $R_{ell}^{gr}(\kk)$ on 
$Ell_{(\mu,\Phi)}(\kk)$ is free. 

Let us now show that this action is transitive. We first observe that 
$R_{ell}^{gr}(\kk)$  can be described as $\{\theta\in
\on{Aut}(\hat\t_{1,2}^\kk) | \exists\tilde\theta\in
\on{Aut}(\hat\t_{1,3}^\kk)$ with $\tilde\theta(t_{ij})=t_{ij}$ 
for $1\leq i\neq j\leq 3$ and $\tilde\theta(x^{i,jk})
= \theta(x)^{i,jk}$ for $\{i,j,k\}=\{1,2,3\}$ and $x\in
\hat\t_{1,2}^\kk\}$. Let $(A_+,A_-)$
and $(\tilde A_+,\tilde A_-)\in Ell_{(\mu,\Phi)}(\kk)$ and let 
$\theta\in\on{Aut}(\hat\t_{1,2}^\kk)$ be the automorphism
such that $\theta(A_\pm)=\tilde A_\pm$. Let us show that 
there exists $\tilde\theta\in\on{Aut}(\hat{\t}_{1,3}^\kk)$, 
such that 
\begin{equation} \label{id:tilde:theta}
\tilde\theta(t_{ij})=t_{ij} \text{\ for\ }1\leq i\neq j\leq 3
\quad\text{and}\quad
\tilde\theta(x^{i,jk})
= \theta(x)^{i,jk} \text{\ for\ }
\{i,j,k\}=\{1,2,3\} \text{\ and\ }x\in
\hat\t_{1,2}^\kk.
\end{equation} 
Let $i_{(\mu,\Phi)}:B_3(\kk)\to\on{exp}(\hat\t_3)\rtimes
S_3$, $i_{(\mu,\Phi,A_\pm)}:B_{1,3}(\kk)\to
\on{exp}(\hat\t_{1,3})\rtimes S_3$ be the isomorphisms
induced by $(\mu,\Phi)$, $(\mu,\Phi,A_\pm)$ and the object
$\bullet(\bullet\bullet)$. We have a commutative diagram 
$$
\xymatrix{
P_3(\kk)\ar[r]_{\tilde i_{(\mu,\Phi)}}\ar[d]
\ar@/^{1pc}/[rr] 
 & \on{exp}(\hat\t_3^\kk)\ar[d]\ar@/^{1pc}/[rr] & 
B_3(\kk) \ar[r]_{i_{(\mu,\Phi)}}\ar[d] 
& \on{exp}(\hat\t_3^\kk)\rtimes S_3\ar[d] \\
P_{1,3}(\kk)\ar[r]^{\tilde i_{(\mu,\Phi,A_\pm)}} \ar@/_{1pc}/[rr]& 
\on{exp}(\hat\t_{1,3}^\kk)\ar@/_{1pc}/[rr] & 
B_{1,3}(\kk)\ar[r]^{i_{(\mu,\Phi,A_\pm)}} & \on{exp}(\hat\t_{1,3}^\kk)
\rtimes S_3 
}
$$
where the maps `$i$' are isomorphisms. Note that for $\sigma
\in B_3$, $i_{(\mu,\Phi)}(\sigma)[\sigma]^{-1}\in\on{exp}(\hat\t_3^\kk)$
(where $\sigma\mapsto [\sigma]$ is the canonical morphism $B_3\to S_3$). 

Then 
\begin{eqnarray*}
&& \tilde i_{(\mu,\Phi,A_\pm)}(X_1^\pm) = A_\pm^{1,23}, \quad 
\tilde i_{(\mu,\Phi,A_\pm)}(X_2^\pm) = \{i_\Phi(\sigma_1^{\pm 1})s_1\}
A_\pm^{2,13} \{s_1 i_\Phi(\sigma_1^{\pm 1})\}, \quad 
\\ && \tilde i_{(\mu,\Phi,A_\pm)}(X_3^\pm) = \{i_\Phi(\sigma_2^{\pm 1}
\sigma_1^{\pm 1})s_1s_2\}
A_\pm^{3,12} \{s_2s_1 i_\Phi(\sigma_1^{\pm 1}\sigma_2^{\pm 1})\}, 
\end{eqnarray*}
where we recall that $x\mapsto \{x\}$ is induced by the canonical 
morphism $\t_3\to\t_{1,3}$. Also 
$\tilde i_{(\mu,\Phi,A_\pm)}(\sigma_i^2) = 
\{\tilde i_{(\mu,\Phi)}(\sigma_i^2)\}$, for $i=1,2$ and 
$\tilde i_{(\mu,\Phi,A_\pm)}(\sigma_1\sigma_2^2\sigma_1) = 
\{\tilde i_{(\mu,\Phi)}(\sigma_1\sigma_2^2\sigma_1)\}$. 

Let $\tilde\theta:= \tilde i_{(\mu,\Phi,\tilde A_\pm)}\circ
\tilde i_{(\mu,\Phi,A_\pm)}^{-1}$. Then $\tilde\theta\in
\on{Aut}(\hat\t_{1,3}^\kk)$, and: 

(a) $\tilde\theta$ leaves $\{\tilde i_{(\mu,\Phi)}(\sigma_i^2)\}$
($i=1,2$) and $\{\tilde i_{(\mu,\Phi)}(\sigma_1\sigma_2^2\sigma_1)\}$
fixed, so it leaves the image of $\hat\t_3\to\hat{\t}_{1,3}$
pointwise fixed; 

(b) $\tilde\theta(A_\pm^{1,23})=\tilde A_\pm^{1,23}$, 
$$
\tilde\theta(\{i_{(\mu,\Phi)}(\sigma_1^{\pm 1})s_1\}A_\pm^{2,13}
\{s_1 i_{(\mu,\Phi)}(\sigma_1^{\pm 1})\})=
\{i_{(\mu,\Phi)}(\sigma_1^{\pm 1})s_1\}\tilde A_\pm^{2,13}
\{s_1 i_{(\mu,\Phi)}(\sigma_1^{\pm 1})\},
$$ which implies, as
$\{i_{(\mu,\Phi)}(\sigma_1^{\pm 1})s_1\}$ and
$\{s_1 i_{(\mu,\Phi)}(\sigma_1^{\pm 1})\}\in\on{im}
(\on{exp}(\hat\t_3^\kk)\to\on{exp}(\hat\t_{1,3}^\kk))$, that 
$\tilde\theta(A_\pm^{2,13})=\tilde A_\pm^{2,13}$; one proves similarly
that $\tilde\theta(A_\pm^{3,12})=\tilde A_\pm^{3,12}$. 

(b) implies that $\tilde\theta(x^{i,jk})=\theta(x)^{i,jk}$
holds for $x=A_\pm$, therefore also for $x$ in the topological 
group generated by $A_\pm$. As $\mu\in\kk^\times$, this group 
is equal to $\on{exp}(\hat\t_{1,2}^\kk)$, so $\tilde\theta$
satisfies (\ref{id:tilde:theta}). So $\theta\in R_{ell}^{gr}(\kk)$. 
\hfill \qed\medskip

\begin{proposition}
The scheme morphisms $\ul{Ell}\to \ul M$ and 
$\ul M\stackrel{\sigma}{\to}\ul{Ell}$ (see Proposition \ref{prop:sigma})
are compatible with the morphisms $\on{GRT}_{ell}(-)
\to\on{GRT}(-)$ and $\on{GRT}(-)\to\on{GRT}_{ell}(-)$ (see Proposition 
\ref{prop26}).  
\end{proposition}

{\em Proof.} We need to prove the second statement only. Let 
$\ul M(\kk)\stackrel{\sigma}{\to}\ul{Ell}(\kk)$ be given by $(\mu,\Phi)
\mapsto (\mu,\Phi,A_\pm(\mu,\Phi))$, then we must show that for
$g\in\on{GRT}_1(\kk)$ and $(\mu,\tilde\Phi)= (\mu,\Phi)*g$, we
have $A_\pm(\mu,\tilde\Phi)=\alpha_g(A(\mu,\Phi))$, where $\alpha_g$
is as in Lemma-Definition \ref{def:alpha:g}. This follows from 
the fact that $\alpha_g$ satisfies $\alpha_g(t_{02}) = 
\on{Ad}(g^{0,2,1})(t_{02})$, $\alpha_g(t_{12})=t_{12}$. It is also 
clear that $\ul M(\kk)\stackrel{\sigma}{\to}\ul{Ell}(\kk)$
is compatible with the action of $\kk^\times$. 
\hfill\qed\medskip 

\begin{remark} In fact, the commutative diagrams 
$\begin{matrix}\ul{Ell} &\to & \ul{M} \\ \downarrow & &\downarrow 
\\ \on{M}_2&\stackrel{\on{det}}{\to} & {\mathbb A}\end{matrix}$
and $\begin{matrix}\ul M &\stackrel{\sigma}{\to} & \ul{Ell} 
\\ \downarrow & &\downarrow 
\\ {\mathbb A}&\stackrel{c\mapsto 
\bigl( \begin{smallmatrix} 0 & -c \\ 1 & 0
\end{smallmatrix}\bigr)}{\to} & \on{M}_2\end{matrix}$
are compatible with the right actions of the diagrams 
$$
\begin{matrix}\on{GRT}_{ell}(-) &\to & \on{GRT}(-) \\ \downarrow & &\downarrow 
\\ \on{GL}_2&\stackrel{\on{det}}{\to} & {\mathbb G}_m\end{matrix}
\text{\ and\ }
\begin{matrix}\on{GRT}(-) &\to & \on{GRT}_{ell}(-) \\ 
\downarrow & &\downarrow 
\\ {\mathbb G}_m &\stackrel{c\mapsto 
\bigl( \begin{smallmatrix} 1 & 0 \\ 0 & c
\end{smallmatrix}\bigr)}{\to} & \on{GL}_2\end{matrix}
$$
\end{remark}

\subsection{Lie algebras} \label{sec:LA}

The graded Grothendieck-Teichm\"uller Lie algebra 
is\footnote{As before, $\f_2=\f_2^\QQ$, etc.} 
\begin{eqnarray*}
\grt_1= \{&& \psi\in\f_2 | 
\psi+\psi^{3,2,1}=0, \psi+\psi^{2,3,1}+\psi^{3,1,2}=0,
[t_{23},\psi^{1,2,3}]+[t_{13},\psi^{2,1,3}]=0,
\\ && \psi^{2,3,4}-\psi^{12,3,4}+\psi^{1,23,4}
-\psi^{1,2,34}+\psi^{1,2,3}=0\},
\end{eqnarray*} 
where we use the inclusion $\f_2\subset \t_3$, 
$A\mapsto t_{12}$, $B\mapsto t_{23}$; it is equipped with the Lie bracket 
$\langle\psi_1,\psi_2\rangle = [\psi_1,\psi_2]+D_{\psi_2}(\psi_1)
-D_{\psi_1}(\psi_2)$, where $D_\psi:A\mapsto [\psi,A]$, $B\mapsto 0$.

The Lie algebra $\QQ$ acts on $\grt_1$ by $[1,\psi]=-(\on{deg}\psi)\psi$
(where $\on{deg}A=\on{deg}B=1$), and we set $\grt:= \grt_1\rtimes\QQ$. 

The Lie algebras $\grt,\grt_{1}$ are $\NN$-graded (where $\on{deg}$ is 
extended to be 0 on $\QQ$), we then have $\on{Lie}\on{GRT}_{(1)}(-)
= \widehat{\grt}_{(1)}$ (the degree completions). 

Let 
\begin{eqnarray*}
\grt_1^{ell}:= \{ && (\psi,\alpha_\pm)\in\f_2\times (\t_{1,2})^2 | 
\psi\in\grt_1, \\
&& \alpha_\pm^{1,23}+\alpha_\pm^{2,31}+\alpha_\pm^{3,12}
+[x^1_\pm,\psi^{1,2,3}]+[x^2_\pm,\psi^{2,1,3}]=0, 
\\ && [x_\pm^1,\alpha_\pm^{3,12}]+[\alpha_\pm^{1,23},x_\pm^3]
-[x^1_\pm,[x^3_\pm,\psi^{1,2,3}]]=0, \\ && 
[x_+^1,\alpha_-^{2,13}]-[x_-^2,\alpha_+^{1,23}] = 
[x_-^2,[x_+^1,\psi^{1,2,3}]] - [x_{+}^{1},[x_{-}^{2},\psi^{2,1,3}]]
\}. 
\end{eqnarray*}
For $\alpha_\pm\in\t_{1,2}$, define $D_{\alpha_\pm}
\in\on{Der}(\t_{1,2})$ by $x_1^\pm\mapsto\alpha_\pm$. 
Then 
$$
[(\psi_1,\alpha_1^\pm),(\psi_2,\alpha_2^\pm)] = 
(\langle\psi_1,\psi_2\rangle,D_{\alpha_2^\pm}(\alpha_1^\pm)
-D_{\alpha_1^\pm}(\alpha_2^\pm))  
$$
defines a Lie bracket on $\grt_1^{ell}$, and 
$$
\grt_{1}^{ell}\subset \grt_{1}\times \on{Der}(\t_{1,2})^{op}. 
$$

The Lie algebra $\QQ e_{22}$ acts on $\grt_1^{ell}$ by 
$$
[e_{22},(\psi,\alpha_+,\alpha_-)] = (-(\on{deg}\psi)\psi,
-(\on{deg}_-\alpha_+)\alpha_+,(1-\on{deg}_-\alpha_-)\alpha_-),
$$ 
where $\on{deg}\psi$ is as above, and $\on{deg}_-\alpha_\pm$ is 
defined by $\on{deg}_- x_1^+=0$, $\on{deg}_- x_1^-=1$. 
We then set $\grt_{ell}:= \grt_1^{ell}\rtimes\QQ e_{22}$. 
 
The Lie algebras $\grt_{(1)}^{ell}$ are $\NN$-graded, 
where $(\psi,\alpha_\pm)$ has degree $n$ if 
$2\on{deg}\psi=\on{deg}\alpha_\pm-1=n$ ($\on{deg}\alpha_\pm$
being defined by $\on{deg}x_1^\pm=1$ and $\on{deg}\psi$ by 
$\on{deg}t_{12}=\on{deg}t_{23}=1$) and $e_{22}$ has degree $0$. 
Then $\on{Lie}\on{GRT}_{(1)}^{ell}(-) = \widehat{\grt}_{(1)}^{ell}$. 

We have a morphism $\grt_1^{ell}\to\SL_2$, $(\psi,\alpha_+,\alpha_-)
\mapsto \bigl( \begin{smallmatrix} a_+ & b_+ \\ a_- & b_-
\end{smallmatrix}\bigr)$, where $\alpha_+\equiv a_\pm x_1+b_\pm y_1$
modulo degree $\geq 2$. It extends to a morphism $\grt_{ell}\to\GL_2$
via $e_{22}\mapsto\bigl( \begin{smallmatrix} 0 & 0 \\ 0 & 1
\end{smallmatrix}\bigr)$. We denote by $\grt_{I_2}^{ell}$ the common
kernel of these morphisms; it coincides with the part of 
$\grt^{ell}$ (or $\grt^{ell}_{1}$) of positive degree. 

These morphisms admit sections $\SL_2\to\grt^{ell}_{1}$ given by 
$\bigl( \begin{smallmatrix} a_+ & b_+ \\ a_- & b_-
\end{smallmatrix}\bigr) \mapsto (0,a_\pm x_1+b_\pm y_1)$
and $\GL_2\to\grt_{ell}$ given by its extension by 
$\bigl( \begin{smallmatrix} 0 & 0 \\ 0 & 1 
\end{smallmatrix}\bigr) \mapsto e_{22}$. We then have 
$\grt_1^{ell} \simeq \grt_{I_2}^{ell}\rtimes\SL_2$, 
$\grt_{ell} \simeq \grt_{I_2}^{ell}\rtimes\GL_2$.  

$\ZZ^2$-gradings may be defined on $\grt^{ell}_{(1)}$ as follows. 
We have a Lie algebra inclusion $\grt_1^{ell}
\subset \grt_1\oplus \on{Der}(\t_{1,2})=:\G$. Recall that 
$\grt_1$ is $\NN$-graded while $\on{Der}(\t_{1,2})$
is $\ZZ^2$-graded by the $\ZZ^2$-grading of $\t_{1,2}$
given by $(\on{deg}_+,\on{deg}_-)(x_1^+) = (1,0)$, 
$(\on{deg}_+,\on{deg}_-)(x_1^-) = (0,1)$.  
We then define a $\ZZ^2$-grading on $\G$ by 
$\G[p,q]:= 
 \left\{ \begin{array}{ll}
         \on{Der}(\t_{1,2})[p,q] & \on{\ if\ } q\neq p\\
        \grt_1[p]\oplus \on{Der}(\t_{1,2})[p,p]  
& \on{\ if\ } q=p \end{array} \right.$. This 
restricts to a $\ZZ^2$-grading $(\on{deg}_+,\on{deg}_-)$ 
of $\grt_1^{ell}$, which extends to $\grt_{ell}$ by 
$(\on{deg}_+,\on{deg}_-)(e_{22})=(0,0)$. 

The $\ZZ^2$-grading of $\grt_{ell}$ is compatible
with the action of the Cartan subalgebra of $\GL_2$: 
we have $[e_{11},x]=-(\on{deg}_+x)x$, $[e_{22},x]=
-(\on{deg}_-x)x$ for 
$e_{11} = \bigl( \begin{smallmatrix} 1 & 0 \\ 0 & 0
\end{smallmatrix}\bigr)\in\GL_2$ and $x\in\grt_{ell}$ homogeneous.   

We have a morphism $\grt_1^{ell}\to\grt_1$, $(\psi,u_\pm)\mapsto\psi$.
It extends to a morphism $\grt_{ell}\to\grt$ by $e_{22}\mapsto 1$. 
Using Proposition \ref{prop26}, sections of these morphisms 
are constructed as follows:  

\begin{proposition} \label{lift:LA}
There is a unique Lie algebra morphism $\grt_1\to\widehat{\grt}_1^{ell}$, 
$\psi\mapsto (\psi,u_+^\psi,u_-^\psi)$, where $u_\pm^\psi:= 
D_\psi(x_1^\pm)$ and $D_\psi\in\on{Der}(\hat\t_{1,2})$
is defined by 
$$
D_\psi(e^{x_1}) = \psi^{0,2,1}e^{x_1}
-e^{x_1}\psi^{0,1,2}, \quad D_\psi(t_{01}) = [\psi^{0,1,2},t_{01}]; 
$$ 
recall that 
$$
\psi^{0,1,2} = \psi(t_{01},t_{12}), \quad 
\psi^{0,2,1} = \psi(t_{02},t_{21}), \quad 
t_{0i} = - {{\on{ad}x_{i}}\over{e^{\on{ad}x_{i}}-1}}(y_{i}), \quad i =1,2. 
$$
It extends to a Lie algebra morphism $\grt\to\widehat{\grt}_{ell}$
by $1\mapsto e_{22}$. It is homogeneous, $\grt$ being equipped 
with its degree and $\widehat{\grt}_{ell}$ with degree $\on{deg}_-$. 
\end{proposition}

Set now $\r^{gr}_{ell}:= \on{Ker}(\grt_{ell}\to\grt)$. 
We have 
\begin{eqnarray*}
\r_{ell}^{gr} = \{  (\alpha_+,\alpha_-)\in (\t_{1,2})^2 | 
&& \alpha_\pm^{1,23}+\alpha_\pm^{2,31}+\alpha_\pm^{3,12}=0, \\ && 
[x_\pm^1,\alpha_\pm^{2,13}]+[\alpha_\pm^{1,23},x^2_\pm]=0, \\ && 
[x_+^1,\alpha_-^{2,13}]+[\alpha_+^{1,23},x_-^2]=0\} \subset 
\on{Der}(\t_{1,2})^{op}. 
\end{eqnarray*}
This is a $\ZZ^2$-graded Lie subalgebra of 
$\grt_{ell}$; it is also $\NN$-graded by $\on{deg}_++
\on{deg}_-$. We have $\r_{ell}^{gr}[0]\simeq \SL_2$
and $\r_{ell}^{gr}\simeq (\oplus_{d>0}\r_{ell}^{gr}[d])\rtimes \SL_2$. 
Its completion for the $\NN$-degree is isomorphic to 
$\on{Lie}R_{ell}^{gr}(-)$. 

Define a partial completion $\hat\r_{ell}^{gr}:= 
\oplus_q (\prod_p \r_{ell}^{gr}[p,q])$. Proposition 
\ref{lift:LA} gives rise to a Lie algebra morphism 
$\grt\to\on{Der}(\hat\r_{ell}^{gr})$. We then have 
$\hat\grt_{ell}\simeq \hat\r_{ell}^{gr}\rtimes\grt$, 
where $\hat\grt_{ell}:= \oplus_q\prod_p\grt_{ell}[p,q]$
is a partial completion. 

Set ${\mathfrak{gt}}_{ell}:=\on{Lie}\on{GT}_{ell}(-)$, 
${\mathfrak{gt}}_{1}^{ell}:=\on{Lie}\on{GT}_{1}^{ell}(-)$, then 
${\mathfrak{gt}}_{ell} = {\mathfrak{gt}}_{1}^{ell} \rtimes \QQ$. 
The Lie algebra ${\mathfrak{gt}}_{1}^{ell}$ admits a description 
as a subspace of $\hat\f_{2}\times(\hat\t_{1,2})^{2}$ similar to that of 
Lemma \ref{descr:R:ell}, and is filtered as follows: 
${\mathfrak{gt}}_{ell}={\mathfrak{gt}}_1^{ell}\rtimes\QQ$, where 
${\mathfrak{gt}}_{1}^{ell} := \on{Lie}\on{GT}_{1}^{ell}(-) \subset 
\hat\f_{2}\times(\hat\t_{1,2})^{2}$. We then set  
$({\mathfrak{gt}}_{1}^{ell})^{\geq n}:= 
 {\mathfrak{gt}}_{1}^{ell} \cap \big( \hat\f_{2}^{\geq n/2}
\times ((\hat\t_{1,2})^{2})^{\geq n+1} \big)$ for $n\geq 0$, 
where the degree in $\hat\f_{2}$ is induced by $\on{deg}(t_{12})
=\on{deg}(t_{23})=1$ and the degree in $\hat\t_{1,2}$
by $\on{deg}(x_{1}^{\pm})=1$. The Lie algebra ${\mathfrak{gt}}_{ell}$
is similarly filtered by $({\mathfrak{gt}}_{ell})^{\geq 0} = 
{\mathfrak{gt}}_{ell}$, $({\mathfrak{gt}}_{ell})^{\geq n} = 
({\mathfrak{gt}}_{1}^{ell})^{\geq n}$ if $n>0$. 
It follows from the form of the conditions under which 
$(\psi,\alpha_{+},\alpha_{-})\in \hat\f_{2}\times(\t_{1,2})^{2}$
belong to ${\mathfrak{gt}}_{1}^{ell}$ that there is a canonical 
morphism $\on{gr}({\mathfrak{gt}}_{ell})\to{\grt}_{ell}$, restricting to 
$\on{gr}({\mathfrak r}_{ell})\to {\mathfrak r}_{ell}^{gr}$ and compatible 
with $\on{gr}(\gt)\to \grt$. In Section \ref{sect:isos}, we will see 
that all these morphisms are isomorphisms. 

\begin{remark}
The relations between Lie groups and algebras are summarized as follows: 
$$
\on{GRT}_{1}(\kk) = \on{exp}(\widehat\grt_1^\kk), \quad 
\on{GRT}(\kk) = \on{exp}(\widehat{\grt}_{1}^{\kk})\rtimes\kk^{\times}, 
$$
$$
\on{GRT}_{1}^{ell}(\kk) = 
\on{exp}(
\widehat\grt^{ell,\kk}_{I_2})\rtimes \on{SL}_2(\kk), \quad 
\on{GRT}_{ell}(\kk) = 
\on{exp}(
\widehat\grt^{ell,\kk}_{I_2})\rtimes \on{GL}_2(\kk),  
$$
$$
R_{ell}^{gr}(\kk)
= \on{exp}(\prod_{d>0}\r_{ell}^{gr}[d]\otimes\kk)\rtimes\on{SL}_2(\kk). 
$$
\end{remark}

\begin{remark}
Any $(\alpha_+,\alpha_-)\in\r_{ell}^{gr}$ satisfies 
$\alpha_\pm+\alpha_\pm^{2,1}=0$, which implies 
that the total degree (in which $x_1^\pm$ have degree 1) 
of $\alpha_\pm$ is odd. So $\r_{ell}^{gr}[d]=0$ unless $d$
is even. 
\end{remark}

\begin{remark} (Relation with the work of H.~Tsunogai.) 
In \cite{Ts}, H.~Tsunogai describes the ``stable derivation 
algebra'' in genus one. This is a graded Lie algebra version of the 
intersection over $n\geq 1$ of the images of the morphisms
$\on{Out}^{*}(P_{1,n})\to \on{Out}^{*}(P_{1,1})$, where 
$\on{Out}^{*}\subset \on{Out}$ are certain subgroups. 
This is a Lie subalgebra ${\mathcal G}_{\on{Ts}}\subset \on{Der}(\t_{1,2})$, 
which may be defined as the set of all $(\alpha_{+},\alpha_{-})\in 
(\t_{1,2})^{2}$, such that there exists $\psi\in\t_{3}$, such that 
$$
\psi^{1,2,3}+\psi^{3,2,1} = [t_{12},\psi^{1,2,3}]+[t_{13},\psi^{2,1,3}]=0, 
\quad
[x_{+}^{1},\alpha_{-}^{1,2}] + [\alpha_{+}^{1,2},x_{-}^{1}]=0, 
$$
$$
[x_{\pm}^{1},\alpha_{\pm}^{3,12}] + [\alpha_{\pm}^{1,23},x_{\pm}^{3}]
= [x_{3}^{\pm},[x_{1}^{\pm},\psi^{1,2,3}]], 
\quad
[x_{+}^{1},\alpha_{-}^{3,12}] + [\alpha_{+}^{1,23},x_{-}^{2}] = 
[t_{13},\psi^{1,3,2}] + [x_{+}^{1},[x_{-}^{2},\psi^{1,3,2}]]
$$
(the relation between the present formalism and that of \cite{Ts} is
as follows: $\t_{3} \leftrightarrow {\mathcal L}_{1}^{(2)\circ}$, 
$\t_{1,2}\leftrightarrow {\mathcal L}_{1}^{(2)}$, $\alpha_{+},\alpha_{-}
\leftrightarrow S,T$, $U^{1,2,3}\leftrightarrow \psi^{2,1,3}$; the present 
relations are obtained from those of \cite{Ts} by some changes of indices). This 
system of conditions is a consequence of the system expressing that 
$(\psi,\alpha_{+},\alpha_{-})\in\grt_{1}^{ell}$; the latter is more restrictive 
as it contains additional conditions, namely the pentagon and hexagon conditions on 
$\psi$, as well as the conditions 
$\alpha_{\pm}^{1,23}+\alpha_{\pm}^{2,31}+\alpha_{\pm}^{3,12}
+[x_{\pm}^{1},\psi^{1,2,3}] + [x_{\pm}^{2},\psi^{2,1,3}] = 0$. It follows that there
is a double inclusion
$$
\on{im}(\grt_{1}^{ell}\to\on{Der}(\t_{1,2})) \subset {\mathcal G}_{\on{Ts}}
\subset \on{Der}(\t_{1,2}). 
$$
\end{remark}

\subsection{A Lie subalgebra $\b_3\subset \r_{ell}^{gr}$}
\label{sec:GC}

\begin{proposition}
For $n\geq 0$, set 
\begin{equation} \label{delta:2n}
\delta_{2n}:= (\alpha_+ = \on{ad}(x_1)^{2n+2}(y_1), \alpha_- = 
{1\over 2}\sum_{\begin{smallmatrix}0\leq p\leq 2n+1,\\ 
p+q=2n+1\end{smallmatrix}}(-1)^p[(\on{ad}x_1)^p(y_1),
(\on{ad}x_1)^q(y_1)]). 
\end{equation}
Then $\delta_{2n}\in\r_{ell}^{gr}[2n+1,1]$. The element
$\delta_0$ is central in $\grt_{ell}$ and coincides with 
$\on{ad}t_{12}$ as an element of $\on{Der}(\t_{1,2})^{op}$.
\end{proposition}

{\em Proof.} In \cite{CEE}, Proposition 3.1, we constructed derivations 
$\dot\delta_{2n}^{(m)}\in\on{Der}(\t_{1,m})$, such that 
$$
\dot\delta_{2n}^{(m)} : x_i\mapsto 0, t_{ij}\mapsto 
[t_{ij},(\on{ad}x_i)^{2n}(t_{ij})], 
y_i\mapsto \sum_{j:j\neq i}{1\over 2}
\sum_{p+q=2n-1} [(\on{ad}x_i)^p(t_{ij}),(-\on{ad}x_i)^q(t_{ij})]. 
$$
Let then $\delta_{2n}^{(m)}:= \dot\delta_{2n}^{(m)} 
+ [\sum_{i<j}(\on{ad}x_i)^{2n}(t_{ij}),-]$. Then 
$$
\delta_{2n}^{(m)}(x_i) = [\sum_{j\neq i}(\on{ad}x_i)^{2n}(t_{ij}),x_i]
= (\on{ad}x_i)^{2n+2}(y_i)
=\alpha_-^{i,1...\check i...n}, \quad \delta_{2n}^{(m)}(t_{ij})=0, 
$$
\begin{eqnarray*}
&& \delta_{2n}^{(m)}(y_i) = \dot\delta_{2n}^{(m)}(y_i)+
[\sum_{j<k}(\on{ad}x_j)^{2n}(t_{jk}),y_i]
\\ && = \dot\delta_{2n}^{(m)}(y_i)
+\sum_{j\neq i}[(\on{ad}x_i)^{2n}(t_{ij}),y_i]
+\sum_{j<k; j,k\neq i}[(\on{ad}x_j)^{2n}(t_{jk}),y_i]
\\ && = \dot\delta_{2n}^{(m)}(y_i)
+\sum_{j\neq i}[(\on{ad}x_i)^{2n}(t_{ij}),y_i]
+\sum_{j<k; j,k\neq i}\sum_{p+q=2n-1}
(\on{ad}x_j)^{p}[t_{ij},(\on{ad}x_j)^{q}
(t_{jk})]
\\ && = \dot\delta_{2n}^{(m)}(y_i)
+\sum_{j\neq i}[(\on{ad}x_i)^{2n}(t_{ij}),y_i]
-\sum_{j<k; j,k\neq i}\sum_{p+q=2n-1}
[(-\on{ad}x_i)^{p}(t_{ij}),(\on{ad}x_i)^{q}(t_{ik})]
\\ && = \dot\delta_{2n}^{(m)}(y_i)
+\sum_{j\neq i}[(\on{ad}x_i)^{2n}(t_{ij}),y_i]
-{1\over 2}\sum_{p+q=2n-1}
[\sum_{j\neq i}(-\on{ad}x_i)^{p}(t_{ij}),\sum_{k\neq i}(\on{ad}x_i)^{q}(t_{ik})]
\\ && +{1\over 2}\sum_{j\neq i}\sum_{p+q=2n-1}
[(-\on{ad}x_i)^{p}(t_{ij}),(\on{ad}x_i)^{q}(t_{ij})]
\\ && = -[(\on{ad}x_i)^{2n+1}(y_i),y_i]
+{1\over 2}\sum_{p+q=2n-1}
[(-\on{ad}x_i)^{p+1}(y_i),(\on{ad}x_i)^{q+1}(y_i)]
\\ && = 
{1\over 2}\sum_{p+q=2n+1}
[(-\on{ad}x_i)^{p}(y_i),(\on{ad}x_i)^{q}(y_i)]
=\alpha_-^{i,1...\check i...n}. 
\end{eqnarray*}
Then $0=\delta_{2n}^{(3)}([x_1^\pm,x_2^\pm])=
[x_1^\pm,\alpha_\pm^{2,13}]+[\alpha_\pm^{1,23},x_2^\pm]$
and $0=\delta_{2n}^{(3)}(t_{12}) = [x_1^+,\alpha_-^{2,13}]
+[\alpha_+^{1,23},x_2^-]$, which implies that $\delta_{2n}
\in\r_{ell}^{gr}$. \hfill \qed\medskip 

We define $\b_3:= \langle\SL_2,\delta_{2n};n\geq 0\rangle \subset 
\r_{ell}^{gr}$ as the Lie subalgebra generated by $\SL_2$
and the $\delta_{2n}$. A basis of $\SL_{2}\subset \b_{3}$ is 
\begin{equation} \label{sl2}
e_{+}:= (\alpha_{+} = 0, \alpha_{-} = x_{1}), \quad 
e_{-} := (\alpha_{+} = y_{1},\alpha_{-} = 0), \quad 
h:= (\alpha_{+} = x_{1}, \alpha_{-} = -y_{-}). 
\end{equation}
The Lie algebra $\b_{3}$ is $\NN$-graded, 
and corresponds to the subgroup $\on{exp}( \hat\b_3^{+,\kk})
\rtimes\on{SL}_2(\kk)\subset R_{ell}^{gr}(\kk)$
(where the hat denotes the degree completion and 
$+$ means the positive degree part). 

\subsection{Isomorphisms of Lie algebras} \label{sect:isos}\label{sec:6:1}

Let $\kk$ be a $\QQ$-ring. As $Ell(\kk)$ is a torsor, each $e\in Ell(\kk)$
gives rise ro an isomorphism $i_{e}:\on{GT}_{ell}(\kk)\to \on{GRT}_{ell}(\kk)$, 
defined by $g*e=e*i_{e}(g)$ for any $g\in\on{GT}_{ell}(\kk)$. Similarly, 
any $\tilde\Phi\in M(\kk)$ gives rise to an isomorphism $i_{\tilde\Phi}:
\on{GT}(\kk)\to \on{GRT}(\kk)$ defined by the same conditions. We then 
have a commutative diagram 
\begin{equation} \label{diag}
\begin{matrix}
\on{GT}_{ell}(\kk) & \stackrel{i_{e}}{\to} & \on{GRT}_{ell}(\kk)\\
 \downarrow  & & \downarrow \\
\on{GT}(\kk)  & \stackrel{i_{\tilde\Phi}}{\to}& \on{GRT}(\kk) 
\end{matrix}\end{equation}
where $\tilde\Phi = \on{im}(e\in Ell(\kk)\to M(\kk))$. In  particular, $i_{e}$
restricts to an isomorphism $i_{e} : R_{ell}(\kk)\to R_{ell}^{gr}(\kk)$. 
When $e\in\on{im}(M(\kk)\stackrel{\sigma}{\to} Ell(\kk))$, the isomorphism 
$R_{ell}(\kk)\stackrel{i_{e}}{\to} R_{ell}^{gr}(\kk)$
is compatible with $i_{\tilde\Phi}$ and the actions of $\on{GT}(\kk)$, $\on{GRT}(\kk)$
on both sides via the lifts $\on{GT}(\kk)\stackrel{\sigma}{\to}
\on{GT}_{ell}(\kk)$, $\on{GRT}(\kk)\stackrel{\sigma}{\to}
\on{GRT}_{ell}(\kk)$.

The isomorphisms $i_{e}$ induce Lie algebra isomorphisms 
${\mathfrak{gt}}_{ell}^{\kk}\to \widehat{\grt}_{ell}^{\kk}$, restricting to 
${\mathfrak r}_{ell}^{\kk}\to\hat{{\mathfrak r}}_{ell}^{gr,\kk}$, 
compatible with the filtrations and whose associated graded isomorphisms
are the canonical morphisms from the end of Subsection \ref{sec:LA}. 
Since $Ell(\QQ)\neq\emptyset$
(e.g., because it contains $\sigma(M(\QQ))$), we obtain: 

\begin{proposition}
There are isomorphisms ${\mathfrak{gt}}_{ell}\simeq \hat{\on{gr}}
({\mathfrak{gt}}_{ell})
= \widehat{\grt}_{ell}$ 
and ${\mathfrak r}_{ell}\simeq \hat{\on{gr}}({\mathfrak r}_{ell})=
\hat{{\mathfrak r}}_{ell}^{gr}$. 
\end{proposition}

\subsection{Actions on prounipotent completions of 
elliptic braid groups}

Let ${\mathbf k}$ be a $\QQ$-ring. We recall that $P_{n}({\mathbf k})$
(resp., $P_{1,n}({\mathbf k})$) is the prounipotent completion of the 
pure (resp., elliptic) braid group $P_{n}$ (resp., $P_{1,n}$), where $n\geq 1$, and 
that $B_{n}({\mathbf k})$ (resp., $B_{1,n}({\mathbf k})$) to be 
the relative completion of the full (resp., elliptic) braid group with $n$ strands 
with respect to the 
canonical morphism to $S_{n}$; it identifes with the pushout  
$B_{n}*_{P_{n}}P_{n}({\mathbf k})$ (resp., 
$B_{1,n}*_{P_{1,n}}P_{1,n}({\mathbf k})$). 

\begin{proposition} \label{prop:act:GT:categs}
1) The action of $\on{GT} = \ZZ/2\ZZ$ on $B_{n}$ via $(-1)\cdot \sigma_{i} 
= \sigma_{i}^{-1}$ extends to the following objects: 

$\bullet$ a morphism $\mu_{O} : \on{GT}({\mathbf k})\to 
\on{Aut}(B_{n}({\mathbf k}))$ for each $O\in{\bf Pa}_{n}$; 

$\bullet$ a map
$$
\on{GT}({\mathbf k})\times {\bf Pa}_{n}\times {\bf Pa}_{n}\to 
P_{n}({\mathbf k}), \quad 
(g,O,O')\mapsto b_{OO'}(g), $$
related by the identities 
\begin{equation} \label{identity:mu:1}
\mu_{O'}(g) =\on{Inn}( b_{OO'}(g) ) \circ  \mu_{O}(g) ,
\end{equation}
\begin{equation} \label{identity:mu:2}
b_{OO'}(gh) = b_{OO'}(g) \cdot \mu_{O}(g)(b_{OO'}(h)), \quad 
b_{OO''}(g) = b_{O'O''}(g)b_{OO'}(g). 
\end{equation}

2) The action of $\on{GT}_{ell} = \tilde B_{3}$ on $B_{1,n}$ given by
(\ref{act:B3:ellbraids}) extends to a collection of morphisms
$$
\mu_{O}^{ell} : \on{GT}_{ell}({\mathbf k})\to \on{Aut}(B_{1,n}({\mathbf k})) 
$$
indexed by $O\in {\bf Pa}_{n}$, related to the morphisms $\mu_{O}$ by the identity 
\begin{equation} \label{identity:mu:3}
\mu_{O}^{ell}(g_{ell})(b_{ell}) = \mu_{O}(g)(b)_{ell}, 
\end{equation}
and satisfying 
\begin{equation} \label{identity:mu:4}
\mu_{O'}^{ell}(g_{ell}) = \on{Inn}( b_{OO'}(g)_{ell} ) \circ \mu_{O}^{ell}(g) ,
\end{equation}
for any $g_{ell}\in \on{GT}_{ell}({\mathbf k})$ and $b\in B_{n}({\mathbf k})$, 
where $g:= \on{im}(g_{ell}\in\on{GT}_{ell}({\mathbf k})
\to \on{GT}({\mathbf k}))$  and 
$b_{ell}:= \on{im}(b\in B_{n}({\mathbf k})\to 
B_{1,n}(\mathbf k))$. 

3) The restriction $\mu^{ell}_{O | R_{ell}({\mathbf k})}$ is 
independent of $O\in{\bf Pa}_{n}$ and will be denoted 
$$
\mu_{ell} : R_{ell}({\mathbf k})\to \on{Aut}(B_{1,n}({\mathbf k})). 
$$ 
If $g_{ell} = (1,1,g_{+},g_{-})\in R_{ell}({\mathbf k})$, where $g_{\pm}
= g_{\pm}(X_{1},Y_{1})\in P_{1,2}({\mathbf k})$, then the action of $g_{ell}$
on $B_{1,n}({\mathbf k})$ induced by $\mu_{ell}$ is such that 
$$ 
g_{ell}\cdot X_{1}^{\pm} = g_{\pm}(X_{1}^{+},X_{1}^{-}), \quad 
g_{ell}\cdot \sigma_{i} = \sigma_{i} \quad \text{for}\quad 
i = 1,\ldots,n-1.  
$$
\end{proposition}

{\em Proof.} 1) Let $\cC := {\bf PaB}_{{\mathbf k}}$ be the 
${\mathbf k}$-prounipotent version of the BMC ${\bf PaB}$, 
$G:= \on{GT}({\mathbf k})$. For $g\in G$, $g*\cC$ is a BMC  
with distinguished object $\bullet$. 
By the universal property of $\cC$, one derives from there a functor 
$\alpha_{g} : \cC\to g*\cC$, uniquely defined by the condition that it is tensor 
and that it induces the identity on objects. As a category, $g*\cC$ canonically 
identifies with $\cC$; let $i_{g} : g*\cC\to \cC$ be this isomorphism. One then 
defines $\beta_{g}:= i_{g}\circ \alpha_{g} : \cC\to\cC$. The identity 
$\beta_{g}\beta_{g'} = \beta_{g'g}$ follows from the commutativity of
$$
\xymatrix{
 &   & \cC\ar[dr]
 \ar@/^{2pc}/[ddrr]^{\beta_{g'}} \\
 & g*\cC\ar[ur]^{\sim}\ar[dr] & & g'*\cC\ar[dr]^{\sim} \\
\cC\ar[rr]\ar[ur]\ar@/^{2pc}/[uurr]^{\beta_{g}}
\ar@/_{1pc}/[rrrr]_{\beta_{gg'}} & & 
gg'*\cC\ar[rr]^{\sim}\ar[ur]^{\sim} 
& & \cC 
}
$$
in which the commutativity of the central square follows from that of
$$\begin{matrix}
g*\cC & \stackrel{g*\varphi}{\to}& g*\cD
\\
 \scriptstyle{\sim}\downarrow & & \downarrow\scriptstyle{\sim}\\
 \cC & \stackrel{\varphi}{\to}& \cD
\end{matrix}$$
for any braided monoidal categories $\cC,\cD$ and any tensor functor
$\varphi:\cC\to\cD$.  

It follows that $g\mapsto\beta_{g^{-1}}$ defines a morphism from $G$ to the group of
autofunctors of $\cC$, i.e., an action of $G$ on $\cC$. 

Let $O,O'\in{\bf Pa}_{n}$. There is a canonical isomorphism
$i_{O} : \on{Aut}_{\cC}(O)\to B_{n}({\mathbf k})$ and 
a canonical element $i_{OO'}\in \on{Iso}_{\cC}(O,O')$ (corresponding to 
the unit in $B_{n}({\mathbf k})$). 
Then for $f\in \on{Aut}_{\cC}(O)$, $i_{O}(f) = i_{O'}(i_{OO'}fi_{OO'}^{-1})$. 

Define the action $\mu_{O}$ of $G$ on $B_{1,n}({\mathbf k})$ as the transport
via $i_{O}$ of its action on $\on{Aut}_{\cC}(O)$, namely 
$\mu_{O}(g)(b):= i_{O}(g*i_{O}^{-1}(b))$. The claimed identities then hold with 
$b_{OO'}(g):= i_{O}(i_{OO'}^{-1}\circ (g*i_{OO'}))$. 

2) The collection of morphisms $\mu_{O}^{ell}$ is then defined in the same way:
$G$ is replaced by $G_{ell}:= \on{GT}_{ell}({\mathbf k})$, $\cC$ by $\cC_{ell}:= 
{\mathbf{PaB}}^{ell}_{{\mathbf k}}$, the isomorphisms $i_{O}$ by $i_{O}^{ell}$
and $i_{OO'}$ by $F(i_{OO'})$, where $F :\cC\to\cC_{ell}$ is the canonical functor. 
The claimed identity follows from $i_{O}^{ell}(F(x)) = i_{O}(F(x))_{ell}$, for 
$x\in \on{Aut}_{\cC}(O)$. 

3) follows from identity (\ref{identity:mu:4}), from the fact that $g=1_{\on{GT}(
{\mathbf k})}$ if $g_{ell}\in R_{ell}({\mathbf k})$, and from 
$b_{OO'}(1_{\on{GT}({\mathbf k})}) =1_{P_{n}({\mathbf k})}$, 
which follows from the first part of (\ref{identity:mu:2}).  
\hfill \qed\medskip 

\begin{proposition} \label{prop:resh}
1) There are morphisms 
$$
\mu_{O}^{gr} : \on{GRT}({\mathbf k})\to 
\on{Aut}(\on{exp}(\hat\t_{n}^{{\mathbf k}})\rtimes S_{n}) 
\quad \text{for each }O\in{\bf Pa}_{n}$$
and a map 
$$
\on{GRT}({\mathbf k})\times {\bf Pa}_{n}\times{\bf Pa}_{n}\to \on{exp}
(\hat{\mathfrak t}_{n}^{{\mathbf k}}), \quad
(g,O,O')\mapsto b^{gr}_{OO'}(g), $$
satisfying the analogues of the identities of Proposition \ref{prop:act:GT:categs}, 1). 

2) There are morphisms 
$$
\mu^{ell,gr}_{O} : \on{GRT}_{ell}({\mathbf k})\to \on{Aut}
(\on{exp}(\hat\t_{1,n}^{{\mathbf k}})\rtimes S_{n}) 
$$
for each $O\in{\mathbf{Pa}}_{n}$, satisfying the analogues of the 
identities of Proposition \ref{prop:act:GT:categs}, 2). 

3) The restriction $\mu_{O | R_{ell}^{gr}({\mathbf k})}^{ell,gr}$
is independent of $O$ and will be denoted 
$$
\mu_{ell}^{gr} : R_{ell}^{gr}({\mathbf k})\to
\on{Aut}(\on{exp}(\hat\t_{1,n}^{{\mathbf k}})
\rtimes S_{n}). 
$$
This morphism factors as $R_{ell}^{gr}({\mathbf k})\to
\on{Aut}(\hat\t_{1,n}^{{\mathbf k}})^{S_{n}}\to 
\on{Aut}(\on{exp}(\hat\t_{1,n}^{{\mathbf k}})
\rtimes S_{n})$. The Lie algebra morphism associated to the 
first factor is 
$$
{\mathfrak r}_{ell}^{gr}\to \on{Der}(\t_{1,n})^{S_{n}}, \quad 
(\alpha_{+},\alpha_{-})\mapsto (x_{i}^{\pm}\mapsto 
\alpha_{\pm}^{i,1\cdots\check i\cdots n}). 
$$
\end{proposition}

{\em Proof.} Similar to that of Proposition \ref{prop:act:GT:categs}. 
\hfill \qed\medskip 

\begin{remark} In \cite{CEE}, we introduced the Lie algebra ${\mathfrak d}
:= {\mathfrak d}_{+}\rtimes \mathfrak{sl}_{2}$, where ${\mathfrak d}_{+}$
is the $\mathfrak{sl}_{2}$-Lie algebra freely generated by the $\tilde\delta_{2m}$, 
$m\geq 0$, which is a highest weight vector for the simple $(2m+1)$-dimensional 
$\mathfrak{sl}_{2}$-module. There is a surjective morphism $\mathfrak{d}\to 
\mathfrak{b}_{3}$, which is the identity on $\mathfrak{sl}_{2}$ and given by 
$\tilde\delta_{2m}\mapsto\delta_{2m}$. In \cite{CEE}, we also constructed a 
morphism 
$$
\mathfrak{d}\to\on{Der}(\mathfrak{t}_{1,n})^{S_{n}}. 
$$
According to Proposition \ref{prop:resh}, 3), this morphism factors as $\mathfrak{d}
\to  \mathfrak{r}_{ell}^{gr}\to\on{Der}(\mathfrak{t}_{1,n})^{S_{n}}$. 
As $\on{im}(\mathfrak{d}\to\mathfrak{r}^{gr}_{ell}) = \mathfrak{b}_{3}$, 
the morphism from \cite{CEE} factors through $\mathfrak{b}_{3}$. 
\hfill \qed\medskip 
\end{remark}

Let us set $B_{n}^{gr}({\mathbf k}):= \on{exp}(\hat\t_{n}^{{\mathbf k}})\rtimes
S_n$, $B_{1,n}^{gr}({\mathbf k}):= \on{exp}(\hat\t_{1,n}^{{\mathbf k}})\rtimes
S_{n}$. We define $P_{n}^{gr}({\mathbf k})$, $P_{1,n}^{gr}({\mathbf k})$ 
as the ``pure'' versions of these groups (i.e., the kernels of their maps to $S_{n}$). 

\begin{proposition}
1) There is a family of isomorphisms $i_{O}^{\tilde \Phi} : 
B_{n}({\mathbf k})\to B_{n}^{gr}(\mathbf k)$ for each $\tilde \Phi := 
(\mu,\Phi)\in M({\mathbf k})$, and a family of maps 
$$
M({\mathbf k})\times {\bf Pa}_{n}\times{\bf Pa}_{n}\to 
P_{n}^{gr}({\mathbf k}), \quad 
(\tilde\Phi,O,O')\mapsto b_{OO'}^{gr}(\tilde\Phi), 
$$
such that 
$$
i_{O'}^{\tilde\Phi} = \on{Inn}(b_{OO'}^{gr}(\tilde\Phi)) \circ 
i_{O}^{\tilde\Phi}, \quad b_{OO''}^{gr}(\tilde\Phi) = 
b_{O'O''}^{gr}(\tilde\Phi)b_{OO'}^{gr}(\tilde\Phi), 
$$
$$
i_{O}^{\tilde\Phi}\circ\mu_{O}(g) = i_{O}^{g^{-1}*\tilde\Phi}, \quad 
\mu_{O}^{gr}(g_{gr}) \circ i_{O}^{\tilde\Phi} = i_{O}^{\tilde\Phi 
* g_{gr}^{-1}}, $$
where $g\in \on{GT}({\mathbf k})$, $g_{gr}\in \on{GRT}({\mathbf k})$. 

2) Each $e\in Ell({\mathbf k})$ gives rise to a family of isomorphisms 
$i_{O}^{ell,e} : B_{1,n}({\mathbf k})\to B_{1,n}^{gr}({\mathbf k})$, 
indexed by $O\in{\bf Pa}_{n}$. They satisfy 
$$
i_{O}^{\tilde\Phi}(b)_{ell}
= i^{ell,e}_{O}(b_{ell}), \quad 
i^{ell,e}_{O'} = \on{Inn}(b_{O,O'}(\tilde\Phi)_{ell}) \circ 
i^{ell,e}_{O} 
$$ for $b\in B_{n}({\mathbf k})$, if $\tilde\Phi := 
\on{im}(e\in Ell({\mathbf k})\to M({\mathbf k}))$, and 
$$
i^{ell,e}_{O}\circ \mu_{O}^{ell}(g) = i^{ell,g^{-1}*e}_{O}, \quad 
\mu^{ell,gr}_{O}(g_{gr})\circ i^{ell,e}_{O} 
= i^{ell,e*g_{gr}^{-1}}_{O},  
$$
for $g\in \on{GT}_{ell}({\mathbf k})$, $g_{gr}\in 
\on{GRT}_{ell}({\mathbf k})$. 

There is a commutative diagram 
$$
\xymatrix{ R_{ell}({\mathbf k}) \ar[r]^{\mu^{gr}_{ell}}\ar[d]_{i_{e}}
& \on{Aut}(B_{1,n}({\mathbf k})) \ar[d]^{(i_{e,O})^{*}}\\
R_{ell}^{gr}({\mathbf k})\ar[r]^{\mu^{ell}}& 
\on{Aut}(B_{1,n}^{gr}({\mathbf k})) }$$
\end{proposition}

{\em Proof.} Let $\cC^{gr}:= {\bf PaCD}_{{\mathbf k}}$, 
$\cC_{ell}^{gr}:= {\bf PaCD}_{ell,{\mathbf k}}$, 
then there are compatible functors $\cC\to \tilde\Phi*\cC_{gr}\simeq 
\cC_{gr}$, $\cC_{ell}\to e*\cC^{gr}_{ell}\simeq\cC^{gr}_{ell}$, where
in each case the first functor arises from universal properties and the second tensor 
forgets about the IBMC (or elliptic IBMC) structures. The statements follow from 
the compatibility of these functors with the actions of $\on{GT}_{ell}({\mathbf k})$, 
$\on{GRT}_{ell}({\mathbf k})$. \hfill \qed\medskip

\section{A family of elliptic associators, $\tau\mapsto e(\tau)$}
\label{sec:5}

In this section, we construct an analytic family of elliptic associators 
$\tau\mapsto e(\tau)$, 
indexed by the Poincar\'e half-plane. This family arises from the KZB connection 
(\cite{CEE}) and may therefore be viewed as an analogue of the KZ associator. 
We study various functional properties of this family: modular properties, behavior
at infinity, differential system. 

\subsection{The KZ associator}

Let $G_{0}(z)$, $G_{1}(z)$ be the analytic solutions of 
$$
G'(z) = ({A\over z}+{B\over{z-1}})G(z)
$$
in $]0,1[$, valued in $\on{exp}(\hat\f_2^\CC)$, with asymptotic behavior
$G_0(z)\sim z^A$ as $z\to 0$ and $G_1(z)\sim (1-z)^B$ as $z\to 1$. 
The KZ associator is defined by
$$
\Phi_{KZ} := G_1(z)^{-1}G_0(z)\in \on{exp}(\hat\f_{2}). 
$$ 
Then\footnote{We set $\on{i}:=\sqrt{-1}$.} 
$(2\pi\on{i},\Phi_{KZ})\in M(\CC)$ (\cite{Dr:Gal}). 

\subsection{Definition of $e(\tau) = (A(\tau),B(\tau))$}

Let $\HH := \{\tau\in\CC|\Im(\tau)>0\}$ be the Poincar\'e half-plane. 
Let $(z,\tau)\mapsto \theta(z|\tau)$ be the holomorphic function on 
$\CC\times\HH$, such that $\theta(z+1|\tau)=-\theta(z|\tau)
=\theta(-z|\tau)$, $\theta(z+\tau|\tau) = - e^{2\pi\on{i}z\tau+\on{i}
\pi\tau}\theta(z|\tau)$, $\{z|\theta(z|\tau)=0\} = \ZZ+\tau\ZZ$, 
$\partial_z\theta(0|\tau)=1$. 

For $\tau\in\HH$, let $F(z|\tau)$ be the holomorphic function on 
$\{z = a+b\tau|a,b\in\RR$, $a$ or $b\in]0,1[\}$, valued in 
$\on{exp}(\hat\t_{1,2}^{\CC})\simeq \on{exp}(\hat\f_{2}^{\CC})$, such that 
$$
\partial_z F(z|\tau) = - {{\theta(z+ 
\on{ad}x|\tau)\on{ad}x}
\over{\theta(z|\tau)\theta( 
\on{ad}x|\tau)}}(y) \cdot 
F(z|\tau) \text{\ and\ }F(z|\tau)
\sim (-2\pi\on{i}z)^{t} \text{\ as\ }z\to 0;   
$$ 
here $x:=x_2^+$, $y:= x_2^-$, $t:= t_{12}$. We then set 
$$
A(\tau):= F(z|\tau)^{-1}F(z+1|\tau), \quad 
B(\tau):= F(z|\tau)^{-1}e^{2\pi\on{i}x} F(z+\tau|\tau).  
$$

\subsection{Algebraic properties of $e(\tau)$}

We set $Ell_{KZ}:= 
Ell(\CC)\times_{M(\CC)}\{(2\pi\on{i},\Phi_{KZ})\}$. 

\begin{proposition}
$\tau\mapsto e(\tau):= (A(\tau),B(\tau))$ is an analytic map 
$\HH\to Ell_{KZ}$. 
\end{proposition}

{\em Proof.} In \cite{CEE}, Section 4.4, we introduced $\tilde A,\tilde B
\in \on{exp}(\hat\t_{1,2}^{\CC})$. We set $\tilde A_+:= \tilde A$, 
$\tilde A_-:= \tilde B$, $A_{+}(\tau):= A(\tau)$, $A_{-}(\tau):=
B(\tau)$, then  
$$
A_\pm(\tau) = \on{Ad}((-2\pi\on{i})^{-t})(\tilde A_{\pm}), 
$$ 
So $(A_+(\tau),A_-(\tau))$ satisfies  
(22), (23), (26) in \cite{CEE}. (22), (23) imply that 
$(A_+(\tau),A_-(\tau))$ satisfies (\ref{def:ell:ass:1}).
(26) implies that 
$$
(A_-(\tau)^{12,3}\{\Phi^{-1}\}(A_-(\tau)^{1,23})^{-1}\{\Phi\},
A_+(\tau)^{12,3}) = \{\Phi^{-1}e^{2\pi\on{i}t_{23}}\Phi\}
$$ and using (23) in 
\cite{CEE}, we rewrite this as    
$$
(\{e^{-\on{i}\pi t_{12}}\Phi^{3,2,1}\}A_-(\tau)^{2,13}
\{\Phi^{2,1,3}e^{-\on{i}\pi t_{12}}\},A_+(\tau)^{12,3})
= \{\Phi^{-1}e^{2\pi\on{i}t_{23}}\Phi\}, 
$$
which as in the proof of Proposition \ref{prop:sigma} implies 
that $(A_+(\tau),A_-(\tau))$ satisfies (\ref{def:ell:ass:2}).
\hfill\qed\medskip 

\subsection{Analytic properties of $e(\tau)$}

\begin{proposition} One has 
$$
2\pi\on{i}{\partial\over{\partial\tau}}e(\tau) = 
e(\tau)*(-e_{-}-\sum_{k\geq 0} (2k+1)G_{2k+2}(\tau)\delta_{2k}), 
$$
where $G_{k}(\tau)$ are the Eisenstein 
series defined by 
$$
G_{k}(\tau) = \sum_{a\in(\ZZ + \tau\ZZ) - \{0\}}a^{-k}\text{ for }k\text{ even } \geq4, 
\quad G_{2}(\tau) = \sum_{m\in\ZZ}(\sum'_{n}(n+m\tau)^{-2}),  
$$
where $\sum'$ means $\sum_{n\in\ZZ}$ if $m\geq 0$ and $\sum_{n\in\ZZ - \{0\}}$
if $m=0$ (notation as in (\ref{delta:2n}), (\ref{sl2})). 
\end{proposition}

{\em Proof.} $R_{ell}(\CC)\subset \on{Aut}(\hat\f_{2}^{\CC})^{op}$
acts from the right on $Ell_{KZ}$ by $(A_{+},A_{-})*(u_{+},u_{-}):= 
(A_{+}(u_{+},u_{-}),A_{-}(u_{+},u_{-}))$. The same formula defines
a left action of $R_{ell}(\CC)^{op}\subset \on{Aut}(\hat\f_{2}^{\CC})$ on 
$Ell_{KZ}$. To prove that 
$$
2\pi\on{i}\partial_{\tau}e(\tau) = e(\tau)*x(\tau)$$
for $x(\tau)\in \hat{\mathfrak{r}}_{ell}^{\CC}\subset 
\on{Der}(\hat\f_{2}^{\CC})^{op}$, it
therefore suffices to prove that 
$$
2\pi\on{i}\partial_{\tau}A(\tau) = x(\tau)(A(\tau)), \quad 
2\pi\on{i}\partial_{\tau}B(\tau) = x(\tau)(B(\tau)), 
$$
where  $x(\tau)$ is now viewed as an element of $\on{Der}(\hat\f_{2}^{\CC})$. 

In \cite{CEE}, Lemma 23, we constructed a function $F^{(2)}(z|\tau)$, defined on 
$\{(z,\tau)\in\CC\times\HH | z = a + b\tau, (a,b)\in
]0,1[\times\RR \cup \RR\times ]0,1[\}$ and valued in 
$\on{exp}(\hat{\mathfrak{f}}_{2}^{\CC})\rtimes
\on{Aut}(\hat{\mathfrak{f}}_{2}^{\CC})$,
such that 
$$
\partial_{z}F^{(2)}(z|\tau) = -{{\theta(z+\on{ad}x|\tau)\on{ad}x}
\over{\theta(z|\tau)\theta(\on{ad}x|\tau)}}(y)\cdot F^{(2)}(z|\tau), 
$$
\begin{align*}
2\pi\on{i} {\partial\over{\partial \tau}} F^{(2)}(z|\tau) & =
-\Big( e_-+\sum_{k\geq 0}
(2k+1)G_{2k+2}(\tau)\dot\delta_{2k}^{(2)} - g(z,\on{ad}x|\tau)(t) \Big) 
\cdot F^{(2)}(z|\tau) \\
 & = -\Big( e_-+\sum_{k\geq 0}
(2k+1)G_{2k+2}(\tau)\delta_{2k}^{(2)} - g(z|\tau)(t) \Big) 
\cdot F^{(2)}(z|\tau) , 
\end{align*}
and 
$F^{(2)}(z|\tau)\sim z^{t}\on{exp}({{-\tau}\over{2\pi\on{i}}}
(e_- +\sum_{k\geq 0} 2(2k+1)\zeta(2k+2)\delta_{2k}^{(2)}))$ 
as $z\to 0$ and $\tau\to\on{i}\infty$. Here $g(z,x|\tau) 
= {{\theta(z+x|\tau)}\over{\theta(z|\tau)\theta(x|\tau)}}
({\theta'\over\theta}(z+x|\tau)-{\theta'\over\theta}(z|\tau))+{1\over{x^{2}}}$,
and $g(z|\tau):= g(z,\on{ad}x|\tau)(t)- g(0,\on{ad}x|\tau)(t)$; in the notation of 
{\it loc.~cit.,} $e_- = \Delta_0$.

These conditions imply that the image of $F^{(2)}(z|\tau)$ in 
$\on{Aut}(\hat{\mathfrak{f}}_{2}^{\CC})$ is independent of $z$. Then 
$$
A_{z_{0}}^{z_{1}}(\tau):= F^{(2)}(z_{1}|\tau)F^{(2)}(z_{0}|\tau)^{-1}
\in\on{exp}(\hat{\mathfrak f}_{2}^{\CC})
$$
and satisfies 
$$
2\pi\on{i}\partial_{\tau}A_{z_{0}}^{z_{1}}(\tau) = 
- \big(e_{-} + \sum_{k\geq 0}(2k+1)G_{2k+2}(\tau)
\delta_{2k}\big)(A_{z_{0}}^{z_{1}}(\tau)) 
+ g(z_{1}|\tau)\cdot A_{z_{0}}^{z_{1}}(\tau) 
- A_{z_{0}}^{z_{1}}(\tau) \cdot g(z_{0}|\tau). 
$$
The function $F(z|\tau)$, basic to the definition of $(A(\tau),B(\tau))$, 
is related to the function $F^{(2)}(z|\tau)$ by 
$F^{(2)}(z|\tau)=F(z|\tau)\varphi(\tau)$, where $\varphi(\tau)$
takes values in $\on{exp}(\hat{\mathfrak{f}}_{2}^{\CC})\rtimes
\on{Aut}(\hat{\mathfrak{f}}_{2}^{\CC})$, as both satisfy the same differential 
equation in $z$. It follows that 
$$
A_{z_{0}}^{z_{1}}(\tau) = F(z_{1}|\tau)F(z_{0}|\tau)^{-1}. 
$$
Therefore $A(\tau) = F(z|\tau)^{-1}A_{z}^{z+1}(\tau)F(z|\tau)$. 
In the limit $z\to 0$, this gives
$$
A(\tau) = \on{lim}_{\epsilon\to 0} (-2\pi\on{i}\epsilon)^{-\on{ad}t}
\big( A_{\epsilon}^{1+\epsilon}(\tau) \big). 
$$
$\varepsilon$ being fixed, $(-2\pi\on{i}\epsilon)^{-\on{ad}t}
\big( A_{\epsilon}^{1+\epsilon}(\tau) \big)$ satisfies the 
same differential equation in $\tau$ as $A_{z_{0}}^{z_{1}}(\tau)$, 
with $g(z_{0}|\tau)$ replaced by $(-2\pi\on{i}\epsilon)^{-\on{ad}(t)}(
g(\epsilon|\tau))$ and $g(z_{1}|\tau)$ replaced by 
$(-2\pi\on{i}\epsilon)^{-\on{ad}(t)}(g(1+\epsilon|\tau))$, which 
both tend to $0$ as $\epsilon\to 0$. It follows that these terms disappear
from the differential equation satisfied by $A(\tau)$, so 
$$
2\pi\on{i}\partial_{\tau}A(\tau) = - (e_{-} + 
\sum_{k\geq 0}(2k+1)G_{2k+2}(\tau)\delta_{2k})(A(\tau)). 
$$
Similarly, $B(\tau) = F(z|\tau)^{-1}e^{2\pi\on{i}x} A_{z}^{z+\tau}(\tau)
F(z|\tau)$, hence 
$$
B(\tau) = \on{lim}_{\epsilon\to 0} (-2\pi\on{i}\epsilon)^{-t}
e^{2\pi\on{i}x} A_{\epsilon}^{\tau+\epsilon}(\tau)
(-2\pi\on{i}\epsilon)^{t}. 
$$
One computes  
\begin{eqnarray*}
&&  \partial_{\tau} (A_{\epsilon}^{\tau+\epsilon}(\tau)) = 
{{-1}\over{2\pi\on{i}}} (e_- +\sum_{k\geq 0}
(2k+1)G_{2k+2}(\tau)\delta_{2k})(A_{\epsilon}^{\tau+\epsilon}(\tau))\\
&& 
 +\Big( 
 {1\over{2\pi\on{i}}} g(\tau+\epsilon|\tau) - {{\theta(\tau+\epsilon+\on{ad}x
 |\tau)\on{ad}x}\over{\theta(\tau+\epsilon|\tau)\theta(\on{ad}x|\tau)}}(y)
 \Big)  A_{\epsilon}^{\tau+\epsilon}(\tau) -  A_{\epsilon}^{
 \tau+\epsilon}(\tau) {1\over {2\pi\on{i}}}g(\epsilon|\tau). 
\end{eqnarray*}
So $X_{\epsilon}(\tau):= (-2\pi\on{i}\epsilon)^{-t}e^{2\pi\on{i}x} A_{\epsilon}^{\tau
+\epsilon}(\tau)(-2\pi\on{i}\epsilon)^{t}$ satisfies ($\epsilon$ being fixed)
\begin{eqnarray*}
&& 2\pi\on{i}\partial_{\tau}(X_{\epsilon}(\tau))
 = - (e_- +\sum_{k\geq 0}
(2k+1)G_{2k+2}(\tau)\delta_{2k})(X_{\epsilon}(\tau) ) - 
X_{\epsilon}(\tau)
\cdot \big(
(-2\pi\on{i}\epsilon)^{-t}g(\epsilon|\tau)
(-2\pi\on{i}\epsilon)^{t}\big)
\\
&&  + \Big( \on{Ad}((-2\pi\on{i}\epsilon)^{-t}
e^{2\pi\on{i}x})\big( 
g(\tau+\epsilon|\tau) - 2\pi\on{i}
{{\theta(\tau+\epsilon+\on{ad}x|\tau)\on{ad}x}
\over{\theta(\tau+\epsilon|\tau)\theta(\on{ad}x|\tau)}}(y)\big)\\
&&  - (-2\pi\on{i}\epsilon)^{-t}e^{2\pi\on{i}x}
(e_- +\sum_{k\geq 0}
(2k+1)G_{2k+2}(\tau)\delta_{2k})(e^{-2\pi
\on{i}x})(-2\pi\on{i}\epsilon)^{t}\Big) \cdot (
X_{\epsilon}(\tau)). 
\end{eqnarray*}
Identity (7) in \cite{CEE} implies that the parenthesis in the two last lines 
equals $\on{Ad}((-2\pi\on{i}\epsilon)^{-t})(g(\epsilon|\tau))$. As before, 
we get in the limit $\epsilon\to 0$ 
$$
2\pi\on{i}\partial_{\tau}B(\tau) = 
-(e_-+\sum_{k\geq 0}(2k+1)G_{2k+2}(\tau)\delta_{2k})(B(\tau)).
$$
\hfill \qed\medskip 

\begin{proposition}
$$
\sigma(\Phi_{KZ})_{ | \begin{smallmatrix}
x\mapsto 2\pi\on{i}x, \\
 y\mapsto (2\pi\on{i})^{-1}y
\end{smallmatrix} } 
= \lim_{\tau\to\on{i}\infty}
e(\tau)* \on{exp}({\tau\over{2\pi\on{i}}}
(e_- +\sum_{k\geq 0}(2k+1)\zeta(2k+2)\delta_{2k})). 
$$
\end{proposition}

{\em Proof.} In \cite{CEE} (proof of Prop. 24 and Lemma 29), is it 
proved that 
$$
A(\tau) = \Phi_{KZ}(\tilde y,t)e^{2\pi\on{i}\tilde y}
 \Phi_{KZ}(\tilde y,t)^{-1} + O(e^{2\pi\on{i}\tau}), 
$$
$$
B(\tau) = e^{\on{i}\pi t}\Phi_{KZ}(-\tilde y-t,t)
e^{2\pi\on{i}x}e^{2\pi\on{i}\tilde y\tau}
 \Phi_{KZ}(\tilde y,t)^{-1} + O(e^{2\pi\on{i}\tau(1-\epsilon)}),$$
for any $\epsilon>0$, where
$$
\tilde y := -{{\on{ad}x}\over{e^{2\pi\on{i}\on{ad}x}-1}}(y).
$$
Let $(A_{pol}(\tau),B_{pol}(\tau))$ be the principal parts of the right 
sides of these equalities; $A_{pol}(\tau)$ is constant in $\tau$, while 
each coordinate of $B_{pol}(\tau)$ in a basis of $U({\mathfrak{f}}_{2})$
is a polynomial in $\tau$. 

It is proved in \cite{CEE} that $\tilde y,t$ are in the kernel of 
$e_{-} + \sum_{k\geq 0}(2k+1)\zeta(2k+2)\delta_{2k}$, while 
$$
\on{exp} \big( {\tau\over{2\pi\on{i}}}(e_{-} + \sum_{k\geq 0}(2k+1)
\zeta(2k+2)\delta_{2k} \big) (e^{2\pi\on{i}x}e^{2\pi\on{i}\tau\tilde y}) = 
e^{2\pi\on{i}x}. 
$$
It follows that 
$$
\on{exp}\big( {\tau\over{2\pi\on{i}}}(e_{-} + \sum_{k\geq 0}(2k+1)
\zeta(2k+2)\delta_{2k} \big)(A_{pol}(\tau),B_{pol}(\tau)) = 
\sigma(\Phi_{KZ})\sigma(\Phi_{KZ})_{ | \begin{smallmatrix}
x\mapsto 2\pi\on{i}x, \\
 y\mapsto (2\pi\on{i})^{-1}y
\end{smallmatrix} }, 
$$
which implies that statement. \hfill \qed\medskip 

Note that the operation $e\mapsto e_{ | \begin{smallmatrix}
x\mapsto 2\pi\on{i}x, \\
 y\mapsto (2\pi\on{i})^{-1}y
\end{smallmatrix} }$ amounts the the action of 
$\on{diag}(2\pi\on{i},(2\pi\on{i})^{-1})\subset 
\on{SL}_{2}(\CC) \subset R_{ell}(\CC)$ on 
$Ell_{KZ}$. 

\subsection{Modularity properties of $e(\tau)$}

We now describe the behavior of the map $\tau\mapsto e(\tau)$
under the action of $\on{SL}_2(\ZZ)$ on $\HH$.

Define $\on{log} : \CC^{\times}\to \CC$ by the condition that 
its image is contained in $\RR + \on{i}[-\pi,\pi[$. We define group 
morphisms $t : \CC^{\times}\to \on{SL}_{2}(\CC)$ and 
$n_{\pm} : \CC\to \on{SL}_{2}(\CC)$ by $t(\lambda):= 
\bigl( \begin{smallmatrix} \lambda^{-1}& 0 \\ 0 & \lambda 
\end{smallmatrix}\bigl)$, $n_{+}(a) := 
\bigl( \begin{smallmatrix} 1 & 0 \\ a & 1\end{smallmatrix}\bigl)$, 
$n_{-}(a) := \bigl( \begin{smallmatrix} 1 & a \\ 0 & 1\end{smallmatrix}\bigl)$. 

\begin{proposition} \label{prop:5:4}
1) There is a unique map 
$$
f : B_{3}\times\HH\to\CC, 
$$
such that 
$$
f(\sigma_{1},\tau) = 0, \quad 
f(\sigma_{2},\tau) = -\on{log}\big({{-1}\over{\tau-1}}\big)
$$
and with the cocycle property $f(gg',\tau) = f(g,\overline{g'}
\cdot \tau) + f(g',\tau)$, 
where $g\mapsto \overline g$ is the morphism $B_{3}\to 
\on{SL}_{2}(\ZZ)$ and the action on $\on{SL}_{2}(\ZZ)$ on $\HH$ is 
$\bigl( \begin{smallmatrix} \alpha& \beta \\ \gamma & \delta 
\end{smallmatrix}\bigl) \cdot\tau = {{\alpha\tau+\beta}\over
{\gamma\tau+\delta}}$. 

2) For any $g\in B_{3}$ and $\tau\in\HH$, one has 
\begin{equation} \label{id:e:g}
e(\overline g\cdot \tau) = \on{Ad}(e^{f(g,\tau)t})
\Big(g*\big(a(\overline g,\tau)\bullet e(\tau)\big)\Big), 
\end{equation}
where: 

$\bullet$ for $\alpha\in\CC$, $\on{Ad}(e^{\alpha t})$ is the self-map of 
$Ell_{KZ}$ given by 
$\on{Ad}(e^{\alpha t})(e) := (e^{\alpha t}Ae^{-\alpha t}, 
e^{\alpha t}Be^{-\alpha t})$ for $e = (A,B)$; 

$\bullet$ $a : \on{SL}_{2}(\ZZ)\times\HH\to\on{SL}_{2}(\CC)$ is 
given by $a(\overline g,\tau) = 
\bigl( \begin{smallmatrix} \gamma\tau+\delta& 0 \\ 2\pi\on{i}\gamma & 
(\gamma\tau+\delta)^{-1} \end{smallmatrix}\bigl) = 
n_{+}({{2\pi\on{i}\gamma}\over{\gamma\tau+\delta}})t((\gamma\tau+\delta
)^{-1})$ if $\overline g = 
\bigl( \begin{smallmatrix} \alpha& \beta \\ \gamma & \delta 
\end{smallmatrix}\bigl)$; 

$\bullet$ $*$ and $\bullet$ are the commuting left actions of 
$B_{3} = R_{ell}\subset R_{ell}(\CC)$ and $\on{SL}_{2}(\CC)
\subset R_{ell}^{gr}(\CC)^{op}$ on $Ell_{KZ}$, given as follows: 

$\bullet$ for $e = (A,B)\in Ell_{KZ}$ and $g\in B_{3}$, 
$g*e:= (\theta_{g}(a)_{|(a,b)\mapsto (A,B)}, \theta_{g}(b)_{|(a,b)\mapsto (A,B)})$, 
where $\theta : B_{3}\to \on{Aut}(F_{2})$ is the action  
of $B_{3}$ on the free group $F_{2}$ generated by $a,b$, and $x\mapsto 
x_{|(a,b)\mapsto (A,B)}$ is the morphism $F_{2}\to\on{exp}(\hat\f_{2}^{\CC})$, 
given by $a,b\mapsto A,B$; 

$\bullet$ for $e = (A,B)\in Ell_{KZ}$ and $a\in \on{SL}_{2}(\CC)$, 
$a\bullet e:= (\alpha_{a}(A),\alpha_{a}(B))$, where 
$\alpha : \on{SL}_{2}(\CC)\to \on{Aut}(\on{exp}(\hat\f_{2}^{\CC}))^{op}$
is induced by 
$\alpha_{a} \bigl( \begin{smallmatrix} x \\ y 
\end{smallmatrix}\bigl) = 
\bigl( \begin{smallmatrix} p & q \\ r & s \end{smallmatrix}\bigl)
\bigl( \begin{smallmatrix} x \\ y 
\end{smallmatrix}\bigl)$ if $a = \bigl( \begin{smallmatrix} p& q \\ r & s 
\end{smallmatrix}\bigl)$. 
\end{proposition}

\begin{remark}
Let $g\in B_{3}$ and  $\overline g = \bigl( \begin{smallmatrix} \alpha & 
\beta \\ \gamma & \delta \end{smallmatrix}\bigl)$ is its image in 
$\on{SL}_{2}(\ZZ)$, then $\on{exp}f(g,\tau) = \gamma\tau+\delta$ 
for any $\tau\in\HH$. 
\end{remark}

\begin{remark}
For $g = (\sigma_{1}\sigma_{2})^{6}$ (a generator of the kernel of $B_{3}\to
\on{SL}_{2}(\ZZ)$), 
$$
g*e = (\on{Ad}(B,A)(A),\on{Ad}(B,A)(B)) = 
(\on{Ad}(e^{2\pi\on{i}t})(A), \on{Ad}(e^{2\pi\on{i}t})(B)),
$$ 
while $f(g,\tau) = -2\pi\on{i}$. One checks this way that the r.h.s. of (\ref{id:e:g})
does not depend of the choice of a lift $g$ of $\overline g$ to $B_{3}$. 
\end{remark}

{\em Proof.} Statement 1) can be checked using the presentation of $B_{3}$. 
It follows from the cocycle identity for $f(g,\tau)$ and from the 
cocycle identity 
$$
a(hh',\tau) = a(h',\tau)a(h,h'\cdot \tau), \quad h,h'\in\on{SL}_{2}(\ZZ), 
\tau\in\HH
$$
that $\Gamma:= \{g\in B_{3} |$ identity $(\ref{id:e:g})$ 
holds for any $\tau\in\HH\}$ is a subgroup of $B_{3}$. So statement 
2) follows from its particular cases $g = \sigma_{1}$, $g = \sigma_{1}\sigma_{2}
\sigma_{1}$. 

Recall that 
$$
A(\tau) = \lim_{\epsilon\to 0^{+}} 
(-2\pi\on{i}\epsilon)^{-t} 
A_{\epsilon}^{1+\epsilon}(\tau)
(-2\pi\on{i}\epsilon)^{t}, \quad 
B(\tau) = \lim_{\epsilon\to 0^{+}} (-2\pi\on{i}\epsilon)^{-t} 
e^{2\pi\on{i}x}A_{\epsilon}^{\tau+\epsilon}(\tau)
(-2\pi\on{i}\epsilon)^{t},
$$ 
where $A_{z_{0}}^{z_{1}}(\tau)$ be the solution of 
$\partial_{z_{1}}A_{z_{0}}^{z_{1}}(\tau) = K(z_{1}|\tau)
A_{z_{0}}^{z_{1}}(\tau)$ such that $A_{z}^{z}(\tau)=1$, where 
$K(z|\tau)= - {{\theta(z+\on{ad}x|\tau)
\on{ad}x}\over{\theta(z|\tau)\theta(\on{ad}x|\tau)}}(y)$ 
and where the chosen branches of $A_{z_{0}}^{z_{1}}(\tau)$
are as in Fig.~1. 

The identity $K(z|\tau)=K(z|\tau+1)$ implies 
$A_{\epsilon}^{1+\epsilon}(\tau+1)=A_{\epsilon}^{1+\epsilon}(\tau)$, 
and using the decomposition of Fig.~2, it also implies 
$A_{\epsilon}^{\tau+1+\epsilon}(\tau+1)
=A_{1+\epsilon}^{\tau+1+\epsilon}(\tau)A_{\epsilon}^{1+\epsilon}(\tau)
=A_{\epsilon}^{\tau+\epsilon}(\tau)A_{\epsilon}^{1+\epsilon}(\tau)$. So 
$A(\tau+1)=A(\tau)$, $B(\tau+1)=B(\tau)A(\tau)$, 
so $e(\tau+1)=\sigma_{1}*e(\tau)$, which shows (\ref{id:e:g}) 
in the case $g=\sigma_{1}$. 

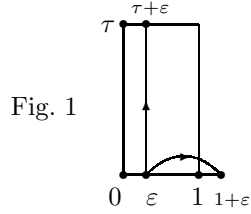
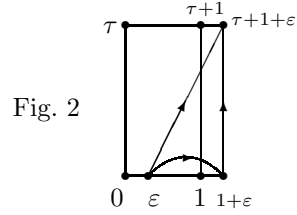
\setlength{\unitlength}{1mm}
\begin{picture}(60,40)(-30,0)
\put(0,10){\line(1,0){13}}
\put(0,10){\line(0,1){20}}
\put(0,30){\line(1,0){10}}
\put(10,10){\line(0,1){20}}
\put(0,10){\circle*{1}}
\put(3,10){\circle*{1}}
\put(10,10){\circle*{1}}
\put(13,10){\circle*{1}}
\put(0,30){\circle*{1}}
\put(3,30){\circle*{1}}
\put(3,10){\line(0,1){20}}
\put(3,10){\vector(0,1){10}}
\qbezier(3,10)(8,15)(13,10)
\put(8,12.4){\vector(1,0){1}}
\put(-2,6){$0$}
\put(3,6){$\eps$}
\put(9,6){$1$}
\put(12,6){${\scriptstyle{1+\eps}}$}
\put(-3,29){$\tau$}
\put(1,31.5){${\scriptstyle{\tau+\eps}}$}
\put(-15,18){Fig.~1}
\put(50,10){\line(1,0){13}}
\put(50,10){\line(0,1){20}}
\put(50,30){\line(1,0){13}}
\put(60,10){\line(0,1){20}}
\put(50,10){\circle*{1}}
\put(53,10){\circle*{1}}
\put(60,10){\circle*{1}}
\put(63,10){\circle*{1}}
\put(50,30){\circle*{1}}
\put(63,30){\circle*{1}}
\put(60,30){\circle*{1}}
\put(63,10){\line(0,1){20}}
\put(63,10){\vector(0,1){10}}
\put(53,10){\vector(1,2){5}}
\put(53,10){\line(1,2){10}}
\qbezier(53,10)(58,15)(63,10)
\put(58,12.4){\vector(1,0){1}}
\put(48,6){$0$}
\put(53,6){$\eps$}
\put(59,6){$1$}
\put(62,6){${\scriptstyle{1+\eps}}$}
\put(47,29){$\tau$}
 \put(58,31.5){${\scriptstyle{\tau+1}}$} 
 \put(64,30){${\scriptstyle{\tau+1+\eps}}$}
\put(35,18){Fig.~2}
\end{picture}
 
Let $w:= -z/\tau$, then 
$$
{\partial\over{\partial w}} - K(w|{{-1}\over\tau}) = 
-\tau e^{-2\pi\on{i}zx}
\Big( {\partial\over{\partial z}}
- \bigl(\begin{smallmatrix} -\tau & 0 \\ 
-2\pi\on{i} & {{-1}\over\tau} \end{smallmatrix}\bigl) 
\bullet K(z|\tau) \Big)
e^{2\pi\on{i}zx}.  
$$
So  
$$
A_{w_{0}}^{w_{1}}({{-1}\over\tau}) = e^{-2\pi\on{i} xz_{1}} \cdot 
\bigl(\begin{smallmatrix} -\tau & 0 \\ 
-2\pi\on{i} & {{-1}\over\tau} \end{smallmatrix}\bigl)
\bullet (A_{z_{0}}^{z_{1}}(\tau)) \cdot 
e^{2\pi\on{i} xz_{0}}.
$$ 

Then 
\begin{eqnarray*}
&& A({{-1}\over\tau}) = \lim_{\eps\to 0^{+}} (-2\pi\on{i}\epsilon)^{-t}
A_{\epsilon}^{1+\epsilon}({{-1}\over\tau})(-2\pi\on{i}\epsilon)^{t} \\
&& = \lim_{\eps\to 0^{+}}
(-2\pi\on{i}\epsilon)^{-t}\cdot 
\bigl(\begin{smallmatrix} -\tau & 0 \\ 
-2\pi\on{i} & {{-1}\over\tau} \end{smallmatrix}\bigl) \bullet
 \Big(e^{-2\pi\on{i}(1+\epsilon)x}
A_{-\eps\tau}^{-\tau-\eps\tau}(\tau)
e^{2\pi\on{i}\epsilon x}\Big) \cdot 
(-2\pi\on{i}\epsilon)^{t}
\\
&& = \on{exp}(-\on{log}({{-1}\over\tau})t)
\cdot 
\bigl(\begin{smallmatrix} -\tau & 0 \\ 
-2\pi\on{i} & {{-1}\over\tau} \end{smallmatrix}\bigl) \bullet 
B(\tau)^{-1} \cdot 
\on{exp}(\on{log}({{-1}\over\tau})t) , 
\end{eqnarray*}
see Fig.~3; and 
\begin{eqnarray*}
&&  B({{-1}\over\tau}) = \lim_{\epsilon\to 0^{+}} (-2\pi\on{i}\epsilon)^{-t}
e^{2\pi\on{i}x}
A_{\epsilon}^{{{-1}\over\tau}+\epsilon}({{-1}\over\tau})
(-2\pi\on{i}\epsilon)^{t} \\
&& = 
\lim_{\eps\to 0^{+}}
(-2\pi\on{i}\epsilon)^{-t}e^{2\pi\on{i}x\tau\epsilon}
\cdot
\bigl(\begin{smallmatrix} -\tau & 0 \\ 
-2\pi\on{i} & {{-1}\over\tau} \end{smallmatrix}\bigl) \bullet 
A_{-\tau\epsilon}^{1-\tau\epsilon}(\tau)
\cdot e^{-2\pi\on{i}x\tau\epsilon}
(-2\pi\on{i}\epsilon)^{-t}
\\
 & & = 
  \\
&& = 
\on{exp}(-\on{log}({{-1}\over\tau})t)
 \cdot 
\bigl(\begin{smallmatrix} -\tau & 0 \\ 
-2\pi\on{i} & {{-1}\over\tau} \end{smallmatrix}\bigl) \bullet
\Big(
\lim_{\eps\to 0^{+}}
(2\pi\on{i}\tau\eps)^{-t} 
A_{-\tau\eps}^{1-\tau\eps}(\tau)
(2\pi\on{i}\tau\eps)^{t} \Big) 
 \cdot    
 \on{exp}(\on{log}({{-1}\over\tau})t)
\\
&& = 
\on{exp}(-\on{log}({{-1}\over\tau})t)
\cdot \bigl(\begin{smallmatrix} -\tau & 0 \\ 
-2\pi\on{i} & {{-1}\over\tau} \end{smallmatrix}\bigl) \bullet
( B(\tau)A(\tau)B(\tau)^{-1})\cdot 
\on{exp}(\on{log}({{-1}\over\tau})t), 
\end{eqnarray*}
see Fig.~4. It follows that 
$$
e({{-1}\over\tau}) = \on{Ad}\Big(
\on{exp}(-\on{log}({{-1}\over\tau})t)\Big) 
\Big( 
\sigma_{1}\sigma_{2}\sigma_{1} * \Big(
\bigl(\begin{smallmatrix} -\tau & 0 \\ 
-2\pi\on{i} & {{-1}\over\tau} \end{smallmatrix}\bigl) 
\bullet e(\tau)\Big) 
\Big) $$
The result for $g=\sigma_{1}\sigma_{2}\sigma_{1}$ then follows. 

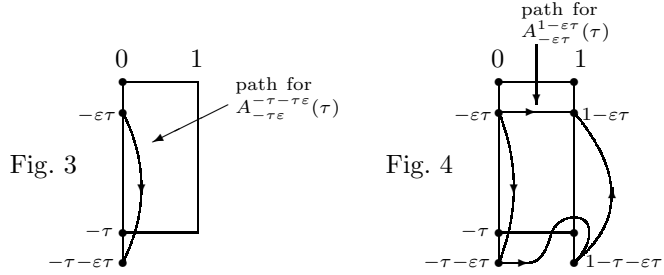
\setlength{\unitlength}{1mm}
\begin{picture}(60,40)(-30,0)
\put(0,10){\line(1,0){10}}
\put(0,6){\line(0,1){24}}
\put(0,30){\line(1,0){10}}
\put(10,10){\line(0,1){20}}
\put(0,6){\circle*{1}}
\put(0,10){\circle*{1}}
\put(0,26){\circle*{1}}
\put(0,30){\circle*{1}}
\qbezier(0,26)(5,16)(0,6)
\put(2.5,16){\vector(0,-1){1}}
\put(-1,32){$0$}
\put(9,32){$1$}
\put(-15,18){Fig.~3}
\put(-6,25){${\scriptstyle{-\eps\tau}}$}
\put(-5,9.5){${\scriptstyle{-\tau}}$}
\put(-10,5.5){${\scriptstyle{-\tau-\eps\tau}}$}
\put(14,27){\vector(-2,-1){10}}
\put(15,29){${\scriptstyle{\on{path\ for}}}$}
\put(15,26){${\scriptstyle{A_{-\tau\eps}^{-\tau-\tau\eps}(\tau)}}$}
\put(50,10){\line(1,0){10}}
\put(50,6){\line(0,1){24}}
\put(50,30){\line(1,0){10}}
\put(60,6){\line(0,1){24}}
\put(50,30){\circle*{1}}
\put(50,26){\circle*{1}}
\put(50,10){\circle*{1}}
\put(50,6){\circle*{1}}
\put(60,30){\circle*{1}}
\put(60,26){\circle*{1}}
\put(60,10){\circle*{1}}
\put(60,6){\circle*{1}}
\put(50,26){\vector(1,0){5}}
\put(50,26){\line(1,0){10}}
\put(55,35){\vector(0,-1){8}}
\qbezier(50,26)(54,16)(50,6)
\put(52,16){\vector(0,-1){1}}
\put(50,6){\vector(1,0){4}}
\qbezier(54,6)(56,6)(57,10)
\qbezier(57,10)(58,12)(60,12)
\qbezier(60,12)(62,12)(62,10)
\qbezier(62,10)(62,8)(60,6)
\qbezier(60,6)(70,14)(60,26)
\put(65.1,15){\vector(0,1){1}}
\put(53,39){${\scriptstyle{\on{path\ for}}}$}
\put(53,36){${\scriptstyle{A_{-\eps\tau}^{1-\eps\tau}(\tau)}}$}
\put(35,18){Fig.~4}
\put(49,32){$0$}
\put(60,32){$1$}
\put(44,25){${\scriptstyle{-\eps\tau}}$}
\put(45,9.5){${\scriptstyle{-\tau}}$}
\put(40,5.5){${\scriptstyle{-\tau-\eps\tau}}$}
\put(61,25){${\scriptstyle{1-\eps\tau}}$}
\put(61.5,5.5){${\scriptstyle{1-\tau-\eps\tau}}$}
\end{picture}

\hfill \qed\medskip 

\section{Computations of Zariski closures}

The action of the mapping class group $B_{3}$ in genus one on the braid 
groups in genus one (see (\ref{act:B3:ellbraids})) restricts to an action on the pure 
braid subgroups. In this section, we compute the Zariski closure of the image 
of $B_{3}$ in the automorphism groups of their prounipotent completions.
This computation relies on the relation between the action of $\on{GT}_{ell}(-)$ on these
prounipotent completions and its the graded counterpart (Section 4), 
and on the properties of the elliptic analogues of the KZ associator (Section 5). 

\subsection{Automorphisms of group schemes}

We will view a $\QQ$-group scheme as a functor $\{\QQ$-algebras$\}\to\{$groups$\}$. 
The Lie algebra of a $\QQ$-group scheme $G(-)$ is then 
$\on{Lie}G(-):= \on{Ker}(G(\QQ[\epsilon]/(\epsilon^{2}))\to G(\QQ))$. 

If $\Gamma$ is a finitely generated group, let $\Gamma(-)$ be its 
$\QQ$-prounipotent completion and let $\on{Lie}\Gamma$ be its Lie algebra
(a pronilpotent $\QQ$-Lie algebra). Let $\ul{\on{Aut}\Gamma}(-)$ be the 
$\QQ$-group scheme defined by 
$\ul{\on{Aut}\Gamma}({\mathbf k}):= \on{Aut}(\on{Lie}\Gamma
\hat\otimes {\mathbf k})$ for ${\mathbf k}$ a $\QQ$-ring, where 
$\on{Lie}\Gamma\hat\otimes {\mathbf k}:= \on{lim}_{\leftarrow}
(\on{Lie}\Gamma/(\on{Lie}\Gamma)^{\geq n})\otimes {\mathbf k}$, and 
$\on{Lie}\Gamma = (\on{Lie}\Gamma)^{\geq 0}\supset  (\on{Lie}\Gamma
)^{\geq 1}\supset\cdots$
is the lower central series filtration of $\on{Lie}\Gamma$. 

Any automorphism of $\Gamma$ gives rise to an automorphism of $\on{Lie}\Gamma$, so 
there are natural morphisms 
$$
\on{Aut}\Gamma\to \ul{\on{Aut}\Gamma}(\QQ)\to \on{Aut}(\Gamma({\mathbf k}))
$$
for any $\QQ$-ring ${\mathbf k}$. One checks that there is a morphism of 
$\QQ$-group schemes 
$$
\mu_{O} : \on{GT}(-)\to \ul{\on{Aut}P_{n}}(-) 
$$
such that the resulting morphism $\on{GT}({\mathbf k})\to 
\on{Aut}(P_{n}({\mathbf k}))$ is compatible with 
$\on{GT}({\mathbf k})\to \on{Aut}(B_{n}({\mathbf k}))$, and 
morphisms 
$$
\mu_{O}^{gr} : \on{GRT}(-)\to \ul{\on{Aut}P^{gr}_{n}}(-), 
\quad \mu_{O}^{ell} : \on{GT}_{ell}(-)\to \ul{\on{Aut}P_{1,n}}(-), 
\quad 
\mu_{O}^{ell,gr} : \on{GRT}_{ell}(-)\to \ul{\on{Aut}P^{gr}_{1,n}}(-), 
$$
$$
\mu_{ell} : R_{ell}(-)\to\ul{\on{Aut}P_{1,n}}(-), \quad 
R_{ell}^{gr}(-)\to\ul{\on{Aut}P_{1,n}^{gr}}(-)
$$
with similar properties. 

\subsection{Results on Zariski closures}

Define the $\QQ$-group scheme $\langle B_{3}\rangle$ to be the Zariski closure of 
the composite group morphism $B_{3}\to \on{Aut}F_{2}\to \ul{\on{Aut}F_{2}}
(\QQ)$; this is a group subscheme of $\ul{\on{Aut}F_{2}}(-)$. 

\begin{theorem} \label{thm:Zar:base}
Any elliptic associator of the form $e(\tau)$, $\tau\in\HH$, or $e_{KZ}$, 
gives rise to an isomorphism of $\CC$-group schemes 
$\langle B_{3}\rangle\otimes \CC \simeq \big( \on{exp}(\hat\b_{3}^{
+})\rtimes\on{SL}_{2} \big) \otimes \CC$. 
Any two isomorphisms arising in this way are related by an inner automorphism. 
There exists an analogous isomorphism for $\QQ$-group schemes. 
\end{theorem}

For $n\geq 1$, define $\langle B_{3}\rangle_{n}$ to be the Zariski 
closure of the composite group morphism $B_{3}\to \on{Aut}P_{1,n}
\to \ul{\on{Aut}P_{1,n}}(\QQ)$; this is a group subscheme of 
$\ul{\on{Aut}P_{1,n}}(-)$. 

\begin{theorem} \label{thm:Zar:n}
For any $n\geq 1$, there is an isomorphism 
$\langle B_{3}\rangle\simeq\langle B_{3}\rangle_{n}$
of $\QQ$-group schemes, which is compatible with the maps from $B_{3}$ 
to both sides.   
\end{theorem}

\subsection{Proof of Theorem \ref{thm:Zar:base}}

Composing (\ref{diag}) with the morphism $\tilde B_{3}\to 
\on{GT}_{ell}({\mathbf k})$, we obtain a commutative diagram 
$\begin{matrix} \tilde B_{3}&
\to & \on{GT}_{ell}({\mathbf k}) & \to  & \on{GRT}_{ell}({\mathbf k})\\
\downarrow & & \downarrow & & \downarrow\\
\{\pm 1\}  & \to & \on{GT}({\mathbf k}) & \to 
& \on{GRT}({\mathbf k})\end{matrix}$ inducing morphisms 
$B_3 \to R_{ell}({\mathbf k}) \to R_{ell}^{gr}({\mathbf k})$. 

Set 
$$
e_{KZ}:= 
\sigma(\Phi_{KZ})_{ | \begin{smallmatrix}
x\mapsto 2\pi\on{i}x, \\
 y\mapsto (2\pi\on{i})^{-1}y
\end{smallmatrix} } . 
$$
When ${\mathbf k}=\CC$ and $e = e(\tau),e_{KZ}$, the morphism 
$B_3\to R_{ell}^{gr}(\CC)$ is computed as follows. 

Define $F(\tau)$ as the map $\HH\to \on{exp}(\hat\b_3^{+,\CC}) 
\rtimes\on{SL}_{2}(\CC)$ such that 
\begin{equation} \label{cnx}
2\pi\on{i}\partial_\tau F(\tau) =  (e_-+\sum_{k\geq 0}
(2k+1)G_{2k+2}(\tau)\delta_{2k})F(\tau)
\end{equation}
and $F(\tau) \sim \on{exp}( {{\tau}\over{2\pi\on{i}}} (e_-
+\sum_{k\geq 0} (2k+1)2\zeta(2k+2) \delta_{2k}))$ as 
$\tau\to\on{i}\infty$. Then the map $\tau\mapsto e(\tau) * F(\tau)$ 
is a constant, and 
\begin{equation} \label{vam}
e_{KZ} = e(\tau) * F(\tau)\; \on{for\ any}\; \tau\in\HH.
\end{equation}

Moreover, for any $\tilde g\in B_3$ with image 
$g = \bigl(\begin{smallmatrix}\alpha & \beta\\ \gamma & \delta
\end{smallmatrix}\bigl)\in\on{SL}_2(\ZZ)$, one has
\begin{align} \label{bro}
e(\tau)*i_{e(\tau)}(\tilde g) & = 
\tilde g*e(\tau)  = 
\on{Ad}(e^{-f(\tilde g,\tau)t})(\bigl(\begin{smallmatrix}\gamma\tau+\delta & 0
\\ 2\pi\on{i}\gamma & (\gamma\tau+\delta)^{-1}
\end{smallmatrix}\bigl)^{-1} \bullet e(g\tau)) 
\\ &  \nonumber = e(\tau)*F(\tau)F(g\tau)^{-1}\bigl(\begin{smallmatrix}
 \gamma\tau+\delta & 0 \\
  2\pi\on{i}\gamma & (\gamma\tau+\delta)^{-1}
\end{smallmatrix}\bigl)^{-1} e^{-f(\tilde g,\tau)\delta_{0}}, 
\end{align}
where the third equality follows from (\ref{vam}) for $\tau$
and $g\tau$. It follows that 
\begin{equation} \label{equality:a}
i_{e(\tau)}(\tilde g) = F(\tau)F(g\tau)^{-1}\bigl(\begin{smallmatrix}\gamma\tau+\delta & 0
\\ 2\pi\on{i} \gamma & (\gamma\tau+\delta)^{-1}
\end{smallmatrix}\bigl)^{-1} e^{-f(\tilde g,\tau)\delta_{0}}. 
\end{equation}
Acting from the right by $F(\tau)$ in the equality between the second and the 
fourth terms of (\ref{bro}),  one gets $\tilde g*e_{KZ} = 
e_{KZ}*F(g\tau)^{-1}\bigl(\begin{smallmatrix}\gamma\tau+\delta & 0
\\  2\pi\on{i}\gamma & (\gamma\tau+\delta)^{-1}
\end{smallmatrix}\bigl)^{-1} F(\tau) e^{-f(\tilde g,\tau)\delta_{0}}$, so 
\begin{equation} \label{equality:b} 
i_{e_{KZ}}(\tilde g) = F(g\tau)^{-1}\bigl(\begin{smallmatrix}
\gamma\tau+\delta & 0
\\ 2\pi\on{i}\gamma & (\gamma\tau+\delta)^{-1}
\end{smallmatrix}\bigl)^{-1} F(\tau) e^{-f(\tilde g,\tau)\delta_{0}}
\end{equation}
for any $\tau\in\HH$. It follows that the images of $i_{e(\tau)},i_{e_{KZ}}$
are contained in $\on{exp}(\hat\b_3^{+,\CC})\rtimes\on{SL}_2(\CC)\subset 
R_{ell}^{gr}(\CC)$. 
The composite morphism $B_3\to \on{exp}(\hat\b_3^{+,\CC})\rtimes\on{SL}_2(\CC)
\to\on{SL}_2(\CC)$ is $\tilde g\mapsto 
\bigl(\begin{smallmatrix}\alpha & -\beta/(2\pi\on{i}) 
\\ -2\pi\on{i}\gamma & \delta \end{smallmatrix}\bigl)$. 

Recall that $B_3\subset \on{Aut}(F_2)$ is generated by 
$\Psi_{+},\Psi_{-}$ given by $\Psi_{+}:X,Y\mapsto X,YX$ and 
$\Psi_{-} : X,Y\mapsto XY^{-1},Y$. Let $\Theta:= (\Psi_{+}\Psi_{-}
\Psi_{+})^{-1}$, then 
$\Psi_{-} = \Theta\Psi_{+}\Theta^{-1}$ and 
$\Theta:X,Y\mapsto XYX^{-1},X^{-1}$. 

Then 
$$
i_{e_{KZ}}(\Psi_{+}) =  F(\tau+1)^{-1}F(\tau) = \on{exp}(-(2\pi\on{i})^{-1}
(e_-+\sum_{k>0}2(2k+1)\zeta(2k+2) \delta_{2k})) 
e^{{2\pi\on{i}\over{12}}\delta_0} =:\psi_{+}, 
$$
and since $i_{e_{KZ}}(\Theta) \in \bigl(\begin{smallmatrix} 0 
& -(2\pi\on{i})^{-1} \\ 2\pi\on{i} & 0
\end{smallmatrix}\bigl) \on{exp}(\hat\b_{3}^{+,\CC})$, 
$$
i_{e_{KZ}}(\Psi_-) = \psi_-, \quad\text{where}\quad 
 \on{log}\psi_- = 2\pi\on{i} (e_+ + \text{ element of }\hat\b_3^{\CC,+}). 
$$

We then prove: 

\begin{proposition} \label{prop:1:3}
For $e=e_{KZ}$, the isomorphism 
$i_{e}:R_{ell}(\CC)\to R_{ell}^{gr}(\CC)$ restricts to an 
isomorphism $\langle B_3\rangle(\CC)
\to \on{exp}(\hat\b_3^{+,\CC})\rtimes\on{SL}_{2}(\CC)$. 
\end{proposition}

{\em Proof.} $i_{e}(\langle B_3\rangle(\CC))$ is the Zariski closure of 
the subgroup of $R_{ell}^{gr}(\CC)$ generated by $\psi_\pm$. These are 
elements of the subgroup 
$\on{exp}(\hat\b_3^{+,\CC})\rtimes\on{SL}_{2}(\CC)$, which is 
Zariski closed, so $i_{e}(\langle B_3\rangle(\CC))$ is contained in this 
group. On the other hand, the Lie algebra of this Zariski closure is the 
topological Lie algebra generated by $\on{log}\psi_\pm$. It then suffices 
to prove  that this Lie algebra coincides with $\hat\b_3^{\CC}$. 

Equip $\on{Aut}(F_{2}(\CC))$ with the topology for which a system of 
neighborhoods of $1$ is $\on{Aut}^{n}(F_{2}(\CC)) = \{\theta|\forall
g\in F_{2}(\CC),\theta(g)\equiv g$ mod $F_{2}^{(n)}(\CC)\} \subset 
\on{Aut}(F_{2}(\CC))$, where $F_{2}^{(1)}(\CC)=F_{2}(\CC)$ and 
$F_{2}^{(n)}(\CC) = (F_{2}^{(n-1)}(\CC),F_{2}(\CC))$. This induces a 
topology on $R_{ell}(\CC)$, which we call the prounipotent topology. 

\begin{lemma}
$\langle B_3\rangle(\CC)\subset R_{ell}(\CC)$ is closed for this topology. 
\end{lemma}

{\em Proof.} We have $\langle B_3\rangle(\CC) = \cap_{G\in{\mathcal G}}G$, 
where ${\mathcal G} = \{G(\CC)|G\subset R_{ell}(-)$ is a subgroup scheme 
such that $G(\QQ)\supset B_3\}$. It then suffices to show that each $G(\CC)$ 
is closed in the prounipotent topology.
Define coordinates on $R_{ell}(\CC)$ as follows: $R_{ell}(\CC)\ni\theta
\leftrightarrow (c_{b},d_{b})_{b}$, where $b$ runs over a homogeneous  
basis of $\f_{2}$ (generated by $\xi=\on{log}X$, $\eta=\on{log}Y$), e.g., 
$\{b\} = \{\xi,\eta,[\xi,\eta],\ldots\}$, and $\on{log}\theta(e^{\xi})
= \sum_{b}c_{b}b$, $\on{log}\theta(e^{\eta})=\sum_{b}d_{b}b$. Then 
$G(\CC)$ is a finite intersection of sets of the form $\{\theta|
P(c_{\xi},c_{\eta},\ldots,d_{\xi},d_{\eta},\ldots)=0\}$, where $P$
is a polynomial in $(c_{b},d_{b})_{b}$, vanishing at the origin. Such a 
$G(\CC)$ necessarily contains $R_{ell}(\CC)\cap \on{Aut}^{n}(F_{2}(\CC))$
for a large enough $n$. \hfill\qed\medskip 

{\em Sequel of proof of Proposition \ref{prop:1:3}.}
It follows that $i_{e}(\langle B_3\rangle(\CC))\subset \on{exp}(
\hat\b_3^{+,\CC})\rtimes\on{SL}_{2}(\CC)$ is closed in the prounipotent 
topology of $R_{ell}^{gr}(\CC)$ (as in the case of $R_{ell}(\CC)$, it is 
defined by the inclusion in $\on{Aut}(\hat\f_{2}^{\CC})$), so 
$\on{Lie} i_{e}(\langle B_3 \rangle(\CC))\subset \hat\b_3^{\CC}$ is closed.  

Recall that $\on{log}\psi = -2\pi\on{i}(e_- + \sum_{k\geq 1}
a_{2k}\delta_{2k}) + {{2\pi\on{i}}\over {12}}\delta_{0}
\in\on{Lie}i_{e}(\langle B_3\rangle(\CC))$, 
while $\log\psi_- \in (2\pi\on{i})^{-1}(e_++\hat\b_3^{+,\CC})$. 
 
\begin{lemma} \label{lemma:subalg}
Let ${\mathfrak{g}}\subset\hat\b_3^\CC$ be a closed (for the total degree 
topology) Lie subalgebra, such that ${\mathfrak{g}}
\ni \tilde e_\pm$, 
where: $\tilde e_+ = e_+ + $ terms of degree $>0$, $\tilde e_- = 
e_-+\sum_{k>0}a_{2k}\delta_{2k}-{1\over{12}}\delta_0+\sum_{p\geq 1,q>1}$
degree $(p,q)$. Then ${\mathfrak{g}}=\hat\b_3^\CC$.  
\end{lemma}

{\em Proof.} Set ${\mathcal G} = \oplus_{k\geq 0}{\mathcal G}_{2k}
:= \b_3^\CC$ (decomposition w.r.t.~the total degree), $\hat{\mathcal G}:= 
\hat\b_3^\CC$. Set $\hat{\mathcal G}_{\geq 2k}:= \prod_{k'\geq k}
{\mathcal G}_{2k'}$, then $\hat{\mathcal G} = \hat{\mathcal G}_{\geq 0}
\supset \hat{\mathcal G}_{\geq 2}\supset \ldots$ is a complete descending 
Lie algebra filtration of $\hat{\mathcal G}$, with associated graded 
Lie algebra $\mathcal G$.
Set ${\mathfrak g}_{\geq 2k}:= {\mathfrak g}\cap {\mathcal G}_{2k}$,
then ${\mathfrak g} = {\mathfrak g}_{\geq 0}\supset {\mathfrak g}_{\geq 2}
\supset \ldots$ is a complete descending filtration of ${\mathfrak g}$. Let $\on{gr}(
{\mathfrak g}):= \oplus_{k\geq 0}{\mathfrak g}_{\geq 2k}$, where $\on{gr}(
{\mathfrak g}):= {\mathfrak g}_{\geq 2k}/{\mathfrak g}_{\geq 2(k+2)}$. 
We then have an
inclusion $\on{gr}({\mathfrak g})\subset {\mathcal G}$ of graded Lie algebras.
We now prove that $\on{gr}({\mathfrak g})= {\mathcal G}$.  

As $\tilde e_\pm\in{\mathfrak g}={\mathfrak g}_{0}$, 
$\on{gr}({\mathfrak g}_{0})$ contains $e_\pm$. Set $h:=[e_+,e_-]$. 
Then $[\tilde e_+,\tilde e_-] = h+\sum_{p,q\geq 1}$ terms of degree $(p,q)$, 
and $[h,e_-]=-2e_-$, $[h,\delta_{2n}]=2n\delta_{2n}$. 
Then ${\mathfrak g}\ni P(\on{ad}[\tilde e_+,\tilde e_-])
(\tilde e_-) = P(-2)\Delta_{0}+\sum_{n\geq 0}a_{2n}P(2n)
\delta_{2n}+\sum_{p\geq 1,q>1}$ terms of degree $(p,q)$ (with 
$a_0=-{{1}\over{12}}$). Taking $P$
such that $P(-2)=P(0)=\ldots=P(2k-2)=0$ and $a_{2k}P(2k)=1$, we see that 
${\mathfrak g}$ contains an element of the form $\delta_{2k}
+\sum_{p\geq 1,q>1}$ terms of degree $(p,q)$. Applying 
$(\on{ad}\tilde e_-)^{2k}$ to this element, and using the fact that 
$(\on{ad}\tilde e_-)^{2k}(x)=0$ for $x\in {\mathcal G}$ of total 
degree $\leq 2(k-1)$, we see that ${\mathfrak g}$ contains an element of 
the form $(\on{ad}e_-)^{2k}(\delta_{2k})+\sum$ terms of total 
degree $\geq 2(k+2)$. As the latter sum belongs to 
${\mathfrak g}_{\geq 2(k+1)}$, we obtain that 
$(\on{ad} e_-)^{2k}(\delta_{2k})
\in \on{gr}({\mathfrak g})_{2(k+1)}$. The Lie subalgebra $\on{gr}(
{\mathfrak g})\subset {\mathcal G}$ then contains $e_\pm$
and $(\on{ad}e_-)^{2k}(\delta_{2k})$, $k\geq 0$. As 
$(\on{ad}e_-)^{2k+1}(\delta_{2k})=0$, 
$(\on{ad}e_+)^{2k}(\on{ad}e_-)^{2k} (\delta_{2k})$ is a nonzero 
multiple of $\delta_{2k}$. So $\on{gr}({\mathfrak g})={\mathcal G}$. 
It follows that ${\mathfrak g}={\mathcal G}$. \hfill \qed\medskip

{\em End of proof of Proposition \ref{prop:1:3}.} 
Applying Lemma \ref{lemma:subalg}
with $\tilde e_+ = 2\pi\on{i}\on{log}\psi_{-}$, 
$\tilde e_- = -(2\pi\on{i})^{-1}\on{log}\psi$, we get 
$i_{e}(\on{Lie}\langle B_3\rangle(\CC))
= \hat\b_3^\CC$, as wanted. \hfill \qed\medskip  

The last part of Theorem \ref{thm:Zar:base} is a consequence of the 
following statement, applied to a torsor of isimorphisms of 
Lie algebras. It was communicated to the author by P.~Etingof; it 
is inspired by the results of \cite{Dr:Gal}. 

\begin{proposition} \label{prop:torsors}
Let $U = \lim_{\leftarrow}  U_{i}$ be a prounipotent $\QQ$-group 
scheme (where $U_{0}=1$) and let $T:= \lim_{\leftarrow}  T_{i}$, where $ T_{i}$
are a compatible system of torsors under $ U_{i}$, defined over $\QQ$. 
If $T(\CC)\neq\emptyset$, then $T(\QQ)\neq \emptyset$.  
\end{proposition}

{\em Proof.} Let $\tilde U_{i}:= \on{im}(U\to U_{i})$, then 
$U=\lim_{\leftarrow}\tilde U_{i}$, 
where $\cdots \to \tilde U_{2}\to \tilde U_{1}\to \tilde U_{0} = 1$ 
is a sequence of epimorphisms of 
unipotent groups. We set $K_{i}:= \on{Ker}(U\to \tilde U_{i})$; then 
$K_{i}\triangleleft U$. 
If we set $\tilde T_{i}:= \on{im}(T\to  T_{i})$, then 
$\tilde T_{i}\simeq T/K_{i}$ is a torsor over $\tilde U_{i}$; $T$ is the 
inverse limit of $\cdots\to \tilde T_{2}\to \tilde T_{1}$, where the morphisms 
are onto.  

We may therefore assume w.l.o.g. that the morphisms $U_{i+1}\to U_{i}$, 
$T_{i+1}\to T_{i}$ are onto; if $K_{i}:= \on{Ker}(U\to U_{i})$, then 
$T_{i} = T/K_{i}$. 

We now show that the projective systems $\cdots\to T_{2}\to T_{1}$, 
$\cdots\to U_{2}\to U_{1}$ may be completed so that for any $i$, 
$\on{Ker}(U_{i+1}\to U_{i})\simeq {\mathbb G}_{a}$. Indeed, 
for $U_{i+1}\to U'\to U_{i}$ a sequence of epimorphisms, we set 
$K':= \on{Ker}(U\to U')$ and $T':= T/K'$. Then $T_{i+1}\to T'\to T_{i}$
is a sequence of epimorphisms, compatible with $U_{i+1}\to U'\to U_{i}$. 

Let $t\in T(\CC)$. We construct a sequence $(k_{i})_{i\geq 0}$, where 
$k_{i}\in K_{i}(\CC)$, such that $\on{im}(k_{i}\cdots k_{0}t\in 
T(\CC)\to T_{i}(\CC))\in T_{i}(\QQ)$. Then $k:= \lim_{i}(k_{i}\cdots k_{0})
\in U(\CC)$ is such that $kt\in T(\QQ)$. 

We first construct $k_{0}$. $T_{1}(\CC)$ is nonempty as it contains $t_{1}
:=\on{im}(t\in
T(\CC)\to T_{1}(\CC))$, hence by Hilbert's Nullstellensatz $T_{1}(\bar\QQ)$
is nonempty. Using then $H^{1}(G_{\QQ},\bar\QQ) = 0$ (\cite{S}), 
we obtain that 
$T_{1}(\QQ)$ is nonempty; let $t'_{1}\in T_{1}(\QQ)$. Let $u_{1}\in U_{1}(\CC)$
be such that $t'_{1}=u_{1}t_{1}$. Let $k_{0}\in U(\CC) = K_{0}(\CC)$ be a 
preimage of $u_{1}$, then $\on{im}(k_{0}t\in T(\CC)\to T_{1}(\CC))\in T_{1}(\QQ)$. 

Assume that $k_{0},\ldots,k_{i-1}$ have been constructed and let us construct 
$k_{i}$. Let $\tilde t:= k_{i-1}\cdots k_{0}t$, then $t_{i-1}:= \on{im}(\tilde t\in 
T(\CC)\to T_{i-1}(\CC)) \in T_{i-1}(\QQ)$. Then 
$T_{i}(\CC)\times_{T_{i-1}(\CC)}\{t_{i-1}\}$ is nonempty 
as it contains $\tau:= \on{im}(\tilde t\in T(\CC)\to T_{i}(\CC))$.
As $t_{i-1}\in T_{i-1}(\QQ)$, we define a functor 
$\{\QQ$-rings$\}\to \{$sets$\}$, ${\mathbf k}\mapsto 
X(\mathbf{k}):= T_{i}(\mathbf{k})\times_{T_{i-1}(\mathbf{k})}\{t_{i-1}\}$; it is a $\QQ$-scheme
and a torsor under $K_{i}/K_{i+1}={\mathbb G}_{a}$. We have seen that $X(\CC)\neq
\emptyset$, from which we derive as above that $X(\QQ)\neq\emptyset$. 
Let $\tau'\in X(\QQ)$ and let $k_{i}\in K_{i}(\CC)$ be such that $\bar k_{i}\tau =
\tau'$, where $\bar k_{i}:= \on{im}(k_{i}\in K_{i}(\CC)\to K_{i}/K_{i-1}(\CC))$; 
then $\on{im}(k_{i}\cdots k_{0}t\in T(\CC)\to T_{i}(\CC)) = 
\on{im}(k_{i}\tilde t\in T(\CC)\to T_{i}(\CC)) = \bar k_{i}\tau
=\tau'\in T_{i}(\QQ)$. \hfill \qed\medskip 

\subsection{Proof of Theorem \ref{thm:Zar:n}}

The morphism $B_{3}\to\ul{\on{Aut}P_{1,n}}(\QQ)$ factors as 
$B_{3}\to R_{ell}(\QQ)\stackrel{\mu_{ell}}{\to}\ul{\on{Aut}P_{1,n}}(\QQ)$. 

The elliptic associator $e_{KZ}$ transports the morphism 
$R_{ell}(-)\to \ul{\on{Aut}P_{1,n}}(-)$ to the morphism 
$R_{ell}^{gr}(-)\to  \ul{\on{Aut}P_{1,n}^{gr}}(-)$, whose Lie 
algebra morphism is 
$$
{\mathfrak{r}}_{ell}^{gr}\to \on{Der}(\t_{1,n}), \quad 
(\alpha_{+},\alpha_{-})\mapsto 
(x_{i}^{\pm}\mapsto \alpha_{\pm}^{i,1\ldots\check i\ldots n}). 
$$
The morphism $\t_{1,n}\to \t_{1,2}$, $x_{i}^{\pm}\mapsto x_{i}^{\pm}$
if $i=1,2$, $x_{i}^{\pm}\mapsto 0$ if $i\in \{3,\ldots,n\}$ can then be used to 
prove that this Lie algebra morphism is injective. It follows that that the group morphism 
$R_{ell}(-)\to \ul{\on{Aut}P_{1,n}}(-)$ is injective. 

One has 
$\langle B_{3}\rangle_{n} = \cap_{H | H \subset  \ul{\on{Aut}P_{1,n}}(-), 
H(\QQ)\supset \on{im}(B_{3})} H$, therefore 
$$
\langle B_{3}\rangle_{n}\cap R_{ell}(-) 
=  \cap_{H | H \subset  \ul{\on{Aut}P_{1,n}}(-), 
H(\QQ)\supset \on{im}(B_{3})} (H\cap R_{ell}(-)). 
$$
The map 
\begin{align*}
& \{H | H \text{ algebraic subgroup of }
\ul{\on{Aut}P_{1,n}}(-), \text{ s.t. }H(\QQ) \supset \on{im}(B_{3})\}
\\ & \to\{ G | G 
\text{ algebraic subgroup of }R_{ell}(-), 
\text{ s.t. }G(\QQ) \supset \on{im}(B_{3})\},
\end{align*} given by $H\mapsto G:= H\cap R_{ell}(-)$, 
is surjective (a preimage of $G$ is $G$ itself). Therefore
$$
\langle B_{3}\rangle_{n}\cap R_{ell}(-)
=  \cap_{G | G \subset R_{ell}(-), G(\QQ)\supset 
\on{im}(B_{3})} G
= \langle B_{3}\rangle. 
$$
The Zariski closure of $\on{im}(B_{3}\to \ul{\on{Aut}P_{1,n}}(\QQ))$
is contained in the Zariski closure of $\on{im}(R_{ell}(\QQ)\to 
\ul{\on{Aut}P_{1,n}}(\QQ))$, which is $R_{ell}(-)$ as the morphism 
$R_{ell}(-)\to \ul{\on{Aut}P_{1,n}}(-)$ is injective. 
So 
$$
\langle B_{3}\rangle_{n}\subset R_{ell}(-)\subset \ul{\on{Aut}P_{1,n}}(-).
$$
All this implies that 
$\langle B_{3}\rangle\to\langle B_{3}\rangle_{n}$ is an isomorphism. 
\hfill \qed\medskip

\section{Iterated integrals of Eisenstein series and MZVs} 

In this section, we define regularized iterated integrals of modular forms. This 
construction generalizes both that of interated integrals of cusp forms ([Ma]) and 
the definition of the Mellin transform of Eisenstein series ([Za]): it is based on a 
truncation procedure and the use of modular properties. We study the relations between 
these numbers arising from modular invariance. We show that the relations 
(26)-(27) from \cite{CEE}, obtained by the study of a monodromy morphism, 
can be recovered from formula (\ref{equality:b}) for the isomorphism 
$i_{e_{KZ}}$. The study of these relations leads to a family 
of algebraic relations between the iterated integrals of Eisenstein series and the MZVs. 

\subsection{Iterated Mellin transforms of modular forms}

Iterated Mellin transforms of cusp modular forms were studied in \cite{Ma}. 
On the other hand, Mellin transforms of non-cusp (e.g.,  Eisenstein)
modular forms were studied in \cite{Za}. In this section, we study iterated 
Mellin transforms of general (i.e., non-necessarily cusp) modular
forms. 

\begin{proposition} Let ${\mathcal E}:= \{f: \on{i}
\RR_{+}^{\times}\to\CC | f$ is smooth and 
$f(\on{i}t) = a + O(e^{-2\pi t})$ as $t\to\infty$ for some $a\in\CC\}$. Set 
$$
F_{t_{0}}^{f_{1},\ldots,f_{n}}(s_{1},\ldots,s_{n}):= 
\int_{t_{0}\leq t_{1}\leq\ldots\leq t_{n}\leq\infty}
f_{1}(\on{i}t_{1})t_{1}^{s_{1}-1}dt_{1}\cdots
f_{n}(\on{i}t_{n})t_{n}^{s_{n}-1}dt_{n}, $$
where $f_{1},\ldots,f_{n}\in{\mathcal E}$ and $t_{0}\in\RR_{+}^{\times}$. 
This function is analytic for $\Re(s_{i})\ll 0$ and admits a meromorphic 
prolongation to $\CC^{n}$, where the only singularities are simple poles at
the hyperplanes $s_{i}+\cdots+s_{j}=0$ ($1\leq i\leq j\leq n$). 
\end{proposition} 

{\em Proof.} Set ${\mathcal E}_{0}:= \{f\in{\mathcal E}| a=0\}$. Then 
${\mathcal E}={\mathcal E}_{0}\oplus\CC 1$. When $f_{1},\ldots,f_{n}\in
{\mathcal E}_{0}$, $F_{t_{0}}^{f_{1},\ldots,f_{n}}$ is analytic on $\CC^{n}$. 
Let now $f_{1},\ldots,f_{n}\in
{\mathcal E}$, and set $f_{i} = \bar f_{i}+a_{i}$, with $\bar f_{i}\in
{\mathcal E}_{0}$. Using 
\begin{eqnarray*}
&& \int_{t\leq t_{1}\leq\ldots\leq t_{n}\leq t'}
t_{1}^{s_{1}-1}dt_{1}\cdots t_{n}^{s_{n}-1}dt_{n}
\\
 && = \sum_{k=0}^{n}(-1)^{k}{{(t')^{s_{k+1}+\cdots+s_{n}}}\over{s_{k+1}
(s_{k+1}+s_{k+2})\cdots(s_{k+1}+\cdots+s_{n})}}{{t^{s_{1}+\cdots+s_{k}}}
\over{s_{k}(s_{k}+s_{k-1})\cdots(s_{k}+\cdots+s_{1})}}, 
\end{eqnarray*}
we get 
\begin{eqnarray*}
&& F_{t_{0}}^{f_{1},\ldots,f_{n}}(s_{1},\ldots,s_{n}) = 
\sum_{k=1}^{n}\sum_{1\leq i_{1}<\ldots<i_{k}\leq n}
(\prod_{j\in\{1,\ldots,n\}-\{i_{1},\ldots,i_{k}\}} a_{j})\\
 && \sum_{\begin{smallmatrix}j_{1}\in\{1,...,j_{1}-1\},\\...,\\
 j_{k}\in\{i_{k-1},...,i_{k}-1\}\end{smallmatrix}}
 {{(-1)^{|A_{1}|+|A_{2}|+...+|A_{k+1}|}}
 \over{\prod}_{i=1}^{k+1}\tilde s_{A_{i}}\prod_{i=1}^{k}
 \tilde{\tilde s}_{B_{i}}} t_{0}^{s_{A_{1}}} F_{t_{0}}^{\bar f_{i_{1}},...,\bar f_{i_{k}}}
 (s_{i_{1}}+s_{B_{1}}+s_{A_{1}},...,s_{i_{k}}+s_{B_{k}}+s_{A_{k}}),
\end{eqnarray*}
where $A_{l} := \{i_{l-1}+1,...,j_{l}\}$, $B_{l} := \{j_{l}+1,...,i_{l}-1\}$
for $l=1,...,k$, and $A_{k+1} := \{i_{k}+1,...,n\}$, $s_{A} := \sum_{\alpha\in 
A}s_{\alpha}$, $\tilde s_{A} := s_{b}(s_{b}+s_{b-1})...(s_{a}+...+s_{b})$,
$\tilde{\tilde s}_{A} := s_{a}(s_{a}+s_{a+1})...(s_{a}+...+s_{b})$, for
$A = \{a,a+1,...,b\}$. This implies the result in general. \hfill \qed\medskip 
 
Note that $F_{t_{0}}^{f_{1},...,f_{n}} = {{(-1)^{n}a_{1}...a_{n}} 
\over{(s_{1}+...+s_{n})...s_{n}}}t_{0}^{s_{1}+...+s_{n}}
+O(t_{0}^{\sigma}e^{-2\pi t_{0}})$ as $t_{0}\to\infty$, where
$\sigma$ depends on the $\Re(s_{i})$. 

Let now $\tilde{\mathcal E}:= \{f\in{\mathcal E}|\exists N\geq 0, 
f(\on{i}t) = O(t^{-N})$ as $t\to 0^{+}\}$. Set 
$$
G_{t_{0}}^{f_{1},...,f_{n}}(s_{1},...,s_{n}):=
\int_{0\leq t_{1}\leq...\leq t_{n}\leq t_{0}}
f_{1}(\on{i}t_{1})t_{1}^{s_{1}-1}dt_{1}
...f_{n}(\on{i}t_{n})t_{n}^{s_{n}-1}dt_{n}
$$
for $f_{1},...,f_{n}\in\tilde{\mathcal E}$. This function is analytic for 
$\Re(s_{i})\gg 0$. 

\begin{proposition} For $f_{1},...,f_{n}\in\tilde{\mathcal E}$, the function 
$$(s_{1},...,s_{n})\mapsto \sum_{k=0}^{n} G_{t_{0}}^{f_{1},...,f_{k}}(s_{1},...,
s_{k})F_{t_{0}}^{f_{k+1},...,f_{n}}(s_{k+1},...,s_{n})$$ is analytic for $\Re(s_{i})
\gg 0$ and independent of $t_{0}$. We denote it $L^{*}_{f_{1},...,f_{n}}
(s_{1},...,s_{n})$. 
\end{proposition}

{\em Proof.} The analyticity follows from the fact that 
$F_{t_{0}}^{f_{k+1},...,f_{n}}(s_{k+1},...,s_{n})$ may be viewed as an
analytic function for $\Re(s_{i})\gg 0$. The independence of $t_{0}$ follows from 
$$\partial_{t_{0}}G_{t_{0}}^{f_{1},...,f_{k}}(s_{1},...,s_{k}) = f_{k}(\on{i}t_{0})
t_{0}^{s_{k}-1}G_{t_{0}}^{f_{1},...,f_{k-1}}(s_{1},...,s_{k-1}),$$ 
$$\partial_{t_{0}}F_{t_{0}}^{f_{k},...,f_{n}}(s_{k},...,s_{n}) = -f_{k}(\on{i}t_{0})
t_{0}^{s_{k}-1}F_{t_{0}}^{f_{k+1},...,f_{n}}(s_{k+1},...,s_{n}),$$ 
where the former identity in valid in the domain $\Re(s_{i})\gg 0$
and the latter is analytically extended from the domain $\Re(s_{i})\ll 0$ to 
$\Re(s_{i})\gg 0$. \hfill \qed\medskip 

Recall that if $f(\tau)$ is a modular form of weight $k$, then $(t\mapsto f(\on{i}t))
\in \tilde{\mathcal E}$, $f(\tau+1)=f(\tau)$ and $f({{-1}\over\tau}) = 
\tau^{k}f(\tau)$. 

\begin{propdef} \label{propdef:13}
Let $f_{i}$ be modular forms of weight $k_{i}$
($i=1,...,n$), then the function $L^{*}_{f_{1},...,f_{n}}(s_{1},...,s_{n})$ extends to a 
meromorphic function on $\CC^{n}$, whose only possible singularities are 
simple poles at the hyperplanes $s_{i}+...+s_{j}=0$ and 
$s_{i}+...+s_{j} =k_{i}+...+k_{j}$ (where $1\leq i\leq j\leq n$). 
We call it the iterated Mellin transform of $f_1,\ldots,f_n$. 
\end{propdef}

{\em Proof.} By modularity, $$G_{t_{0}}^{f_{1},...,f_{l}}(s_{1},...,s_{l})
=(-1)^{(k_{1}+\cdots+k_{l})/2}F_{1/t_{0}}^{f_{l},...,f_{1}}(k_{l}-s_{l},...,
k_{1}-s_{1}).$$ 
Plugging this equality in the definition of $L^{*}_{f_{1},...,f_{n}}(s_{1},...,s_{n})$
and using the poles structures of the functions $F_{1/t_{0}}^{f_{l},...,f_{1}}$, 
$F_{t_{0}}^{f_{l+1},...,f_{n}}$, we obtain the result. \hfill \qed\medskip 

When $n=1$, we now relate $L^{*}_{f}(s)$ with the Mellin transform $L^{*}(f,s)$
defined in \cite{Za}. Let $f$ be a modular form with $f(\tau)\to a$ as $\tau\to
\on{i}\infty$. Then $L^{*}(f,s)$ is defined for $\Re(s)\gg 0$ by 
$L^{*}(f,s) = \int_{0}^{\infty} (f(\on{i}t)-a)t^{s-1}dt$. Then: 

\begin{proposition} $L^{*}(f,s) = L^{*}_{f}(s)$. 
\end{proposition} 

{\em Proof.} $L^{*}(f,s) = \int_{0}^{t_{0}} f(\on{i}t)t^{s-1}dt 
- a{{t_{0}^{s}}\over s} + \int_{t_{0}}^{\infty} (f(\on{i}t)-a)t^{s-1}dt$
for $\Re(s)\gg 0$. 

On the other hand, $G_{t_{0}}^{f}(s) = \int_{0}^{t_{0}} f(\on{i}t)t^{s-1}dt$
for $\Re(s)\gg 0$, while 
$F_{t_{0}}^{f}(s) = \int_{t_{0}}^{\infty} f(\on{i}t)t^{s-1}dt =
\int_{t_{0}}^{\infty} (f(\on{i}t)-a)t^{s-1}dt -a {{t_{0}^{s}}\over s}$
for $\Re(s)\ll 0$. The second expression of $F_{t_{0}}^{f}(s)$ is meromorphic 
on $\CC$ with as its only possible singularity, a simple pole at $s=0$; in 
particular, this expression coincides with $F_{t_0}^f(s)$ for $\Re(s)\gg 0$. 
Then for $\Re(s)\gg 0$, $L^{*}_{f}(s) = G_{t_{0}}^{f}(s) + F_{t_{0}}^{f}(s) 
= L^{*}(f,s)$. \hfill \qed\medskip 

For $s_{1},\ldots,s_{n}\in\ZZ$, one sets 
$$ 
L^{\sharp}_{f_{1},\ldots,f_{n}}(s_{1},\ldots,s_{n})
:= \on{i}^{s_{1}+\cdots+s_{n}} L^{*}_{f_{1},\ldots,f_{n}}(s_{1},\ldots,
s_{n}). 
$$
According to Proposition-Definition \ref{propdef:13}, the numbers
\begin{equation} \label{def:L:star}
L^{\sharp}_{k_{1},\ldots,k_{n}}(l_{1},\ldots,l_{n}):= 
L^{\sharp}_{G_{k_{1}},\ldots,G_{k_{n}}}(l_{1},\ldots,l_{n}),  
\end{equation}
for $k_{1},\ldots,k_{n}$ even integers $\geq 4$, $l_{i}\in \{1,\ldots,k_{i}-1\}$, 
are well-defined. One can prove that $L^{\sharp}_{k_{1},\ldots,k_{n}}
(b_{1},\ldots,b_{n}) \in \on{i}^{l_{1}+\cdots + l_{n}}{\mathbb R}$.


\subsection{Monodromy relations and the isomorphism $i_{e_{KZ}}$}

(\ref{equality:b}) defines a morphism 
$$
i_{e_{KZ}} : B_{3}\to \on{exp}(\hat\b_{3}^{+})\rtimes \on{SL}_{2}(\CC)
(\subset \on{Aut}(\hat\f_{2}^{\CC})^{op}), 
$$
such that 
\begin{equation} \label{rel:KZ}
\forall \tilde g_{3}\in B_{3}, \quad 
\tilde g * e_{KZ} = e_{KZ} * i_{e_{KZ}}(\tilde g) = 
(i_{e_{KZ}}(\tilde g)(A_{KZ}), i_{e_{KZ}}(\tilde g)(B_{KZ})), 
\end{equation}
where $e_{KZ} = (A_{KZ},B_{KZ})$. 

Specializing to $\tilde g = \Psi_{+}$, this gives
$$
i_{e_{KZ}}(\Psi_{+}) : A_{KZ}\mapsto A_{KZ}, \quad 
B_{KZ}\mapsto B_{KZ}A_{KZ}, 
$$
and for $\tilde g = \Theta$, this gives
$$
i_{e_{KZ}}(\Theta) : A_{KZ}\mapsto B_{KZ}^{-1}, \quad 
B_{KZ}\mapsto B_{KZ}A_{KZ}B_{KZ}^{-1}. 
$$
In \cite{CEE}, we introduced $\tilde A,\tilde B\in\on{exp}(\hat\t_{1,2})$
related to $A_{KZ},B_{KZ}$ by 
$$
\tilde A = (2\pi/\on{i})^{t}A_{KZ}  (2\pi/\on{i})^{-t}, \quad 
\tilde B = (2\pi/\on{i})^{t}B_{KZ}  (2\pi/\on{i})^{-t}, 
$$ 
and elements $[\Psi],[\Theta]\in \on{exp}(\hat\b_{3}^{+})\rtimes
\on{SL}_{2}(\CC)$, and studying a monodromy morphism, 
showed relations (numbered (26), (27) in \cite{CEE}) 
$$
[\Psi]e^{\on{i}{\pi\over 6}\on{ad}t} : A_{KZ}\mapsto A_{KZ}, \quad
B_{KZ}\mapsto B_{KZ}A_{KZ}, $$
$$
[\Theta]e^{\on{i}{\pi\over 2}\on{ad}t} : A_{KZ}\mapsto B_{KZ}^{-1}, 
\quad B_{KZ}\mapsto B_{KZ}A_{KZ}B_{KZ}^{-1}. $$
One checks that $[\Psi] e^{\on{i}{\pi\over 6}\on{Ad}t} = i_{e_{KZ}}(\Psi)$, 
$[\Theta]e^{\on{i}{\pi\over 2}\on{ad}t} = i_{e_{KZ}}(\Theta)$. So 
(\ref{equality:b}) allows to recover relations (26), (27) from \cite{CEE}. 

\subsection{Relations between iterated Mellin transforms and MZVs} \label{sect:IIMZV}

Another consequence of (\ref{rel:KZ}) is the behavior of the automorphism 
$i_{e_{KZ}}(\Psi_{-})$, namely 
\begin{equation} \label{behav:Psi-}
i_{e_{KZ}}(\Psi_{-}) : A_{KZ}\mapsto A_{KZ}B_{KZ}^{-1}, \quad 
B_{KZ}\mapsto B_{KZ}. 
\end{equation}
Notice that $\Psi_{-} = \Theta\Psi_{+}\Theta^{-1}$, and that 
$\on{log}i_{e_{KZ}}(\Psi_{+})$ is a well-defined derivation of 
$\hat\f_{2}^{\CC}$. Set 
$$
x_{KZ}:= \on{log}A_{KZ|x\mapsto (2\pi\on{i})^{-1}x, y\mapsto 2\pi\on{i}y}
\in \hat\f_{2}^{\CC}, \quad 
y_{KZ}:= \on{log} B_{KZ|x\mapsto (2\pi\on{i})^{-1}x, y\mapsto 2\pi\on{i}y}
\in\hat\f_{2}^{\CC}, 
$$
so $\sigma(\Phi_{KZ}) = (e^{x_{KZ}},e^{y_{KZ}})$. Then 
(\ref{behav:Psi-}) is equivalent to the statement that the derivation 
$$D:= \on{Ad}(\on{diag}((2\pi\on{i})^{-1},2\pi\on{i})
\circ i_{e_{KZ}}(\Theta))(\on{log}i_{e_{KZ}}(\Psi_{+}))
\in \on{Der}(\hat\f_{2}^{\CC})$$ 
acts as follows 
$$
D : 
x_{KZ}\mapsto 
- {{\on{ad} x_{KZ}}
\over{1 - e^{- \on{ad} x_{KZ}}}}(y_{KZ}), \quad 
y_{KZ}\mapsto 0, 
$$
where $t(2\pi\on{i}) = \on{diag}((2\pi\on{i})^{-1},2\pi\on{i})\in \on{SL}_{2}(\CC)$ is viewed
as an automorphism of $\hat\f_{2}^{\CC}$ (see Proposition \ref{prop:5:4}). 

There is a decomposition $\on{Der}(\hat\f_{2}^{\CC}) = \prod_{k,l\in\ZZ}
\on{Der}(\f_{2}^{\CC})[k,l]$, where the bracket indicates the bidegree in 
$x,y$. Let $D = \sum_{k,l}D[k,l]$ be the corresponding decomposition of $D$. 
One has $\on{Der}(\f_{2}^{\CC})[k,l] = \on{Der}(\f_{2}^{\QQ})[k,l]
\otimes\CC$. 

Set ${\mathcal Z}_{0}:= \QQ$ and for $l\geq 1$, set 
$$
{\mathcal Z}_{l}:= \on{Span}_{\QQ}\{\zeta(l_{1},\ldots,l_{s}) |s \geq 1,  
l_{1}\geq 1,\ldots, l_{s-1}\geq 1, l_{s}\geq 2, l_{1}+\cdots + l_{s} = l\}
\subset \CC, 
$$
where 
$$
\zeta(l_{1},\ldots,l_{s}) = \sum_{1\leq k_{1}\leq\ldots\leq k_{s}} 
k_{1}^{-l_{1}}\cdots k_{s}^{-l_{s}}. 
$$
For $V\subset \CC$ a $\QQ$-vector subspace and $k\in\ZZ$, 
set $V(k):= (2\pi\on{i})^{k}V$. 

\begin{proposition} \label{form:D}
$D[k,l]$ has the following properties: 

$\bullet$ it lies in $\b_{3}^{\QQ}[k,l] \otimes \QQ(l)$ if $k=1$, $l\geq -1$; 

$\bullet$ it lies in $\b_{3}^{\QQ}[k,l] \otimes ({\mathcal Z}_{l}(0) 
+ {\mathcal Z}_{l+1}(-1))$ if $k\geq 2$, $l\geq 1$; 

$\bullet$ it is equal to zero in all the other cases. 
\end{proposition}

{\em Proof.} $\on{log}i_{e_{KZ}}(\Psi_{+}) \in 
\sum_{k\geq -1}\QQ(k) \otimes \b_{3}^{\QQ}[k,1]$, and 
$\on{diag}((2\pi\on{i})^{-1},2\pi\on{i})\circ i_{e_{KZ}}(\Theta) = 
($element of $\on{exp}(\hat\b_{3}^{+}))\times 
(x\mapsto y, y\mapsto x)$, and the support of $\hat\b_{3}^{+}$ is contained 
in $\{1,2,\ldots\}^{2}$. All this implies that 
$$D\in \prod_{l\geq -1}\QQ(l)\otimes \b_{3}^{\QQ}[1,l] \oplus 
\prod_{k\geq 2,l\geq 0}\b_{3}^{\CC}[k,l].$$ 
Since $D$ lies in $\hat\b_{3}^{\CC}$, 
whose support is contained in $\{(1,-1),(0,0),(-1,1)\}\cup \{1,2,\ldots\}^{2}$, 
this statement can be improved by changing the second product into 
$\prod_{k\geq 2,l\geq 1}\b_{3}^{\CC}[k,l]$. This implies the first and the 
last statement of the proposition. 

Recall that 
$$
x_{KZ} = \on{Ad}\Big(\Phi_{KZ}(- {{\on{ad}x}\over{e^{\on{ad}x}-1}}(y),t)
\Big)\Big(2\pi\on{i}{{- \on{ad}x}\over{e^{\on{ad}x}-1}}(y)\Big), 
$$
$$
y_{KZ} = \on{i}\pi t * \on{log}\Phi_{KZ}(
- {{\on{ad}x}\over{e^{\on{ad}x}-1}}(y),t) * x * \on{log}\Phi_{KZ}(
{{\on{ad}x}\over{e^{\on{ad}x}-1}}(y)+t,t), $$
where $*$ is the CBH product $a*b:= \on{log}e^{a}e^{b}$. 

There exists a unique derivation $\tilde D$ of $\hat\f_{2}^{\CC}$, such that 
$$
\tilde D : x\mapsto 0, \quad 
y\mapsto -{1\over{2\pi\on{i}}} {{e^{\on{ad}x}-1}\over{\on{ad}x}}
\varphi\Big(\on{ad}\Big(-2\pi\on{i}{{\on{ad}x}\over{e^{\on{ad}x}-1}}
(y)\Big)\Big)(x), 
$$
where $\varphi(t) = (-t)/(1-e^{-t})$, and a unique automorphism $\theta$ of the 
same Lie algebra, such that 
$$
\theta : x \mapsto y_{KZ}, \quad y\mapsto - {1\over{2\pi\on{i}}}
{{e^{\on{ad}y_{KZ}}-1}\over{\on{ad}y_{KZ}}}(x_{KZ}) ;  
$$
then $D = \theta\tilde D\theta^{-1}$. One has 
\begin{equation} \label{supp:tilde:D}
\tilde D \in \QQ(-1)\otimes \on{Der}(\f_{2})[1,-1]
+ \prod_{k\geq 1,l\geq 0}\QQ(l)\otimes \on{Der}(\f_{2})[k,l]. 
\end{equation}
One computes 
$$
\on{log}\Phi_{KZ}( - {{\on{ad}x}\over{e^{\on{ad}x}-1}}(y),t), 
\on{log}\Phi_{KZ}({{\on{ad}x}\over{e^{\on{ad}x}-1}}(y)+t,t)
\in \prod_{k\geq 1,l\geq 2}{\mathcal Z}_{l}\otimes \f_{2}^{\QQ}[k,l], \quad 
\on{i}\pi t\in \QQ(1)\otimes \f_{2}[1,1],  
$$
which implies 
\begin{equation} \label{supp:y:KZ}
y_{KZ} \in x + \prod_{k\geq 1,l\geq 1}
({\mathcal Z}_{l} + {\mathcal Z}_{l-1}(1)) \otimes \f_{2}^{\QQ}[k,l]. 
\end{equation} 
It also implies 
$$
-{1\over{2\pi\on{i}}}x_{KZ} = y 
+ \prod_{k,l\geq 1} {\mathcal Z}_{l-1}\otimes \on{Der}\f_{2}^{\QQ}[k,l], 
$$
which then implies 
\begin{equation} \label{supp:x:KZ}
-{1\over{2\pi\on{i}}}
{{e^{\on{ad}y_{KZ}}-1}\over{\on{ad}y_{KZ}}}
(x_{KZ}) \in y + \prod_{k,l\geq 1}
 ({\mathcal Z}_{l} + {\mathcal Z}_{l-1}(1)) \otimes \f_{2}^{\QQ}[k,l]. 
\end{equation}
(\ref{supp:y:KZ}) and (\ref{supp:x:KZ}) imply that 
$\theta - \on{id}\in \prod_{k\geq 0,l\geq 1} 
({\mathcal Z}_{l} + {\mathcal Z}_{l-1}(1))\otimes \on{End}(\f_{2}^{\QQ})[k,l]$, 
so that $\on{log}\theta$ belongs to the same space. Together with the estimate on 
$\tilde D$, this implies that for any $k\geq 1$, 
$$
\on{ad}(\on{log}\theta)^{k}(\tilde D)\in 
\prod_{k\geq 1,l\geq 0} \on{Der}(\f_{2}^{\QQ})[k,l]
\otimes \big({\mathcal Z}_{l} + {\mathcal Z}_{l+1}(-1)\big). 
$$
Combining this with the estimate on $\tilde D$, one obtains
that $D = \theta\tilde D\theta^{-1}$ belongs to the direct sum of 
$\on{Der}(\f_{2}^{\QQ})[1,-1]\otimes \QQ(-1)$ with this space, 
which together with the first and third statements of the proposition, 
and the  fact that $D\in\hat\b_{3}^{\CC}$, implies the second statement of the 
proposition. 
\hfill \qed \medskip 

For $\lambda\in\CC^{\times}$, set $w(\lambda):= \bigl(\begin{smallmatrix}
0 & -\lambda^{-1} \\ \lambda & 0 \end{smallmatrix}\bigl)\in\on{SL}_{2}(\CC)$. 

\begin{lemma}
\begin{align*}
i_{e_{KZ}}(\Theta) \equiv w(2\pi\on{i}) \cdot 
\sum_{n\geq 0} \sum_{ {k_{1},\ldots,k_{n}\geq 1} \atop 
{ l_{i}\in\{0,\ldots, 2k_{i}\} } } & ({-1\over{2\pi\on{i}}})^{l_{1}+1}\cdots 
({-1\over{2\pi\on{i}}})^{l_{n}+1}L^{\sharp}_{2k_{1}+2,\ldots,2k_{n}+2}(l_{1}+1,\ldots,l_{n}+1) 
\\
 & {{2k_{1}+1}\over{l_{1}!}}\on{ad}(e_{-})^{l_{1}}(\delta_{2k_{1}})
 \cdots {{2k_{n}+1}\over{l_{n}!}}\on{ad}(e_{-})^{l_{n}}(\delta_{2k_{n}}) 
\end{align*}
in $\on{exp}(\hat\b_{3}^{+,\CC})\rtimes \on{SL}_{2}(\CC)$, 
up to multiplication by an element of $\on{exp}(\CC\delta_{0})$. 
\end{lemma}

{\em Proof.} $i_{e_{KZ}}(\Theta) = \tilde F({-1\over\tau})^{-1}w(2\pi\on{i})
\tilde F(\tau) e^{\on{log}({-1\over\tau})\delta_{0}}$, where $w = 
\bigl(\begin{smallmatrix} 0 & -1 \\
 1 & 0 \end{smallmatrix}\bigl)$ and $\tilde F(\tau):= 
 n_{-}({\tau\over{2\pi\on{i}}})^{-1}F(\tau)$.
 As $\tilde F(\tau)$ satisfies 
$$
\partial_{\tau}\tilde F(\tau) = - \Big(\sum_{k\geq 0}\sum_{l=0}^{2k}
\tau^{l}G_{2k+2}(\tau) ({-1\over{2\pi\on{i}}})^{l+1}
{{2k+1}\over{l!}}\on{ad}(e_{-})^{l}(\delta_{2k})\Big)\tilde F(\tau), $$ 
and taking into account the behavior of $\tilde F(\tau)$ at $\tau\to\on{i}\infty$, 
one obtains 
\begin{align*}
\tilde F(\tau) \equiv \sum_{n\geq 0}\sum_{k_{1},\ldots,k_{n}\geq 1}
\sum_{l_{i}\in\{0,\ldots,2k_{i}\}}
& \phi_{\tau}^{G_{2k_{1}+2},\ldots,G_{2k_{n}+2}}(l_{1}+1,\ldots,l_{n}+1)
({-1\over{2\pi\on{i}}})^{l_{1}+1}\cdots
({-1\over{2\pi\on{i}}})^{l_{n}+1}
\\
 & {{2k_{1}+1}\over{l_{1}!}}\on{ad}(e_{-})^{l_{1}}(\delta_{2k_{1}})\cdots 
{{2k_{n}+1}\over{l_{n}!}}\on{ad}(e_{-})^{l_{n}}(\delta_{2k_{n}}),  
\end{align*}
where $\phi_{\on{i}t_{0}}^{f_{1},\ldots,f_{n}}(s_{1},\ldots,s_{n}):= 
\on{i}^{s_{1}+\cdots+s_{n}}
F_{t_{0}}^{f_{1},\ldots,f_{n}}(s_{1},\ldots,s_{n})$. 
Combining this with the similar formula for $\tilde F({-1\over\tau})^{-1}$, 
one obtains the result. \hfill \qed\medskip 

Set $w:= w(1)\in \on{SL}_{2}(\CC)$. 
\begin{lemma} \label{form:D:bis}
\begin{align*}
& \on{Ad}(w^{-1}) \circ D = 
- {1\over{2\pi\on{i}}}e_{-} + {{2\pi\on{i}}\over{12}}\delta_{0}
+ \sum_{{n>0}\atop{{k_{1},\ldots,k_{n}\geq 1}
\atop {l_{i}\in\{0,\ldots,2k_{i}\}}}}
({{-1}\over{2\pi\on{i}}})^{l_{1}+1} \cdots 
({{-1}\over{2\pi\on{i}}})^{l_{n}+1} \times
\\
 & \times \Big\{\begin{matrix} 
- L_{2k_{1}+2,\ldots,2k_{n}+2}^{\sharp}(l_{1}+1,\ldots,l_{n}+1)\cdot 
l_{n} \quad\text{ if }l_{n}\neq 0 \\ 
L_{2k_{1}+2,\ldots,2k_{n-1}+2}^{\sharp}(l_{1}+1,\ldots,l_{n-1}+1)\cdot
2\zeta(2k_{n}+2) \text{ if }l_{n}=0 \end{matrix}\Big\}\times \\
 & \times 
 [{{2k_{1}+1}\over{l_{1}!}}
 (\on{ad}e_{-})^{l_{1}}(\delta_{2k_{1}}), \cdots, 
 {{2k_{n}+1}\over{l_{n}!}}
 (\on{ad}e_{-})^{l_{n}}(\delta_{2k_{n}})],  
\end{align*}
where $L^{\sharp}_{\emptyset}(\emptyset)=1$ by convention, and 
$[a_{1},\ldots,a_{n}]:= \on{ad}a_{1}\circ \cdots \circ \on{ad}a_{n-1}(a_{n})$. 
\end{lemma}

{\em Proof.} One has $t(2\pi\on{i}) \circ w(2\pi\on{i}) = w(1)=w$
in $\on{Aut}(\hat\f_{2})$. One also has 
$$
\on{log}i_{e_{KZ}}(\Psi_{+}) 
= {2\pi\on{i}\over{12}}\delta_{0} - {1\over{2\pi\on{i}}}(e_{-} + 
\sum_{k>0} 2(2k+1)\zeta(2k+2)\delta_{2k}).
$$ 
The result then follows from 
the expansion of $i_{e_{KZ}}(\Theta)$ and from the identity 
$\on{Ad}(g)(y) = \sum_{n\geq 0} a_{i_{1},\ldots,i_{n}}
[x_{i_{1}},\ldots,x_{i_{n}},y]$ for 
$g = \sum_{n\geq 0}\sum_{i_{1},\ldots,i_{n}\in I}
a_{i_{1},\ldots,i_{n}}x_{i_{1}}\cdots x_{i_{n}}$ a group-like element of 
$U{\mathfrak{g}}$ and $y\in {\mathfrak{g}}$, where ${\mathfrak{g}}$ 
is a topological Lie algebra and $x_{i}$, $i\in I$ are positive 
degree elements of ${\mathfrak{g}}$. \hfill \qed\medskip 

Combining Proposition \ref{form:D} and Lemma \ref{form:D:bis}, one 
obtains the following family of relations between iterated integrals of Eisenstein 
series and MZVs: 

\begin{proposition} \label{prop:II:MZV}
Let $I:= \{(a,b) | a,b\geq 1, a+b\text{ is even}\}$. For $(a,b)\in I$, 
let 
$$
e_{a,b}:= {{a+b-1}\over{(b-1)!}}(\on{ad}e_{-})^{b-1}(\delta_{a+b-2})
\in \b_{3}^{\QQ}[a,b].
$$ 
Let $A\geq 2,B\geq 1$. Any $\xi\in\b_{3}^{\QQ}[A,B]^{*}$ gives rise to a 
relation
\begin{align*}
& \sum_{n>0}\sum_{{(a_{1},b_{1}),\ldots,(a_{n},b_{n}) | }\atop{
(a_{1},b_{1}) + \cdots + (a_{n},b_{n}) = (A,B)}}
\langle \xi, [e_{a_{1},b_{1}},\ldots,e_{a_{n},b_{n}}] \rangle \times 
\\
 & \times 
\Big\lbrace \begin{matrix} 
- L^{\sharp}_{a_{1}+b_{1},\ldots,a_{n}+b_{n}}(b_{1},\ldots,b_{n})\cdot 
(b_{n}-1) \text{ if }b_{n}\neq 1
\\
L^{\sharp}_{a_{1}+b_{1},\ldots,a_{n-1}+b_{n-1}}(b_{1},\ldots,b_{n-1})
\cdot 2\zeta(a_{n}+b_{n}) \text{ if }b_{n}=1 
\end{matrix}\Big\rbrace \in {\mathcal Z}_{A}(B) + {\mathcal Z}_{A+1}(B-1). 
\end{align*} 
\end{proposition}

\subsection{Modular and shuffle relations} \label{sect:modular}

The numbers $L^{\sharp}_{k_{1},\ldots,k_{n}}(b_{1},\ldots,b_{n})$ are 
subject to other relations:

(a) the shuffle relations
\begin{equation} \label{rels:shuffle}
L^{\sharp}_{k_{1},\ldots,k_{n}}(b_{1},\ldots,b_{n})
L^{\sharp}_{k_{n+1},\ldots,k_{n+m}}(b_{n+1},\ldots,b_{n+m})
= \sum_{\sigma\in S_{n,m}}
L^{\sharp}_{k_{\sigma(1)},\ldots,k_{\sigma(n+m)}}(b_{\sigma(1)},\ldots,
b_{\sigma(n+m)}), 
\end{equation}
where $S_{n,m} = \{\sigma\in S_{n+m} | \sigma(i)<\sigma(j)$ if $i<j\leq n$ or 
$n+1\leq i<j\}$, which can be reexpressed as the following statement: 
let ${\mathfrak M}:= 
\oplus_{k\geq 4}\CC G_{k}\otimes \CC[t]_{\leq k-2}$, then the linear map 
$I : T({\mathfrak M})\to \CC$ such that $I(G_{k_{1}}(t)t^{b_{1}},\ldots,
G_{k_{n}}(t)t^{b_{n}}):= L^{\sharp}_{k_{1},\ldots,k_{n}}(b_{1}+1,\ldots,
b_{n}+1)$ is an algebra morphism, $T({\mathfrak M})$ being equipped with the 
shuffle algebra product ; 

(b) the modular relations 
\begin{equation} \label{modular}
I^{\otimes 2}\circ (\on{id}\otimes S)\circ \Delta = 
J^{\otimes 3}\circ (\on{id}\otimes U \otimes U^{2})\circ \Delta^{(2)} 
= \varepsilon 
\end{equation}
(equalities in $\on{Hom}_{alg}(T({\mathfrak M}),{\mathbb C})$, 
$T({\mathfrak M})$ being equipped with the shuffle product), where: 

$\bullet$ $J : T({\mathfrak M})\to\CC$ is defined by $J:= (I\otimes\psi)\circ 
\Delta$, $\psi : T({\mathfrak M})\to \CC$ being defined by 
$$
\psi(G_{k_{1}}t^{b_{1}-1}\otimes \cdots \otimes 
G_{k_{n}}t^{b_{n}-1}):= 
{{2\zeta(k_{1})\cdots 2\zeta(k_{n})}
\over{b_{1}(b_{1}+b_{2})\cdots (b_{1}+\cdots+b_{n})}} ;  
$$ 
$\bullet$ $S = \bigl(\begin{smallmatrix} 0 & -1 \\
 1 & 0 \end{smallmatrix}\bigl)$, $U = \bigl(\begin{smallmatrix} 1 & -1 \\
 1 & 0 \end{smallmatrix}\bigl)\in \on{PSL}_{2}(\ZZ)$ act on ${\mathfrak M}$ by 
 $S\cdot t^{b-1}G_{k}:= t^{k-2}({{-1}\over t})^{b-1}G_{k}$, 
 $U\cdot t^{b-1}G_{k}:= t^{k-2}(1-{1\over t})^{b-1}G_{k}$; 
 
 $\bullet$ $\varepsilon : 
 T({\mathfrak M})\to\CC$ is the augmentation morphism,  $\Delta : T({\mathfrak M})
\to T({\mathfrak M})^{\otimes 2}$ is the shuffle coproduct morphism 
$x_{1}\otimes \cdots \otimes x_{n}\mapsto \sum_{k=0}^{n}
(x_{1}\otimes \cdots \otimes x_{k})\otimes (x_{k+1}\otimes \cdots
\otimes x_{n})$, $\Delta^{(2)} := (\Delta\otimes\on{id}) \circ \Delta$. 

The relations (\ref{modular}) are proved as follows. Let 
$\tilde{\mathfrak b}_{3}:= \on{Lie}(\oplus_{a,b\geq 1, a+b\text{\ even}}
{\mathbb C}\tilde e_{a,b})\rtimes {\mathfrak{sl}}_{2}$, where $\on{Lie}(-)$
means the free Lie algebra generated by a vector space, and $\oplus_{a,b\geq 1, a+b\text{\ even}}
{\mathbb C}\tilde e_{a,b} = \oplus_{k\text{\ even}}(\oplus_{a,b\geq 1, 
a+b = k}\CC\tilde e_{a,b})$ is the direct sum of all the odd-dimensional 
simple ${\mathfrak{sl}}_{2}$-modules, the action being normalized by 
$e_{-}\cdot \tilde e_{a,b} = b \tilde e_{a-1,b+1}$. There is a unique morphism 
$\tilde{\mathfrak b}_{3}\to {\mathfrak b}_{3}$, such that it induces the 
identity on $\mathfrak{sl}_{2}$ and such that $\tilde e_{a,b}\mapsto 
e_{a,b}$. The $\on{SL}_{2}(\ZZ)$-equivariant connection on $\HH$ with values in the 
trivial principal bundle with group $\on{exp}(\hat\b_{3}^{+,\CC})\rtimes
\on{SL}_{2}(\CC)$ defined by (\ref{cnx}) admits a lift to a similar connection, 
where this group is replaced by its analogue with $\tilde{\mathfrak b}_{3}$
replacing ${\mathfrak b}_{3}$. This connection therefore gives rise to a morphism 
$\on{SL}_{2}(\ZZ)\to \on{exp}(\hat{\tilde{\mathfrak{b}}}_{3}^{+,\CC})
\rtimes \on{SL}_{2}(\CC)$ to this group. The relations (\ref{modular}) express the 
fact that the relations between the usual generators of $\on{SL}_{2}(\ZZ)$ are 
satisfied by their images. 

\begin{remark} Relations (\ref{modular}) are generalizations of the modular relations 
satisfied by the period polynomials of Eisenstein series (\cite{Za}, Proposition p. 453).
The contribution of $\psi$ to $J$ is the analogue of the contributions of the values
at cusps to the period polynomials of Eisenstein series as defined in \cite{Za}, (9). 
\end{remark}

\begin{remark}
Let ${\mathfrak Z} = \oplus_{k\geq 0}{\mathfrak Z}_{k}$ be the $\QQ$-algebra
of formal MZVs, i.e., the $\QQ$-algebra generated by formal versions of $2\pi\on{i}$
and of the $\zeta(k_{1},\ldots,k_{s})$, subject to the associator relations. 
Define ${\mathfrak Z}^{*}$ as the ${\mathfrak Z}$-algebra generated by 
formal analogues of the $L^{\sharp}_{k_{1},\ldots,k_{n}}(b_{1},\ldots,b_{n})$, 
$k_{1},\ldots,k_{n}$ even $\geq 4$, $b_{i}\in \{1,\ldots,k_{i}-1\}$, 
modulo the shuffle relations (\ref{rels:shuffle}), the modular relations (\ref{modular}), 
and the relations from Proposition \ref{prop:II:MZV}, in which the right hand side is 
replaced by any lift in ${\mathfrak Z}_{A}(B) + {\mathfrak Z}_{A+1}(B-1)$. 
Then ${\mathfrak Z}^{*}$ is $\NN$-graded, with the degree of 
$L^{\sharp}_{k_{1},\ldots,k_{n}}(b_{1},\ldots,b_{n})$ being equal to
$k_{1}+\cdots + k_{n}$. 
\end{remark}

\subsection{Computation of some regularized iterated integrals}

Denote by $\on{Sh}({\mathfrak M})$ the vector space $T({\mathfrak M})$, 
equipped with its (commutative)  shuffle  algebra structure. Let 
$\on{Lie}({\mathfrak M}) \subset T({\mathfrak M})$ be the (free) Lie 
subalgebra of $T({\mathfrak M})$ generated by ${\mathfrak M}$, 
$T({\mathfrak M})$ being equipped with its tensor algebra structure. 
This inclusion gives rise to a commutative algebra morphism 
$S(\on{Lie}({\mathfrak M}))\to \on{Sh}({\mathfrak M})$, 
which can be shown to be an isomorphism. As $I : \on{Sh}({\mathfrak M})
\to \CC$ is an algebra morphism, it is uniquely determined by its restriction
$$
I : \on{Lie}({\mathfrak M})\to \CC. 
$$
$\on{Lie}({\mathfrak M})$ decomposes as ${\mathfrak M}\oplus 
\on{Lie}_{2}({\mathfrak M})\oplus\cdots$. The restriction of $I$
to ${\mathfrak M}$ has been determined in \cite{Za}: for $k$ even $\geq 4$, 
\begin{equation} \label{I:on:M}
I(t^{k-2}G_{k}) = -I(G_{k}) = {{2\pi\on{i}}\over{k-1}}\zeta(k-1), 
\quad 
I(t^{a}G_{k}) = {{(-1)^{a+1}}\over{(k-1)!}}{{B_{a+1}}\over{a+1}}
{{B_{k-a-1}}\over{k-a-1}}(2\pi\on{i})^{k} \text{\ for\ }
a = 1,\ldots,k-3. 
\end{equation}
The grading ${\mathfrak M} = \oplus_{k\geq 4}{\mathfrak M}_{k}$, 
where ${\mathfrak M}_{k} = \CC G_{k}\otimes \on{Span}_{\CC}(1,t,\ldots,
t^{k-2})$, induces a grading $\on{Lie}({\mathfrak M}) = \oplus_{k\geq 4, 
k\text{\ even}}\on{Lie}({\mathfrak M})_{k}$. The restriction of $I$ to 
$\on{Lie}({\mathfrak M})_{k}$ for the first values of $k$ can be carried out 
as follows. 

$\bullet$ $k = 4,6$. In these cases, $\on{Lie}({\mathfrak M})_{k} 
= {\mathfrak M}_{k}$, so (\ref{I:on:M}) determines the restriction of $I$
to $\on{Lie}({\mathfrak M})_{k}$. 

$\bullet$ $k = 8$. $\on{Lie}({\mathfrak M})_{8} 
 = {\mathfrak M}_{8} \oplus \on{Lie}_{2}({\mathfrak M}_{4})$. 
 (\ref{I:on:M}) determines the restriction of $I$ to ${\mathfrak M}_{8}$, so it 
 remains to compute its restriction to $\on{Lie}_{2}({\mathfrak M}_{4})
 = \on{Span}_{\CC}([G_{4},tG_{4}],[G_{4},t^{2}G_{4}],[tG_{4},t^{2}G_{4}])$. 
 The modular relations imply that 
$$
I([G_{4},t^{2}G_{4}]) = - \big({{2\pi\on{i}}\over 3}\zeta(3)\big)^{2}
- {{418}\over{45}}\zeta(4)^{2},  
$$
and that $I([G_{4},tG_{4}]) + I([tG_{4},t^{2}G_{4}]) = 0$. 
Proposition \ref{prop:II:MZV} for $(A,B) = (4,4)$, together with the 
fact that the restriction of the morphism $\tilde\b_{3}\to \b_{3}$ to degree 8 
is an isomorphism, then implies that 
$I([G_{4},tG_{4}])\in {\mathcal Z}_{5}(3) + {\mathcal Z}_{4}(4)
= \on{Span}_{\QQ}((2\pi\on{i})^{3}\zeta(5),(2\pi\on{i})^{5}\zeta(3),
(2\pi\on{i})^{8})$. As $I([G_{4},tG_{4}])$ is pure imaginary, one even has 
$$
I([G_{4},tG_{4}]) = - I([tG_{4},t^{2}G_{4}])\in 
\QQ (2\pi\on{i})^{3}\zeta(5) + \QQ (2\pi\on{i})^{5}\zeta(3)
$$
(the rational coefficients can be determined from the expression of the 
components of the derivation $D$ in a generating family of MZVs). 

$\bullet$ $k=10$. $\on{Lie}({\mathfrak M})_{10} = {\mathfrak M}_{10}
\oplus {\mathfrak M}_{4}\otimes {\mathfrak M}_{6}$, and as a
$\on{SL}_{2}(\CC)$-module, ${\mathfrak M}_{4}\otimes{\mathfrak M}_{6}$
decomposes as a direct sum $V_{7}\oplus V_{5}\oplus V_{3}$ of irreducible modules
of the indicated dimensions, generated by the highest weight vectors 
\begin{equation} \label{hw:vect}
[G_{4},G_{6}], \quad [tG_{4},G_{6}] - [G_{4},tG_{6}], \quad 
[t^{2}G_{4},G_{6}]-2[tG_{4},tG_{6}]+[G_{4},t^{2}G_{6}]. 
\end{equation}
The modular relations determine the restriction of $I$ to $1$-codimensional subspaces
of $V_{i}$ ($i=3,5,7$), for which the highest weight vectors 
(\ref{hw:vect}) span supplementary subspaces. 

On the other hand, the expansion of $\on{log}(w(2\pi\on{i})^{-1}i_{e_{KZ}}(\Theta))$
up to degree 10 yields the identity 
\begin{align*}
& \on{Ad}(w^{-1})(D)\\
 &  = \on{Ad} \on{exp}\Big(
\sum_{a,b} ({{-1}\over{2\pi\on{i}}})^{b}
I(t^{b-1}G_{a+b})e_{a,b} + {1\over 4}\sum_{a,b,a',b'}
({{-1}\over{2\pi\on{i}}})^{b+b'} 
I([t^{b-1}G_{a+b},t^{b'-1}G_{a'+b'}])[e_{a,b},e_{a',b'}]\Big) \cdot \\
 & \cdot \Big({{-1}\over{2\pi\on{i}}}(e_{-} + \sum_{k\geq 1}
 2\zeta(2k+2)\delta_{2k})\Big) 
\end{align*}
modulo degree $\geq 12$, from where one derives the expression in terms of MZVs of 
\begin{equation} \label{known}
[e_{-},\sum_{a,b,a',b' | a+b+a'+b' = 10}
({{-1}\over{2\pi\on{i}}})^{b+b'}
I([t^{b-1}G_{a+b},t^{b'-1}G_{a'+b'}])[e_{a,b},e_{a',b'}]]. 
\end{equation}
On the other hand, let $\tilde V_{k}:= \on{Span}_{\CC}(\tilde e_{2k,1},\ldots,
\tilde e_{1,2k}) \subset \tilde{\mathfrak b}_{3}$; the degree 10 part of 
$\tilde\b_{3}$
decomposes as $(\tilde\b_{3})_{10} = \tilde V_{10} \oplus \tilde V_{4}
\otimes \tilde V_{6}$. This $\on{SL}_{2}(\CC)$-module is dual to 
$\on{Lie}({\mathfrak M})_{10}$, in particular 
\begin{equation} \label{can:elt}
\sum_{a,b,a',b' | a+b+a'+b'=10}
({{-1}\over{2\pi\on{i}}})^{b+b'}[t^{b-1}G_{a+b},t^{b'-1}G_{a'+b'}]
\otimes [\tilde e_{ab},\tilde e_{a'b'}]
\end{equation}
is the canonical element of $({\mathfrak M}_{4}\otimes{\mathfrak M}_{6})
\otimes (\tilde V_{4}\otimes \tilde V_{6})$. 
Decompose $\tilde V_{4}\otimes\tilde V_{6}$ as a direct sum 
$\tilde W_{7}\oplus \tilde W_{5}\oplus\tilde W_{3}$ of irreducible
$\on{SL}_{2}(\CC)$-modules of the indicated dimensions, then (\ref{can:elt})
is the sum of the canonical elements in each summand of 
$(V_{7}\otimes\tilde W_{7}) \oplus 
(V_{5}\otimes\tilde W_{5}) \oplus (V_{3}\otimes\tilde W_{3})$; these
canonical elements have the form 
$$
[G_{4},G_{6}]\otimes (\text{lowest weight vector of }\tilde W_{7}) + 
\text{ a sum of tensors of different weights}, 
$$
$$
([tG_{4},G_{6}]-[G_{4},tG_{6}])\otimes(\text{lowest weight vector of }
\tilde W_{5}) + \text{ a sum of tensors of different weights}, 
$$ 
$$
([t^{2}G_{4},G_{6}]-2[tG_{4},tG_{6}]+[G_{4},t^{2}G_{6}])\otimes 
(\text{l.w.v. of }\tilde W_{3}) + \text{tensors of different weights}.
$$
\begin{lemma} 
The composite maps $\tilde W_{7}\subset (\tilde\b_{3})_{10}
\to (\b_{3})_{10}$, $\tilde W_{5}\subset (\tilde\b_{3})_{10}\to 
(\b_{3})_{10}$ and $\tilde W_{3}\to(\tilde\b_{3})_{10}\to(\b_{3})_{10}$
are injective. 
\end{lemma}

{\em Proof.} The images of the highest weight vectors of $\tilde W_{7}$, 
$\tilde W_{5}$ in $\b_{3}\subset \on{Der}_{t}(\f_{2})$ can 
be partially computed (here $t = -[x,y]$ and $\on{Der}_{t}$ means the derivations
taking $t$ to zero) as follows. The commutator of derivations induces a map 
$\on{Der}_{t}(\f_{2},\f_{2}')^{\otimes 2} \to \on{Der}_{t}(\f_{2},\f_{2}'')$
(where $\f_{2}'= [\f_{2},\f_{2}]$, $\f_{2}'' = [\f_{2}',\f_{2}']$), which 
in its turn induces a map $D_{1}^{\otimes 2} \to D_{2}$, where 
$D_{1} := \on{Der}_{t}(\f_{2},\f_{2}')/ \on{Der}_{t}(\f_{2},\f_{2}'')$, 
$D_{2}:= \on{Der}_{t}(\f_{2},\f_{2}'')½/\on{Der}_{t}(\f_{2},\f_{2}''')$
(where $\f_{2}''':= [\f_{2}',\f_{2}'']$). 
There is a natural map $D_{1}\to \f_{2}'/\f_{2}''$ 
induced by $\on{Der}_{t}(\f_{2},\f_{2}')\to \f_{2}'/\f_{2}''$, $D\mapsto (\text{the 
class of an element }a\in\f_{2}'\text{ such that }D - \on{ad}a\in 
\on{Der}(\f_{2},\f_{2}''))$ and $D_{2}\to \f_{2}''/\f_{2}'''$ defined similarly. 
There are isomorphisms $\CC[u,v]\simeq \f_{2}'/\f_{2}''$, defined by $u^{n}v^{m}
\mapsto (\text{the class of }(\on{ad}x)^{n}(\on{ad}y)^{m}[x,y])$, and 
$\bigwedge^{2}\CC[u,v]\simeq \f_{2}''/\f_{2}'''$, induced by the Lie bracket
$\bigwedge^{2}\f_{2}'\to \f_{2}''$. The map $D_{1}^{\otimes 2}\to D_{2}$
is then compatible with an explicit map $\CC[u,v]^{\otimes 2}\to \bigwedge^{2}
\CC[u,v]$. The images in $\b_{3}$ of the highest weight vectors of $\tilde W_{7},
\tilde W_{5}$ in fact lie in $\on{Der}_{t}(\f_{2},\f_{2}'')$, and their images 
in $\bigwedge^{2}\CC[u,v]$ can be computed using the above map 
$\CC[u,v]^{\otimes 2}\to \bigwedge^{2}\CC[u,v]$ and shown to be nonzero. 
On the other hand, the image in $ \bigwedge^{2}\CC[u,v]$ of the highest 
weight vector of $\tilde W_{3}$ is zero, so the image of this highest weight 
vector in $\b_{3}$ lies in $\on{Der}_{t}(\f_{2},\f_{2}''')$. This derivation 
can be computed explicitly (by computer) and shown to be nonzero (this can also be
derived from \cite{Po}, Thm. 3, where $\on{Ker}(D_{1}^{\otimes 2}\to D_{2})$
is computed). Note that 
$[\on{Der}_{t}(\f_{2},\f_{2}'')_{10} : \underline{3}] = 1$,
where $\underline{3}$ is the irreducible 3-dimension representation of 
$\on{SL}_{2}(\CC)$, so this multiplicity space is spanned by the image of the 
highest weight vector of $\tilde W_{3}$. \hfill \qed\medskip 

The expression of (\ref{known}) in terms of MZVs therefore allows one to express
$I([G_{4},G_{6}])$, $I([tG_{4},G_{6}]-[G_{4},tG_{6}])$ 
and $I([t^{2}G_{4},G_{6}]-2[tG_{4},tG_{6}]+[G_{4},t^{2}G_{6}])$ 

in terms of MZVs, thereby completing the computation of 
the restriction of $I$ to $V_{7}$, $V_{5}$ and $V_{3}$. 
To summarize, the results of Sections \ref{sect:IIMZV}, \ref{sect:modular}
allow one to determine the restriction of $I$ to  
$\on{Lie}({\mathfrak M})_{10}$ in terms of MZVs of weight 10. 

$\bullet$ $k=14$. It has been shown in \cite{Po} that 
$[\delta_{2},\delta_{8}] =3 [\delta_{4},\delta_{6}]$. 
Using the same techniques as for $k=10$, one can prove that  
$81\cdot I([G_{4},G_{10}]) + 35 \cdot I([G_{6},G_{8}])$ is a MZV of 
weight 14. These techniques do not give any information on the individual values 
of $I([G_{4},G_{10}])$ and $I([G_{6},G_{8}])$.

\section{Galois aspects} 

In this section, we recall the links between $G_{\QQ}$, 
$\widehat{\on{GT}}$ and the Teichm\"uller groupoids in genus zero. We then 
establish the analogous results in genus one: they relate the arithmetic
fundamental group $\pi_{1}(M_{1,1}^{\QQ})$, $\widehat{\on{GT}}_{ell}$
and the Teichm\"uller groupoid in genus one.

\subsection{Galois groups and Teichm\"uller groupoids in genus zero}

\subsubsection{Profinite Galois representations}

Let $n\geq 3$ and $M_{0,n}^\QQ$ be the moduli stack over $\QQ$ of genus 
zero smooth projective curves with $n$ marked points and 
$\overline{M}_{0,n}^{\QQ}$ its Deligne-Mumford compactification. 
Maximally degenerate curves are rational points of this stack, and 
correspond bijectively to planar unrooted trivalent trees with 
leaves indexed bijectively by $\{1,\ldots,n\}$, modulo `mirror' symmetry. 
For $T$ such a tree, let $X_{T}^{0}$ the corresponding curve. 
The formal neighborhood of $X_{T}^{0}$ is a fibration 
$X_{T}\to\on{Spec}\QQ[[q_{e},e$ inner edge of $T]]$. 
Then the pull-back $X_{T}\otimes_{\QQ[[\{q_{e}\}_{e}]]}\QQ[[q]]$ 
corresponding to the morphism given by $q_{e}\mapsto q$ is a 
rational tangential base point of $M_{0,n}^{\QQ}$ (recall that a rational 
tangential base point of a scheme $X$ is a morphism $X\to \on{Spec}\QQ((q))$); 
see \cite{Mum,Ih:Nak}. 

Let $S$ be this set of rational tangential base points. The fundamental groupoid
$\widehat T_{0,n}:= \pi_{1}^{geom}(M_{0,n}^{\QQ},S)$ relative to this base 
set is the profinite completion of the groupoid $T_{0,n}$ described in \cite{Sch}. 
There is a split exact sequence 
$$
1\to \widehat{T}_{0,n}\to \pi_{1}(M_{0,n}^{\QQ},S)
\stackrel{\curvearrowleft}{\to} G_{\QQ}\to 1
$$
with section induced by $S$. It results in a group morphism 
 $G_\QQ = \on{Gal}(\bar\QQ/\QQ)\to 
\on{Aut}(\wh T_{0,n})$ (see \cite{G:Esq,Dr:Gal}). 

\begin{theorem} (\cite{Dr:Gal,Sch}) \label{thm:GT:1}This morphism factors as 
$G_\QQ\to \wh{\on{GT}}\to \on{Aut}(\wh T_{0,n})$. 
\end{theorem}

\subsubsection{Pro-$l$ and prounipotent completions}

Let $\pi$ be a finitely generated group, and let $\pi_{\QQ}(-)$ denote its
$\QQ$-prounipotent completion. It has the following properties:
$\pi_{\QQ}(-)$ is a prounipotent $\QQ$-group scheme; there is a group 
morphism $\pi\to\pi_{\QQ}(\QQ)$; any morphism 
$\pi\to U(\QQ)$, where $U(-)$ is a unipotent $\QQ$-group scheme, 
induces a $\QQ$-group scheme morphism $\pi_{\QQ}(-)\to U(-)$, such that 
$(\pi\to U(\QQ)) = (\pi\to \pi_{\QQ}(\QQ)\to U(\QQ))$. 

If ${\mathbf k}$ is a $\QQ$-ring, then $\pi_{{\mathbf k}}(-):= 
\pi_{\QQ}(-)\otimes{\mathbf k}$ is a prounipotent ${\mathbf k}$-group scheme
(it is the functor $\{{\mathbf k}$-rings$\}\to\{$groups$\}$, $K\mapsto 
\pi_{\QQ}(K)$).  
There is a morphism $(\pi\to\pi_{{\mathbf k}}({\mathbf k})):= 
(\pi\to \pi_{\QQ}(\QQ)\to \pi_{\QQ}({\mathbf k}) = 
\pi_{{\mathbf k}}({\mathbf k}))$. Any morphism $\pi\to U({\mathbf k})$, 
where $U(-)$ is a prounipotent ${\mathbf k}$-group scheme, gives rise to a 
morphism $\pi_{{\mathbf k}}(-)\to U(-)$, such that $(\pi\to U({\mathbf k})) 
= (\pi\to \pi_{{\mathbf k}}({\mathbf k})\to U({\mathbf k}))$ 
(\cite{H}, Section 3).  

Let $l$ be a prime number, and let $\pi_{l}$ be the prounipotent completion of 
$\pi$. According to \cite{HM}, Lemma A.7, there exists a morphism $\pi_{l}\to
\pi_{\QQ}(\QQ_{l})$, compatible with the maps from $\pi$.  

If $\pi,\pi'$ are finitely generated groups, then a continuous morphism 
$\pi_{l}\to\pi'_{l}$ gives rise to the 
morphism $\pi\to\pi_{l}\to \pi'_{l}\to \pi'_{\QQ}(\QQ_{l})$, and hence 
to a $\QQ_{l}$-group scheme morphism $\pi_{\QQ_{l}}(-)\to \pi'_{\QQ_{l}}(-)$, 
such that $(\pi\to \pi'_{\QQ}(\QQ_{l})) = (\pi\to \pi_{\QQ_{l}}(\QQ_{l})
\to \pi'_{\QQ_{l}}(\QQ_{l}))$. The resulting map $\on{Hom}(\pi_{l},\pi'_{l})
\to \on{Hom}_{\QQ_{l}\text{-group}\atop \text{schemes}}(\pi_{\QQ_{l}},\pi'_{\QQ_{l}})$
is compatible with compositions, hence gives rise to a group morphism 
\begin{equation} \label{morph:1}
\on{Aut}(\pi_{l})\to \on{Aut}_{\QQ_{l}\text{-group}\atop \text{schemes}}(\pi_{\QQ_{l}}).
\end{equation} 

Let $U(-)$ be a prounipotent $\QQ$-group scheme. Let $\ul{\on{Aut}U}$
be the $\QQ$-group scheme defined as the functor $\{\QQ$-rings$\}\to 
\{$groups$\}$, ${\mathbf k}\mapsto \ul{\on{Aut}U}({\mathbf k}) := 
\on{Aut}_{{\mathbf k}\text{-group}
\atop\text{schemes}}
(U\otimes {\mathbf k}) = \on{Aut}_{{\mathbf k}\text{-Lie}\atop
\text{algebras}}({\mathfrak u}\otimes {\mathbf k})$, 
where ${\mathfrak u} = \on{Lie}U$. Then $\ul{\on{Aut}U}$ is an extension
of a group $\QQ$-subscheme $G \subset \on{GL}({\mathfrak u}^{ab})$ by a 
prounipotent $\QQ$-group scheme, explicitly 
$$
1\to \underline{\on{Aut}^{+}U}\to \underline{\on{Aut}U}\to G \to 1.$$ 
Namely, $G$ is the intersection of the decreasing sequence of group schemes 
$\on{Im}(\underline{\on{Aut}U/U^{(n)}}\to \on{GL}({\mathfrak u}^{ab}))$, 
which is stationary. 

The morphism (\ref{morph:1}) may therefore be interpreted as a morphism 
$$
\on{Aut}(\pi_{l})\to \underline{\on{Aut}\pi}(\QQ_{l}). 
$$

Let $\cG\rightrightarrows B$ be a groupoid where for any $b\in B$, 
$\cG_{b}:= \cG_{bb}$ is finitely generated. We denote by 
$\cG_l\rightrightarrows B$, $\cG_{\QQ}(-)\rightrightarrows B$ 
its pro-$l$ and $\QQ$-prounipotent completions, given by 
$(\cG_l)_{bc}:= (\cG_{b})_l\times_{\cG_{b}} \cG_{bc}$ and  
$\cG_{\QQ}({\mathbf k})_{bc}:= \cG_{b}({\mathbf k})\times_{\cG_{b}} \cG_{bc}$. 

Assume that $\cG$ is connected (i.e., for any $b,c\in B$, $\cG_{bc}\neq\emptyset$).
Define the group scheme $\underline{\text{Aut}{\mathcal G}}$ by 
$\underline{\text{Aut}{\mathcal{G}}}({\mathbf{k}}):= 
\text{Aut}({\mathcal G}(\mathbf k))$. If 
${\mathcal G}_{ab}\times {\mathcal G}_{bc}\to 
{\mathcal G}_{ac}$, $(g_{ab},g_{bc})\mapsto g_{bc}g_{ab}$
is the composition of $\mathcal G$, then 
$\on{Aut}({\mathcal G}({\mathbf k})) = 
\{\theta_{ab} : {\mathcal G}_{ab}\to {\mathcal G}_{ab}({\mathbf k}) |
\forall a,b,c, \forall g_{ab}, g_{bc}, \theta_{ac}(g_{bc}g_{ab}) = 
\theta_{bc}(g_{bc})\theta_{ab}(g_{ab})\}$. The choice of $b\in B$
and of particular elements $g_{ab}^0\in {\mathcal G}_{ab}$ for any 
$a\in B - \{b\}$ gives rise to an isomorphism 
$\on{Aut}({\mathcal G}({\mathbf k}))\simeq 
{\mathcal G}_b({\mathbf k})^{B - \{b\}}\rtimes 
\underline{\text{Aut}{\mathcal G}_b}({\mathbf k})$, the inverse 
isomorphism taking $((X_a)_a,\theta)$ to the automorphism 
such that ${\mathcal G}_b({\mathbf k})\ni g_b\mapsto \theta(g_b)\in 
{\mathcal G}_b({\mathbf k})$, and ${\mathcal G}_{ab}\ni g_{ab}^0\mapsto 
X_a g_{ab}^0 \in {\mathcal G}_{ab}({\mathbf k})$. 
The morphisms $\pi_{l}\to\pi_{\QQ}(\QQ_{l})$ and 
$\on{Aut}(\pi_{l})\to \underline{\on{Aut}\pi}(\QQ_{l})$, 
where $\pi = {\mathcal G}_{b}$, give rise to a morphism 
$$
\on{Aut}(\cG_{l})\to \underline{\on{Aut}\cG}(\QQ_{l}).
$$ 

\subsubsection{Pro-$l$ Galois representations}
 
The following statement can be derived from \cite{Dr:Gal,Sch}. 

\begin{proposition} 
There exist morphisms $\on{GT}_{l}\to \on{Aut}(T_{0,n}^{l})$, 
$\on{GT}(-)\to \ul{\on{Aut}T_{0,n}}(-)$, such that 
the squares in the following diagram commute 
$$
\xymatrix{
\wh{\on{GT}} \ar[r]\ar[d]& \on{GT}_{l} \ar[r]\ar[d]& \on{GT}(\QQ_{l})\ar[d]\\
\on{Aut}(\wh T_{0,n}) \ar[r]& \on{Aut}(T_{0,n}^{l}) \ar[r]& 
\ul{\on{Aut}T_{0,n}}(\QQ_{l})
}$$
\end{proposition}

%
%
%

\subsection{Arithmetic fundamental groups and Teichm\"uller groupoids in 
genus one} \label{sec:main:results}

The Galois theoretic counterpart of the theory of elliptic associators
is the action of the arithmetic fundamental group $\pi_1(M_{1,\vec{1}}^\QQ)$ 
on the completions of elliptic braid groups, based on the fibration 
$M_{1,n}^\QQ\to M_{1,1}^\QQ$, as studied in \cite{G:LM,O}. 
We first recall the main points of this study.

\subsubsection{Arithmetic fundamental groups of moduli spaces}

Let $M_{1,1}^{\QQ}$ (resp., $M_{1,\vec{1}}^{\QQ}$, $\tilde M_{1,1}^{\QQ}$)
be the moduli space of elliptic curves with one puncture (resp., with one puncture
and a nonzero tangent vector at the puncture, with one puncture and a formal coordinate
at the puncture). 

A rational tangential base point $\xi$ of $M_{1,1}^{\QQ}$ is defined as follows. The 
Deligne-Mumford compactification $\overline M_{1,1}^{\QQ}$ of $M_{1,1}^{\QQ}$ 
contains a unique curve $X^{0}$, which corresponds to the tadpole graph. 
A formal neighborhood of $X^{0}$ in $\overline M_{1,1}^{\QQ}$ is 
a curve $X\to \on{Spec}\QQ[[q]]$, whose generic fiber
is the Tate elliptic curve ${\mathbb G}_{m}/q^{\ZZ}$ with marked 
point $[1] = q^{\ZZ}$. This may be viewed as a morphism
$\on{Spec}\QQ[[q]]\to \overline M_{1,1}^{\QQ}$, which 
restricts to $\xi : \on{Spec}\QQ((q))\to  M_{1,1}^{\QQ}$. 

A lift $\tilde\xi$ of $\xi$ to $\tilde M_{1,1}^{\QQ}$ is defined 
by choosing the local coordinate $\on{log}z$ at $[1] = q^{\ZZ}$, 
$z$ being the canonical coordinate on ${\mathbb G}_{m}$ (such that the 
function ring is $\QQ[z,z^{-1}]$). Let $\vec\xi$ be the lift of 
$\xi$ to $M_{1,\vec1}$ given by the expansion of the local coordinate
of $\tilde\xi$ at order one. 

The isomorphism $\pi_{1}^{geom}(M_{1,\vec{1}}^{\QQ},\vec{\xi})
\simeq \wh B_{3}$ gives rise to a split exact sequence 
\begin{equation} \label{ES}
1\to \widehat B_3 \to \pi_1(M_{1,\vec{1}}^\QQ,\vec\xi)
\stackrel{\curvearrowleft}{\to} G_\QQ\to 1, 
\end{equation}
where the section is provided by the base point $\vec\xi$; the induced morphism 
$G_{\QQ}\to \on{Aut}(\widehat{B}_{3})$ has been computed
explicitly in \cite{Na}, Cor. 4.15 (it is recalled in Subsection \ref{sec:morphism}). 

The result of \cite{Na} can be complemented as follows. 

\begin{proposition} \label{prop:CD}
There is a morphism from (\ref{ES}) to the split exact sequence 
\begin{equation} \label{ES'}
1\to \on{SL}_{2}(\widehat\ZZ)\to \on{GL}_{2}(\widehat\ZZ) 
\stackrel{\curvearrowleft}{\to}\ZZ^{\times}\to 1, 
\end{equation}
where the second morphism is the determinant $\on{det}$ and the 
section is the morphism $\lambda\mapsto 
\bigl(\begin{smallmatrix} \lambda & 0 \\ 0 & 1
\end{smallmatrix}\bigl)$.  The rightmost morphism in 
$(\ref{ES})\to (\ref{ES'})$ is the cyclotomic character $G_{\QQ}\to
\widehat\ZZ^{\times}$. 
\end{proposition}

The proof will be carried out in Section \ref{sec:morphism}. 

\subsubsection{Profinite representations}

Let $\tilde M_{1,n}^{\QQ}$ be the moduli space of elliptic curves with $n$
punctures. There is a fibration $M_{1,n}^{\QQ}\to M_{1,1}^{\QQ}$ defined by 
forgetting all the punctures except the first one. One sets $\tilde M_{1,n}^{\QQ}
:= \tilde M_{1,1}^{\QQ}\times_{M_{1,1}^{\QQ}}M_{1,n}^{\QQ}$. 

A {\it tangential section} of a morphism $X\to Y$ of $\QQ$-schemes is defined 
to be a morphism $Y\times \on{Spec}\QQ((t))\to X$, such that 
its composition with $X\to Y$ is the canonical projection. 

A {\it $n$-tree} $T$ is defined to be a rooted trivalent planar tree, equipped with a 
bijection $i_{T} : \{\text{leaves}\} \to \{1,\ldots,n\}$ (the root is not a leaf), such 
that the leftmost leaf is labeled $1$. 
Such a tree gives rise to the assignment, to each $i\in\{1,\ldots,n\}$, of a pair
$(d_{i},s_{i})$, where  $d_{i}$ is an integer $\geq 1$ (the distance between 
the leaf labeled $i$ and the root), and of a map $s_{i}\in \{l,r\}^{d_{i}}$
describing the path from the root to the leaf labeled $i$ ($s_{i}(k) = l$ or $r$
according to whether the $k$th interval of the path is a left or right descendant). 
It also gives rise to a permutation $s_{T}\in S_{n}$ such that $s_{T}(1)=1$: 
$s_{T}$ is the composite map $\{1,\ldots,n\}\to \{\on{leaves}\}
\stackrel{i_{T}}{\to} \{1,\ldots,n\}$, where the first map is the inverse of the 
lexicographic (according to the order $\text{left} < \text{right}$) 
indexation of the leaves.  

A tangential section $\sigma_{T}$ of the morphism $\tilde M_{1,n}^{\QQ}\to 
\tilde M_{1,1}^{\QQ}$ may be associated to each $n$-tree $T$ as follows: 
$\sigma_{T}$ is the morphism $\tilde M_{1,1}^{\QQ}\times
\on{Spec}\QQ((t))\to \tilde M_{1,n}$, taking a pair 
$((E,p,z),t)$ to $(E,p_{1},\ldots,p_{n},z)$, where 
$p_{i}:= z^{-1}(\sum_{k\in s_{i}^{-1}(r)} t^{k})$. 

Let ${\mathcal F}_{\xi}$ be the fiber over $\xi$ of $M_{1,n}^{\QQ}\to 
M_{1,1}^{\QQ}$. There is a split exact sequence of groupoids
$$
1\to \pi_{1}^{geom}({\mathcal F}_{\xi},\{\sigma_{T}(\xi)\})
\to \pi_{1}(\tilde M_{1,n}^{\QQ},\{\sigma_{T}(\tilde\xi)\})
\stackrel{\curvearrowleft}{\to} \pi_{1}(\tilde M_{1,1}^{\QQ},
\tilde\xi)\to 1. 
$$
(see \cite{G:LM,O} and also \cite{Na}, Section 4.1), which gives rise to a morphism 
\begin{equation} \label{morphism}
\pi_{1}(\tilde M_{1,1}^{\QQ},\tilde\xi) \to 
\on{Aut}(\pi_{1}^{geom}({\mathcal F}_{\xi},\{\sigma_{T}(\xi)\})). 
\end{equation}
The fiber at $(E,p)$ of $M_{1,n}^{\QQ}\to M_{1,1}^{\QQ}$ is $(E - \{p\})^{n-1}
- (\on{diagonals})$, whose geometric fundamental group is the profinite completion of 
$\overline P_{1,n}$ (the quotient of the elliptic braid group with $n$ strands 
$P_{1,n}$ by the central $\ZZ^{2}$). The geometric fundamental groupoid
$\pi_{1}^{geom}({\mathcal F}_{\xi},\{\sigma_{T}(\xi)\})$ is the profinite
completion of the groupoid $T_{ell,n}$ where objects are $n$-trees and the set of 
morphisms from $T$ to $T'$ is $\overline P_{1,n}\times_{S_{n}}
\{s_{T'}s_{T}^{-1}\}$, equipped with the composition of morphisms 
induced from the product in $\overline P_{1,n}$. On the other hand, there is an 
isomorphism $\pi_{1}(\tilde M_{1,1}^{\QQ},\tilde\xi)\simeq 
\pi_{1}(M_{1,\vec{1}}^{\QQ},\vec{\xi})$. (\ref{morphism}) 
therefore gives rise to a morphism 
\begin{equation} \label{ell:morph}
\pi_{1}(M_{1,\vec{1}}^{\QQ},\vec{\xi}) \to 
\on{Aut}(\widehat{T}_{ell,n}). 
\end{equation}

\begin{theorem} \label{thm:GTell:LM}
There exists a morphism $\pi_1(M_{1,\vec{1}}^{\QQ},\vec{\xi})
\to\widehat{\on{GT}}_{ell}$ and an action of 
$\widehat{\on{GT}}_{ell}$ on $\widehat T_{ell,n}$, such that: 

(a) the morphism (\ref{ell:morph}) factors as 
$\pi_1(M_{1,\vec{1}}^{\QQ},\vec{\xi})\to\widehat{\on{GT}}_{ell}\to 
\on{Aut}(\widehat T_{ell,n})$; 

(b) the morphism of split morphisms induced by $(\ref{ES})\to (\ref{ES'})$ 
factors as $(\pi_{1}(M_{1,\vec1}^{\QQ},\vec\xi)
\stackrel{\curvearrowleft}{\to} G_{\QQ})
\to (\widehat{\on{GT}}_{ell} \stackrel{\curvearrowleft}{\to}
\widehat{\on{GT}})\to(
\on{GL}_{2}(\wh\ZZ)\stackrel{\curvearrowleft}{\to}
\wh\ZZ^{\times})$. 
\end{theorem}

The proof will be carried out in Section \ref{sec:pf:LM}. 

\subsubsection{Pro-$l$ representations}


\begin{proposition} \label{thm:details:GTell}

There exist morphisms $\on{GT}_{ell}^{l}\to \on{Aut}T_{ell,n}^{l}$, 
$\on{GT}_{ell}(-)\to \ul{\on{Aut}T_{ell,n}}(-)$, such that the squares
in the following diagram commute
$$
\xymatrix{
\wh{\on{GT}}_{ell} \ar[r]\ar[d]& \on{GT}_{ell}^{l} 
\ar[r]\ar[d]& \on{GT}_{ell}(\QQ_{l})\ar[d]\\
\on{Aut}(\wh T_{ell,n}) \ar[r]& \on{Aut}(T_{ell,n}^{l}) \ar[r]& 
\ul{\on{Aut}T_{ell,n}}(\QQ_{l})
}$$
\end{proposition}

%

\subsection{Proof of Proposition \ref{prop:CD}} \label{sec:morphism}

As in \cite{Na}, let ${\mathfrak s}_{0} : G_{\QQ}\to\pi_{1}
(M_{1,\vec{1}}^{\QQ},\vec{\xi})$
be the section induced by $\vec{\xi}$. The diagram 
$\pi_{1}(M_{1,\vec{1}}^{\QQ},\vec{\xi})\stackrel{\curvearrowleft}{\to} 
G_{\QQ}$ gives rise to the semidirect product decomposition
$\pi_{1}(M_{1,\vec{1}}^{\QQ},\vec{\xi})\simeq \widehat B_{3}\rtimes G_{\QQ}$, 
where the action $G_{\QQ}\to \on{Aut}(\widehat B_{3})$ is $g*x:= 
{\mathfrak s}_{0}(g)x {\mathfrak s}_{0}(g)^{-1}$. On the other hand, 
the diagram $\on{GL}_{2}(\ZZ)\stackrel{\curvearrowleft}{\to} 
\widehat\ZZ^{\times}$ gives rise to the semidirect product 
decomposition $\on{GL}_{2}(\widehat\ZZ)\simeq \on{SL}_{2}(\ZZ)
\rtimes\widehat\ZZ^{\times}$, where the action $\widehat\ZZ^{\times}
\to\on{Aut}(\on{SL}_{2}(\widehat\ZZ))$ is $\lambda \bullet m:= 
\bigl(\begin{smallmatrix}\lambda & 0 \\
 0 & 1\end{smallmatrix}\bigl) m 
  \bigl(\begin{smallmatrix}\lambda & 0 \\
 0 & 1\end{smallmatrix}\bigl)^{-1}$. 
 
 Let $\sigma_{1},\sigma_{2}$ be the Artin generators of $B_{3}$
 (denoted $\bar a_{1},\bar a_{2}$ in \cite{Na}). As $\on{SL}_{2}(\widehat\ZZ)$
 is profinite, there is a unique morphism 
\begin{equation} \label{mor:B3}
 \widehat B_{3}\to \on{SL}_{2}(\widehat\ZZ), \quad x\mapsto
 \overline{x}
\end{equation}
 extending the quotient morphism $B_{3}\to B_{3}/Z(B_{3})\simeq 
 \on{SL}_{2}(\ZZ)$, $\sigma_{1}\mapsto  
\bigl(\begin{smallmatrix} 1 & 1 \\
 0 & 1\end{smallmatrix}\bigl)$,  $\sigma_{2}\mapsto  
\bigl(\begin{smallmatrix} 1 & 0 \\
 -1 & 1\end{smallmatrix}\bigl)$. 
 
The action of $G_{\QQ}$ on $\widehat B_{3}$ can be made explicit as follows. 
Denote the map $G_{\QQ}\to \widehat{\on{GT}} \subset 
 \widehat\ZZ^{\times}\times \widehat F_{2}$ by $g\mapsto 
 (\chi(g),f_{g})$. Using the formula $\beta_{0}(g)=\sigma_{1}^{8\rho_{2}(g)}
 {\mathfrak s}_{0}(g)$ in \cite{Na} before Proposition 4.12, and Corollary 4.15 in the 
 same paper, one obtains 
 $$
 g*\sigma_{1}=\sigma_{1}^{\chi(g)}, \quad
 g*\sigma_{2} = \on{Ad}_{\sigma_{1}^{-8\rho_{2}(g)}f_{g}(\sigma_{1}^{2},
 \sigma_{2}^{2})^{-1}} (\sigma_{2}^{\chi(g)}) 
 $$
 (here $\rho_{2} : G_{\QQ}\to\widehat\ZZ$ is the Kummer cocycle related to the 
 roots of 2).
 
 Then 
$$
\overline{g*\sigma_{1}} = \overline{\sigma_{1}^{\chi(g)}} = 
\bigl(\begin{smallmatrix} 1 & \chi(g) \\
0 & 1 \end{smallmatrix}\bigl) = 
\bigl(\begin{smallmatrix} \chi(g) & 0 \\
0 & 1 \end{smallmatrix}\bigl) 
\bigl(\begin{smallmatrix} 1 & 1 \\
0 & 1 \end{smallmatrix}\bigl) 
\bigl(\begin{smallmatrix} \chi(g) & 0 \\
0 & 1 \end{smallmatrix}\bigl) ^{-1} = \chi(g)\bullet
\bigl(\begin{smallmatrix} 1 & 1 \\
0 & 1 \end{smallmatrix}\bigl)  = \chi(g)\bullet \overline{\sigma}_{1}; 
$$
on the other hand, Corollary 4.15 in \cite{Na} says that 
$$
f_{g}(\bigl(\begin{smallmatrix} 1 & 2 \\
 0 & 1 \end{smallmatrix}\bigl), 
 \bigl(\begin{smallmatrix} 1 & 0 \\
 -2 & 1 \end{smallmatrix}\bigl)) = 
 \pm \bigl(\begin{smallmatrix} 1 & 0 \\
 -8\rho_{2}(g) & 1 \end{smallmatrix}\bigl)
  \bigl(\begin{smallmatrix} \chi(g)^{-1} & 0 \\
 0 & \chi(g) \end{smallmatrix}\bigl)
 \bigl(\begin{smallmatrix} 1 & -8\rho_{2}(g)  \\
 0 & 1 \end{smallmatrix}\bigl)
 $$
 (identity in $\on{SL}_{2}(\widehat\ZZ)$), 
therefore 
\begin{eqnarray*}
& \overline{g*\sigma_{2}} = 
\overline{\on{Ad}_{\sigma_{1}^{-8\rho_{2}(g)}
f_{g}(\sigma_{1}^{2},\sigma_{2}^{2})^{-1}}\sigma_{2}^{\chi(g)}}
= \on{Ad}_{\pm 
\bigl(\begin{smallmatrix}\chi(g) & 0 \\
0 & \chi(g)^{-1}\end{smallmatrix}\bigl)
\bigl(\begin{smallmatrix} 1 & 0 \\
8\rho_{2}(g) & 1\end{smallmatrix}\bigl)}
\bigl(\begin{smallmatrix} 1 & 0 \\
-\chi(g) & 1\end{smallmatrix}\bigl)\\
& 
=\bigl(\begin{smallmatrix} 1 & 0 \\
-\chi(g)^{-1} & 1\end{smallmatrix}\bigl)  
=\bigl(\begin{smallmatrix} \chi(g) & 0 \\
0 & 1\end{smallmatrix}\bigl)
\bigl(\begin{smallmatrix} 1 & 0 \\
-1 & 1\end{smallmatrix}\bigl)
\bigl(\begin{smallmatrix} \chi(g) & 0 \\
0 & 1\end{smallmatrix}\bigl)^{-1}
= \chi(g)\bullet \overline{\sigma}_{2}. \\
\end{eqnarray*}
It follows that (\ref{mor:B3})
intertwines the actions of $G_{\QQ}$ and $\widehat\ZZ^{\times}$, 
which proves Proposition \ref{prop:CD}. \hfill \qed\medskip 

\subsection{Proof of Theorem \ref{thm:GTell:LM}} \label{sec:pf:LM}

Theorem \ref{thm:GTell:LM} states the existence of a morphism 
$\pi_1(M_{1,\vec{1}}^\QQ,\vec\xi)\to\wh{\on{GT}}_{ell}$, 
which will now be constructed. 

\begin{proposition} \label{prop:mor:SD}

Set $\wh R_{ell}:= \on{Ker}(\wh{\on{GT}}_{ell}\to\wh{\on{GT}})$. 

(a) There is a unique morphism $\wh B_{3}\to \wh R_{ell}$, 
extending the canonical morphism $B_{3}\simeq R_{ell}\to \wh R_{ell}
(\subset \on{Aut}(\wh F_{2}))$. 

(b) There is a unique morphism 
$\pi_1(M_{1,\vec{1}}^\QQ,\vec\xi)\to\wh{\on{GT}}_{ell}$, 
such that the diagram 
$$
\begin{matrix}
1& \to & \wh B_{3}& \to & \pi_{1}(M_{1,\vec{1}}^\QQ,\vec\xi) & 
 \stackrel{\stackrel{{\mathfrak s}_0}{\curvearrowleft}}{\to} 
& G_{\QQ}  & \to & 1 \\
&  & \downarrow&  &  \downarrow & & \downarrow & &   \\
1& \to & \wh R_{ell}& \to &\widehat{\on{GT}}_{ell} & 
  \stackrel{\stackrel{\sigma}{\curvearrowleft}}{\to} 
 & \widehat{\on{GT}}     & \to & 1
\end{matrix}
$$
commutes. 
\end{proposition}

{\em Proof.} (a) Recall that $\widehat{R}_{ell}$ is a subgroup of 
$\on{Aut}(\widehat F_{2})$. $\on{Aut}(\widehat F_{2})$ is profinite 
(\cite{DDMS}, Thm.~5.3), and the map $\on{Aut}(\widehat F_{2})\to 
\widehat F_{2}^{2}$, $\theta\mapsto (\theta(X),\theta(Y))$ is continuous (loc.~cit., 
Ex.~2, p.~96). As $\widehat{R}_{ell}$ is the preimage of $1$ by a
continuous map $(\widehat F_{2})^{2}\to 
(\widehat{B}_{1,3})^{2}\times \widehat F_{2}$, it is closed, 
so $\widehat{R}_{ell}$ is a closed subgroup of $\on{Aut}(\widehat F_{2})$, 
hence is profinite. The morphism 
$B_3\to \widehat R_{ell}$ therefore extends to a morphism 
$\wh B_3 \to \widehat R_{ell}$. 

Statement (b) is equivalent of the compatibility of 
the morphism $\wh B_3 \to \widehat R_{ell}$ with the actions of 
$G_\QQ$ and $\wh{\on{GT}}$ on both sides via ${\mathfrak s}_0$ and $\sigma$
and the morphism $G_{\QQ}\to \widehat{\on{GT}}$, i.e., to the 
commutativity of 
\begin{equation} \label{diag:gal}
\begin{matrix}
G_{\QQ}\times \wh B_3 & \to & \wh B_3 \\
\downarrow  & & \downarrow \\
\widehat{\on{GT}}\times \widehat{R}_{ell}   & \to & 
\widehat{R}_{ell}\end{matrix}
\end{equation}
Consider the following cubic diagram 
$$
\xymatrix{
G_{\QQ}\times\pi_{1}^{geom}(M_{1,\vec1},\vec\xi)
 \ar[rr]\ar[dd]\ar[dr] & & \pi_{1}^{geom}(M_{1,\vec1},\vec\xi)
  \ar[dd]\ar[dr] & \\
& G_{\QQ}\times \on{Aut}\pi_{1}^{geom}(C_{\xi},\xi_{C})\ar[rr]
\ar[dd] & &  \on{Aut}\pi_{1}^{geom}(C_{\xi},\xi_{C}) \ar[dd] \\
\wh{\on{GT}}\times \wh R_{ell} \ar[rr]\ar[dr] & & \wh R_{ell}\ar[dr] & \\
& \wh{\on{GT}}\times \on{Aut}\wh F_{2}\ar[rr] & & \on{Aut}\wh F_{2} }
$$
where $C_{\xi}$ is the fiber of $M_{1,2}^{\QQ}\to M_{1,1}^{\QQ}$ at $\xi$
(which identifies with the fiber of $M_{1,\vec 2}^{\QQ}\to M_{1,\vec 1}^{\QQ}$, 
where $M_{1,\vec 2}^{\QQ} = M_{1,\vec 1}^{\QQ}\times_{M_{1,1}^{\QQ}}
M_{1,2}^{\QQ}$), $\xi_{C}$ is a tangential base point of $C_{\xi}$ supported
at the marked point, and the maps are defined as follows : 

$\bullet$ the upper horizontal maps are the Galois actions; the map 
$\wh{\on{GT}}\times \wh R_{ell}\to\wh R_{ell}$ is the action induced by the section
of $\wh{\on{GT}}_{ell}\to \wh{\on{GT}}$ defined in Proposition \ref{prop:lift}; the 
map $\wh{\on{GT}}\times \on{Aut}\wh F_{2}\to \on{Aut}\wh F_{2}$
is induced by the composite map $\wh{\on{GT}}\to \on{Aut}\wh F_{2}
\to \on{Aut}(\on{Aut}\wh F_{2})$, where the second map is the inner action of 
$\on{Aut}\wh F_{2}$ on itself, and the first map is the composite morphism 
$\wh{\on{GT}}\to \wh{\on{GT}}_{ell} \subset \wh{\on{GT}}\times
\on{Aut}\wh F_{2}\to \on{Aut}\wh F_{2}$, where 
$\wh{\on{GT}}_{ell}\to \wh{\on{GT}}$ is the same morphism as above and 
$\wh{\on{GT}}\times\on{Aut}\wh F_{2}\to \on{Aut}\wh F_{2}$ is the second 
projection; 

$\bullet$ the vertical maps are induced by the morphisms $G_{\QQ}\to \wh{\on{GT}}$,
$\pi_{1}^{geom}(M_{1,\vec1}^{\QQ},\vec\xi)\simeq\wh B_{3}\to \wh R_{ell}$, 
$\pi_{1}^{geom}(C_{\xi},\xi_{C})\stackrel{\sim}{\to}\wh F_{2}$; 

$\bullet$ the diagonal maps are induced by the canonical inclusion $\wh R_{ell}
\to \on{Aut}\wh F_{2}$, and by the action of $\pi_{1}^{geom}(M_{1,\vec1},
\vec\xi)$ on $\pi_{1}^{geom}(C_{\xi},\xi_{C})$ induced by the fibration 
$M_{1,\vec2}^{\QQ}\to M^{\QQ}_{1,\vec1}$. 

The square corresponding to the upper face of the cube commutes because the
action of $\pi_{1}^{geom}(M_{1,\vec1},\vec\xi)$ on $\pi_{1}^{geom}(C_{\xi},
\xi_{C})$ is compatible with the Galois action. 

The square corresponding to the sides of the cube commute because this action 
identifies with the profinite completion of the action of $B_{3}$ on $F_{2}$. 

The square corresponding to the lower face of the cube commutes by construction of the 
map $\wh{\on{GT}}\times \on{Aut}\wh F_{2}\to \on{Aut}\wh F_{2}$. 

The square corresponding to the lower front face commutes for the following 
reason. According to \cite{Na}, Corollary 4.5, the action of $G_{\QQ}$ on 
$\pi_{1}^{geom}(C_\xi,\xi_{C})$ may be described as follows. 
$\pi_{1}^{geom}(C_\xi,\xi_{C})$ is topologically free, generated 
by $x_{1},x_{2}$. The action of $\gamma\in G_{\QQ}$ on this group is
\begin{equation} \label{act:nak:1}
\gamma^{*}(x_{1}) = f_{\gamma}(x_{1},z_{1})x_{1}^{\chi(\gamma)}
f_{\gamma}(x_{1},z_{1})^{-1}, 
\end{equation}
\begin{equation} \label{act:nak:2}
\gamma^{*}(x_{2}) =  f_{\gamma}(x_{1},z_{1})x_{1}^{1-\chi(\gamma)}
 f_{\gamma}^{\stackrel{\to}{\infty 1}}(x_{1},x_{1}^{-1}z_{1}x_{1})^{-1}
x_{2} x_{1}^{\chi(\gamma)-1} f_{\gamma}(x_{1},z_{1})^{-1},
\end{equation}
where $z_{1} = (x_{2},x_{1}) = x_{2}x_{1}x_{2}^{-1}x_{1}^{-1}$,  
$\gamma\mapsto 
(\chi(\gamma),f_{\gamma})$ is the map $G_{\QQ}\to 
\widehat{\on{GT}}$, and $f_{\gamma}^{\stackrel{\to}{\infty1}}(a,b) 
= f_{\gamma}(b,c)b^{(\chi(\gamma)-1)/2}f_{\gamma}(a,b)$ for $abc=1$. 

Under the identification $x_{1}\mapsto X$, $x_{2}\mapsto Y$, formula 
(\ref{act:nak:1}) corresponds to the expression of $g_{+}$ in Proposition 
\ref{prop:lift}. It follows from the hexagon and duality identities that any 
$(\lambda,f)\in\wh{\on{GT}}$ satisfies the octagon identity 
$$
f(X^{-1}Z^{-1},Z)(ZX)^{-\lambda}f(Z,X^{-1}Z^{-1})Z^{(\lambda+1)/2}
f(X,Z)X^{\lambda}f(Z,X)Z^{(\lambda-1)/2}=1, 
$$
where $Z:= (Y,X)$.  This identity implies 
\begin{eqnarray*}
& f(X,(Y,X))X^{1-\lambda}f^{\stackrel{\to}{\infty1}}(X,(X^{-1},Y))^{-1}
YX^{\lambda-1}f(X,(Y,X))^{-1} \\
 & = Z^{(\lambda-1)/2}f(X^{-1}Z^{-1},Z)Y
f(X,(Y,X))^{-1} 
\end{eqnarray*}
so that (\ref{act:nak:2}) corresponds to $g_{-}$ in  Proposition 
\ref{prop:lift}. All this implies the commutativity of 
$$
\begin{matrix}
G_{\QQ} & \to & 
\on{Aut} \pi_{1}^{geom}(C_\xi,\xi_{C})\\
\downarrow   & & \downarrow \\
\widehat{\on{GT}}_{ell}   & \to & \on{Aut}\widehat{F}_{2}
\end{matrix}$$
Composing this square with the commutative square 
$$
\begin{matrix}
\on{Aut}\pi_{1}^{geom}(C_{\xi},\xi_{C}) & \to & 
\on{Aut}(\on{Aut}\pi_{1}^{geom}(C_{\xi},\xi_{C}))\\
\downarrow & & \downarrow \\
\on{Aut}\wh F_{2} & \to & \on{Aut}(\on{Aut}\wh F_{2})
\end{matrix}
$$
where the horizontal maps are inner action morphisms, one obtains the commutativity of 
the square corresponding to the lower front face. 

The commutativity of all these squares implies that the two composite maps 
$$
G_{\QQ}\times\pi_{1}^{geom}(M_{1,\vec1},\vec\xi)\to\wh R_{ell}
\to\on{Aut}\wh F_{2}
$$ 
coincide, where the maps $G_{\QQ}\times \pi_{1}^{geom}
(M_{1,\vec1},\vec\xi)\to\wh R_{ell}$ are the two composite maps which can be 
obtained from the upper front face. As $\wh R_{ell}\to \on{Aut}\wh F_{2}$
is injective, this implies the commutativity of the square corresponding to the upper
front face, and therefore of (\ref{diag:gal}). 
\hfill \qed\medskip 

The next statement of Theorem \ref{thm:GTell:LM} is the existence of an action of 
$\wh{\on{GT}}_{ell}$ on $\wh T_{ell,n}$, which will now be constructed
(Definition \ref{def:act}).

 If $\cC$ is a category, let $\on{Aut}(\cC)$ be its group of automorphisms 
 (as a category, even if $\cC$ has a braided monoidal structure).    

For $(\lambda,f)\in\wh{\ul{\on{GT}}}$, let $i_{\lambda,f}$ be the composite functor 
$\wh{\bf PaB}\stackrel{\alpha_{(\lambda,f)}}{\to} (\lambda,f)*\wh{\bf PaB}
\stackrel{\sim}{\to}
\wh{\bf PaB}$, where the first functor is the unique tensor functor which 
induces the identity on objects, and the second functor is the identity functor (which
is not tensor). $i_{\lambda,f}$ is then an endofunctor of $\wh{{\bf PaB}}$. 

\begin{lemma} \label{prop:aut:BMC}
$(\lambda,f)\mapsto i_{\lambda,f}^{-1}$ is a 
morphism $\wh{\on{GT}}\to\on{Aut}(\wh{{\bf PaB}})$. 
\end{lemma}

{\em Proof.} 
The identity $i_{(\lambda',f')}i_{(\lambda,f)} = 
i_{(\lambda,f)(\lambda',f')}$ follows from the commutativity of the 
diagram 
$$
\xymatrix{
 &   & \cC\ar[dr]
 \ar@/^{2pc}/[ddrr]^{i_{(\lambda',f')}} \\
 & (\lambda,f)*\cC\ar[ur]^{\sim}\ar[dr] & & (\lambda',f')*\cC\ar[dr]^{\sim} \\
\cC\ar[rr]\ar[ur]\ar@/^{2pc}/[uurr]^{i_{(\lambda,f)}}
\ar@/_{1pc}/[rrrr]_{i_{(\lambda,f)
(\lambda',f')}} & & (\lambda,f)(\lambda',f')
*\cC\ar[rr]^{\sim}\ar[ur]^{\sim} 
& & \cC 
}
$$
in which the commutativity of the central square follows from that of
$$\begin{matrix}
(\lambda,f)*\cC & \stackrel{(\lambda,f)*\varphi}{\to}& (\lambda,f)*\cD
\\
 \scriptstyle{\sim}\downarrow & & \downarrow\scriptstyle{\sim}\\
 \cC & \stackrel{\varphi}{\to}& \cD
\end{matrix}$$
for any braided monoidal categories $\cC,\cD$ and any tensor functor
$\varphi:\cC\to\cD$.  \hfill\qed\medskip 

One constructs in the same way a morphism 
\begin{equation} \label{GTell:to:AutPaB}
\wh{\on{GT}}_{ell}\to \on{Aut}(\wh{\bf PaB}_{ell}).
\end{equation} 
 
If $\cC_{0}$ is a braided monoidal category, then $\on{Ob}\cC_{0}$ is a 
magma (i.e., a set equipped with a composition map and a unit). Let $\phi : M\to
\on{Ob}\cC_{0}$ be a magma morphism, then a braided monoidal 
category $\phi*\cC_{0}$ can be constructed by $\on{Ob}\phi^{*}\cC_{0} = M$, 
$\phi^{*}\cC_{0}(m,m'):= \cC_{0}(\phi(m),\phi(m'))$ and by the 
condition that the obvious functor $\phi^{*}\cC_{0}\to\cC_{0}$ is tensor. 
If $\cC_{0}\to\cC$ is an elliptic structure over $\cC_{0}$, then one 
defines an elliptic structure $\phi^{*}\cC_{0}\to\phi^{*}\cC$ over
$\phi^{*}\cC_{0}$ in the same way. Then there are natural group morphisms
\begin{equation} \label{gp:morph}
\on{Aut}\cC_{0}\to \on{Aut}\phi^{*}\cC_{0}, \quad 
\on{Aut}\cC\to \on{Aut}\phi^{*}\cC. 
\end{equation}

Let $\mu(S)$ be the free magma generated by a set $S$. 
The unique map $S\to\{\bullet\}$ induces a magma morphism $\phi:
\mu(S)\to \mu(\{\bullet\})\simeq \on{Ob}\wh{\bf PaB}$. Set 
$\wh{\bf PaB}_{S}:= \phi^{*}\wh{\bf PaB}$, 
$\wh{\bf PaB}_{ell,S}:= \phi^{*}\wh{\bf PaB}_{ell}$. 
The morphisms (\ref{gp:morph}) then specialize to morphisms 
\begin{equation} \label{gp:morph'}
\on{Aut}(\wh{\bf PaB}) \to \on{Aut}(\wh{\bf PaB}_{S}), \quad 
\on{Aut}(\wh{\bf PaB}_{ell}) \to \on{Aut}(\wh{\bf PaB}_{ell,S}). 
\end{equation}

$\wh T_{ell,n}$ may be viewed as the full subcategory of 
$\wh{\bf PaB}_{ell,[n]}$, where the objects are the preimages of 
$1\dot+\cdots\dot+ n$ under the map $\mu([n])\stackrel{\psi}{\to}\NN[n]$, 
extending the identity on $[n]$, where $\NN[n]$ is the free abelian semigroup 
generated by $[n] = \{1,\ldots,n\}$ (in which the addition is denoted $\dot+$). 
If $\cC$ is any category and $\cC'$ is any full subcategory, then there is a natural 
morphism $\on{Aut}(\cC)\to \on{Aut}(\cC')$. It specializes to a group morphism 
\begin{equation} \label{Aut:to:AutT}
\on{Aut}(\wh{\bf PaB}_{ell,[n]})\to \on{Aut}(\wh T_{ell,n}). 
\end{equation}

\begin{definition} \label{def:act}
The action of $\wh{\on{GT}}_{ell}$ on $\wh T_{ell,n}$
is given by the composite morphism 
$$
\wh{\on{GT}}_{ell}\to  \on{Aut}\wh{\bf{PaB}}_{ell}\to
\on{Aut}\wh{\bf{PaB}}_{ell,[n]}\to 
\on{Aut}(\wh T_{ell,n}). 
$$
obtained from (\ref{GTell:to:AutPaB}), (\ref{gp:morph'}) and 
(\ref{Aut:to:AutT}). 
\end{definition}

Theorem \ref{thm:GTell:LM} next states the compatibility of the `arithmetic' action 
$$
\pi_{1}(\tilde M_{1,1}^{\QQ},\tilde\xi) \to 
\on{Aut}(\pi_{1}^{geom}({\mathcal F}_{\xi},\{\sigma_{T}(\xi)\}))
$$
(see (\ref{morphism})) with its `algebraic model' $\wh{\on{GT}}_{ell}\to 
\on{Aut}(\wh T_{ell,n})$ (see Definition \ref{def:act}), namely the 
commutativity of  
\begin{equation} \label{comm:diag:pi1}
\xymatrix{ \pi_{1}(\tilde M_{1,1}^{\QQ},\tilde\xi) \ar[d]\ar[r]& 
\on{Aut}(\pi_{1}^{geom}({\mathcal F}_{\xi},\{\sigma_{T}(\xi)\})) \ar[d]\\
\wh{\on{GT}}_{ell} \ar[r]& \on{Aut}(\wh{T}_{ell,n})}
\end{equation}

The commutativity of the restriction of (\ref{comm:diag:pi1}) to 
$\wh B_{3}\subset \pi_{1}(\tilde M_{1,1}^{\QQ},\tilde\xi)$ can be proved as
follows. Let $\wh{\bf B}_{ell}$ be the category with $\on{Ob}\wh{\bf B}_{ell}=\NN$, 
$\wh{\bf B}_{ell}(n,m) = \emptyset$ if $m\neq n$, and 
$\wh{\bf B}_{ell}(n,n) = \wh B_{1,n}$. There is a natural functor 
$\wh{\bf PaB}_{ell}\to\wh{\bf B}_{ell}$, defined as the length map 
$l:\mu(\{\bullet\})\to\NN$ at the level of objects and as the identity at
the level of morphisms; actually $\wh{\bf PaB}_{ell}\simeq l^{*}
\wh{\bf B}_{ell}$. As $\wh R_{ell}\subset \wh{\on{GT}}_{ell}$ acts trivially
on the images of the associativity constraints under $\wh{\bf PaB}\to
\wh{\bf PaB}_{ell}$, its action on $\wh{\bf PaB}_{ell}$ is the lift of an 
action of $\wh R_{ell}$ on $\wh{\bf B}_{ell}$. One checks explicitly that the 
composition of this action with the morphism $\wh B_{3}\to\wh R_{ell}$
coincides with the action of $\wh B_{3}$  $\wh{\bf B}_{ell}$, which arises from 
its geometric action on the various groups $\wh B_{1,n}$.

The commutativity of the composition of (\ref{comm:diag:pi1}) with 
$G_{\QQ}\stackrel{\sigma}{\to}\pi_{1}(M_{1,\vec1},\vec\xi)$ 
can be shown as follows. As the diagram 
$$
\xymatrix{ G_{\QQ}\ar[r]\ar[d]& \pi_{1}(M_{1,\vec 1}^{\QQ},\vec\xi)\ar[d]\\
\wh{\on{GT}}\ar[r] & \wh{\on{GT}}_{ell}}
$$
commutes, it suffices to proves that its composition 
with (\ref{comm:diag:pi1}) commutes, i.e., that 
$$
\xymatrix{ G_{\QQ}\ar[r]\ar[d]& \on{Aut}(\cF_{\xi},\{\sigma_{T}(\xi)\})
\ar[d]\\
\wh{\on{GT}}\ar[r] & \on{Aut}(\wh T_{ell,n})}
$$
commutes. According to  \cite{Mt}, the morphism 
$G_{\QQ}\to \on{Aut}(\cF_{\xi},\{\sigma_{T}(\xi)\})$
can be derived explicitly from the actions of $G_{\QQ}$ on 
$\pi_{1}(C_{\xi},\xi_{C})$ and on the profinite braid groups in 
genus zero. The former action has been computed in \cite{Na}, 
Cor. 4.5. The resulting formulas for the action of $G_{\QQ}$ can be shown 
to match those for the action of $\wh{\on{GT}}$ on $\wh T_{ell,n}$. 

The last statement of Theorem \ref{thm:GTell:LM} says that the morphism 
$(\pi_{1}(M_{1,\vec1}^{\QQ},\vec\xi)
\stackrel{\curvearrowleft}{\to} G_{\QQ})\to 
(G_{\QQ}\stackrel{\curvearrowleft}{\to} \wh\ZZ^{\times})$
factors through $(\wh{\on{GT}}_{ell}
\stackrel{\curvearrowleft}{\to} \wh{\on{GT}})$. 
This can be proved as follows. First one checks that 
the morphism $\wh B_{3}\to \on{SL}_{2}(\wh\ZZ)$ factors through 
$\wh R_{ell}$. The three morphisms between $\wh B_{3}$, $\wh R_{ell}$
and $\on{SL}_{2}(\wh\ZZ)$ are compatible with the actions of 
$G_{\QQ}$, $\wh{\on{GT}}$ and $\wh\ZZ^{\times}$; and the morphism 
$G_{\QQ}\to \wh\ZZ^{\times}$ factors through $\wh{\on{GT}}$. 
This ends the proof of Theorem \ref{thm:GTell:LM}. 

\subsection{Proof of Proposition \ref{thm:details:GTell}} \label{sec:pf:details}
This statement follows from the form taken by the action of $\wh{\on{GT}}_{ell}$
on $\wh{T}_{ell,n}$. \hfill \qed\medskip

\section{A question} \label{sect:conjs}

In this section, we ask whether ${\mathfrak{r}}_{ell}$ is generated by 
the elements $\delta_{2n}$ arising from \cite{CEE}. This question is analogous to the 
problem of whether ${\mathfrak{grt}}_{1}$ is generated by its Drinfeld generators, which 
is also open. We give an indication in favor of a positive answer: such an answer would imply 
a statement which is also implied by a transcendence conjecture about MZVs; this last
conjecture would follow from Grothendieck's transcendence conjecture for the category 
of mixed Tate motives and the equality of the motivic Galois group $G_{MTM}(-)$ 
with $\on{GT}(-)$ (see \cite{An2}).  
We record that in contrast with the fact that the Drinfeld generators of 
${\mathfrak{grt}}_{1}$ generate a free Lie algebra (Brown), several families of 
relations between the $\delta_{2n}$ have been found (see \cite{Po}). 

\subsection{A generation conjecture (GC)}

The Drinfeld generators of $\grt_1$ are obtained from the homogeneous 
decomposition of the logarithm of
$\on{im}(-1\in\ul{\on{GT}}\stackrel{j_{\Phi_{KZ}}}{\to}
\on{GRT}(\CC))\cdot \on{can}(-1)$, where $\on{can}:\CC^\times
\to\on{GRT}(\CC)$ is the canonical morphism. 
The analogue of the conjecture that these elements 
generate $\grt_1$ is then: 

\begin{conjecture} (Generation Conjecture)
$\b_3\subset \r_{ell}^{gr}$ is an equality, i.e., $\r_{ell}^{gr}$
is generated by $\SL_2$ and the $\delta_{2n}$, $n\geq 0$. 
\end{conjecture}

This conjecture is equivalent to the inclusion $\on{exp}( \hat\b_3^{+,\kk})
\rtimes\on{SL}_2(\kk)\subset R_{ell}^{gr}(\kk)$ being an equality. 

\begin{proposition}
GC is equivalent to the Zariski density of $B_3\subset R_{ell}(-)$, 
i.e., $\langle B_3 \rangle = R_{ell}(-)$.  
\end{proposition}

{\em Proof.} According to Lemma \ref{lemma:Levi}, $\langle B_3\rangle$ 
is uniquely determined by its Lie algebra. This fact and Proposition 
\ref{prop:Lie<SL_2>} immediately imply that $\langle B_3 \rangle = R_{ell}(-)$ 
iff the inclusion $\on{Lie}(u_+,u_-) = \langle \on{log}\psi,\on{log}
\psi_-\rangle \subset \on{Lie}R_{ell}(-)$ is actually an equality. Tensoring 
with $\CC$, this holds iff $$\langle\on{log}\psi,\on{log}\psi_{-}
\rangle^{\CC}\subset \hat{\mathfrak{r}}_{ell}^{\CC}$$ is an equality. 
On the other hand, GC holds iff $\hat\b_3^{+,\CC}\subset 
\hat{\mathfrak r}_{ell}^{gr,\CC}$ is an equality. Now $i_{e_{KZ}}$ sets up 
a diagram 
$$
\begin{matrix}
\langle\on{log}\psi,\on{log}\psi_{-}\rangle^{\CC}
 & \hookrightarrow & {\mathfrak{r}}_{ell}^{\CC}\\
\scriptstyle{\simeq} \downarrow & & \downarrow \scriptstyle{\simeq}\\
\hat\b_3^{+,\CC} & \hookrightarrow & \hat{\mathfrak{r}}_{ell}^{gr,\CC}
\end{matrix}
$$ 
It follows that the upper inclusion is an equality iff the lower is. 
\hfill \qed\medskip

\subsection{Relation with a transcendence conjecture}

We first present the transcendence conjecture. 

\subsubsection{The coordinate ring of associators}

The functors $\{\QQ$-rings$\}\to\{$sets$\}$, $\kk\mapsto \ul{M}(\kk),M(\kk)$
may be represented  as follows. Let $\on{pent}_{\kk} : \on{exp}(\hat\f_{2}^{\kk})
\to \on{exp}(\hat\t_{3}^{\kk})$, $\on{hex}:\kk\times\on{exp}(\hat\f_{2}^{\kk})
\to \on{exp}(\hat\t_{3}^{\kk})$, $\on{dual}:\on{exp}(\hat\f_{2}^{\kk})
\to \on{exp}(\hat\t_{3}^{\kk})$ be the maps $\on{pent}(\Phi):= 
\on{lhs}((\ref{pent}))^{-1} \on{rhs}((\ref{pent}))$,  $\on{hex}(\mu,\Phi):= 
\on{lhs}((\ref{hex:Phi}))\on{rhs}((\ref{hex:Phi}))^{-1}$, 
$\on{dual}(\Phi):= \Phi\Phi^{3,2,1}$. 

Let $B,B',C$ be homogeneous bases of $\f_{2},\t_{3},\t_{4}$. Let 
$\mu,\phi_{b}$ ($b\in B$) be free commutative variables and set 
$\kk_{0}:= \QQ[\varphi_{b},b\in B]$, $\kk_{1}:= \QQ[\mu,\varphi_{b},b\in B]$. 
Then $\kk_{0}\subset \kk_{1}$. Let $\Phi:= \on{exp}(\sum_{b\in B}
\varphi_{b}b)\in\on{exp}(\hat\f_{2}^{\kk_{0}})\subset 
\on{exp}(\hat\f_{2}^{\kk_{1}})$. For $b'\in B'$, $c\in C$, define
$\on{pent}_{c},\on{dual}_{b'}\in \kk_{0}\subset\kk_{1}$ by 
$\sum_{c\in C}\on{pent}_{c}c = \on{log}(\on{pent}(\Phi))$, 
$\sum_{b'\in B'}\on{dual}_{b'}b' = \on{log}(\on{dual}(\Phi))$, 
$\sum_{b'\in B'}\on{hex}_{b'}b' = \on{log}(\on{hex}(\mu,\Phi))$. 

We then set $\QQ[\underline{M}]:= \kk_{1}/(\on{pent}_{c},\on{dual}_{b'},
\on{hex}_{b'}, b'\in B',c\in C)$ and $\QQ[M]:= \QQ[\underline{M}][\mu^{-1}]
= \QQ[\mu^{\pm1},\varphi_{b},b\in B]/\{$ideal with the same generators$\}$. 
Then for any $\QQ$-ring $\kk$, we have (functorial in $\kk$) bijections 
$\begin{matrix} M(\kk)& \simeq & \on{Hom}_{\QQ\text{-rings}}(\QQ[M],\kk)\\
\scriptstyle{\subset}\downarrow & & \downarrow\scriptstyle{\subset}\\
\underline{M}(\kk)  &\simeq  & \on{Hom}_{\QQ\text{-rings}}
(\QQ[\underline{M}],\kk)\end{matrix}$

\subsubsection{The Transcendence Conjecture}

The KZ associator $(2\pi\on{i},\Phi_{KZ})\in M(\CC)$ gives to a morphism 
$\varphi_{KZ}:\QQ[M]\to\CC$. 

\begin{conjecture} (Transcendence Conjecture) $\varphi_{KZ}$ is injective. 
\end{conjecture}

Let $\kk_{MZV}:= \on{im}(\QQ[M]\stackrel{\varphi_{KZ}}{\to}
\CC)$. This is a subring of $\CC$ (according to \cite{LM}, 
this is the subring generated
by $(2\pi\on{i})^{\pm1}$ and the MZVs (multizeta values)). 

\subsection{Consequences of GC}

\begin{proposition} \label{prop:77}
The inclusion
\begin{equation} \label{incl}
i_{e}(B_3) \subset \on{exp}(\hat\b_3^{+,\CC}) \rtimes\on{SL}_{2}(\CC)
\end{equation}
holds: 

(a) for any $e\in Ell(\CC)\times_{M(\CC)}\{\Phi_{KZ}\}$ iff 
$\b_3 \triangleleft {\mathfrak{r}}_{ell}^{gr}$; 

(b) for any $e\in \sigma(M(\CC))$ iff $[\sigma({\mathfrak{grt}}),
\b_3] \subset \hat\b_3$; 

(c) for any $e\in Ell(\CC)$ iff $\b_3 \triangleleft {\mathfrak{grt}}_{ell}$, 
i.e., iff the two above-mentioned conditions are realized. 

Moreover, GC implies that (\ref{incl}) holds for any $e\in Ell(\CC)$. 
\end{proposition}

We do not know whether the Lie algebraic statements in (a), (b), (c)
hold, so they may be viewed as conjectures implied by GC. 
\medskip 

{\em Proof.} Note first that for any $e\in Ell(\CC)$, and by Zariski 
density, (\ref{incl}) $\Leftrightarrow (  i_{e}(\langle
B_3\rangle(\CC)) = \on{exp}(\hat\b_3^{+,\CC})
\rtimes\on{SL}_{2}(\CC)  )$. 

(a) is proved as follows. $e\in Ell(\CC)\times_{M(\CC)}\{\Phi_{KZ}\}$ iff 
$e=\sigma(\Phi_{KZ})*g$ for some $g\in R_{ell}^{gr}(\CC)$. So 
\begin{eqnarray*}
&& ((\ref{incl}) \text{\ holds\ for\ any\ }e\in
Ell(\CC)\times_{M(\CC)}\{\Phi_{KZ}\})\\
&& \Leftrightarrow (g(i_{\sigma(\Phi_{KZ})}(\langle
B_3 \rangle(\CC))) =  \on{exp}(\hat\b_3^{+,\CC})
\rtimes\on{SL}_{2}(\CC)
\text{\ for\ any\ }g\in R_{ell}^{gr}(\CC))\\
&& \Leftrightarrow (g(\Gamma)=\Gamma \text{\ for\ any\ }
g\in R_{ell}^{gr}(\CC), \text{\ where\ }
\Gamma = \on{exp}(\hat\b_3^{+,\CC})
\rtimes\on{SL}_{2}(\CC))\\
&& \Leftrightarrow (\b_3 \triangleleft {\mathfrak{r}}_{ell}^{gr}). 
\end{eqnarray*}
Here the second equivalence follows from Proposition \ref{prop:1:3}. 

(b), (c) are then proved in the same way, using 
$$
(e\in\sigma(M(\CC)))\Leftrightarrow (e=\sigma(\Phi_{KZ})*\sigma(g)
\text{\ for\ some\ }g\in \on{GRT}(\CC)), $$
$$
(e\in Ell(\CC))\Leftrightarrow (e=\sigma(\Phi_{KZ})*g
\text{\ for\ some\ }g\in \on{GRT}_{ell}(\CC)). $$
The equivalence (c) $\Leftrightarrow$ ((a) and (b)) follows from
${\mathfrak{grt}}_{ell} = \sigma({\mathfrak{grt}})\oplus 
{\mathfrak{r}}_{ell}^{gr}$. Finally, GC means that 
$\langle{\mathfrak{sl}}_{2},\delta_{2k}\rangle
={\mathfrak{r}}_{ell}^{gr}$, which immediately implies 
(a), (b) and (c) as ${\mathfrak{r}}_{ell}^{gr}\triangleleft
{\mathfrak{grt}}_{ell}$. \hfill \qed\medskip 

\subsection{Consequences of the Transcendence Conjecture (TC)}

\begin{proposition} \label{prop:78}
If TC holds, then for any $\QQ$-ring $\kk$ and any $\Phi\in M(\kk)$, 
$i_{\sigma(\Phi)}(B_3) \subset 
\on{exp}(\hat\b_3^{+,\kk})\rtimes\on{SL}_{2}(\kk)$. 
\end{proposition}

{\em Proof.} Recall that $\langle B_3 \rangle(-)
\hookrightarrow R_{ell}(-)$, $\on{exp}(\hat\b_3^+)\rtimes
\on{SL}_{2}(-)\hookrightarrow R_{ell}^{gr}(-)$ are inclusions
of $\QQ$-group schemes, and $Ell\to M$, $M\stackrel{\sigma}{\to}M$ are 
morphisms of $\QQ$-group schemes.

In the notation of Definition \ref{def:torsors}, any $x\in X(\kk)$
gives rise to a morphism $i_{x}:G(\kk)\to H(\kk)$, defined by $g*x = 
x*i_{x}(g)$ for any $g\in G(\kk)$. The assignment $x\mapsto i_{x}$ is functorial  
in the following sense: if $\kk\to\kk'$ is a morphism of $\QQ$-rings and 
$x':= \on{im}(x\in X(\kk)\to X(\kk'))$, then 
$$
\begin{matrix}
G(\kk) & \stackrel{i_{x}}{\to}& H(\kk)\\
\downarrow  & & \downarrow \\
G(\kk') & \stackrel{i_{x'}}{\to}& H(\kk')  \end{matrix}
$$
commutes. 

For any $\QQ$-scheme $X$ and any $\QQ$-ring $\kk$, let $X\otimes\kk$
be the $\kk$-scheme $(X\otimes\kk)(\kk'):= X(\kk')$ for any $\kk'\in
\{\kk$-rings$\}$. Again with the notation of Definition \ref{def:torsors}, 
a torsor even gives rise to an assignment $X(\kk)\ni x\mapsto (G\otimes\kk
\stackrel{i^\kk_{x}}{\to} H\otimes\kk)$, where $i_{x}^{\kk}$ is a
morphism of $\kk$-group schemes, defined by: $\forall \kk'\in\{\kk$-rings$\}$, 
$g*\bar x = \bar x*i_{x}^{\kk}(g)$ for any $g\in (G\otimes \kk)(\kk')
= G(\kk')$, where $\bar x:= \on{im}(x\in X(\kk)\to X(\kk'))$.

In particular, $\Phi_{KZ}\in M(\kk_{MZV})$ gives rise to an isomorphism 
$i_{\sigma(\Phi_{KZ})}:
R_{ell}(-)\otimes \kk_{MZV}\stackrel{\sim}{\to} R_{ell}^{gr}(-)
\otimes \kk_{MZV}$, and therefore to a Lie algebra isomorphism 
$\on{Lie}i_{\sigma(\Phi_{KZ})}:
{\mathfrak r}_{ell}\otimes \kk_{MZV}\stackrel{\sim}{\to} 
({\mathfrak r}_{ell}^{gr}\otimes \kk_{MZV})^{\wedge}$, 
whose $\otimes_{\kk_{MZV}}
\CC$ is the infinitesimal of the isomorphism of Proposition \ref{prop:1:3}. 

The group scheme inclusions $\langle B_3\rangle(-)
\subset R_{ell}(-)$ and $\on{exp}(\hat\b_3^+)\rtimes\on{SL}_{2} 
\subset R_{ell}^{gr}(-)$ give rise to Lie algebra inclusions 
$\on{Lie}\langle B_3\rangle(-)\subset {\mathfrak{r}}_{ell}$
and $\hat\b_3^+\subset \hat{\mathfrak{r}}_{ell}^{gr}$ and Proposition 
\ref{prop:1:3} implies that $\on{Lie}i_{\sigma(\Phi_{KZ})}
\otimes_{\kk_{MZV}}\CC$ restricts to an isomorphism 
$\on{Lie}\langle B_3\rangle(-)\otimes_{\QQ}\CC\to (\b_3
\otimes_{\QQ}\CC)^{\wedge}$. This implies that 
$\on{Lie}i_{\sigma(\Phi_{KZ})}$ restricts to a Lie algebra isomorphism
$$\on{Lie}\langle B_3\rangle(-)\otimes_{\QQ}\kk_{MZV}\to 
(\b_3 \otimes_{\QQ}\kk_{MZV})^{\wedge}.$$ 

There are Lie subalgebras $\QQ\on{log}\psi_{\pm}\otimes \kk_{MZV}$
in the l.h.s., mapping to $(\QQ e_{\pm}+$ terms of degree $>0)\otimes\kk_{MZV}$
in the r.h.s. (where $e_{+}=e$, $e_{-}=f$, $\psi_{+}=\psi$). This induces a diagram 
$$
\begin{matrix}
\langle B_3\rangle(-)\otimes\kk_{MZV} 
& \stackrel{i_{\sigma(\Phi_{KZ})}}{\to}& (
\on{exp}(\hat\b_3^+)
\rtimes \on{SL}_{2})\otimes\kk_{MZV}\\
\scriptstyle{\subset}\uparrow  & & \uparrow\scriptstyle{\subset}\\
{\mathbb G}_{a} \otimes\kk_{MZV}  & = & {\mathbb G}_{a} \otimes\kk_{MZV} 
\end{matrix}$$ 
If now $\Phi\in M(\kk)$, the transcendence conjecture says that 
there exists a $\QQ$-ring morphism $\kk_{MZV}\stackrel{\varphi}{\to}\kk$, 
such that $\Phi = \varphi_{*}(\Phi_{KZ})$. Applying this morphism to the above
diagram, one gets 
$$
\begin{matrix}
\langle B_3\rangle(-)\otimes\kk 
& \stackrel{i_{\sigma(\Phi)}}{\to}& (
\on{exp}(\hat\b_3^+)
\rtimes \on{SL}_{2})\otimes\kk\\
\scriptstyle{\subset}\uparrow  & & \uparrow\scriptstyle{\subset}\\
{\mathbb G}_{a} \otimes\kk  & = & {\mathbb G}_{a} \otimes\kk 
\end{matrix}$$ 
Taking $\kk$-points, one obtains a commutative diagram
$$
\begin{matrix}
\langle B_3 \rangle(\kk) 
& \stackrel{i_{\sigma(\Phi)}}{\to}& 
\on{exp}(\hat\b_3^{+,\kk})
\rtimes \on{SL}_{2}(\kk) \\
\scriptstyle{\subset}\uparrow  & & \uparrow\scriptstyle{\subset}\\
\kk  & = & \kk 
\end{matrix}$$ 
The image of $1\in\kk$ is $\Psi_{\pm} \subset \langle B_3\rangle(\kk)$; 
then $i_{\sigma(\Phi)}(\Psi_{\pm})\in \on{exp}(\hat\b_3^{+,\kk})
\rtimes \on{SL}_{2}(\kk) \subset \on{exp}(\hat\b_3^{+,\CC})
\rtimes \on{SL}_{2}(\CC)$. \hfill \qed


\medskip 










\begin{thebibliography}{DDMS}

\bibitem[An1]{An} Y.~Andr\'e, Une introduction aux motifs (motifs purs, 
motifs mixtes, p\'eriodes). Panoramas et Synth\`eses, no. 17, Soci\'et\'e 
Math\'ematique de France, Paris, 2004. 

\bibitem[An2]{An2} Y.~Andr\'e, Galois theory, motives and transcendental numbers. 
Preprint arXiv/0805.2569; in `Renormalization and Galois Theories' (A. Connes, 
F. Fauvet, J.-P. Ramis, eds.), IRMA Lectures in Math. and Theor. Phys., no. 15, 
European Math. Soc., Z\"urich, 2009. 

\bibitem[Ba]{Ba} D.~Bar-Natan, On associators and the Grothendieck-Teichm\"uller 
group. I. Selecta Math. (N.S.) 4 (1998), no. 2, 183--212. 


\bibitem[Bi]{Bi} J.~Birman, On braid groups, Commun. Pure Appl. Math. 
12 (1969), 41-79. 

\bibitem[BS]{BS} A.~Borel, J.-P.~Serre, Th\'eor\`emes de finitude en 
cohomologie galoisienne. Comment. Math. Helv. 39 1964 111--164. 

\bibitem[BL]{BL} F.~Brown, A.~Levin, Multiple elliptic polylogarithms. 
Preprint arXiv/1009.2652. 

\bibitem[CEE]{CEE} D.~Calaque, B.~Enriquez, P.~Etingof,  Universal KZB 
equations I: the elliptic case. Preprint arXiv:math/0702670. In 
'Algebra, arithmetic, and geometry: in honor of Yu. I. Manin' (Yu. Tschinkel, Yu. Zarhin, eds.). 
Vol. I, 165Ð266,
Progr. Math., 269, Birkh\"auser Boston, Inc., Boston, MA, 2009. 


\bibitem[DDMS]{DDMS} J.~Dixon, M.~du Sautoy, A.~Mann, D.~Segal, 
Analytic pro-$p$ groups. Second edition. Cambridge Studies in 
Advanced Mathematics, 61. Cambridge University Press, Cambridge, 1999. 

\bibitem[Dr]{Dr:Gal} V.~Drinfeld, On quasitriangular quasi-Hopf algebras and a 
group closely connected with $\on{Gal}(\bar\QQ/\QQ)$, Leningrad Math. J. 2 
(1991), 829-860.


\bibitem[Gr1]{G:Esq} A.~Grothendieck, Esquisse d'un programme. 
London Math. Soc. Lecture Note Ser., 242,  Geometric Galois actions, 1, 
5--48, Cambridge Univ. Press, Cambridge, 1997. 

\bibitem[Gr2]{G:LM} A.~Grothendieck, La longue marche \`a travers la th\'eorie 
de Galois. Tome 1. Transcription d'un manuscrit in\'edit. \'Edit\'e par Jean 
Malgoire. Universit\'e Montpellier II, D\'epartement des Sciences Math\'ematiques, 
Montpellier, 1995. 

\bibitem[Ha]{H} R.~Hain, Completions of mapping class groups and the 
cycle $C - C^{-}$, Contemp. Math., vol. 150 (1993). 

\bibitem[HM1]{HM} R.~Hain, M.~Matsumoto, Weighted completion of Galois groups 
and Galois actions on the fundamental group of ${\mathbb P}^1-\{0,1,\infty\}$.  
Compositio Math.  139  (2003),  no. 2, 119--167. 
 
\bibitem[HM2]{HM:ell} R.~Hain, M.~Matsumoto, Universal elliptic motives. Talk 
given at Heilbronn Institute, Bristol (2011). 

\bibitem[IN]{Ih:Nak} Y.~Ihara, H.~Nakamura, On deformation of maximally degenerate 
stable marked curves and Oda's problem. J. reine. angew. Math., 487 (1997), 125-151. 

\bibitem[JS]{JS} A.~Joyal, R.~Street, Braided tensor categories. Adv. Math.  
102  (1993),  no. 1, 20--78. 

\bibitem[Ka]{Ka} C.~Kassel, Quantum groups. Graduate Texts in Mathematics, 155. 
Springer-Verlag, New York, 1995. 

\bibitem[LM]{LM} T.T.Q.~Le, J.~Murakami, Kontsevich's integral for the 
Kauffman polynomial. Nagoya Math. J. 142 (1996), 39--65. 

\bibitem[LR]{LR} A.~Levin, G.~Racinet, Towards multiple 
elliptic polylogarithms, arXiv:math/0703237. 

\bibitem[Ma]{Ma} Yu.~Manin, Iterated integrals of modular forms 
and noncommutative modular symbols.  Algebraic geometry and number 
theory,  565--597, Progr. Math., 253, Birkh\"auser Boston, Boston, MA, 2006. 

\bibitem[Mt]{Mt} M.~Matsumoto, Galois representations on profinite 
braid groups on curves. J. reine angew. Math. 474 (1996), 169--219. 

\bibitem[Mm]{Mum} D.~Mumford, An analytic construction of degenerating curves 
over complete local rings. Compositio Math. 24 (1972), 129-174. 

\bibitem[Mo]{Mo} G.~Mostow, Fully reducible subgroups of algebraic groups.
Amer. J. Math. 78 (1956), 200--221. 

\bibitem[Na]{Na} H.~Nakamura, Limits of Galois representations 
in fundamental groups along maximal degeneration of marked curves. 
I.  Amer. J. Math.  121  (1999),  no. 2, 315--358. 

\bibitem[Ni1]{N} J.~Nielsen, Die Isomorphismen der allgemeinen, 
unendlichen Gruppe mit zwei Erzeugenden. Math. Ann. 78 (1918/1964), 
no. 1, 385--397. 

\bibitem[Ni2]{N:danish:1921} J.~Nielsen, Om regning me ikke-kommutative 
faktorer og dens anvendelse i gruppeteorie (Danish), Math. Tidsskrift B, 1921, 
78-94 (1921).

\bibitem[Oda]{O} T.~Oda, Etale homotopy type of the moduli spaces of 
algebraic curves.  Geometric Galois actions, 1,  85--95, London 
Math. Soc. Lecture Note Ser., 242, Cambridge Univ. Press, Cambridge, 
1997. 

\bibitem[Po]{Po} A.~Pollack, Relations between derivations arising from 
modular forms, Ph.D. thesis, Duke University, 2009. 

\bibitem[Sch]{Sch} L.~Schneps, Groupoides 
fondamentaux des espaces de modules en genre 0 et cat\'egories tensorielles 
tress\'ees, pp. 71--116, in Espaces de modules des courbes, groupes modulaires 
et th\'eorie des champs.  Panoramas et Synth\`eses, 7. Soci\'et\'e Math\'ematique de 
France, Paris, 1999. 

\bibitem[Se]{S} J.-P.~Serre, Corps locaux. Deuxi\`eme \'edition. Publications de 
l'Universit\'e de Nancago, No. VIII. Hermann, Paris, 1968. 

\bibitem[TY]{TY} P.~Tauvel, R.~Yu, Lie algebras and algebraic groups. 
Springer Monographs in Mathematics. Springer-Verlag, Berlin, 2005. 

\bibitem[Ts]{Ts} H.~Tsunogai, The stable derivation algebra for higher genera. 
Israel J. Math., 136 (2003), 221--250. 

\bibitem[Za]{Za} D.~Zagier, Periods of modular forms and Jacobi theta functions.
Invent. Math. 104 (1991), no. 3, 449--465. 

\end{thebibliography}
\end{document}